%% file: AJL-nongeneric.tex
\documentclass{memo-l}

\usepackage[english]{babel}
\usepackage{amssymb}
\usepackage{graphicx}
\usepackage{ltxtable} 
\usepackage{setspace}
\usepackage{float}
\usepackage{amsrefs}

\floatstyle{plaintop}
\restylefloat{table}
\numberwithin{table}{chapter}
\numberwithin{section}{chapter}

\newtheorem{maintheorem}{Theorem}
\newtheorem{maincorollary}[maintheorem]{Corollary}
\newtheorem{theorem}{Theorem}[chapter]
\newtheorem{lemma}[theorem]{Lemma}
\newtheorem{proposition}[theorem]{Proposition}
\newtheorem{corollary}[theorem]{Corollary}

\theoremstyle{definition}
\newtheorem{definition}[theorem]{Definition}
\newtheorem{remark}[theorem]{Remark}
\newtheorem*{remark*}{Remark}

\newcommand{\nongcr}{\textbf{N}}
\newcommand{\possprim}{\textbf{P}}

\DeclareMathOperator{\Alt}{Alt}
\DeclareMathOperator{\Sym}{Sym}

\makeindex

\begin{document}

\frontmatter

\title{On Non-Generic Finite Subgroups of Exceptional Algebraic Groups}

\author{Alastair J. Litterick}
\address{University of Auckland, Auckland, New Zealand}
\email{ajlitterick@gmail.com}
\thanks{The author acknowledges support from the University of Auckland, as well as a Doctoral Training Award from the UK's Engineering and Physical Sciences Research Council. A large portion of the results presented here were derived during the author's doctoral study at Imperial College London, under the supervision of Professor Martin Liebeck. The author would like to thank Professor Liebeck for guidance and support received. Finally, the author would like to thank the anonymous referee for their numerous helpful suggestions.}



\begin{abstract}
The study of finite subgroups of a simple algebraic group $G$ reduces in a sense to those which are almost simple. If an almost simple subgroup of $G$ has a socle which is not isomorphic to a group of Lie type in the underlying characteristic of $G$, then the subgroup is called \emph{non-generic}. This paper considers non-generic subgroups of simple algebraic groups of exceptional type in arbitrary characteristic.

A finite subgroup is called Lie primitive if it lies in no proper subgroup of positive dimension. We prove here that many non-generic subgroup types, including the alternating and symmetric groups $\Alt_{n}$, $\Sym_{n}$ for $n \ge 10$, do not occur as Lie primitive subgroups of an exceptional algebraic group.

A subgroup of $G$ is called $G$-completely reducible if, whenever it lies in a parabolic subgroup of $G$, it lies in a conjugate of the corresponding Levi factor. Here, we derive a fairly short list of possible isomorphism types of non-$G$-completely reducible, non-generic simple subgroups.

As an intermediate result, for each simply connected $G$ of exceptional type, and each non-generic finite simple group $H$ which embeds into $G/Z(G)$, we derive a set of \emph{feasible characters}, which restrict the possible composition factors of $V \downarrow S$, whenever $S$ is a subgroup of $G$ with image $H$ in $G/Z(G)$, and $V$ is either the Lie algebra of $G$ or a non-trivial Weyl module for $G$ of least dimension.

This has implications for the subgroup structure of the finite groups of exceptional Lie type. For instance, we show that for $n \ge 10$, $\Alt_n$ and $\Sym_n$, as well as numerous other almost simple groups, cannot occur as a maximal subgroup of an almost simple group whose socle is a finite simple group of exceptional Lie type.
\end{abstract}

\maketitle

\tableofcontents

\mainmatter

\include{AJL-intro}
\include{AJL-background}

\include{AJL-disproving}
\include{AJL-stab}

\include{AJL-gcr}

\include{AJL-thetables}
\include{AJL-auxiliary}

\backmatter

\bibliography{AJL-biblio.bib}

\printindex

\end{document}

%% file: AJL-intro.tex

\chapter{Introduction and Results} \label{chap:intro}

This paper concerns the closed subgroup structure of simple algebraic groups of exceptional type over an algebraically closed field. Along with the corresponding study for classical-type groups, this topic has been studied extensively ever since Chevalley's classification in the 1950s of reductive groups via their root systems.

The subgroup structure of a simple algebraic group divides naturally according to properties of the subgroups under consideration. Some of the first major results were produced by Dynkin \cites{MR0049903,MR0047629}, who classified maximal \emph{connected} subgroups in characteristic zero. The corresponding study over fields of positive characteristic was initiated by Seitz \cites{MR888704,MR1048074} and extended by Liebeck and Seitz \cites{MR1066572,MR2044850}, giving a classification of maximal connected closed subgroups, and more generally maximal closed subgroups of positive dimension.

Here, we consider the opposite extreme, the case of finite subgroups. Let $G$ be a simple algebraic group of exceptional type, over an algebraically closed field $K$ of characteristic $p \ge 0$. When considering finite subgroups of $G$, it is natural to restrict attention to those not lying in any intermediate closed subgroup of positive dimension. Such a finite subgroup is called \emph{Lie primitive}, and a result due to Borovik (Theorem \ref{thm:boro}) reduces the study of Lie primitive subgroups to those which are \emph{almost simple}, that is, those groups $H$ such that $H_{0} \le H \le \textup{Aut}(H_{0})$ for some non-abelian finite simple group $H_{0}$. The study then splits naturally according to whether or not $H_{0}$ is isomorphic to a member of $\textup{Lie}(p)$, the collection of finite simple groups of Lie type in characteristic $p$ (with the convention $\textup{Lie}(0) = \varnothing$).

Subgroups isomorphic to a member of $\textup{Lie}(p)$ are generally well-understood. If $G$ has a semisimple subgroup $X$ and $p > 0$, then $G$ has a subgroup isomorphic to a central extension of the corresponding finite simple group $X(q)$ for each power $q$ of $p$. In \cite{MR1458329} it is shown that, for $q$ above some explicit bound depending on the root system of $G$, any such subgroup of $G$ arises in this manner.

Turning to those simple groups not isomorphic to a member of $\textup{Lie}(p)$, the so-called \emph{non-generic} case, we have two problems: which of these admit an embedding into $G$, and how do they embed? For example, are there any Lie primitive such subgroups? If not, what can we say about the positive-dimensional overgroups which occur?

The question of which simple groups admit an embedding now has a complete answer, due to the sustained efforts of many authors \cites{MR1220771,MR1087228,MR1721818,MR1767573,MR689418,MR1416728,MR933426,MR1653177,MR1717629}. The end product of these efforts is Theorem 2 of \cite{MR1717629}, where Liebeck and Seitz present tables detailing precisely which non-generic simple groups admit an embedding into an exceptional algebraic group, and for which characteristics this occurs. The following is a summary of this information.

\addtocounter{maintheorem}{-1}
\addtocounter{table}{-1}

\begin{maintheorem}[{\cite{MR1717629}}] \label{thm:subtypes}
Let $G$ be an adjoint exceptional algebraic group, over an algebraically closed field of characteristic $p \ge 0$, and let $H$ be a non-abelian finite simple group, not isomorphic to a member of $\textup{Lie}(p)$. Then $G$ has a subgroup isomorphic to $H$ if and only if $(G,H,p)$ appears in Table \ref{tab:subtypes}.
\end{maintheorem}

\begin{table}[htbp]
\small \onehalfspacing
\caption{Non-generic subgroup types of adjoint exceptional $G$.}
\label{tab:subtypes}
\begin{tabularx}{\linewidth}{l|>{\centering\arraybackslash}X}
\hline $G$ & $H$ \\ \hline
$G_2$ & $\Alt_5$, $L_2(q) \ (q = 7,8,13)$, $U_3(3)$, \\
        & $J_1 \ (p = 11)$, $J_2 \ (p = 2)$ \vspace{3pt} \\ \hline
$F_4$ & $\Alt_{5-6}$, $L_2(q) \ (q = 7,8,13,17,25,27)$, $L_3(3)$, $U_3(3)$, $^{3}D_4(2)$, \\
        & $\Alt_{7} \ (p = 2, 5)$, $\Alt_{9-10} \ (p = 2)$, $M_{11} \ (p = 11)$, $J_1 \ (p = 11)$, $J_2 \ (p = 2)$, $L_4(3) \ (p = 2)$ \vspace{3pt} \\ \hline
$E_6$ & $\Alt_{5-7}$, $M_{11}$, $L_2(q) \ (q = 7,8,11,13,17,19,25,27)$, $L_3(3)$, $U_3(3)$, $U_4(2)$, $^{3}D_4(2)$, $^{2}F_4(2)'$, \\
        & $\Alt_{9-12} \ (p = 2)$, $M_{12} \ (p = 2, 5)$, $M_{22} \ (p = 2)$, $J_1 \ (p = 11)$, $J_2 \ (p = 2)$, $J_3 \ (p = 2)$, $Fi_{22} \ (p = 2)$, $L_4(3) \ (p = 2)$, $U_4(3) \ (p = 2)$, $\Omega_7(3) \ (p = 2)$, $G_2(3) \ (p = 2)$ \vspace{3pt} \\ \hline
$E_7$ & $\Alt_{5-9}$, $M_{11}$, $M_{12}$, $J_2$, $L_2(q) \ (q = 7,8,11,13,17,19,25,27,29,37)$, $L_3(3)$, $L_3(4)$, $U_3(3)$, $U_3(8)$, $U_4(2)$, $Sp_6(2)$, $\Omega_8^{+}(2)$, $^{3}D_4(2)$, $^{2}F_4(2)'$, \\
        & $\Alt_{10} \ (p = 2,5)$, $\Alt_{11-13} \ (p = 2)$, $M_{22} \ (p = 5)$, $J_1 \ (p = 11)$, $Ru \ (p = 5)$, $HS \ (p = 5)$, $L_4(3) \ (p = 2)$, ${}^{2}B_{2}(8)$ $(p = 5)$\footnotemark \vspace{3pt} \\ \hline
$E_8$ & $\Alt_{5-10}$, $M_{11}$, $L_2(q) \ (q = 7,8,11,13,16,17,19,25,27,29,31,32,41,49,61)$, $L_3(3)$, $L_3(5)$, $U_3(3)$, $U_3(8)$, $U_4(2)$, $Sp_6(2)$, $\Omega_8^{+}(2)$, $G_2(3)$, $^{3}D_4(2)$, $^{2}F_4(2)'$, $^{2}B_2(8)$, \\
        & $\Alt_{11} \ (p = 2,11)$, $\Alt_{12-17} \ (p = 2)$, $M_{12} \ (p = 2, 5)$, $J_1 \ (p = 11)$, $J_2 \ (p = 2)$, $J_3 \ (p = 2)$, $Th \ (p = 3)$, $L_2(37) \ (p = 2)$, $L_4(3) \ (p = 2)$, $L_4(5) \ (p = 2)$, $PSp_4(5) \ (p = 2)$, $^{2}B_2(32) \ (p = 5)$ \\ \hline
\end{tabularx}
\end{table}
\footnotetext{This subgroup was erroneously omitted from \cite{MR1717629}, however
such a subgroup is easily seen to exist, since an $8$-dimensional module for
the double-cover of ${}^{2}B_{2}(8)$ gives an embedding into $A_7$.}

Note that we have the following isomorphisms:
\[
\begin{array}{c}
\Alt_5 \cong L_2(4) \cong L_2(5), \quad \Alt_6 \cong L_2(9) \cong Sp_4(2)', \quad \Alt_8 \cong L_4(2),\\
L_2(7) \cong L_3(2), \quad U_4(2) \cong PSp_4(3), \quad U_3(3) \cong G_2(2)',
\end{array}
\]
and we therefore consider these groups to be of Lie type in each such characteristic.

\addtocontents{toc}{\setcounter{tocdepth}{-1}}
\section*{Lie Primitivity}
\addtocontents{toc}{\setcounter{tocdepth}{1}}

At this point it remains to determine, for each pair $(G,H)$ above, properties of the possible embeddings $H \to G$. Here, our aim is to determine whether or not $H$ occurs as a Lie primitive subgroup of $G$, and if not, to determine useful information regarding the intermediate subgroups of positive dimension.

Our main result is the following. Let $\textup{Aut}(G)$ be the abstract group generated by all inner automorphisms of $G$, as well as graph and field morphisms.

\begin{maintheorem} \label{THM:MAIN}
Let $G$ be an adjoint exceptional simple algebraic group, over an algebraically closed field of characteristic $p \ge 0$, and let $S$ be a subgroup of $G$, isomorphic to the finite simple group $H \notin \textup{Lie}(p)$.

If $G$, $H$, $p$ appear in Table \ref{tab:main} then $S$ lies in a proper, closed, connected subgroup of $G$ which is stable under all automorphisms in $N_{\textup{Aut}(G)}(S)$.
\end{maintheorem}

\begin{table}[htbp]
\centering \small \onehalfspacing
\caption{Imprimitive subgroup types.} \label{tab:main}
\begin{tabularx}{\linewidth}{l|>{\centering\arraybackslash}X}
\hline $G$ & $H$ \\ \hline
$F_4$ & $\Alt_5$, $\Alt_{7}$, $\Alt_{9-10}$, $M_{11}$, $J_1$, $J_2$, \\ & $\Alt_6 \ (p = 5)$, $L_2(7) \ (p = 3)$, $L_2(17) \ (p = 2)$, $U_3(3) \ (p \neq 7)$ \\ \hline
$E_6$ & $\Alt_{9-12}$, $M_{22}$, $L_2(25)$, $L_2(27)$, $L_4(3)$, $U_4(2)$, $U_4(3)$, $^{3}D_4(2)$, \\ & $\Alt_{5} \ (p \neq 3)$, $\Alt_7 \ (p \neq 3, 5)$, $M_{11} \ (p \neq 3, 5)$, $M_{12} \ (p = 2)$, $L_2(7) \ (p = 3)$, $L_{2}(8) \ (p = 7)$, $L_{2}(11) \ (p = 5)$, $L_{2}(13) \ (p = 7)$, $L_2(17) \ (p = 2)$, $L_3(3) \ (p = 2)$, $U_3(3) \ (p = 7)$ \\ \hline
$E_7$ & $\Alt_{10-13}$, $M_{11}$, $J_1$, $L_2(17)$, $L_2(25)$, $L_3(3)$, $L_4(3)$, $U_4(2)$, $Sp_{6}(2)$, $^{3}D_4(2)$, $^{2}F_4(2)'$, \\ & $\Alt_{9} \ (p \neq 3)$, $\Alt_{8} \ (p \neq 3,5)$, $\Alt_{7} \ (p \neq 5)$, $M_{12} \ (p \neq 5)$, $J_2 \ (p \neq 2)$, $L_2(8) \ (p \neq 3,7)$ \\ \hline
$E_8$ & $\Alt_{8}$, $\Alt_{10-17}$, $M_{12}$, $J_1$, $J_2$, $L_2(27)$, $L_2(37)$, $L_4(3)$, $U_3(8)$, $Sp_{6}(2)$, $\Omega_{8}^{+}(2)$, $G_2(3)$, \\ & $\Alt_{9} \ (p \neq 2,3)$, $M_{11} \ (p \neq 3,11)$, $U_{3}(3) \ (p \neq 2,3,7)$, $^{2}F_4(2)' \ (p \neq 3)$ \\ \hline
\end{tabularx}
\end{table}

The `normaliser stability' condition here is important for applications to the subgroup structure of finite groups of Lie type. For instance, if $G = F_{4}(K)$ where $K$ has characteristic $5$, then each subgroup $S \cong \Alt_{6}$ of $G$ lies in a proper, closed connected subgroup of $G$, which is stable under $N_{G}(S)$ and under every Frobenius morphism of $G$ stabilising $S$. This can be used to deduce that no finite almost simple group with socle $F_{4}(5^{r})$ $(r \ge 1)$ contains a maximal subgroup with socle $\Alt_{6}$, as in Theorem \ref{THM:FIN} below.

The proof of Theorem \ref{THM:MAIN} proceeds by considering the action of the finite groups in question on certain modules for the algebraic group. We consider in particular the adjoint module $L(G)$, and a Weyl module of least dimension, denoted $V_{\textup{min}}$, for the simply connected cover $\tilde{G}$ of $G$. Each such module has dimension at most 248. For $G$ of type $E_{6}$ and $E_{7}$, the group $\tilde{G}$ has a centre of order $3/(3,p)$ and $2/(2,p)$ respectively, acting by scalars on $V_{\textup{min}}$. Thus to make use of $V_{\textup{min}}$, we must consider embeddings $\tilde{H} \to \tilde{G}$, where $\tilde{H}$ is a perfect central extension of the simple group $H$, and where the image of $Z(\tilde{H})$ lies in $Z(\tilde{G})$.

For each $(G,H,p)$ in Table \ref{tab:subtypes}, we calculate \emph{compatible feasible characters} of $\tilde{H}$ on $L(G)$ and $V_{\textup{min}}$, which are Brauer characters of $\tilde{H}$ that agree with potential restrictions of $L(G)$ and $V_{\textup{min}}$ to a subgroup $\tilde{S} \cong \tilde{H}$ of $\tilde{G}$ with $Z(\tilde{S}) \le Z(\tilde{G})$ (see Definition \ref{def:feasible}). This requires knowledge of the irreducible Brauer characters of $\tilde{H}$ of degree at most $\textup{dim}(L(G))$ for each subgroup type $H$, as well as knowledge of the eigenvalues of various semisimple elements of $\tilde{G}$ on the relevant $K\tilde{G}$-modules. The necessary theory is already well-developed, and we give an outline in Chapter \ref{chap:disproving}.

\begin{maintheorem} \label{THM:FEASIBLES}
Let $G$ be a simple algebraic group of type $F_4$, $E_6$, $E_7$ or $E_8$ over an algebraically closed field of characteristic $p \ge 0$. Let $S$ be a non-abelian finite simple subgroup of $G$, not isomorphic to a member of $\textup{Lie}(p)$, and let $\tilde{S}$ be a minimal preimage of $S$ in the simply connected cover of $G$. Then the composition factors of $L(G) \downarrow \tilde{S}$ and $V_{\textup{min}} \downarrow \tilde{S}$ are given by a line of the appropriate table in Chapter \ref{chap:thetables}.
\end{maintheorem}

Once these composition factors are known, the representation theory of the finite subgroup can be used to determine further information on the structure of $L(G)$ and $V_{\textup{min}}$ as an $\tilde{S}$-module, which in turn allows us to determine information on the inclusion of $S$ into $G$. In particular, if $\tilde{S}$ fixes a nonzero vector on a non-trivial $G$-composition factor of $L(G)$ or $V_{\textup{min}}$, then $\tilde{S}$ lies in the full stabiliser of this vector, which is positive-dimensional and proper. In Chapter \ref{chap:disproving}, we present several techniques for determining the existence of such a fixed vector.

Theorem \ref{THM:MAIN} complements a number of existing results. Frey \cites{Fre1,MR1617620,MR1839999,MR1423302} and Lusztig \cite{MR1976697} give much information on embeddings of alternating groups $\Alt_{n}$ and their proper covers into exceptional groups in characteristic zero. Lifting and reduction results such as \cite{MR1369427}*{Lemme 4 and Proposition 8} and \cite{MR1320515}*{Corollary 2.2 and Theorem 3.4} can then be used to pass results between characteristic zero and positive characteristic not dividing the order of the subgroup in question.

For $G$ of type $G_2$, embeddings of non-generic finite simple groups are well understood, and hence we omit these from our study. In particular, by \cite{MR1717629}*{Corollary 4}, the only non-generic finite simple subgroups of $G_2(K)$ which are not Lie primitive are isomorphic to $\Alt_{5}$ or $L_2(7)$. From \cite{MR898346}*{Theorems 8, 9}, it follows that $\Alt_{5}$ does not occur as a Lie primitive subgroup of $G_{2}(K)$, while $L_2(7)$ occurs both as a Lie primitive subgroup and as a non-Lie primitive subgroup. In addition, Magaard \cite{MR2638705} gives a necessary condition for a given non-generic simple group to occur as a maximal subgroup of $F_4(F)$, where $F$ is a finite or algebraically closed field of characteristic $\neq 2,3$. This can be used to limit the possible isomorphism types of Lie primitive finite subgroups. The methods used here differ from those of the above references, as well as from relevant papers of Aschbacher \cites{MR892190,MR928524,MR986684,MR1054997,MR898346}. There, results focus on the geometry of certain low-dimensional modules supporting a multilinear form or algebra product. Here, however, we primarily use techniques from the representation theory of finite groups and reductive algebraic groups.

For each adjoint exceptional group $G$, Corollary 4 of \cite{MR1717629} lists those non-generic finite simple groups $H$ which embed into $G$ and occur only as Lie primitive subgroups there. Since we will refer to this in Section \ref{sec:reps}, we record this information below in Table \ref{tab:onlyprim}. Combining this with Theorem \ref{THM:MAIN} gives the following:

\begin{maincorollary} \label{cor:prim}
Let $G$ be an adjoint simple algebraic group of type $F_4$, $E_6$, $E_7$ or $E_8$, over an algebraically closed field of characteristic $p \ge 0$, and let $H \notin \textup{Lie}(p)$ be a finite non-abelian simple group which embeds into $G$. Then exactly one of the following holds:
\begin{itemize}
\item[\textup{(i)}] $G$, $H$ appear in Table \ref{tab:main} and $G$ has no Lie primitive subgroups $\cong H$.
\item[\textup{(ii)}] $G$, $H$ appear in Table \ref{tab:onlyprim} and every subgroup $S \cong H$ of $G$ is Lie primitive.
\item[\textup{(iii)}] $G$, $H$ appear in Table \ref{tab:unclassified}.
\end{itemize}
\end{maincorollary} 

\begin{table}[ht]
\centering \small \onehalfspacing
\caption{Subgroup types occurring only as Lie primitive subgroups of $G$} \label{tab:onlyprim}
\begin{tabularx}{.99\linewidth}{c|>{\centering\arraybackslash}X}
\hline $G$ & $H$ \\ \hline
$F_4$ & $L_2(25)$, $L_2(27)$, $L_3(3)$, $^{3}D_4(2)$, $L_4(3) \ (p = 2)$ \\ \hline
$E_6$ & $L_2(19)$, $^{2}F_4(2)'$, $\Omega_7(3) \ (p = 2)$, $G_2(3) \ (p = 2)$, $M_{12} \ (p = 5)$, $J_3 \ (p = 2)$, $Fi_{22} \ (p = 2)$ \\ \hline
$E_7$ & $L_2(29)$, $L_2(37)$, $U_3(8)$, $M_{22} \ (p = 5)$, $Ru \ (p = 5)$, $HS \ (p = 5)$ \\ \hline
$E_8$ & $L_2(q) \ (q = 31,32,41,49,61)$, $L_3(5)$, $L_4(5) \ (p = 2)$, $^{2}B_2(8) \ (p \neq 5, 13)$, $^{2}B_2(32) \ (p = 5)$, $Th \ (p = 3)$ \\ \hline
\end{tabularx}
\end{table}

\begin{table}[htbp]
\centering \small \onehalfspacing
\caption{Subgroup types $H \notin \textup{Lie}(p)$ possibly occurring both as a Lie primitive subgroup and as a Lie imprimitive subgroup of $G$}
\begin{tabularx}{.99\linewidth}{c|>{\centering\arraybackslash}X}
\hline $G$ & $H$ \\ \hline
$F_4$ & $L_2(q) \ (q = 8,13)$, $\Alt_{6} \ (p \neq 5)$, $L_2(7) \ (p \neq 3)$, $L_2(17) \ (p \neq 2)$, $U_3(3) \ (p = 7)$ \vspace{3pt} \\ \hline
$E_6$ & $\Alt_{6}$, $\Alt_5 \ (p = 3)$, $\Alt_7 \ (p = 3,5)$, $M_{11} \ (p = 3,5)$, $J_1 \ (p = 11)$, $J_2 \ (p = 2)$, $L_2(7) \ (p \neq 3)$,  $L_2(8) \ (p \neq 7)$, $L_{2}(11) \ (p \neq 5)$, $L_{2}(13) \ (p \neq 7)$, $L_2(17) \ (p \neq 2)$, $L_3(3) \ (p \neq 2)$, $U_3(3) \ (p \neq 7)$ \vspace{3pt} \\ \hline
$E_7$ & $\Alt_{5}$, $\Alt_{6}$, $L_2(q) \ (q = 7,11,13,19,27)$, $L_3(4)$, $U_3(3)$, $\Omega_8^{+}(2)$, $\Alt_{7} \ (p = 5)$, $\Alt_8 \ (p = 3, 5)$, $\Alt_9 \ (p = 3)$, $M_{12} \ (p = 5)$, $J_2 \ (p = 2)$,  $L_2(8) \ (p = 3,7)$, \\ & ${}^{2}B_{2}(8)$ $(p = 5)$ \vspace{3pt} \\ \hline
$E_8$ & $\Alt_{5-7}$, $L_2(q) \ (q = 7,8,11,13,16,17,19,25,29)$, $L_3(3)$, $U_4(2)$, $^{3}D_4(2)$, $\Alt_9 \ (p = 2,3)$, $M_{11} \ (p = 3,11)$, $J_3 \ (p = 2)$, $U_3(3) \ (p = 7)$, $PSp_4(5) \ (p = 2)$, $^{2}B_2(8) \ (p = 5, 13)$, $^{2}F_4(2)' \ (p = 3)$ \\ \hline
\end{tabularx}
\label{tab:unclassified}
\end{table}

We remark that many of the groups appearing in Table \ref{tab:onlyprim} admit only a single feasible character on $L(G)$ and $V_{\textup{min}}$ in each characteristic where they embed. The exceptions to this occur primarily for subgroups $L_{2}(q)$, and these tend to have elements of large order, whose Brauer character values we have not considered (further on this is given at the start of Chapter \ref{chap:thetables}). It is therefore likely that a small number of the feasible characters for these groups are not realised by an embedding, so that these groups also have very few possible actions on $L(G)$ and $V_{\textup{min}}$.

\addtocontents{toc}{\setcounter{tocdepth}{-1}}
\section*{Complete Reducibility}
\addtocontents{toc}{\setcounter{tocdepth}{1}}

In the course of proving Theorem \ref{THM:MAIN}, we derive much information on the intermediate positive-dimensional subgroups which occur when a finite group is shown not to appear as a Lie primitive subgroup of $G$. Recall that a subgroup $X$ of a reductive algebraic group $G$ is called \emph{$G$-completely reducible} (in the sense of Serre \cite{Ser3}) if, whenever $X$ is contained in a parabolic subgroup $P$ of $G$, it lies in a Levi factor of $P$. In case $G = GL_n(K)$ or $SL_n(K)$, this coincides with $X$ acting completely reducibly on the natural $n$-dimensional module.

The concept of $G$-complete reducibility has seen much attention recently, for instance in the context of classifying reductive subgroups of $G$, see for example \cites{MR2604850,MR3075783}. Similar techniques can be brought to bear in the case of finite simple subgroups.

\begin{maintheorem}
Let $G$ be an adjoint exceptional simple algebraic group in characteristic $p > 0$, and let $H$ be a non-abelian finite simple group, not isomorphic to a member of $\textup{Lie}(p)$, which admits an embedding into $G$.

If $G$ has a subgroup isomorphic to $H$, which is not $G$-completely reducible, then $(G,H)$ appears in Table \ref{tab:nongcr}, with $p$ one of the primes given there.
\label{THM:NONGCR}
\end{maintheorem}

\begin{table}[htbp]
\centering \small \onehalfspacing
\caption{Isomorphism types of potential non-$G$-cr subgroups}
\begin{tabularx}{.99\linewidth}{c|c|>{\raggedright\arraybackslash}X}
\hline $G$ & $p$ & Subgroup types $H \notin \textup{Lie}(p)$ \\ \hline
$F_4$ & $2$ & $J_{2}$, $L_{2}(13)$ \\
& $3$ & $\Alt_{5}$, $L_{2}(7)$, $L_{2}(8)$ \\
& $7$ & $L_{2}(8)$ \\ \hline

$E_6$ & $2$ & $M_{22}$, $L_{2}(q) \ (q = 11,13,17)$, $U_{4}(3)$ \\
& $3$ & $\Alt_{5}$, $\Alt_{7}$, $M_{11}$, $L_{2}(q) \ (q = 7,8,11,17)$ \\
& $5$ & $\Alt_{6}$, $M_{11}$, $L_{2}(11)$ \\
& $7$ & $\Alt_{7}$, $L_{2}(8)$ \\ \hline

$E_7$ & $2$ & $\Alt_{n} \ (n = 7,8,10,12)$, $M_{11}$, $M_{12}$, $J_{2}$, $L_{2}(q) \ (q = 11,13,17,19,27)$, $L_{3}(3)$ \\
& $3$ & $\Alt_{n} \ (n = 5,7,8,9)$, $M_{11}$, $M_{12}$, $L_{2}(q) \ (q = 7,8,11,13,17,25)$, $^{3}D_4(2)$ \\
& $5$ & $\Alt_{6}$, $\Alt_{7}$, $M_{11}$, $L_{2}(11)$, $^{2}F_4(2)'$ \\
& $7$ & $\Alt_{7}$, $L_{2}(8)$, $L_{2}(13)$, $U_{3}(3)$ \\
& $11$ & $M_{11}$ \\
& $13$ & $L_{3}(3)$ \\ \hline

$E_8$ & $2$ & $\Alt_{n} \ (n = 7,9,10,12,16)$, $M_{11}$, $M_{12}$, $J_{2}$, $L_{2}(q) \ (q = 11,13,17,19,25,27,37)$, $L_{3}(3)$, $PSp_{4}(5)$ \\
& $3$ & $\Alt_{n} \ (n = 5,7,8,9)$, $M_{11}$, $L_{2}(q) \ (q = 7,8,11,13,17,19,25)$, $U_{3}(8)$, $^{3}D_4(2)$ \\
& $5$ & $\Alt_{n} \ (n = 6,7,10)$, $M_{11}$, $L_{2}(11)$, $L_{2}(19)$, $U_{4}(2)$, $^{2}F_{4}(2)'$ \\
& $7$ & $\Alt_{7}$, $L_{2}(8)$, $L_{2}(13)$, $U_{3}(3)$ \\
& $11$ & $M_{11}$ \\
& $13$ & $L_{3}(3)$ \\ \hline
\end{tabularx}
\label{tab:nongcr}
\end{table}

\begin{maincorollary} \label{cor:nongcr}
If $G$ is an exceptional simple algebraic group over an algebraically closed field of characteristic $p = 0$ or $p > 13$, then all non-generic finite simple subgroups of $G$ are $G$-completely reducible.
\end{maincorollary}

Theorem \ref{THM:NONGCR} and Corollary \ref{cor:nongcr} mirror results of Guralnick \cite{MR1717357}*{Theorem A}, which state that if $G = GL_{n}(K)$, where $K$ has characteristic $p > n + 1$, then all finite subgroups of $G$ having no non-trivial normal $p$-subgroup are $G$-completely reducible. Another general result in this direction appears in \cite{MR2178661}*{Remark 3.43(ii)}, which tells us that if $p > 3$ for $G$ of type $G_{2}$, $F_{4}$, $E_{6}$ or $E_{7}$, or if $p > 5$ for $G$ of type $E_{8}$, then a subgroup $S$ of $G$ is $G$-completely reducible if $L(G)$ is a completely reducible $S$-module. All of these results are in the spirit of those of Serre in \cite{MR2167207}*{\S 4, 5}.

At this stage we do not attempt to classify those non-generic subgroups which are not $G$-completely reducible. For most triples $(G,H,p)$ in Table \ref{tab:nongcr} it is straightforward to show the existence of such a subgroup. However, some cases are quite involved; we comment briefly on this in Section \ref{rem:existence} on page \pageref{rem:existence}.

\addtocontents{toc}{\setcounter{tocdepth}{-1}}
\section*{Almost Simple Subgroups}
\addtocontents{toc}{\setcounter{tocdepth}{1}}

The `normaliser stability' of the connected subgroups in Theorem \ref{THM:MAIN} also allows us to extend our results to almost simple finite subgroups. Let $S_0$ be almost simple, with simple socle $S$, and suppose that the isomorphism type of $S$ appears in Table \ref{tab:main}. If $\bar{S}$ is the $N_{\textup{Aut}(G)}(S)$-stable connected subgroup given by Theorem \ref{THM:MAIN}, then $S_0 < N_G(\bar{S})$, and this latter group is positive-dimensional and proper since $G$ is simple.

\begin{maincorollary} \label{cor:almostsimple}
If $G$, $H$, $p$ appear in Table \ref{tab:main}, then $G$ has no Lie primitive finite almost simple subgroup whose socle is isomorphic to $H$. In particular, no group $\Alt_{n}$ or $\Sym_{n}$ $(n \ge 10)$ occurs as a Lie primitive subgroup of an exceptional simple algebraic group in any characteristic.
\end{maincorollary}

In the case that the field of definition has characteristic $0$, we obtain the following. For $\Alt_n$, this is proved independently in forthcoming work of D. Frey \cite{Fre1}.
\begin{maincorollary} \label{cor:almostsimple_zero}
No exceptional simple algebraic group over an algebraically closed field of characteristic $0$ has a Lie primitive subgroup isomorphic to $\Alt_n$ or $\Sym_n$ for $n \ge 8$.
\end{maincorollary}

\addtocontents{toc}{\setcounter{tocdepth}{-1}}
\section*{Application: Finite Groups of Lie Type}
\addtocontents{toc}{\setcounter{tocdepth}{1}}

Ever since the classification of the finite simple groups, a question of primary importance in finite group theory has been to understand their subgroups, and in particular their maximal subgroups. Recall that a group of Lie type arises as a subquotient of the group of fixed points of a simple algebraic group under a Frobenius morphism (more details are given in Section \ref{sec:notation} of Chapter \ref{chap:background}). If $\sigma$ is a Frobenius morphism of an adjoint group $G$ over a field of characteristic $p > 0$, then the group $O^{p'}(G_\sigma)$ is usually simple. Simple groups of Lie type make up all the non-abelian finite simple groups besides the alternating and sporadic groups.

For the alternating groups and classical groups of Lie type, the O'Nan-Scott Theorem \cite{MR1409812}*{Theorem 4.1A} and Aschbacher's Theorem \cite{MR746539} reduce the study of maximal subgroups to understanding primitive permutation actions and modular representations of almost simple groups. For the exceptional groups of Lie type, understanding maximal subgroups again reduces naturally to embeddings of almost simple groups by an analogue of Borovik's Theorem \ref{thm:boro} (see \cite{MR1066315} for the full statement). As with the corresponding algebraic groups, we expect more explicit results here than in the classical case.

In Section \ref{sec:pfcorfin} we use Theorem \ref{THM:MAIN} to prove the following; note that again we have not considered the groups of type $G_2$ since a complete description of their maximal subgroups is already available in \cite{MR618376}, \cite{MR955589} and \cite{MR898346}.

\begin{maintheorem} \label{THM:FIN}
Let $G$ be an adjoint exceptional simple algebraic group over an algebraically closed field of positive characteristic $p$. Let $\sigma$ be a Frobenius morphism of $G$ such that $L = O^{p'}(G_\sigma)$ is simple, and let $L \le L_1 \le \textup{Aut}(L)$. 

If $(G,H,p)$ appears in Table \ref{tab:main}, then $L_{1}$ has no maximal subgroup with socle isomorphic to $H$.
\end{maintheorem}

Finally we remark that while the proof of Theorems \ref{THM:MAIN} and \ref{THM:FIN} are independent of the classification of finite simple groups, Theorem \ref{thm:subtypes} does rely on this, and hence so do Theorems \ref{THM:FEASIBLES} and \ref{THM:NONGCR}, as well as Corollaries \ref{cor:prim} and \ref{cor:nongcr}.

\section*{Layout}

Our analysis requires an amount of background on the representation theory of reductive algebraic groups and finite groups, and the subgroup structure of reductive groups. Chapter \ref{chap:background} gives a concise overview of the necessary theory. In Chapter \ref{chap:disproving} we describe the calculations used to deduce Theorem \ref{THM:FEASIBLES}, illustrating with a typical example. The feasible characters themselves appear in Chapter \ref{chap:thetables}.

For most triples $(G,H,p)$ in Table \ref{tab:main}, the proof that a subgroup $S \cong H$ of $G$ lies in a proper connected subgroup is simply a matter of inspecting the feasible characters and applying an elementary bound on the number of trivial composition factors (Proposition \ref{prop:substab}). In Chapter \ref{chap:thetables} we have marked with `\possprim' those feasible characters which do \emph{not} satisfy the bound; thus every Lie primitive subgroup of $G$ gives rise to one of these feasible characters. The groups with no characters marked `\possprim' are collected in Table \ref{tab:substabprop}. Table \ref{tab:main} consists of these together with some additional cases, which are considered in Proposition \ref{prop:algorithm}.

Chapter \ref{chap:stab} begins by assuming that $(G,H,p)$ appears in Table \ref{tab:main} and that if $S \cong H$ is a subgroup of $G$, then $S$ lies in a proper connected subgroup of $G$. A variety of techniques are then applied in order to show that $S$ in fact lies in a proper, connected, $N_{\textup{Aut}(G)}(S)$-stable subgroup. Theorem \ref{THM:FIN} follows from Theorem \ref{THM:MAIN} through a short argument which is independent of the remainder of the paper, which we give in Section \ref{sec:pfcorfin}. Finally, Chapter \ref{chap:gcr} considers $G$-complete reducibility of the finite subgroups occurring, and Appendix \ref{chap:auxiliary} contains auxiliary data on representations of finite simple groups used in Chapters \ref{chap:disproving}, \ref{chap:stab} and \ref{chap:gcr}.

\section*{Notation}

Throughout, unless stated otherwise, all algebraic groups are affine and defined over $K$, a fixed algebraically closed field of characteristic $p \ge 0$. Subgroups are assumed to be closed, and modules are assumed to be rational and of finite dimension. A homomorphism between algebraic groups is assumed to be a morphism of varieties, unless stated otherwise.

For a group $X$, the derived subgroup is denoted by $X'$. The image of $X$ under a homomorphism $\phi$ is denoted $X^{\phi}$, and if $\phi$ is an endomorphism of $X$, then $X_{\phi}$ denotes the fixed points of $\phi$ in $X$. If $X$ is finite then $O^{p'}(X)$ denotes the smallest normal subgroup of $X$ with index coprime to $p$. If $X_{1}$, $X_{2}$, $\ldots$, $X_{r}$ are simple algebraic groups or tori, then $X_{1} X_{2} \ldots X_{n}$ denotes a commuting product with pairwise finite intersections.

If $V$ is a $KX$-module, then $V^{*}$ denotes the dual module $\textup{Hom}(V,K)$. If $Y$ is a subgroup of $X$, then $V \downarrow Y$ denotes the restriction of $V$ to $Y$. If $M_{1}$, $M_{2}$, $\ldots$ $M_{r}$ are $KX$-modules and $n_{1}$, $\ldots$, $n_{r}$ are integers then we write
\[ M_{1}^{n_{1}}/M_{2}^{n_{2}}/\ldots/M_{r}^{n_{r}} \]
to denote a module having the same composition factors as the direct sum $M_{1}^{n_{1}} \oplus M_{2}^{n_{2}} \oplus \ldots \oplus M_{r}^{n_{r}}$. Finally, the notation
\[ M_{1}|M_{2}|\ldots|M_{r} \]
denotes a $KX$-module $V$ with a series $V = V_{1} > V_{2} > \ldots > V_{r+1} = 0$ of submodules such that $\textup{soc}(V/V_{i+1}) = V_{i}/V_{i+1} \cong M_{i}$ for $1 \le i \le r$.

%% file: AJL-background.tex

\chapter{Background} \label{chap:background}

Good general references for the theory of algebraic groups covered in this chapter are \cite{MR1102012}, \cite{MR0396773}, \cite{MR2850737} and \cite{MR2015057} and we will also give references to specific results when appropriate.

\section{Affine Algebraic groups} \label{sec:notation}

\subsection{Reductive and semisimple groups; root systems} For an algebraic group $G$, the connected component containing the identity is denoted $G^{\circ}$. The soluble and unipotent radicals of $G$ are respectively denoted $R(G)$ and $R_u(G)$; they are the maximal normal connected soluble (respectively unipotent) subgroups of $G$. Then $G$ is called \emph{reductive} if $R_u(G) = 1$ and \emph{semisimple} if $R(G) = 1$. An algebraic group is \emph{simple} if it has no non-trivial, proper, connected normal subgroups, and a semisimple algebraic group is a commuting product of simple subgroups whose pairwise intersections are finite. \label{term:radicals}

Let $G$ be reductive. We let $T$ denote a fixed maximal torus of $G$, and $B$ a Borel subgroup (maximal connected soluble subgroup) containing $T$. The \emph{rank} of $G$ is $r = \textup{dim}(T)$. The \emph{character group} $X(T) \stackrel{\textup{def}}{=} \textup{Hom}(T,K^{*})$ is a free abelian group of rank $r$, written additively. Any $KG$-module $V$ restricts to $T$ as a direct sum of \emph{weight spaces}
\[ V_{\lambda} \stackrel{\textup{def}}{=} \{ v \in V \ : \ t.v = \lambda(t)v \textup{ for all }t \in T \} \]
where $\lambda$ runs over elements of $X(T)$. Those $\lambda$ with $V_{\lambda} \neq \{0\}$ are the \emph{weights} of $V$. The Lie algebra of $G$, denoted $L(G)$\label{term:lg}, is a $KG$-module under the adjoint action, and the nonzero weights of $L(G)$ are the \emph{roots} of $G$; these form an abstract root system in the sense of \cite{MR0396773}*{Appendix}. We let $\Phi$ denote the set of roots. The choice of $B$ gives a base of simple roots $\Pi = \{\alpha_1,\ldots,\alpha_r\}$ and a partition of $\Phi$ into positive roots $\Phi^{+}$ and negative roots $\Phi^{-}$, which are those roots expressible as a positive (respectively negative) integer sum of simple roots.

Abstract root systems are classified by their \emph{Dynkin diagram}; nodes correspond to simple roots, with bonds indicating relative lengths and angle of incidence between non-orthogonal roots. A semisimple group is simple if and only if its Dynkin diagram is connected, and the well-known classification of possible connected diagrams, as well as our chosen numbering of the nodes, are as follows:
\begin{center}
\includegraphics{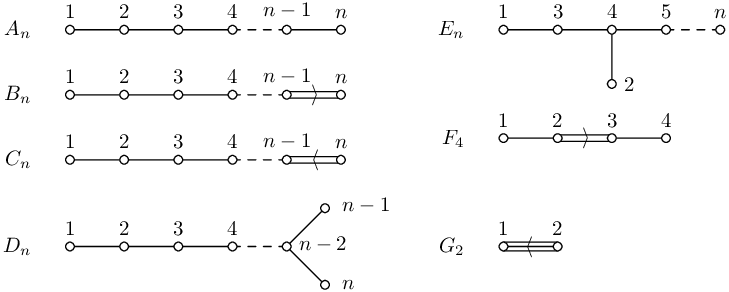}
\end{center}
where type $E_n$ occurs for $n \in \{6,7,8\}$.

For a simple algebraic group $G$, the corresponding label ($A_n$, etc.) is called the (Lie) \emph{type} of $G$. An \emph{isogeny} is a surjective homomorphism of algebraic groups with finite kernel, and for each Lie type there exists a \emph{simply connected} group $G_{\textup{sc}}$ and an \emph{adjoint} group $G_{\textup{ad}}$, each unique up to isomorphism, such that there exist isogenies $G_{\textup{sc}} \to G$ and $G \to G_{\textup{ad}}$ for each simple group $G$ over $K$ of the same type.

\subsection{Weights and rational modules} The group $\mathbb{Z}\Phi$ is naturally a lattice in the Euclidean space $E \stackrel{\textup{def}}{=} X(T) \otimes_{\mathbb{Z}} \mathbb{R}$. In this space, an \emph{abstract weight} is a point $\lambda$ satisfying
\[ \left<\lambda,\alpha\right> \stackrel{\textup{def}}{=} 2\frac{(\lambda,\alpha)}{(\alpha,\alpha)} \in \mathbb{Z} \]
for each simple root $\alpha$. The abstract weights form a lattice $\Lambda$ (the \emph{weight lattice}); the quotient $\Lambda/\mathbb{Z}\Phi$ is finite, and $\mathbb{Z}\Phi \le X(T) \le \Lambda$. The possible subgroups between $\mathbb{Z}\Phi$ and $\Lambda$ correspond to isogeny types of simple group; $G$ is simply connected if and only if $X(T) = \Lambda$, and adjoint if and only if $X(T) = \mathbb{Z}\Phi$.

A weight $\lambda \in \Lambda$ is called \emph{dominant} if $\left<\lambda,\alpha\right> \ge 0$ for all $\alpha \in \Pi$, and $\Lambda$ is free abelian on the set $\{\lambda_1,\ldots,\lambda_r\}$ of \emph{fundamental dominant weights}, which are defined by $\left<\lambda_i,\alpha_j\right> = \delta_{ij}$ for all $i$ and $j$. If $\lambda$, $\mu$ are weights, we say $\mu \le \lambda$ if $\lambda - \mu$ is a non-negative integer combination of positive roots.

If $V$ is a (rational, finite-dimensional) $KG$-module, the Lie-Kolchin theorem \cite{MR0396773}*{\S 17.6} implies that $B$ stabilises a 1-space of $V$; an element spanning a $B$-stable 1-space is called a \emph{maximal vector} of $V$. If $V$ is generated as a $KG$-module by a maximal vector of weight $\lambda$, we say that $V$ is a \emph{highest weight module} for $G$, of \emph{high weight} $\lambda$. Then $\lambda$ is dominant, and all other weights of $V$ are strictly less than $\lambda$ under the above ordering.

A $B$-stable 1-space of a $KG$-module generates a $KG$-submodule under the action of $KG$, and thus an irreducible $KG$-module has a unique $B$-stable 1-space and is a highest weight module. Conversely, for any dominant $\lambda \in X(T)$ there exists an irreducible module of highest weight $\lambda$, which we denote $V_G(\lambda)$\label{term:vglambda}. Then $V_G(\lambda) \cong V_G(\mu)$ implies $\lambda \le \mu$ and $\mu \le \lambda$, hence $\mu = \lambda$.

The \emph{Weyl group} $W = N_G(T)/T$ acts naturally on $X(T)$, inducing an action on $E = X(T) \otimes_{\mathbb{Z}} \mathbb{R}$ which preserves the set of roots, as well as the root and weight lattices. This identifies $W$ with the group of isometries of $E$ generated by \emph{simple reflections} $s_i \, : \, x \mapsto x - \left<x,\alpha_i\right>\alpha_i$ for each $i$. The \emph{length} of an element $w \in W$ is the smallest $k$ such that $w$ is expressible as a product of $k$ simple reflections. The \emph{longest element} of the Weyl group is the unique element of maximal length; it is also characterised as the unique element of $W$ sending every positive root to a negative one. Hence if $w_{\circ}$ is the longest element, then the linear transformation $-w_{\circ}$ of $E$ permutes the dominant weights. Since the weights of a module are precisely the negatives of those of its dual, the simple modules $V_G(-w_{\circ}\lambda)$ and $V_G(\lambda)^{*}$ have identical weight spaces, and thus $V_G(-w_{\circ}\lambda) \cong V_G(\lambda)^{*}$. When no confusion can arise, we abbreviate $V_G(\lambda)$ to simply $\lambda$.

For each dominant $\lambda \in X(T)$, we also have a \emph{Weyl module}, denoted $W_G(\lambda)$\label{term:wglambda}, which is universal in the sense that if $V$ is a rational module of finite dimension with $V/\textup{Rad}(V) = V_{G}(\lambda)$, then $V$ is a quotient of $W_G(\lambda)$. The Weyl module can be defined in terms of the 1-dimensional $B$-module $k_{\lambda}$, with $R_u(B)$ acting trivially and $T \cong B/R_u(B)$ acting with weight $\lambda$. Then $H^{0}(G/B,\lambda) = \textup{Ind}_{B}^{G}(k_{\lambda})$ has socle $\cong V_G(\lambda)$, and we set $W_G(\lambda) \stackrel{\textup{def}}{=} H^{0}(G/B,-w_{\circ}\lambda)^{*}$.

Recall that a rational module is called \emph{tilting} if it has a filtration by Weyl modules and a filtration by duals of Weyl modules. It is well-known (for instance, see \cite{MR1200163}) that the class of tilting modules is closed under taking direct sums and summands. An indecomposable tilting module has a unique highest weight, occurring with multiplicity $1$, and moreover for each dominant $\lambda \in X(T)$ there exists a unique such indecomposable tilting module, denoted $T_{G}(\lambda)$ or $T(\lambda)$ when $G$ is understood.

If the ambient field $K$ has characteristic zero then $T_{G}(\lambda)$, $W_G(\lambda)$ and $V_G(\lambda)$ coincide for each $\lambda$; in general, the module structure of $W_G(\lambda)$ can be quite complicated, although much can still be said without too much effort. For instance, the weights and multiplicities of $W_G(\lambda)$ are independent of the field $K$, and can be calculated using the `Weyl character formula' (see \cite{MR1153249}*{\S 24}).

We will make frequent use of the module $V_{\textup{min}}$, a non-trivial Weyl module for $G$ of least dimension, when $G$ is simply connected of exceptional type, with high weight and dimension as follows:
\begin{center}
\begin{tabular}{c|ccccc}
$G$ & $F_4$ & $E_6$ & $E_7$ & $E_8$ \\ \hline
$V_{\textup{min}}$ & $W_G(\lambda_4)$ & $W_G(\lambda_1)$ & $W_G(\lambda_7)$ & $W_G(\lambda_8) = L(G)$\\
$\textup{dim }V_{\textup{min}}$ & $26$ & $27$ & $56$ & $248$
\end{tabular}
\end{center}
We will give further information on the structure of Weyl modules when discussing representations of reductive groups in Section \ref{sec:algreps}.

\subsection{Automorphisms, endomorphisms, groups of Lie type} With $X$ a semisimple algebraic group over an $K$, let $\phi \, : \, X \to X$ be an abstract group automorphism of $X$. If $\phi$ is a morphism of varieties, then by \cite{MR0230728}*{10.3} exactly one of the following holds:
\begin{itemize}
\item The map $\phi$ is an automorphism of algebraic groups, that is, $\phi^{-1}$ is also a morphism, or
\item the fixed-point subgroup $X_\phi$ is finite. In this case, $\phi$ is called a \emph{Frobenius morphism}.
\end{itemize}

We reserve the term `automorphism' for maps of the first type. As well as inner automorphisms, a simple group of type $A_n$ $(n \ge 2)$, $D_n$ $(n \ge 3)$ or $E_6$ additionally has \emph{graph automorphisms}, induced from symmetries of the Dynkin diagram.

Turning to Frobenius morphisms, when $K$ has positive characteristic $p$, for each power $q$ of $p$, the $q$-power field automorphism there exists a field automorphism sending $x \in K$ to $x^q$. This induces a Frobenius morphism $F_q$ of $GL_n(K)$ for each $n > 0$. A surjective endomorphism $\sigma \, : \, X \to X$ is called a \emph{field morphism} or \emph{standard Frobenius morphism} if there exists some $n > 0$, some power $q$ of $p$ and and injective morphism $i \, : \, G \to GL_n(K)$ such that $i(\sigma(x)) = F_q(i(x))$ for each $x \in X$. If $G$ is simple of type $B_2$ or $F_4$ with $p = 2$, or type $G_2$ with $p = 3$, then $G$ additionally has \emph{exceptional graph morphisms}, which are Frobenius morphisms corresponding to symmetries of the Dynkin diagram when the root lengths are ignored. For every Frobenius morphism $\sigma$, there is an integer $m$ such that $\sigma^{m}$ is a standard Frobenius morphism \cite{MR794307}*{p.\ 31}.

If $\sigma$ is a Frobenius morphism of $G$, the fixed-point set $G_{\sigma}$ is a \emph{finite group of Lie type}, and if $G$ is simple of adjoint type, the group $O^{p'}(G_\sigma)$ generated by elements of order a power of $p$ is usually simple, and is then called a \emph{simple group of Lie type}. For example, if $G = {\rm PGL}_n(K)$ is adjoint of type $A_{n-1}$ and $\sigma = \sigma_q$ is a $q$-power field morphism, then $G_{\sigma} = {\rm PGL}_n(q)$ and $O^{p'}(G_\sigma) = L_n(q)$. When $\sigma$ involves an exceptional graph morphism of $G$, the corresponding finite simple group is called a \emph{Suzuki-Ree group}. Further Frobenius morphisms and groups of Lie type arise through composition with an automorphism of $G$; for instance, if $G = SL_n(K)$ (type $A_{n-1}$), then composition of a $q$-power Frobenius automorphism with a graph automorphism of $G$ gives rise to the finite group $SU_n(q)$, hence to the simple unitary group $U_n(q)$. We let $\textup{Lie}(p)$ denote the collection of finite simple groups of Lie type in characteristic $p$, adopting the convention $\textup{Lie}(0) = \varnothing$.

\section{Finite Subgroups}

\label{sec:generic}
\label{sec:nongeneric}

When considering finite subgroups of an algebraic group, Lie primitivity is a natural maximality condition. If $X$ is an algebraic group of positive dimension, and if $S$ is a finite subgroup of $X$, by dimension considerations there exists a positive-dimensional subgroup $Y$ containing $S$ such that $S$ is Lie primitive in $Y$, and then $S$ normalises $Y^{\circ}$. In this way, understanding embeddings of finite groups into $X$ reduces to understanding Lie primitive subgroups of normalisers of connected subgroups; we can then appeal to the broad array of available results concerning connected subgroups.

A natural idea when studying finite subgroups is to make use of the \emph{socle} (product of minimal normal subgroups). A minimal normal subgroup of a finite group is characteristically simple, hence is a product of isomorphic finite simple groups. We then have a dichotomy between the case that the finite subgroup is \emph{local}, that is, normalises a non-trivial elementary abelian subgroup, and the case that the socle is a product of non-abelian simple groups. Local maximal subgroups of exceptional algebraic groups and finite groups of Lie type are well understood; see \cite{MR1132853} for more details.

This dichotomy is typified by the theorem of Borovik below, which reduces the study of Lie primitive finite subgroups of $G$ to \emph{almost simple} groups, that is, those finite groups $H$ satisfying $H_0 \le H \le \textup{Aut}(H_0)$ for some non-abelian simple group $H_0$. In this case, the socle is $\textup{soc}(H) = H_{0}$.

Recall from \cite{MR0379748} that a \emph{Jordan subgroup} of a simple algebraic group $G$ in characteristic $p$ is an elementary abelian $p_0$-subgroup $E$ (where $p_0$ is a prime not equal to $p$) such that:
\begin{itemize}
\item $N_G(E)$ is finite,
\item $E$ is a minimal normal subgroup of $N_G(E)$, and
\item $N_G(E) \ge N_G(A)$ for any abelian $p_0$-subgroup $A \unlhd N_G(E)$ containing $E$.
\end{itemize}

\begin{theorem}[{\cite{MR1066315}}] \label{thm:boro}
Let $G$ be an adjoint simple algebraic group over an algebraically closed field $K$ of characteristic $p \ge 0$, and let $S$ be a Lie primitive finite subgroup of $G$. Then one of the following holds:
\begin{itemize}
\item[\textup{(i)}] $S \le N_G(E)$ where $E$ is a Jordan subgroup of $G$.
\item[\textup{(ii)}] $G = E_8(K)$, $p \neq 2,3,5$ and $\textup{soc}(S) \cong \Alt_5 \times \Alt_6$.
\item[\textup{(iii)}] $S$ is almost simple.
\end{itemize}
\end{theorem}
Cases (i) and (ii) here are well understood. For each adjoint simple algebraic group $G$, the Jordan subgroups have been classified up to $\textup{Aut}(G)$-conjugacy by Alekseevskii \cite{MR0379748} when $p = 0$, and by Borovik \cite{MR1065645} when $p > 0$. Case (ii), sometimes known as `the Borovik subgroup', is unique up to conjugacy and is described by Borovik in \cite{MR1066315}. This theorem thus focuses our attention on the almost simple groups in (iii).

If $G$ is a simple algebraic group in positive characteristic $p$ and $X$ is a connected simple subgroup of $G$, then a Frobenius endomorphism of $X$ gives rise to a finite subgroup as in Section \ref{sec:notation}, which will be a (possibly trivial) central extension of an almost simple group. Due to the abundance of such subgroups, we say that an almost simple subgroup $S$ of $G$ is \emph{generic} if $\textup{soc}(S)$ is isomorphic to a member of $\textup{Lie}(p)$ and \emph{non-generic} otherwise.

Generic almost simple subgroups of exceptional simple algebraic groups are mostly well understood. In \cite{MR1458329}*{Theorem 1}, Liebeck and Seitz prove that if $S = S(q)$ is a generic simple subgroup of $G$, then $S$ arises in the manner above if $q$ is larger than an explicit bound depending on the Lie type of $S$ and the root system of $G$ ($q > 9$ is usually sufficient, and $q > 2624$ always suffices). Subject to this bound on $q$, they construct a closed, connected subgroup $\bar{S}$ containing $S$, such that $S$ and $\bar{S}$ fix precisely the same submodules on $L(G)$, and if $S < G_{\sigma}$ for a Frobenius morphism $\sigma$, then $\bar{S}$ is also $\sigma$-stable and $N_G(S)$-stable.

In the same article \cite{MR1458329}*{Theorem 6}, this latter statement is used to deduce information on the subgroup structure of the finite groups of Lie type; roughly, an almost simple subgroup of $G_{\sigma}$ with generic socle, is either of the same Lie type as $G$, or arises via the fixed points $X_{\sigma}$ of a maximal closed, connected, reductive $\sigma$-stable subgroup $X$ of $G$. We will apply a similar argument in Section \ref{sec:pfcorfin} to subgroups with non-generic socle.

As mentioned in the introduction, it is now known precisely which non-generic simple groups admit an embedding to each exceptional simple algebraic group, and it remains to classify these embeddings. Theorems \ref{THM:MAIN} and \ref{THM:FIN} are an analogue for non-generic subgroups of the aforementioned results of Liebeck and Seitz. If $S$ is a non-generic finite simple subgroup of the exceptional simple algebraic group $G$, we wish to construct a `sufficiently small' connected subgroup $\bar{S}$ of positive dimension containing $S$. In the generic case, `sufficiently small' was taken to mean that each $S$-submodule of $L(G)$ is $\bar{S}$-invariant. This captures the idea that the representation theory of $S$ is in some sense `close' to that of $\bar{S}$ when considering the action on $L(G)$. In the non-generic setting, we do not expect such a strong result to hold in general. On the other hand, as we shall see in Chapter \ref{chap:stab}, it is often possible to find a subgroup $\bar{S}$ which stabilises `sufficiently many' $S$-submodules, for instance, all those of a particular dimension. This is sufficient to prove Theorem \ref{THM:MAIN}.

\section{Subgroups of Positive Dimension} \label{sec:posdim}

Once a given finite simple group is shown not to occur as a Lie primitive subgroup of $G$, proving the remainder of Theorem \ref{THM:MAIN} and subsequent results requires knowledge of the possible intermediate subgroups of positive dimension which can occur.

\subsection{Maximal connected subgroups} In characteristic zero, connected subgroups of an (affine) algebraic group $G$ are in 1-1 correspondence with Lie subalgebras of $L(G)$ (see \cite{MR0396773}*{\S 13}). Thus Dynkin's classification of the maximal Lie subalgebras of simple Lie algebras over $\mathbb{C}$ \cites{MR0047629,MR0049903}, gives also a classification of maximal connected subgroups of a simple algebraic group. This has been extended into positive characteristic by Seitz \cite{MR1048074} and by Liebeck and Seitz \cites{MR1066572,MR2044850}. The following result enumerates the maximal connected subgroups of the exceptional simple algebraic groups, in arbitrary characteristic. Thus a finite subgroup of $G$ which lies in a proper, connected subgroup, lies in a conjugate of one of these.

\begin{theorem}[\cite{MR2044850}*{Corollary 2}] \label{thm:connected}
Let $G$ be an adjoint exceptional simple algebraic group over an algebraically closed field of characteristic $p \ge 0$, and let $X$ be maximal among connected closed subgroups of $G$. Then either $X$ is parabolic or semisimple of maximal rank, or $X$ appears in the table below, and is given up to ${\rm Aut}(G)$-conjugacy.
\end{theorem}

\begin{table}[htbp]
\centering
\begin{tabular}{|l|l|l|}
\hline $G$ & $X$ simple & $X$ not simple \\ \hline \hline
$G_2$ & $A_1$ $(p \ge 7)$ & \\ \hline
$F_4$ & $A_1$ $(p \ge 13)$, $G_2$ $(p = 7)$ & $A_1 G_2$ $(p \neq 2)$ \\ \hline
$E_6$ & $A_2$ $(p \neq 2,3)$, $G_2$ $(p \neq 7)$, & $A_2G_2$ \\
& $C_4$ $(p \neq 2)$, $F_4$ & \\ \hline
$E_7$ & $A_1$ (2 classes, $p \ge 17, 19$ resp.), & $A_1 A_1$ $(p \neq 2,3)$, $A_1 G_2$ $(p \neq 2)$, \\
& $A_2$ $(p \ge 5)$ & $A_1F_4$, $G_2 C_3$ \\ \hline
$E_8$ & $A_1$ (3 classes, $p \ge 23,29,31$ resp.), & $A_1A_2$ $(p \neq 2,3)$, \\
& $B_2$ $(p \ge 5)$ & $G_2 F_4$ \\ \hline
\end{tabular}
\end{table}

Since semisimple subgroups of maximal rank are examples of \emph{subsystem subgroups}, i.e.\ semisimple subgroups normalised by a maximal torus, they can be enumerated by an algorithm of Borel and de Siebenthal \cite{MR0032659}. Similarly, up to conjugacy a parabolic subgroup corresponds to a choice of nodes in the Dynkin diagram; then the parabolic has a non-trivial unipotent radical, whose structure is discussed in Section \ref{sec:radfilt}, and the Dynkin diagram generated by the chosen nodes determines the (reductive) Levi complement.

\subsection{Parabolic subgroups and complete reducibility} \label{sec:parabs}

Let $G$ be a reductive algebraic group, with maximal torus $T$ and set of roots $\Phi$. Then for each $\alpha \in \Phi$, there is a homomorphism $x_{\alpha} \, : \, (K,+) \to G$ of algebraic groups, whose image $U_{\alpha}$ is $T$-stable; this $U_{\alpha}$ is the \emph{root subgroup} corresponding to $\alpha$. We have $tx_{\alpha}(c)t^{-1} = x_{\alpha}(\alpha(t)c)$ for all $t \in T$ and all $c \in K$. 

Let $\Pi$ be a base of simple roots, with corresponding positive roots $\Phi^{+}$. Then a subset $I \subseteq \Pi$ generates a root subsystem $\Phi_I \subseteq \Phi$ in a natural way, and we have a corresponding \emph{standard parabolic subgroup}
\[ P_I \stackrel{\textup{def}} = \left<T, U_{\alpha},U_{\pm \beta} \ : \ \alpha \in \Phi^{+},\ \beta \in I \right>. \]
The unipotent radical of $P_I$ is then
\[ Q_I \stackrel{\textup{def}}{=} \left<U_{\alpha} \ : \ \alpha \in \Phi^{+},\ \alpha \notin \Phi_I \right>\]
and in $P_I$, the unipotent radical has a reductive closed complement, the \emph{Levi factor}
\[ L_I \stackrel{\textup{def}}{=} \left<T, U_{\pm \beta} \ : \ \beta \in I \right> \]
whose root system is precisely $\Phi_I$. A \emph{parabolic subgroup} is any $G$-conjugate of a standard parabolic subgroup, and each is conjugate to precisely one standard parabolic subgroup. A \emph{Levi subgroup} is any conjugate of some $L_I$; since each Levi subgroup $L$ contains a maximal torus of $G$, the derived subgroup $L'$ is a subsystem subgroup of $G$.

Recall that a subgroup $X$ of a reductive group $G$ is \emph{$G$-completely reducible} ($G$-cr) if, whenever $X$ lies in a parabolic subgroup $P = QL$ of $G$, $X$ then lies in a conjugate of the Levi factor $L$. Similarly, $X$ is \emph{$G$-irreducible} ($G$-irr) if there is no parabolic subgroup of $G$ containing $X$, and \emph{$G$-reducible} otherwise. A finite subgroup which is Lie primitive in $G$ is necessarily $G$-irreducible.

There now exists much literature on the concept of $G$-complete reducibility, particularly with regard to classifying connected $G$-completely reducible subgroups. Such a subgroup is necessarily reductive, by the following theorem of Borel and Tits:

\begin{theorem}[\cite{MR2850737}*{Theorem 17.10}] \label{thm:boreltits}
Let $U$ be a unipotent subgroup, not necessarily closed or connected, of a reductive algebraic group $G$. Then there exists a parabolic subgroup $P$ of $G$ with $U \le R_u(P)$ and $N_G(U) \le P$.
\end{theorem}
In particular, if $X$ is a connected subgroup with non-trivial unipotent radical, then $X$ lies in some parabolic $P$ with $X \cap R_u(P) = R_u(X) \neq 1$, so $X$ is not contained in a Levi factor of $P$.

In \cite{MR1329942}, Liebeck and Seitz produce a constant $N(X,G) \le 7$, for each exceptional simple algebraic group $G$ and each connected simple subgroup type $X$, such that if $p > N(X,G)$ then all subgroups of $G$ of type $X$ are $G$-completely reducible, and this has been further refined by D.\ Stewart \cites{MR3282997}. Theorem \ref{THM:NONGCR} can be viewed as an analogue of this result for finite simple subgroups of $G$.

Knowledge of reductive subgroups, and $G$-cr subgroups in particular, is of use in proving Theorem \ref{THM:MAIN} since the precise action of these subgroups on $L(G)$ and $V_{\textup{min}}$ has either already been determined or is straightforward to determine. For example \cite{MR1329942}*{Tables 8.1-8.7} gives the composition factors of all connected simple subgroups of $G$ of rank $\ge 2$ on $L(G)$ or $V_{\textup{min}}$, provided $p$ is greater than the largest entry of the above table. The action is usually completely reducible, and if not, it is usually straightforward to determine the precise action from the composition factors, using techniques outlined shortly.

\subsection{Irreducibility in classical groups} If $G$ is a simple algebraic group of classical type $A_n$, $B_n$, $C_n$ or $D_n$, then $G$ is closely related to $SL(V)$, $SO(V)$ or $Sp(V)$ for some vector space $V$ with an appropriate bilinear or quadratic form. Then $V$ is isomorphic to the Weyl module $W_G(\lambda_1)$, and parabolic subgroups of $G$ have a straightforward characterisation in terms of the irreducible module $V_G(\lambda_1)$ \cite{MR1329942}*{pp.\ 32-33}:

\begin{lemma} \label{lem:classicalparabs}
Let $X$ be a $G$-irreducible subgroup of a simple algebraic group $G$ of classical type, and let $V = V_G(\lambda_1)$. Then one of the following holds:
\begin{itemize}
\item $G = A_n$ and $X$ is irreducible on $V$,
\item $G = B_n$, $C_n$ or $D_n$ and $V \downarrow S = V_1 \perp \dots \perp V_k$ with the $V_i$ all nondegenerate, irreducible and inequivalent as $X$-modules,
\item $G = D_n$, $p = 2$ and $X$ fixes a non-singular vector $v \in V$, such that $X$ is then $C_G(v)$-irreducible, where $C_G(v)$ is simple of type $B_{n-1}$.
\end{itemize}
\end{lemma}
Thus determining when a finite simple group $S$ admits a $G$-irreducible embedding into some classical $G$ is straightforward if we have sufficient information about the representation theory of $KS$. For us, the necessary information is available in the literature, and we give a summary in Section \ref{sec:reps}.

\subsection{Internal modules of parabolic subgroups} \label{sec:radfilt}
We recall some information from \cite{MR1047327}. With $G$ a semisimple algebraic group, let $P$ be a parabolic subgroup of $G$, and without loss of generality assume that $P = P_I$ is a standard parabolic subgroup for some $I \subseteq \Pi$, with unipotent radical $Q = Q_I$ and Levi factor $L = L_I$ as before.

An arbitrary root $\beta \in \Phi$ can be written as $\beta = \beta_I + \beta_{I'}$ where $\beta_I = \sum_{\alpha_i \in I}c_i\alpha_i$ and $\beta_{I'}=\sum_{\alpha_j \in \Pi - I}d_j\alpha_j$. We then define
\begin{align*}
\textup{height}(\beta) &= \sum c_i + \sum d_j, \\
\textup{level}(\beta) &= \sum d_j, \\
\textup{shape}(\beta) &= \beta_{I'},
\end{align*}
 and for each $i \ge 1$ we define subgroup $Q(i) = \left<U_{\beta} \ : \ \textup{level}(\beta) \ge i \right>$.

Call $(G,\textup{char }K)$ \emph{special} if it is one of $(B_n,2)$, $(C_n,2)$, $(F_4,2)$, $(G_2,2)$ or $(G_2,3)$.

\begin{lemma}[\cite{MR1047327}]
The subgroups $Q(i)$ are each normal in $P$. There is a natural $KL$-module structure on the quotient groups $Q(i)/Q(i+1)$, with decomposition $Q(i)/Q(i+1) = \prod V_S$, the product being over all shapes $S$ of level $i$. Each $V_S$ is an indecomposable $KL$-module of highest weight $\beta$ where $\beta$ is the unique root of maximal height and shape $S$. If $(G,p)$ is not special, then $V_S$ is irreducible.
\end{lemma}

Given $G$ and the subset $I$ of simple roots corresponding to $P$, it is a matter of straightforward combinatorics to calculate the modules occurring in this filtration, as laid out in \cite{MR1047327}. A quick summary for the exceptional groups is provided by the following lemma, which is \cite{MR1329942}*{Lemma 3.1}.

\begin{lemma} \label{lem:radfilts}
If $P = QL$ is a parabolic subgroup of an exceptional algebraic group $G$ and $L_0$ is a simple factor of $L$, then the possible high weights $\lambda$ of non-trivial $L_0$-composition factors occurring in the module filtration of $Q$ are as follows:
\begin{itemize}
\item $L_0 = A_n$: $\lambda = \lambda_j$ or $\lambda_{n+1-j}$ $(j = 1,2,3)$;
\item $L_0 = A_1$ or $A_2$ ($G = F_4$ only ): $\lambda = 2\lambda_1$ or $2\lambda_2$;
\item $L_0 = A_1$ $(G = G_2 \textup{ only})$: $\lambda = 3\lambda_1$;
\item $L_0 = B_n$, $C_n$ $(G = F_4$, $n = 2$ or $3)$: $\lambda = \lambda_1$, $\lambda_2$ or $\lambda_3$;
\item $L_0 = D_n$: $\lambda_1$, $\lambda_{n-1}$ or $\lambda_n$;
\item $L_0 = E_6$: $\lambda = \lambda_1$ or $\lambda_6$;
\item $L_0 = E_7$: $\lambda = \lambda_7$.
\end{itemize}
\end{lemma}

In Section \ref{sec:complements} we will use this information and some basic cohomology theory to parametrise complements to $Q$ in the subgroup $QX$, for $X$ a finite simple subgroup of $G$ contained in $P$.

\section{Representation Theory of Semisimple Groups}

\subsection{Structure of Weyl Modules} \label{sec:algreps}

The following proposition summarises information on the module structure of Weyl modules which we will use throughout subsequent chapters. This information is well-known (for example, see \cite{MR1901354}), and can quickly be verified with computer calculations, for example using the \emph{Weyl Modules} GAP package of S. Doty \cite{Dot1}.

\begin{proposition} \label{prop:weyls}
Let $G$ be a simple algebraic group in characteristic $p \ge 0$, and let $\lambda$ be a dominant weight for $G$.
\begin{itemize}
\item If $(G,\lambda)$ appear in Table \ref{tab:irredweyls}, then $W_G(\lambda) = V_G(\lambda)$ is irreducible in all characteristics, with the given dimension.
\item If $(G,\lambda)$ appear in Table \ref{tab:redweyls}, then $W_G(\lambda)$ is reducible in each characteristic given there, with the stated factors, and is irreducible in all other characteristics.
\end{itemize}
\end{proposition}

\begin{table}[H]\small
\centering
\onehalfspacing
\caption{Irreducible Weyl modules}
\begin{tabular}{c|c|c}
$G$ & $\lambda$ & Dimension \\ \hline
$A_n$ & $\lambda_i$ $(i = 1,\ldots,n)$ & $\binom{n+1}{i}$ \\
$B_n$ & $\lambda_n$ & $2^{n}$ \\
$C_n$ & $\lambda_1$ & $2n$ \\
$D_n$ & $\lambda_1$, $\lambda_{n-1}$, $\lambda_n$ & $2n$, $2^{n-1}$, $2^{n-1}$  \\
$E_6$ & $\lambda_1$, $\lambda_6$ & $27$ \\
$E_7$ & $\lambda_7$ & $56$\\
$E_8$ & $\lambda_8$ & $248$
\end{tabular}
\label{tab:irredweyls}
\end{table}

Note that for the classical types, the module $W_G(\lambda_1)$ is the natural module. For type $A_n$ we have $W_G(\lambda_i) = V_G(\lambda_i) = \bigwedge^{i} V_G(\lambda_1)$ for $i = 1,\ldots,n$.

The composition factors of the adjoint module $L(G)$ are the same as those of $W_G(\lambda)$ as in the following table (for a proof, see \cite{MR1458329}*{Proposition 1.10}):
\begin{center}
\begin{tabular}{c|*{9}{c}}
Type of $G$ & $A_n$ & $B_n$ & $C_n$ & $D_n$ & $E_6$ & $E_7$ & $E_8$ & $F_4$ & $G_2$ \\ \hline
$\lambda$ & $\lambda_1 + \lambda_n$ & $\lambda_2$ & $2\lambda_1$ & $\lambda_2$ & $\lambda_2$ & $\lambda_1$ & $\lambda_8$ & $\lambda_1$ & $\lambda_2$
\end{tabular}
\end{center}

\begin{table}[p]
\centering
\onehalfspacing
\caption{Composition Factors of Reducible Weyl modules}
\begin{tabular}{c|c|c|l|l}
$G$ & $\lambda$ & $p$ & High Weights & Factor dimensions \\ \hline
$A_n$   & $\lambda_1 + \lambda_n$ & $p \mid n+1$ & $0$, $\lambda_1 + \lambda_n$ & $1$, $n^{2} + 2n - 1$ \\
        & $2\lambda_1$ $(n > 1)$ & 2 & $\lambda_2$, $2\lambda_1$ & $\binom{n+1}{2}$, $n+1$\\
$B_n$   & $\lambda_1$ & 2 & $0$, $\lambda_1$ & $1$, $2n$ \\
        & $\lambda_2$ $(n > 2)$ & 2 & $0^{(n,2)}$, $\lambda_1$, $\lambda_2$ & $1$, $2n$, \\
        & & & & $2n^{2} - n - (2,n)$ \\
$D_n$   & $\lambda_2$ & 2 & $0^{(n,2)}$, $\lambda_2$ & $1$, $\binom{2n}{2}-(2,n)$ \\
$A_3$   & $2\lambda_2$ & 2 & $2\lambda_2$, $\lambda_1+\lambda_3$ & $6$, $14$ \\
		& & 3 & $0$, $2\lambda_2$ & $1$, $19$ \\
$B_3$	& $2\lambda_1$ & 2 & $0$, $\lambda_1$, $\lambda_2$, $2\lambda_1$ & $1$, $6$, $14$, $6$ \\
		& & 7 & $0$, $2\lambda_1$ & $1$, $26$ \\
		& $2\lambda_3$ & 2 & $0$, $\lambda_1^{2}$, $\lambda_2$, $2\lambda_3$ & $1$, $6$, $14$, $8$ \\
$B_4$   & $\lambda_3$ & 2 & $0^{2}$, $\lambda_1$, $\lambda_2$, $\lambda_3$ & $1$, $8$, $26$, $48$ \\
        & $\lambda_1 + \lambda_4$ & 3 & $\lambda_4$, $\lambda_1+\lambda_4$ & $16$, $112$ \\
        & $2\lambda_1$ & 2 & $0^{2}$, $\lambda_1$, $\lambda_2$, $2\lambda_1$ & $1$, $8$, $26$, $8$ \\
        & & 3 & $0$, $2\lambda_1$ & $1$, $43$ \\
$C_3$   & $\lambda_3$ & 2 & $\lambda_1$, $\lambda_3$ & $6$, $8$ \\
$C_4$   & $2\lambda_1$ & 2 & $0^{2}$, $2\lambda_1$, $\lambda_2$ & $1$, $8$, $26$ \\
		& $\lambda_2$ & 2 & $0$, $\lambda_2$ & $1$, $26$ \\
        & $\lambda_4$ & 2 & $\lambda_2$, $\lambda_4$ & $26$, $16$ \\
        & & 3 & $0$, $\lambda_4$ & $1$, $41$ \\
$D_4$   & $2\lambda_i$ $(i = 1,3,4)$ & 2 & $0$, $\lambda_2$, $2\lambda_i$ & $1$, $26$, $8$ \\
		& $\lambda_i + \lambda_j$, where & 2 & $\lambda_i+\lambda_j$, $\lambda_k$ & $48$, $8$ \\
		& $\{i,j,k\} = \{1,3,4\}$ & & \\
$E_6$   & $\lambda_2$ & 3 & $0$, $\lambda_2$ & $1$, $77$ \\
$E_7$   & $\lambda_1$ & 2 & $0$, $\lambda_1$ & $1$, $132$ \\
$F_4$   & $\lambda_1$ & 2 & $\lambda_1$, $\lambda_4$ & $26$, $26$ \\
        & $\lambda_4$ & 3 & $0$, $\lambda_4$ & $1$, $25$ \\
$G_2$   & $\lambda_1$ & 2 & $0$, $\lambda_1$ & $1$, $6$ \\
		& $\lambda_2$ & 3 & $\lambda_1$, $\lambda_2$ & $7$, $7$ \\
        & $\lambda_1 + \lambda_2$ & 3 & $0$, $\lambda_1$, $\lambda_2$, $\lambda_1+\lambda_2$ & $1$, $7$, $7$, $49$ \\
        & & 7 & $\lambda_1+\lambda_2$, $2\lambda_1$ & $38$, $26$ \\
        & $2\lambda_1$ & 2 & $0$, $\lambda_1$, $\lambda_2$, $2\lambda_1$ & $1$, $6$, $14$, $6$ \\
        & & 7 & $0$, $2\lambda_1$ & $1$, $26$ \\
        & $3\lambda_{1}$ & 2 & $0^{3}$, $\lambda_1^{2}$, $\lambda_2$, $2\lambda_{1}^{2}$, $3\lambda_1$ & $1$, $6$, $14$, $6$, $36$ \\
        & & 3 & $\lambda_1$, $\lambda_2^{2}$, $\lambda_1 + \lambda_2$, $3\lambda_1$ & $7$, $7$, $49$, $7$
\end{tabular}
\label{tab:redweyls}
\end{table}

Note also that the simply connected groups form a chain
\[ F_4(K) < E_6(K) < E_7(K) < E_8(K) \]
and up to composition factors, we have the following well-known restrictions (see, for example \cite{MR1329942}*{Tables 8.1-8.4}):
\begin{align*}
L(E_8) \downarrow E_7 = L(&E_7)/V_{E_7}(\lambda_7)^{2}/0^{3},\\
L(E_7) \downarrow E_6 = L(E_6)/V_{E_6}(\lambda_1)/V_{E_6}(\lambda_6)/0, \quad&\quad V_{E_7}(\lambda_7) \downarrow E_6 = V_{E_6}(\lambda_1)/V_{E_6}(\lambda_6)/0^{2},\\
L(E_6) \downarrow F_4 = L(F_4)/W_{F_4}(\lambda_4), \quad&\quad V_{E_6}(\lambda_1) \downarrow F_4 = W_{F_4}(\lambda_4)/0.
\end{align*}

In addition, the longest element $w_{\circ}$ of the Weyl group of $G$ induces the scalar transformation $-1$ on $X(T) \otimes_{\mathbb{Z}} \mathbb{R}$ if $G$ is simple of type $B_n$, $C_n$, $D_{2n}$, $G_2$, $F_4$, $E_7$ or $E_8$. In this case, every irreducible $KG$-module is self-dual, as $V_G(\lambda)^{*} \cong V_G(-w_{\circ}\lambda) = V_G(\lambda)$. If $G$ is instead simple of type $A_n$, $D_{2n+1}$ or $E_6$, then $w_{\circ}$ induces $-\tau$, where $\tau$ corresponds to a non-trivial symmetry of the Dynkin diagram. This gives rise to isomorphisms:
\begin{itemize}
\item $V_{G}(\lambda_i) \cong V_{G}(\lambda_{n+1-i})^{*}$ for $G$ of type $A_n$ and any $1 \le i \le n$,
\item $V_G(\lambda_1) \cong V_G(\lambda_6)^{*}$ for $G$ of type $E_6$,
\item $V_G(\lambda_{2n}) \cong V_G(\lambda_{2n+1})^{*}$ for $G$ of type $D_{2n+1}$,
\end{itemize}
For $G$ of any type, if $\lambda$ is fixed by each graph automorphism of $G$, the module $V_G(\lambda)$ is self-dual.

Finally, recall that a \emph{spin module} for $X$ is the irreducible module $V_X(\lambda_n)$, of dimension $2^{n}$, for $X$ of type $B_n$, or one of the modules $V_X(\lambda_{n-1})$ or $V_X(\lambda_{n})$, of dimension $2^{n-1}$, for $X$ of type $D_n$.

\begin{lemma}[\cite{MR1329942}*{Lemma 2.7}] \label{lem:spins}
Let $X = B_{n}$ $(n \ge 3)$ or $D_{n+1}$ ($n \ge 4)$, and let $Y$ be either a Levi subgroup of type $B_{r}$ $(r \ge 1)$ or $D_{r}$ $(r \ge 3)$ of $X$, or a subgroup $B_{n}$ of $X = D_{n+1}$. If $V$ is a spin module for $X$, then all composition factors of $V \downarrow Y$ are spin modules for $Y$.
\end{lemma}

\subsection{Tilting Modules}

Tilting modules enjoy several properties that will be of use here. For instance, $L(G)$ and $V_{\textup{min}}$ are often tilting; this only fails for an exceptional simple algebraic group $G$ in the following cases:
\begin{itemize}
\item $G = E_{7}$, $p = 2$, $L(G) = V_{G}(\lambda_1)/0$, $T(\lambda_1) = 0|V_{G}(\lambda_1)|0$;
\item $G = E_{6}$, $p = 3$, $L(G) = V_{G}(\lambda_2)/0$, $T(\lambda_2) = 0|V_{G}(\lambda_2)|0$;
\item $G = F_{4}$, $p = 2$, $L(G) = V_{G}(\lambda_1)/V_{G}(\lambda_4)$, $T(\lambda_1) = V_{G}(\lambda_4)|V_{G}(\lambda_1)|V_{G}(\lambda_4)$;
\item $G = F_{4}$, $p = 3$, $V_{\textup{min}} = V_{G}(\lambda_4)/0$, $T(\lambda_4) = 0|V_{G}(\lambda_4)|0$.
\end{itemize}

In every scenario encountered here, it is a trivial task to determine the structure of a tilting module $T_{G}(\lambda)$ from the (known) structure of the Weyl modules $W_{G}(\mu)$ for $\mu \le \lambda$. Moreover, the following will be useful in determining the action of subgroups of $G$ on the various $G$-modules. This is parts (i)--(iii) of \cite{MR1200163}*{Proposition 1.2} (cf.\ also \cite{MR1072820}).

\begin{lemma}
Let $X$ be a reductive algebraic group.
\begin{itemize}
\item[\textup{(i)}]If $M$ and $N$ are tilting modules for $X$, then so is $M \otimes N$;
\item[\textup{(ii)}] If $V$ is a tilting module for $X$ and $L$ is a Levi subgroup of $X$, then $V \downarrow L$ is a tilting module for $L$;
\item[\textup{(iii)}] $T_{X}(\lambda)^{*} = T_{X}(-w_{\circ}\lambda)$;
\end{itemize}
\label{lem:tilting}
\end{lemma}

\subsection{Extensions of Rational Modules}
Suppose that a finite subgroup $S$ of the simple algebraic group $G$ is known to lie in some proper, connected subgroup of $G$. Then $S$ lies in either a reductive subgroup of $G$, and then in the semisimple derived subgroup, or in a parabolic subgroup of $G$. In the latter case, $S$ can be studied using its image under the projection to a Levi factor. Thus, as we shall see in Chapter \ref{chap:stab}, useful information on embeddings $S \to G$ can be procured by comparing the possible actions of $S$ on $L(G)$ and $V_{\textup{min}}$ with those of the various semisimple subgroups of $G$ admitting an embedding of $S$.

We therefore require some results on the representation theory of semisimple groups. A comprehensive reference for the material of this section is \cite{MR2015057}.

Let $X$ be a semisimple algebraic group over the algebraically closed field $K$. Recall that for rational $KX$-modules $V$ and $W$, we denote by $\textup{Ext}_X^{1}(V,W)$ the set of all \emph{rational} extensions of $V$ by $W$ up to equivalence.

\begin{lemma}[\cite{MR2015057}*{p.183, Proposition}] \label{lem:exthom}
If $\lambda$, $\mu$ are dominant weights for $X$ with $\lambda$ not less than $\mu$, then
\[ \textup{Ext}_{X}^{1}(V_X(\lambda),V_X(\mu)) \cong \textup{Hom}_X(\textup{rad}(W_X(\lambda)),V_X(\mu)). \]
\end{lemma}

This lemma is particularly useful in light of the information on the structure of Weyl modules above. Another result of use in this direction is the following.
\begin{lemma}[\cite{MR1048074}*{Lemma 1.6}] \label{lem:2step}
In a short exact sequence of $KX$-modules
\[ 0 \to V_X(\lambda) \to M \to V_X(\mu) \to 0, \]
one of the following occurs:
\begin{itemize}
\item[\textup{(i)}] The sequence splits, so $M \cong V_X(\lambda) \oplus V_X(\mu)$,
\item[\textup{(ii)}] $\lambda < \mu$ and $M$ is a quotient of the Weyl module $W_X(\mu)$,
\item[\textup{(iii)}] $\mu < \lambda$ and $M^{*}$ is a quotient of the Weyl module $W_X(-w_{\circ}\lambda)$.
\end{itemize}
\end{lemma}

\begin{corollary} \label{cor:2step}
If $V$ is a rational $KX$-module with high weights $\{\mu_1, \ldots, \mu_t\}$, and if $W_X(\mu_i)$ has no high weight $\mu_j$ for all $i \neq j$, then $V$ is completely reducible.

In particular, if each $W_X(\mu_i)$ is irreducible, then $V$ is completely reducible.
\end{corollary}

Each of these will be useful in narrowing down the possibilities for intermediate subgroups $S < X < G$ when $S$ is known not to be Lie primitive in $G$.

%% file: AJL-disproving.tex

\chapter{Calculating and Utilising Feasible Characters} \label{chap:disproving}

With $G$ an adjoint exceptional simple algebraic group over $K$, our strategy for studying embeddings of finite groups into $G$ is to determine the possible composition factors of restrictions of $L(G)$ and $V_{\textup{min}}$. Here we describe the tools which allow us to achieve this.

\section{Feasible Characters} \label{sec:brauer}

\subsection{Definitions} Let $H$ be a finite group and let the field $K$ have characteristic $p \ge 0$. Recall that to a $KH$-module $V$ we can assign a \emph{Brauer character}, a map from $H \to \mathbb{C}$ which encodes information about $V$ in much the same way as an ordinary character encodes information about a module in characteristic zero. Let $n$ be the exponent of $H$ if $p = 0$, or the $p'$-part of the exponent if $p > 0$. Then the eigenvalues of elements of $H$ on $V$ are $n$-th roots of unity in $K$. Fix an isomorphism $\phi$ between the group of $n$-th roots of unity in $K$ and in $\mathbb{C}$. The Brauer character of $V$ is then defined by mapping $h \in H$ to $\sum \phi(\zeta)$, the sum over eigenvalues $\zeta$ of $h$ on $V$.

When $K$ has characteristic zero, a Brauer character is simply an ordinary character (a Galois conjugate of the usual character). In general, a Brauer character has a unique expression as a sum of irreducible Brauer characters (those arising from irreducible modules), though the decomposition of the Brauer character into irreducibles only determines the corresponding module up to composition factors, not up to isomorphism.

We now give some crucial definitions, the first of which is found in Frey \cite{MR1617620}.
\begin{definition}
A \emph{fusion pattern} from $H$ to $G$ is a map $f$ from the $p'$-conjugacy classes of $H$ to the conjugacy classes of $G$, which preserves element orders and is compatible with power maps, i.e.\ for each $i \in \mathbb{Z}$, $f$ maps the class $(x^{i})^{H}$ to the class of $i$-th powers of elements in $f(x^H)$.
\end{definition}

\begin{definition} \label{def:feasible}
A \emph{feasible decomposition} of $H$ on a finite-dimensional $G$-module $V$ is a $KH$-module $V_0$ such that for some fusion pattern $f$, the Brauer character of any $x \in H$ on $V_0$ is equal to the trace of elements in $f(x^{H})$ on $V$. The Brauer character of $V_0$ is then called a \emph{feasible character}. We say that a collection of feasible decompositions or feasible characters of $H$ on various $G$-modules is \emph{compatible} if they all correspond to the same fusion pattern.
\end{definition}

Any subgroup $S$ of $G$ gives rise to a fusion pattern (map the $S$-conjugacy class to the $G$-conjugacy class of its elements), and the restriction of any set of finite-dimensional $G$-modules gives a compatible collection of feasible decompositions. Not all feasible characters (or fusion patterns) are necessarily realised by an embedding; however, determining all feasible characters places a strong restriction on possible embeddings. Note that our definition of a feasible character is more restrictive than that given by Cohen and Wales \cite{MR1416728}*{p.\ 113}, since the definition there does not take power maps into account.

\subsection{Determining feasible characters}

Given $G$, a finite group $H$ and a collection $\left\{V_i\right\}$ of rational $KG$-modules, determining the compatible collections of feasible characters of $H$ on $\left\{V_i\right\}$ is a three-step process.
\begin{itemize}
\item Firstly, we need the Brauer character values of all irreducible $KH$-modules of dimension at most $\textup{Max}(\textup{dim}(V_i)$). The $\{V_i\}$ used here each have dimension at most $248 = \textup{dim } L(E_8)$. The necessary information on Brauer characters either exists in the literature or can be calculated directly. We give more details on this in Section \ref{sec:reps}.
\item Secondly, for each $m$ coprime to $p$ such that $H$ has elements of order $m$, we will need to know the eigenvalues of elements of $G$ of order $m$ on each module $V_i$. These can be determined using the weight theory of $G$, which we outline in Section \ref{sec:sselts}. We then take the corresponding sum in $\mathbb{C}$ under the bijection determining the Brauer characters.
\item Finally, determining feasible characters becomes a matter of enumerating non-negative integer solutions to simultaneous equations, one equation for each class in $H$. Each solution gives the irreducible character multiplicities in a feasible character. This step is entirely routine, and we give an illustrative example in Section \ref{sec:example}.
\end{itemize}

In Chapter \ref{chap:thetables}, proceeding as above we give the compatible feasible characters of each finite simple group $H \notin \textup{Lie}(p)$ on the $KG$-modules $L(G)$ and $V_{\textup{min}}$ defined in Section \ref{sec:algreps}.

\subsection{Irreducible Modules for Finite Quasisimple Groups} \label{sec:reps}

Let $S$ be a finite simple subgroup of a semisimple algebraic group $G$, let $\tilde{G}$ be the simply connected cover of $G$ and let $\tilde{S}$ be a minimal preimage of $S$ under the natural projection $\tilde{G} \twoheadrightarrow G$, so that $\tilde{S}$ is a \emph{cover} (perfect central extension) of $S$. If $\tilde{S} \cong S$ then we have an induced action of $S$ on each $K\tilde{G}$ module and no isogeny issues are encountered. However, if $\tilde{S}$ has non-trivial centre, then in order to make use of any faithful $\tilde{G}$-modules (in particular, the module $V_{\textup{min}}$ for $G = E_6$ and $E_7$), we will need to consider the action of $\tilde{S}$ rather than $S$. Note that in this case, we have $Z(\tilde{S}) \le \text{Ker}(\tilde{G} \twoheadrightarrow G) \le Z(\tilde{G})$.

A perfect finite group has, up to isomorphism, a unique covering group of maximal order, the \emph{universal cover}. All covers are then quotients of this, and the centre of the universal cover is called the \emph{Schur multiplier} of $S$. The Schur multipliers of the finite simple groups are all known.

In order to carry out the calculations described above, for each finite simple group $H$ appearing in Table \ref{tab:subtypes}, we will need to know information on the irreducible $K\tilde{H}$-modules of dimension $\le 248$, where $\tilde{H}$ is a cover of $H$ with $|Z(H)| \le 3$. For most such groups encountered here, the necessary Brauer characters appear in the Atlas \cite{MR827219} or the Modular Atlas \cite{MR1367961}. In addition, Hiss and Malle \cite{MR1942256} give a list of each dimension at most $250$ in which each finite quasisimple group has an absolutely irreducible module (not in the ambient characteristic of the group is of Lie type).

The groups from Table \ref{tab:subtypes} whose Brauer characters are not given explicitly in the Modular Atlas are $\Alt_{n}$ ($13 \le n \le 17$), $L_2(q)$ ($q = 37$, $41$, $49$, $61$), $Fi_{22}$, $Ru$, $Th$, $L_4(5)$ and $\Omega_7(3)$. The modular character tables of $L_{2}(q)$ for $q$ and $p$ coprime are well-known, cf. \cite{MR0480710}. For the alternating groups, the following reciprocity result is of use in determining all irreducible Brauer characters of a given degree.
\begin{lemma}[\cite{MR860771}*{p.58}]
Let $Y \le X$ be finite groups, $V$ a finite-dimensional $KY$-module and $U$ a finite-dimensional $KX$-module. Then
\[
\textup{Hom}_{KX}(U, V\uparrow X) \cong \textup{Hom}_{KY}(U \downarrow Y, V).
\]
In particular, if $U$ is simple, of dimension $N$, and if $V$ is any $KY$-module quotient of $U \downarrow Y$, then $U$ is a submodule of $V \uparrow X$.
\end{lemma}

In particular, the irreducible Brauer characters of degree $\le 248$ for $\Alt_{n}$ can all be found by inducing and decomposing the irreducible Brauer characters of degree $\le 248$ for $\Alt_{n-1}$. For completeness, in Appendix \ref{chap:auxiliary} we include the Brauer character values used for $\Alt_{13}$ to $\Alt_{17}$ in characteristic $2$.

According to Corollary 4 of \cite{MR1717629} (cf. Table \ref{tab:onlyprim}), the remaining groups from the above list occur only as Lie primitive subgroups of an exceptional algebraic group. Then according to \cite{MR1942256}, the faithful irreducible modules of dimension at most 248 for these groups and their double or triple covers are as follows:

\begin{center}
\begin{tabular}{r|l|c}
\multicolumn{1}{c|}{$H$} & \multicolumn{1}{c|}{Module Dimensions} & Adjoint $G$ containing $H$ \\ \hline
$Fi_{22}$ $(p = 2)$ & $78$ & $E_6$ \\
$3.Fi_{22}$ $(p = 2)$ & $27$ (not self-dual) \\
$Ru$ $(p = 5)$ & $133$ & $E_7$ \\
$2.Ru$ $(p = 5)$ & $28$ (not self-dual) \\
$Th$ $(p = 3)$ & $248$ & $E_8$ \\
$L_4(5)$ $(p = 2)$ & $154$, $248$ & $E_8$ \\
$\Omega_7(3)$ $(p = 2)$ & $78$, $90$, $104$ & $E_6$ \\
$3.\Omega_7(3)$ $(p = 2)$ & $27$ (not self-dual)
\end{tabular}
\end{center}

We see immediately that each embedding of $Fi_{22}$, $\Omega_7(3)$ and $Th$ must be irreducible on both $L(G)$ and $V_{\textup{min}}$. An embedding of $Ru$ into $E_7$ must arise from an embedding of $2\cdot Ru$ into the simply connected group, acting irreducibly on the Lie algebra, and with composition factors $28/28^{*}$ on the 56-dimensional self-dual module $V_{\textup{min}}$. Finally, since $L_4(5)$ has no non-trivial modules of dimension less than $154$, it must be irreducible on $L(E_8)$ in any embedding, otherwise it would fix a nonzero vector on $L(G)$, and would not be Lie primitive in $G$ by Lemma \ref{lem:properembed} below. This determines all the feasible characters for these groups on $L(G)$ and $V_{\textup{min}}$.

\subsection{Frobenius-Schur indicators} The (Frobenius-Schur) indicator of a $KH$-module $V$ encodes whether the image of $H$ in $GL(V)$ lies in an orthogonal or symplectic group. An irreducible $KH$-module $V$ supports a nondegenerate $H$-invariant bilinear form if and only if $V \cong V^{*}$, and the form is then symmetric or alternating. If $K$ has characteristic $\neq 2$, then $H$ preserves a symmetric form if and only if it preserves a quadratic form. In characteristic $2$, any $H$-invariant bilinear form on $V$ is symmetric, and a nondegenerate $H$-invariant quadratic form on $V$ gives rise to a nondegenerate $H$-invariant bilinear form, but not conversely.

The indicator of $V$ is then defined as
\[ \textup{ind}(V) = \left\{
\begin{array}{rll}
0 & (\textup{or }\circ) & \textup{ if } V \ncong V^{*},\\
1 & (\textup{or } +) & \textup{ if } H \textup{ preserves a nondegenerate quadratic form on } V,\\
-1 & (\textup{or } -) & \textup{ otherwise.}
\end{array}\right.
\]
Thus $V$ has indicator $-$ if and only if $V$ is self-dual, supports an $H$-invariant alternating bilinear form, and does not support an $H$-invariant quadratic form.

\subsection{First cohomology groups} The last piece of information we will use is the \emph{cohomology group} $H^{1}(H,V)$ for various $KH$-modules $V$. Recall that this is the quotient of the additive group $Z^{1}(H,V)$ of 1-cocycles (maps $\phi\, : \, H \to V$ satisfying $\phi(xy) = \phi(x) + x.\phi(y)$) by the subgroup $B^{1}(H,V)$ of 1-coboundaries (cocycles $\phi$ such that $\phi(x) = x.v - v$ for some $v \in V$).

The cohomology group is a $K$-vector space which parametrises conjugacy classes of complements to $H$ in the semidirect product $VH$, via
\[ \phi \mapsto \{\phi(h)h \ : \ h \in H\}, \]
and also parametrises short exact sequences of $KH$-modules
\[ 0 \to V \to E \to K \to 0 \]
under equivalence, where an equivalence is an isomorphism $E \to E'$ of $KH$-modules inducing the identity map on $V$ and $K$.

Knowledge of cohomology groups will be useful in determining the existence of fixed points in group actions (see Proposition \ref{prop:substab} for a good example). Computational routines exist for determining the dimension of $H^{1}(H,V)$ (for instance, Magma implements such routines). The information we have made use of here is summarised as follows:

\begin{lemma} \label{lem:reps}
Let $K$ be an algebraically closed field of characteristic $p$ such that the simple group $H \notin \textup{Lie}(p)$ embeds into an adjoint exceptional simple algebraic group over $K$ (i.e.\ $H$ appears in Table \ref{tab:subtypes}). Then Tables \ref{tab:altreps} to \ref{tab:ccreps} give every non-trivial irreducible $KH$-module of dimension $\le 248$.
\end{lemma}
We also give there the Frobenius-Schur indicator $\textup{ind}(V)$ of each module, as well as $\textup{dim}(H^{1}(H,V))$ when this has been used.

Note that the computational packages used to calculate these cohomology group dimensions do not perform calculations over an algebraically closed field, but rather over a \emph{finite} field. This is sufficient for us if we ensure that the field being used is a \emph{splitting field} for the group $H$ (see for example \cite{MR1110581}*{\S 1}). Such fields always exist; for instance if $|H| = p^{e}r$ with $r$ coprime to $p$, then a field containing a primitive $r$-th root of unity suffices.

\subsection{Semisimple Elements of Exceptional Groups} \label{sec:sselts}

In order to calculate feasible characters, we need to know the eigenvalues of semisimple elements of small order in $F_4(K)$, $E_6(K)$, $E_7(K)$ and $E_8(K)$. A semisimple element necessarily has order coprime to the characteristic $p$ of $K$.

Let $H$ be a group such that $(G,H)$ appears in Table \ref{tab:subtypes}. Let $n$ be the $p'$ part of the exponent of $H$. Let $\omega_n$ be an $n$-th root of unity in $K$, let $\zeta_n$ be an $n$-th root of unity in $\mathbb{C}$, and let $\phi$ be the isomorphism of cyclic groups sending $\omega_n$ to $\zeta_n$. For each $m$ dividing $n$, let $\omega_m = \omega_n^{(n/m)}$, and similarly for $\zeta_m$.

If $m > 0$ is coprime to $p$, then an element of order $m$ in $G$ is semisimple, hence lies in some maximal torus $T$. From the existence and uniqueness of the Bruhat decomposition for elements in $G$, it follows (see for example \cite{MR794307}*{Section 3.7}) that two elements of order $m$ in $T$ are conjugate in $G$ if and only if they are in the same orbit under the action of the Weyl group $W = N_G(T)/T$.

Now, let $\{\chi_i\}$ be a free basis of the character group $X(T)$. If $t \in T$ has order $m$, then each $\chi_i(t)$ is a power of $\omega_m$, say $\chi_i(t) = \omega_m^{n_i}$, where $0 \le n_{i} < m$ and gcd$(\{n_{i}\}) = 1$. Conversely, for any $r$-tuple of integers $(n_1,\ldots,n_r)$ satisfying these conditions, there exists $t \in T$ with $\chi_i(t) = \omega_{m}^{n_i}$ (see \cite{MR0396773}*{Lemma 16.2C}), which then has order $m$.

Thus, for a fixed basis of $X(T)$, elements of $T$ of order $m$ correspond to $r$-tuples of integers $(n_1,\ldots,n_r)$ as above, where $r = \textup{rank}(G)$, and the action of $W$ on $T$ induces an action on these. This latter action makes no reference to $K$, only to $X(T)$ and the basis. In particular, if $G_1$ is a simple algebraic group over $\mathbb{C}$ with an isomorphism of root systems $\Phi(G) \to \Phi(G_1)$ identifying the character groups and Weyl groups, then classes of elements of order $m$ in $G$ and in $G_1$ are each in 1-1 correspondence with orbits of the Weyl group on these $r$-tuples, hence are in 1-1 correspondence with each other. 

This correspondence respects Brauer character values; suppose $V$ and $W$ are respectively modules for $G$ and $G_1$, whose weight spaces correspond under the isomorphism of root systems. If we express a weight $\lambda$ as a linear combination of basis elements $\lambda = \sum m_i(\lambda)\chi_i$, then an element $g \in G$ which is represented by $(n_1,\ldots,n_r)$ has a corresponding eigenvalue $\prod_{i=1}^{r} \omega_m^{n_i m_i(\lambda)}$, and thus the Brauer character value of $g$ on $V$ is
\[ \sum_{\lambda \textup{ a weight of }V} \left(\prod_{i=1}^{r}\zeta_m^{n_i m_i(\lambda)}\right) \]
which is equal to the (true) character value of the corresponding elements in $G_1$.

There now exists a well-established theory of elements of finite order in simple complex Lie groups, and in \cite{MR752042}, Moody and Patera provide an algorithmic approach to enumerating elements of finite order, as well as determining their eigenvalues on rational modules. Classes of semisimple elements in a simply connected simple group $G$ over $\mathbb{C}$ are in 1-1 correspondence with $(r+1)$-tuples $(s_0,s_1,\ldots,s_r)$ with $\textup{gcd}(s_0,s_1,\ldots,s_r) = 1$, where $r = \textup{rank}(G)$. Under projection to $G/Z(G)$, elements represented by $(s_0,s_1,\ldots,s_r)$ have order $\sum_{i = 1}^{r} n_is_i$, where $\alpha_0 = \sum n_i \alpha_i$ is the highest root. The full order of the element is also determined by the $s_i$ (see \cite{MR752042}*{Section 4}), and the algorithm additionally tells us precisely when elements represented by distinct $(r+1)$-tuples are conjugate in the adjoint group.

By implementing the procedure above in Magma, we have calculated the eigenvalues of all semisimple elements of order at most 37 in the simply connected groups $F_4(\mathbb{C})$, $E_6(\mathbb{C})$, $E_7(\mathbb{C})$ and $E_8(\mathbb{C})$ on the adjoint and minimal modules. There are a total of 2,098,586 conjugacy classes of such elements between these groups, and it is impractical to give representatives of each. However, all but a few tables of Chapter \ref{chap:thetables} can be verified using only elements of order at most 7, and the eigenvalues of such elements are already known. In particular, our calculations have been checked against \cite{MR933426} (elements of $E_8$ and simply connected $E_7$) and \cite{MR1416728} (elements of $F_4$ and simply connected $E_6$).

\section{Deriving Theorem \ref{THM:FEASIBLES}}

\subsection{Example: Feasible characters of $\Alt_{17}$ on $L(E_8)$, $p = 2$} \label{sec:example}

We now have sufficient information to derive the feasible characters for each $(G,H,p)$ in Table \ref{tab:subtypes}, which is the content of Theorem \ref{THM:FEASIBLES}. We illustrate with the case $H \cong \Alt_{17}$ and $G = E_8(K)$, where $p = \textup{char}(K) = 2$, on the 248-dimensional adjoint module $L(G)$. It transpires that there is a unique feasible character, up to a permutation of irreducible $H$-modules corresponding to an outer automorphism of $H$. Note that $G$ does indeed have a subgroup isomorphic to $H$, since $G$ has a simple subgroup of type $D_8$, which is abstractly isomorphic to $SO_{16}(K)$ since $p = 2$, and it is well-known that $\Alt_{17}$ has a 16-dimensional irreducible module which supports a nondegenerate quadratic form.

As given in the Appendix (Table \ref{tab:altreps}), in characteristic 2 there is a unique $KH$-module of each dimension 1, 16 and 118 up to isomorphism, and two irreducible modules of dimension 128, which are interchanged by an outer automorphism of $H$. There are no other irreducible $KH$-modules of dimension $\le 248$. The 16-dimensional module is a quotient of the natural 17-dimensional permutation module, the 118-dimensional module is a section of $\bigwedge^{2} 16$ and the 128-dimensional `spin' modules arise from embeddings $\Alt_{17} \le \textup{SO}_{16}^{+}(2) \le {\rm SL}_{2^7}(2)$.

In this case it suffices to consider \emph{rational} elements of $H$, that is, elements which are $H$-conjugate to all their proper powers of the same order. All Brauer character values of rational elements are integers, and the character values of such elements of $H$ can be calculated by hand. It is well-known that the Brauer character of $h \in H$ on the 16-dimensional deleted permutation module is $|\textup{fix}(h)| - 1$, giving also a formula for the Brauer character of the alternating square. The eigenvalues of $h$ on a spin module can be inferred directly from the eigenvalues of $H$ on the 16-dimensional module; this is done, for example, in \cite{MR1057341}*{pp.\ 195-196}. For the elements of orders $3$, $5$ and $7$, we obtain the following:
\begin{center}\begin{tabular}{c|ccccccccccc}
$\chi$ & $e$ & $3$ & $3^2$ & $3^3$ & $3^4$ & $3^5$ & 5 & $5^2$ & $5^3$ & $7$ & $7^2$\\ \hline
$\chi_{16}$ & $16$ & $13$ & $10$ & $7$ & $4$ & $1$ & $11$ & $6$ & $1$ & $9$ & $2$ \\
$\chi_{118}$ & $118$ & $76$ & $43$ & $19$ & $4$ & $-2$ & $53$ & $13$ & $-2$ & $34$ & $-1$\\
$\chi_{128_a}$, $\chi_{128_b}$ & $128$ & $-64$ & $32$ & $-16$ & $8$ & $-4$ & $-32$ & $8$ & $-2$ & $16$ & $2$
\end{tabular}\end{center}

Finding the possible Brauer characters of $L(G) \downarrow H$ then involves finding non-negative integers $a$, $b$, $c$, $d_1$, $d_2$ such that
\[ L(G) \downarrow H = 1^a/16^b/118^c/128_a^{d_1}/128_b^{d_2}\]
where we are denoting $KH$-modules by their degree. Let $d = d_1+d_2$ be the total number of 128-dimensional factors in the feasible character.

From the calculations described in \ref{sec:sselts}, we find that $G$ has two classes of rational elements of order 5; elements of these have traces $-2$ and $23$ on $L(G)$ (see also \cite{MR933426}*{Table 1}). There is also a class of elements with trace $3$ on $L(G)$, however such elements are not rational (they have non-integer trace on the 3875-dimensional $G$-module). As all elements of $H$ of order 5 are rational, each feasible character of $H$ on $L(G)$ must take a value in $\{-2,23\}$ on each such class.

Thus, evaluating the character of $L(G) \downarrow H$ on the classes $e$, $5$ and $5^{2}$ gives the following equations:
\begin{align}
248 &= a + 16b + 118c + 128d \\
-2 {\rm\ or\ } 23 &= a + 11b + 53c - 32d \\
-2 {\rm\ or\ } 23 &= a + 6b + 13c + 8d
\end{align}
By (1) we have $c \le 2$, $d \le 1$. The third line must equal 23, since coefficients are non-negative. Subtracting (1) from (3) then gives:
\[ 45 = 2b + 21c + 24d \]
and so $c$ must be odd, hence $c = 1$. Assuming $d = 0$ forces $b = 12$, making (1) inconsistent. Therefore $d = 1$, which therefore means $b = 0$, $a = 2$. Thus we have
\[ L(G) \downarrow H = 1^2/118/128_a \quad \textup{or}\quad 1^{2}/118/128_b \]
which are the same up to the action of an outer automorphism of $H$. This also determines a fusion pattern from $H$ to $G$, which is unique up to interchanging classes according to an outer automorphism of $H$. In similar calculations with $G$ not of type $E_8$, it may be necessary to calculate a feasible character on both $L(G)$ and another non-trivial $KG$-module before a fusion pattern is determined. Here, we have
\begin{center}\begin{tabular}{c|cccccccccccc}
$x \in H$ & $e$ & $3$ & $3^2$ & $3^3$ & $3^4$ & $3^5$ & $5$ & $5^2$ & $5^3$ & $7$ & $7^2$\\ \hline
$\chi_{L(G) \downarrow H}(x)$ & $248$ & $14$ & $77$ & $5$ & $14$ & $-4$ & $23$ & $23$ & $-2$ & $52$ & $3$ \\
Class in $G$ & 1A & 3C & 3D & 3B & 3C & 3A & 5G & 5G & 5C & 7N & 7H
\end{tabular}\end{center}
where the conjugacy class labels are taken from \cite{MR933426}*{Table 1}.

Given knowledge of the necessary Brauer characters and semisimple elements, identical calculations to the above are possible for each pair $(G,H)$ in Table \ref{tab:subtypes}, as well as their double and triple covers when $G$ is respectively of type $E_7$ and $E_6$. This is entirely routine, and we have used Magma to facilitate calculations and help avoid errors. The results are the tables of feasible characters of Chapter \ref{chap:thetables}.

\section{Finding Fixed Vectors} \label{sec:notprim}

Let $G$ be an exceptional simple algebraic group and let $S$ be a finite quasisimple subgroup of $G$ with $Z(S) \le Z(G)$. Our approach to studying the embedding of $S$ into $G$ is to use the representation theory of $S$ to find a nonzero fixed vector in the action of $S$ on $L(G)$ or $V_{\textup{min}}$, since then $S < C_G(v)$, a closed subgroup, and
\[ \textup{dim}(C_G(v)) = \textup{dim}(G) - \textup{dim}(G.v). \]
In particular, if $\textup{dim}(V) < \textup{dim}(G)$ then $C_G(v)$ is of positive dimension. If also $v$ is not fixed by $G$, then $C_G(v)$ is proper, and thus $S$ is not Lie primitive in $G$. Note that $G$ can only fix a nonzero vector on $V$ if $(G,V,p) = (E_6,L(G),3)$, $(E_{7},L(G),2)$ or $(F_{4},V_{G}(\lambda_4),3)$.

If $V = L(G)$, then although $\textup{dim}(V) = \textup{dim}(G)$, a similar conclusion nevertheless holds. By the uniqueness of the Jordan decomposition, any endomorphism fixing $v \in L(G)$ must also fix the semisimple and nilpotent parts of $v$. Thus if $S$ fixes a nonzero vector in its action on $L(G)$, then it fixes a nonzero vector which is semisimple or nilpotent. We then appeal to the following result; note that the proof given in \cite{MR1048074} is valid for an arbitrary reductive group $G$.
\begin{lemma}[\cite{MR1048074}*{Lemma 1.3}] \label{lem:properembed}
Let $0 \neq v \in L(G)$.
\begin{itemize}
\item[\textup{(i)}] If $v$ is semisimple then $C_G(v)$ contains a maximal torus of $G$.
\item[\textup{(ii)}] If $v$ is nilpotent, then $R_u(C_G(v)) \neq 1$ and hence $C_G(v)$ is contained in a proper parabolic subgroup of $G$.
\end{itemize}
\end{lemma}
We are therefore interested in conditions on a feasible character which will guarantee the existence of a fixed vector.

\subsection{Group cohomology}
 To begin, recall that if $V$ and $W$ are $KS$-modules then we denote by Ext$_{S}^{1}(V,W)$ the set of equivalence classes of short exact sequences of $S$-modules:
\[ 0 \to W \to E \to V \to 0, \]
and we have isomorphisms (see \cite{MR2015057}*{1, Chapter 4}):
\[ \textup{Ext}_{S}^{1}(V,W) \cong \textup{Ext}_{S}^{1}(K,V^{*} \otimes W) \cong H^{1}(S,V^{*}\otimes W), \]
where $H^{1}(S,V^{*}\otimes W)$ is the first cohomology group.

The following result is based on \cite{MR1367085}*{Lemma 1.2}, and is a highly useful tool for deducing the existence of a fixed vector in the action of $S$ on a module $V$.

\begin{proposition} \label{prop:substab}
Let $S$ be a finite group and $M$ a finite-dimensional $KS$-module, with composition factors $W_1,\ldots,W_r$, of which $m$ are trivial. Set $n = \sum {\rm dim\ } H^{1}(S,W_i)$, and assume $H^{1}(S,K)=\{0\}$.
\begin{itemize}
\item[\textup{(i)}] If $n < m$ then $M$ contains a trivial submodule of dimension at least $m - n$. \label{substabi}
\item[\textup{(ii)}] If $m = n$ and $M$ has no nonzero trivial submodules, then $H^{1}(S,M) = \{0\}$. \label{substabii}
\item[\textup{(iii)}] Suppose that $m = n > 0$, and that for each $i$ we have
\[ H^{1}(S,W_i) = \left\{0\right\} \Longleftrightarrow H^{1}(S,W_i^{*}) = \left\{0\right\}. \]
\end{itemize}
Then $M$ has a nonzero trivial submodule or quotient. \label{substabiii}
\end{proposition}

\proof In each case, we proceed by induction on the number $r$ of composition factors of $M$.

(i) If $r = 1$ or if $M$ has only trivial composition factors, the result is immediate. So let $r > 1$ and assume that $M$ has a non-trivial composition factor. Let $W \subseteq M$ be a submodule which is maximal such that $M/W$ has a non-trivial composition factor. Let $n' = \sum \textup{dim}(H^1(S,W_i))$, the sum being over composition factors of $W$, and let $m'$ be the number of trivial composition factors of $W$. If $n - n' < m - m'$, then by induction $M/W$ would have a trivial submodule, contradicting the choice of $W$. Thus $n - n' \ge m - m'$ and so $n' \le m' + (n - m) < m'$; by induction, $W$ contains a trivial submodule of dimension $m' - n' \ge m - n$.

(ii) Suppose that $n = m$, that $M$ has no trivial submodule and that $H^{1}(S,M)$ is nonzero, so that there exists a non-split extension $0 \rightarrow M \rightarrow N \rightarrow K \rightarrow 0.$ Since $M$ contains no trivial submodule, neither does $N$. This contradicts (i), since $N$ has $m + 1$ trivial composition factors, while the sum $\sum{\rm dim\ }H^1(S,W_i)$ over composition factors $W_i$ of $N$ is equal to $m$.

(iii) Now suppose $n = m > 0$. Assume that $M$ has no trivial submodules. We will show that $M$ has a nonzero trivial quotient. Let $N$ be a maximal submodule of $M$. Since $M$ has no nonzero trivial submodules, neither does $N$. Hence $H^{1}(S,M/N) = \{0\}$, otherwise $N$ would have a nonzero trivial submodule by part (i). By induction on $r$, we deduce that $N$ has a nonzero trivial quotient. Let $Q$ be a maximal submodule of $N$ such that $N/Q$ is trivial. Then $M/Q$ is an extension of a trivial module by the irreducible module $M/N$. By our hypothesis on cohomology groups, we have
\[ \textup{Ext}_S^{1}(M/N,K) \cong \textup{Ext}_S^{1}(K,(M/N)^{*}) \cong H^{1}(S,(M/N)^{*}) = \{0\}. \]
Thus the extension splits and $M$ has a trivial quotient, as required. \qed

\subsection{Some representation theory of finite groups} \label{sec:finreps}

While Proposition \ref{prop:substab} is widely applicable, it is sometimes possible to infer the existence of a fixed vector even when it does not apply. In Section \ref{sec:thealg} we detail a representation-theoretic approach to determining whether a module can exist with a prescribed set of composition factors and no nonzero trivial submodules. To describe this approach properly, we give here a survey of preliminary results from the representation theory of finite groups. A good reference is \cite{MR860771}.

Let $S$ be a finite group and $K$ be an algebraically closed field of characteristic $p$. The group algebra $KS$ admits a $K$-algebra decomposition into indecomposable \emph{block algebras} $B_i$, giving also a decomposition of the identity element:
\begin{align*}
KS &= B_1 \oplus B_2 \oplus \dots \oplus B_n,\\
e &= e_1 + e_2 + \ldots + e_n.
\end{align*}
In turn, this gives a canonical direct-sum decomposition of any $KS$-module $M$:
\[ M = e_1 M + e_2 M + \dots + e_n M. \]
We say that a module $M$ \emph{belongs to the block} $B_i$ if $M = e_i M$ (in which case $e_j M = 0$ for all $j \neq i$). It is immediate that any indecomposable module lies in a unique block, and that if $M$ lies in the block $B_i$, then so do all submodules and quotients of $M$. Hence if we know \emph{a priori} to which block each irreducible $KS$-module belongs, and if we know the composition factors of some $KS$-module, we know that it must split into direct summands accordingly.

Next, recall that the \emph{Jacobson Radical} $\textup{rad}(V)$ of a $KS$-module $V$ is the intersection of all maximal submodules of $V$, or equivalently, is the smallest submodule $J$ of $V$ such that $V/J$ is completely reducible; we call $V/\textup{rad}(V)$ the \emph{head} of $V$. Dually, the \emph{socle} $\textup{soc}(V)$ of $V$ is the sum of all irreducible submodules of $V$, or equivalently, is the unique maximal semisimple submodule of $V$.

A $KS$-module is called \emph{projective} if it is a direct summand of a free module, or equivalently, if any surjection onto it must split. We have the following basic facts (see \cite{MR860771}*{Chapter II}):
\begin{lemma} \label{lem:projcoverdef}
Let $V$ be a finite-dimensional $KS$-module. Then there exists a finite-dimensional projective $KS$-module $P$ such that
\begin{itemize}
\item[\textup{(i)}] $V/\textup{rad}(V) \cong P/\textup{rad}(P)$.
\item[\textup{(ii)}] $P$ is determined up to isomorphism by $V$.
\item[\textup{(iii)}] $V$ is a homomorphic image of $P$.
\item[\textup{(iv)}] $P/\textup{rad}(P) \cong \textup{soc}(P)$.
\item[\textup{(v)}] $\textup{dim}(P)$ is divisible by the order of a Sylow $p$-subgroup of $S$.
\end{itemize}
\end{lemma}
Such a $P$ is then called the \emph{projective cover} of $V$. Projective covers provide a highly useful computational tool for studying the submodule structure of $KS$-modules with known composition factors. Every finite-dimensional projective module is a direct sum of indecomposable projective modules, and these are in 1-1 correspondence with the irreducible $KS$-modules $\{S_i\}$ via $P_i/\textup{rad}(P_i) \cong \textup{soc}(P_i) \cong S_i$. The projective indecomposable modules are the indecomposable direct summands of the free module $KS$, and the projective module $P_i$ occurs precisely $\textup{dim}(S_i)$ times in a direct-sum decomposition of $KS$. We thus obtain the formula
\[ \sum_{i = 1}^{r} \textup{dim}(P_i)\ \textup{dim}(S_i) = \textup{dim}(KS) = |S|. \]

It is well-known that the number of isomorphism types of irreducible $KS$-modules is equal to the number of conjugacy classes of $S$ of elements of order coprime to $p$. Let $m$ be the number of such classes. Then if a $KS$-module $V$ has composition factors $S_1^{r_1}/S_2^{r_2}/.../S_m^{r_m}$ where $r_i \ge 0$, then by the above lemma the projective cover of $V$ has the form
\[ P = P_1^{n_1} + P_2^{n_2} + \ldots + P_m^{n_m} \]
with $n_i \le r_i$ for each $i$.

Additionally, the Brauer characters of the irreducible $KS$-modules can be used to determine the composition factors of the indecomposable projective modules (cf.\ \cite{MR661045}*{Chapter IV}). If the Brauer character of the irreducible $KS$-module $S_i$ is $\chi_i$, we extend each $\chi_i$ to a class function on $S$ by setting $\chi_i(x) = 0$ whenever $x$ has order divisible by $p$. We can then define the inner product $\left<\chi_i,\chi_j\right> = \frac{1}{|S|}\sum_{x \in S}\chi_i(x)\chi_j(x^{-1})$, as in the case of ordinary characters. The $(m \times m)$ matrix $( \left<\chi_i,\chi_j\right>_{i,j})$ is invertible, and the $(i,j)$-entry of its inverse is the multiplicity of $S_i$ as a composition factor of $P_j$, and also to the multiplicity of $S_j$ as a composition factor of $P_i$.

Finally, a result involving the defect group of a $KS$-module tells us that if $p^{e}$ is the order of a Sylow $p$-subgroup of $S$, where $p$ is the characteristic of $K$, and if $p^{e} \mid \textup{dim }S_i$, then $S_i \cong P_i$ is its own projective cover \cite{MR661045}*{Theorem IV.4.5}.

As a quick example, let us determine the structure of the projective indecomposables for $S \cong \Alt_{5}$ with $p = 3$. From \cite{MR1942256}, we know that $S$ has four irreducible modules $S_1$, ..., $S_4$, of dimension $1$, $3$, $3$ and $4$, respectively, each of which is self-dual. Immediately, the 3-dimensional modules are projective since their dimension is divisible by 3 (the order of a Sylow 3-subgroup), and the other two modules are not, since their dimension is not divisible by 3.

If $P_1$ and $P_4$ are the projective covers of the 1- and 4-dimensional $KS$-modules, we have
\begin{align*}
60 &= \textup{dim}(P_1) + 9 + 9 + 4\,\textup{dim}(P_4)
\end{align*}
Since $3 | \textup{dim }P_i$, we have $\textup{dim }P_4 < 10$. It follows that $P_4$ is uniserial (has a unique composition series) with composition series $S_4 | S_1 | S_4$. This in turn implies that $P_1$ has a single 4-dimensional composition factor, and is thus uniserial of shape $S_1 | S_4 | S_1$.

Methods for constructing and manipulating projective indecomposable modules have now been implemented in various computational algebra packages. As with calculating cohomology groups, these implementations are designed to work over finite extensions of the prime field $\mathbb{Q}$ or $\mathbb{F}_p$. This is sufficient for determining the submodule structure of the projective indecomposables, in particular their socle series and radical series, since each projective $KS$-module can be obtained by extending scalars from projective $kS$-module, whenever $k \subset K$ is a splitting field for $S$.

\subsection{Projective covers and fixed vectors} \label{sec:thealg}

As indicated above, the submodule structure of projective indecomposable modules for many of the simple groups in Table \ref{tab:subtypes} can be determined either by hand or using computational techniques. They thus provide a powerful tool for studying the possible structure of a module with known composition factors.

Let $S$ be a finite group, $K$ an algebraically closed field, $\{S_i\}$ the irreducible $KS$-modules and $\{P_i\}$ the corresponding projective indecomposable modules. Let $V$ be a $KS$-module with composition factors $S_1^{r_1}/\ldots/S_m^{r_m}$. As above, the projective cover of $V$ has the form
\[ P = P_1^{n_1} + P_2^{n_2} + \ldots + P_m^{n_m} \]
with $n_i \le r_i$ for each $i$. To simplify calculations, it will be useful to find smaller upper bounds for the $n_i$ such that $V$ must still be a quotient of $P$.

\begin{lemma} \label{lem:projcover}
Let $S_1$, $\ldots$, $S_r$ be the irreducible $KS$-modules and $P_1$, $\ldots$, $P_r$ the corresponding projective indecomposables. Let $V = S_{1}^{r_1}/S_{2}^{r_2}/\ldots/S_{m}^{r_m}$ be a self-dual $KS$-module with no irreducible direct summands, and let $P = \bigoplus P_i^{n(P_i)}$ be the projective cover of $V$.

Then $n(P_i) + n(P_i^{*}) \le r_i$ for all $i$. In particular, $n(P_i) \le r_i/2$ when $S_i$ is self-dual.
\end{lemma}

\proof If $\textup{soc}(V) \nsubseteq \textup{rad}(V)$, then we may pick an irreducible submodule $W \subseteq \textup{soc}(V)$ such that $W \,\cap\, \textup{rad}(V) = \left\{0\right\}$. Since $V/\textup{rad}(V)$ is semisimple, we then have the composition of surjective maps
\[ V \twoheadrightarrow V/\textup{rad}(V) \twoheadrightarrow W \]
whose kernel does not intersect $W$, hence $W$ is an irreducible direct summand, contrary to hypothesis. Therefore we have $\textup{soc}(V) \subseteq \textup{rad}(V)$. As $V$ is self-dual, we have $V/\textup{rad}(V) \cong \textup{soc}(V)^{*}$. Hence we have
\[ P/\textup{rad}(P) \cong V/\textup{rad}(V) \cong \textup{soc}(V)^{*} \cong \bigoplus_{i=1}^{r}S_i^{n(P_i)} \]
and the result follows as the multiplicity of $S_i$ as a composition factor of $V$ is at least the sum of multiplicities of $S_i$ in $V/\textup{rad}(V)$ and $\textup{soc}(V)$. \qed

Now, suppose that $V = S_1^{r_1}/S_2^{r_2}/.../S_m^{r_m}$ is a self-dual $KS$-module (in practice $V$ will usually be the restriction to $S$ of a self-dual $KG$-module for an algebraic group $G$). Suppose that $V$ has trivial composition factors, but does not necessarily satisfy the hypotheses of Proposition \ref{prop:substab}, and we want to deduce that $V$ nevertheless contains a nonzero trivial $KS$-submodule.

Let $W$ be a direct summand of $V$ which is minimal subject to being self-dual and containing all trivial composition factors of $V$. Then $W$ lies in the principal block (that is, the block to which the trivial irreducible module belongs), and has no irreducible direct summands. In addition, since $W$ is self-dual and has no nonzero trivial submodules, it has no nonzero trivial quotients, and hence the projective cover of $W$ will have no projective indecomposable summand corresponding to the trivial module. Applying Lemma \ref{lem:projcover}, we deduce that $W$ is an image of $P = \bigoplus P_i^{m(P_i)}$, where
\[ m(P_i) = \left\{
\begin{array}{cl}
0 & : S_i \text{ is trivial or does not lie in the principal block,} \\
\lfloor r_i/2 \rfloor & : S_i \text{ is non-trivial}, S \cong S^{*} \text{ and $S$ lies in the principal block,} \\ 
r_i & : S_i \text{ is non-trivial}, S \ncong S^{*} \text{ and $S$ lies in the principal block.}
\end{array} \right.
\]

We therefore proceed by taking this module $P$, and looking for quotients which:
\begin{itemize}
\item are self-dual,
\item have composition factor multiplicities bounded above by those of $V$,
\item have precisely as many trivial composition factors as $V$, and
\item have no trivial submodules.
\end{itemize}
If no such quotients exist, then $V$ must contain a nonzero trivial submodule. 

We adopt this approach in Proposition \ref{prop:algorithm}, considering pairs $(G,H)$ not satisfying Proposition \ref{prop:substab}(i) or (iii).

\begin{remark}
In the course of proving Proposition \ref{prop:algorithm}, we refer on occasion to calculations performed over a finite splitting field for $S$, say $k \subset K$. Since $KS = kS \otimes_{k} K$, it follows that a projective $KS$-module $P$ is equal to $P_{0} \otimes_{k} K$, where $P_{0}$ is a projective $kS$-module. However, it need not be the case that every $KS$-module quotient of $P$ is obtained by extending scalars from a $kS$-module quotient of $P_{0}$. Here, the only calculations performed over a finite field are the determination of the socle and radical series of $P$, and of quotients $P/M$, where $M$ is the smallest submodule of $P$ such that $P/M$ has no composition factors of a given isomorphism type. Such an $M$ is equal to $M_{0} \otimes_{k} K$, where $M_{0}$ is the smallest submodule of $P_{0}$ such that $(P_{0}/M_{0}) \otimes_{k} K$ has no composition factors of the given isomorphism type. It suffices to work over a finite splitting field for these calculations.
\end{remark}

\subsection{Connectedness of Proper Overgroups}

Once we have deduced that a finite simple group $H$ does not occur as a Lie primitive subgroup of the adjoint exceptional simple algebraic group $G$, we must next show that each subgroup $S \cong H$ of $G$ lies in a proper \emph{connected} subgroup of $G$. Proposition \ref{prop:consub} below guarantees this in all but a few cases. We begin with a preliminary lemma.

\begin{lemma} \label{lem:maxtorus}
Let $S$ be a non-abelian finite simple subgroup of an adjoint exceptional simple algebraic group $G$. Suppose that $S$ normalises a maximal torus $T$ of $G$, so that $S$ is isomorphic to a subgroup of $W(G) \cong N_G(T)/T$. Then either $S$ lies in a proper subsystem subgroup of $G$, or $G = E_6$ and $S \cong U_{4}(2)$, or $G = E_7$ and $S \cong L_2(8)$, $U_3(3)$ or $Sp_6(2)$.
\end{lemma}

\proof The exceptional Weyl groups and their subgroup structure are well-known (for instance, see \cite{MR2562037}*{\S\S 2.8.4, 3.12.4}). The Weyl groups of type $G_2$ and $F_4$ are soluble, and hence $G$ is not one of these types. The remaining groups are $W(E_6)$, which has a subgroup of index 2 isomorphic to $U_{4}(2)$; $W(E_7) \cong 2 \times Sp_6(2)$; and $W(E_8) \cong 2\cdot\Omega_8^{+}(2).2$. The maximal subgroups of each classical group here appear in the Altas \cite{MR827219}.

If $G$ is of type $E_6$, then besides the subgroup $U_{4}(2)$ as in the statement of the lemma, there are three conjugacy classes of non-abelian simple subgroups of $W = W(E_6)$. There are two classes of subgroups isomorphic to $\Alt_5$, and each such subgroup lies in a subgroup isomorphic to $\Alt_6$, which is unique up to $W$-conjugacy. On the other hand, $G$ has a subsystem subgroup of type $A_5$, giving rise to a class of subgroups of $W$ which are isomorphic to $W(A_5) \cong \textup{Sym}_6$. Thus each simple alternating subgroup of $W$ lies in one of these, and if $S \cong \Alt_5$ or $\Alt_6$ then $TS$ lies in a Levi subgroup of $G$ of type $A_5$, as required.

Similarly, for $G = E_7$ we have $W \cong 2 \times Sp_{6}(2)$, and the non-abelian simple subgroups of $W$ are either isomorphic to one of $L_2(8)$, $U_3(3)$ or $Sp_6(2)$, which appear in the statement of the lemma, or lie in the Weyl group of a subsystem subgroup, call it $X$. The possibilities are $S \cong \Alt_5$ or $\Alt_6$ with $X$ of type $D_6$; $L_{2}(7)$, $\Alt_{7}$, or $\Alt_8$ with $X$ of type $A_7$, or $U_{4}(2)$ with $X$ of type $E_6$.

Similarly, for $G = E_8$, $W$ is a double cover of $\Omega_{8}^{+}(2):2$. The non-abelian simple subgroups of $W$ each lie in either a subgroup $L_2(7) \le W(D_8)$ (2 classes), $Sp_6(2) \le W(E_7)$, $\Alt_9 \le W(A_8)$, or $\Alt_5 \le W(A_4 A_4)$ (2 classes), and in each case $S$ lies in a proper subsystem subgroup. \qed

\begin{proposition}
Let $S$ be a non-abelian finite simple subgroup of an adjoint exceptional simple algebraic group $G$ which is not Lie primitive in $G$. Then either $S$ lies in a proper connected subgroup of $G$, or the type of $G$ and the isomorphism type of $S$ appear in the table below.
\label{prop:consub}
\end{proposition}
\begin{center}
\begin{tabular}{c|c}
$G$ & Type of $S$ \\ \hline
$E_{6}$ & $U_{4}(2)$ \\
$E_{7}$ & $L_{2}(7)$, $L_{2}(8)$, $U_{3}(3)$, $Sp_{6}(2)$ \\
$E_{8}$ & $\Alt_{5}$, $L_{2}(7)$
\end{tabular}
\end{center}

\proof Since $S$ is not Lie primitive in $G$, it lies in some maximal subgroup $X$ of positive dimension and normalises the identity component. If $S$ does not lie in $X^{\circ}$, then the image of $S$ under $X \to X/X^{\circ}$ is isomorphic to $S$. The possible maximal closed subgroups $X$, as well as $N_{G}(X)/N_{G}(X)^{\circ}$, are given by \cite{MR2044850}*{Corollary 2(i)}. In particular, the subgroups such that $N_{G}(X)/N_{G}(X)^{\circ}$ contains a finite simple group are either maximal tori, in which case $S$ is isomorphic to a subgroup of the Weyl group as in Lemma \ref{lem:maxtorus}, or $(G,X,S) = (E_{7},A_{1}^{7},L_{2}(7))$, $(E_{8},A_{1}^{8},L_{2}(7))$ or $(E_{8},A_{1},\Alt_{5})$, where this latter subgroup exists only if $p \neq 2,3,5$. \qed

\section{Lie Imprimitivity of Subgroups in Theorem \ref{THM:MAIN}: Standard Cases}

In this section, we prove that if a triple $(G,H,p)$ appears in Table \ref{tab:main}, then any subgroup $S \cong H$ of $G$ must lie in a proper, connected subgroup of $G$. We split the proof into two general propositions, as well as some cases requiring ad-hoc arguments. Recall that $\tilde{G}$ denotes the simply connected cover of $G$, and that $V_{\textup{min}}$ is a Weyl module for $\tilde{G}$ of least dimension; this has highest weight $\lambda_4$, $\lambda_1$, $\lambda_7$ or $\lambda_8$, and dimension $26$, $27$, $56$ or $248$, for $G$ respectively of type $F_4$, $E_6$, $E_7$ or $E_8$.

\begin{proposition} \label{prop:notprim}
Let $H$ be a finite simple group, not isomorphic to a member of $\textup{Lie}(p)$, and let $S \cong H$ be a finite simple subgroup of the adjoint exceptional simple algebraic group $G$, in characteristic $p \ge 0$. Let $\tilde{S}$ be a minimal preimage of $S$ in the simply connected cover $\tilde{G}$ of $G$, and let $L$ and $V$ respectively denote the quotient of $L(G)$ and $V_{\textup{min}}$ by any trivial $\tilde{G}$-submodules.

If $(G,H,p)$ appears in Table \ref{tab:substabprop}, then $\tilde{S}$ fixes a nonzero vector on $L$, $V$ or $V^{*}$, and is therefore not Lie primitive in $\tilde{G}$.

Furthermore, $\tilde{S}$ lies in a proper connected subgroup of $G$.
\end{proposition}

\begin{table}[htbp]
\caption{{Subgroup types satisfying Proposition \ref{prop:substab}}} \label{tab:substabprop}
\centering \small \onehalfspacing
\begin{tabularx}{\linewidth}{l|>{\centering\arraybackslash}X}
\hline $G$ & $H$ \\ \hline
$F_4$ & $\Alt_{7-10}$, $M_{11}$, $J_1$, $J_2$,\\ & $L_2(17) \ (p = 2)$, $U_3(3) \ (p \neq 7)$ \\ \hline
$E_6$ & $\Alt_{10-12}$, $M_{22}$, $L_2(25)$, $L_4(3)$, $U_4(2)$, $U_4(3)$, ${}^{3}D_4(2)$,\\ & $\Alt_{5} \ (p \neq 3)$, $\Alt_7 \ (p \neq 3, 5)$, $M_{11} \ (p \neq 3, 5)$, $M_{12} \ (p = 2)$, $L_2(8) \ (p = 7)$, $L_2(11) \ (p = 5)$, $L_2(13) \ (p = 7)$, $L_2(17) \ (p = 2)$, $L_2(27) \ (p \neq 2)$, $L_3(3) \ (p = 2)$, $U_3(3) \ (p = 7)$ \\ \hline
$E_7$ & $\Alt_{11-13}$, $M_{11}$, $J_1$, $L_2(17)$, $L_2(25)$, $L_3(3)$, $L_4(3)$, $U_4(2)$, $Sp_{6}(2)$, $^{3}D_4(2)$, $^{2}F_{4}(2)'$,\\ & $\Alt_{10}$ $(p = 2)$, $\Alt_{9} \ (p = 2)$, $\Alt_{8} \ (p \neq 3,5)$, $\Alt_{7} \ (p \neq 5)$, $M_{12}$ $(p \neq 5)$, $J_2$ $(p \neq 2)$, $L_2(8) \ (p \neq 3,7)$ \\ \hline
$E_8$ & $\Alt_{17}$, $\Alt_{12-15}$, $\Alt_{8}$, $M_{12}$, $J_1$, $J_2$, $L_2(27)$, $L_2(37)$, $L_4(3)$, $U_{3}(8)$, $Sp_{6}(2)$, $\Omega_{8}^{+}(2)$, $G_{2}(3)$,\\ & $\Alt_{11}$ $(p = 2)$, $\Alt_{9} \ (p \neq 2,3)$, $M_{11} \ (p \neq 3,11)$, $U_{3}(3)$ $(p \neq 7)$, $^{2}F_{4}(2)'$ $(p \neq 3)$ \\ \hline
\end{tabularx}
\end{table}

\proof Since the composition factors of $L \downarrow \tilde{S}$ and $V \downarrow \tilde{S}$ appear in the appropriate table in Chapter \ref{chap:thetables}, proving that $\tilde{S}$ is not Lie primitive in $\tilde{G}$ comes down to inspecting the corresponding table and comparing this with the information in Appendix \ref{chap:auxiliary} to decide whether the conditions in Proposition \ref{prop:substab}(i) or (iii) hold for $L \downarrow \tilde{S}$ or $V \downarrow \tilde{S}$. For convenience, in Chapter \ref{chap:thetables} we have labelled with `\possprim' those feasible characters for which we \emph{cannot} infer the existence of a fixed vector using Proposition \ref{prop:substab}. Thus $(G,H,p)$ appears in Table \ref{tab:substabprop} if and only if the corresponding table in Chapter \ref{chap:thetables} has no rows labelled `\possprim'.

As a typical example, take $G = E_6$, $H = U_3(3)$, $p = 7$. Here, we have four pairs of compatible feasible characters on $L(G)$ and $V_{\textup{min}}$, given by Table \ref{u33e6p7} on page \ref{u33e6p7}. As stated in Table \ref{tab:sporadicreps} of the Appendix, we know that $H^{1}(U_3(3),26)$ is 1-dimensional, while the corresponding group for other composition factors vanishes. Thus any subgroup $S \cong H$ of $G$, having composition factors as in Cases 2), 3) or 4) of Table \ref{u33e6p7}, satisfies Proposition \ref{prop:substab}(i) in its action on $L = L(G)$ and fixes a nonzero vector.

In Case 1), the feasible character has no trivial composition factors on $L(G)$. On the other hand, the corresponding composition factors of $V = V_{27} = V_{\textup{min}}$ are `$1$' and `$26$'. Hence by Proposition \ref{prop:substab}(iii), if $S \cong H$ gives rise to these feasible characters, a preimage $\tilde{S}$ of $S$ in $\tilde{G}$ must fix a nonzero vector on either $V$ or its dual, as required. 

Now, with the exception of $(G,H) = (E_{7},Sp_{6}(2))$ or $(E_{7},L_{2}(8))$, Proposition \ref{prop:consub} applies and so $\tilde{S}$ lies in a proper, connected subgroup of $\tilde{G}$. To show that the same holds for these two cases, for a contradiction assume that $S$ lies in no proper connected subgroup of $G$. Inspecting Tables \ref{l28e7p0}, \ref{sp62e7p0}, \ref{sp62e7p7}, \ref{sp62e7p5}, \ref{sp62e7p3}, we see that $S$ fixes a nonzero vector $v \in L(G)$. Since $S$ lies in no parabolic subgroup of $G$ by assumption, by Lemma \ref{lem:properembed} it follows that $v$ is semisimple and $C_G(v)$ contains a maximal torus, say $T$. Then $S$ normalises $C_{G}(v)^{\circ}$ and moreover, since $H$ does not occur as a subgroup of $\textup{Sym}_{7}$, the proof of Proposition \ref{prop:consub} shows that $S$ cannot normalise a non-trivial connected semisimple subgroup of $G$. It follows that $C_{G}(v)^{\circ} = T$ is a maximal torus of $G$. Now since $H = Sp_{6}(2)$ or $L_2(8)$ and $H \notin \textup{Lie}(p)$, the ambient characteristic $p$ is not $2$, hence a non-trivial $KH$-module of least dimension is $7$-dimensional. Since $S$ normalises $T$ it follows that $S$ acts irreducibly on $L(T)$. Now, we have a well-known decomposition
\[ L(G) \downarrow T = L(T) \oplus \bigoplus_{\alpha \in \Phi} L_{\alpha} \]
where $\Phi$ is the set of roots corresponding to $T$ and each $L_{\alpha}$ is a non-trivial $1$-dimensional $T$-module. Since $S$ acts irreducibly on $L(T)$ it follows that $TS$ cannot fix a nonzero vector on $L(G)$, which contradicts $TS \le C_G(v)$. \qed

\begin{proposition} \label{prop:algorithm}
With the notation of the previous Proposition, if $(G,H,p)$ appears in Table \ref{tab:algprop}, then $\tilde{S}$ lies in a proper, connected subgroup of $G$.
\end{proposition}

\begin{table}[htbp]
\caption{Further simple groups not arising as Lie primitive subgroups of $G$} \label{tab:algprop}
\begin{tabular}{c|l}
\hline $G$ & $(H,p)$ \\ \hline
$F_4$ & $(\Alt_{6},5)$, $(\Alt_{5},p \neq 2,5)$, $(L_2(7),3)$ \\
$E_6$ & $(\Alt_{9},2)$ $(L_2(7),3)$, $(L_2(27),2)$ \\
$E_8$ & $(\Alt_{10},3)$, $(\Alt_{10},2)$ \\ \hline
\end{tabular}
\end{table}

\proof For each $(G,H,p)$, we let $S$ be a hypothetical Lie primitive subgroup of $G$, and derive a contradiction by showing that a minimal pre-image $\tilde{S}$ of $S$ in $\tilde{G}$ must fix a nonzero vector on some non-trivial $\tilde{G}$-composition factor of $L(G)$ or $V_{\textup{min}}$, using the approach described in Section \ref{sec:thealg}. Since Proposition \ref{prop:consub} applies to each $(G,H)$ in Table \ref{tab:algprop}, the conclusion follows. Since $S$ is Lie primitive, recall that the composition factors of $L(G) \downarrow \tilde{S}$ and $V_{\textup{min}} \downarrow \tilde{S}$ are given by a row in the appropriate table of Chapter \ref{chap:thetables} which is marked with `\possprim'.

\subsection*{Case: $(G,H,p) = (F_4,\Alt_{6},5)$} As given in Table \ref{tab:altreps}, there are four non-trivial irreducible $KH$-modules, of dimensions 5, 5, 8 and 10. Since $|H| = 2^{3}.3^{2}.5$, the 5- and 10-dimensional modules are projective. If $P_1$ and $P_8$ are respectively the projective covers of the trivial and 8-dimensional irreducible modules, we have
\[ |H| = 360 = \textup{dim}(P_1) + 25 + 25 + 8\,\textup{dim}(P_8) + 100 \]
and thus $\textup{dim}(P_8) \le 26$. Since $5 \mid \textup{dim}(P_{8})$, this has at least two 8-dimensional composition factors, hence has composition factors $1^{4}/8^{2}$ or $1/8^{3}$. The former implies that $P_1$ has four 8-dimensional factors, contradicting the above equation. Denoting modules by their dimension, we therefore deduce that $P_8 = 8|(1+8)|8$, $P_1 = 1|8|1$.

Suppose that $S \cong H$ is a Lie primitive subgroup of $G$, so its composition factors on $L(G)$ and $V_{G}(\lambda_1)$ are given by Case 1) of Table \ref{a6f4p5}, and $V_G(\lambda_1) \downarrow S$ has no nonzero trivial submodules. If $V_{G}(\lambda_1) \downarrow S$ had an 8-dimensional $S$-direct summand, then its complement would satisfy Proposition \ref{prop:substab}(iii), and $S$ would fix a nonzero vector. Thus $V_G(\lambda_1) \downarrow S$ must be indecomposable, and is therefore an image of the projective module $P_8$ by Lemma \ref{lem:projcover}. But this has only a single trivial composition factor; a contradiction.

\subsection*{Case: $(G,H,p) = (F_4,\Alt_{5},p \neq 2,3,5)$} Here there is a unique compatible pair of fixed-point free feasible characters of $H$ on $L(G)$ and $V_{\textup{min}}$. It is proved in \cite{MR2638705}*{pp.\ 117--118}, that a subgroup $\Alt_{5}$ having these composition factors on these modules is not Lie primitive (more precisely, the subgroup stabilises a certain subalgebra on a 27-dimensional `Jordan algebra' on which $G$ acts; the stabiliser of such a structure is a positive-dimensional subgroup of $G$).

\subsection*{Case: $(G,H,p) = (F_4,\Alt_{5},3)$} The projective $KH$-modules are $3_{a}$, $3_{b}$, $P_{1} = 1|4|1$ and $P_{4} = 4|1|4$, where the latter two are uniserial. A Lie primitive subgroup $S \cong H$ of $G$ must act on $L(G)$ with composition factors as in Case 2) or 3) of Table \ref{a5f4p3}, fixing no nonzero vector. Now, all the trivial composition factors of $S$ on $L(G)$ and all the composition factors `4' must occur in a single indecomposable $S$-direct summand $W$, otherwise Proposition \ref{prop:substab}(iii) would apply to an $S$-direct summand of $L(G)$. Furthermore $W$ is self-dual since $W^{*}$ is also a direct summand of $L(G)$ with trivial composition factors. By Lemma \ref{lem:projcover} therefore, $W$ is an image of $P_4^{3}$. But this has only three trivial composition factors; a contradiction.

\subsection*{Case: $(G,H,p) = (F_4\ \textup{or}\ E_{6},L_2(7),3)$} The projective $KH$-modules here are $3_{a}$, $3_{b}$, $6$, $P_{1} = 1|7|1$ and $P_7 = 7|1|7$. Let $L = L(G)$ if $G$ has type $F_{4}$, or $V_{G}(\lambda_2)$, of dimension 77, if $G$ has type $E_{6}$. Using Table \ref{l27f4p3} and \ref{l27e6p3}, we see that $L \downarrow S$ has four trivial composition factors, and a similar argument to the above shows that if $S \cong H$ is a subgroup of $G$ fixing no nonzero vectors on $L$, then $L \downarrow S$ has an indecomposable self-dual direct summand $W$ containing all trivial composition factors and at least four 7-dimensional composition factors. This summand is then an image of $P_{7}^{3}$ by Lemma \ref{lem:projcover}. But this only has three trivial composition factors; a contradiction.

\subsection*{Case: $(G,H,p) = (E_6,\Alt_{9},2)$}
Here a Lie primitive subgroup $S \cong H$ of $G$ must give rise to Case 1) of Table \ref{a9e6p2}. Since $V_{27} \downarrow \tilde{S}$ then only has two composition factors, one of which is trivial, it follows that either $V_G(\lambda_1)$ or its dual $V_G(\lambda_6)$ has a nonzero trivial $\tilde{S}$-submodule; a contradiction.

\subsection*{Case: $(G,H,p) = (E_6,L_2(27),2)$}
Here a Lie primitive subgroup $S \cong H$ of $G$ must give rise to Case 1) of Table \ref{l227e6p2}. Since $\tilde{S}$ fixes no nonzero vector on $V_{G}(\lambda_1)$ or its dual, we see that $V_{G}(\lambda_1) \downarrow \tilde{S} = 13|1|13^{*}$ or $13^{*}|1|13$. This is therefore an image of $P_{13}$ or $P_{13}^{*}$.

Now, the six 28-dimensional irreducible $KH$-modules are projective. Calculations with dimensions, and the fact that $P/\textup{rad}(P) \cong \textup{soc}(P)$ for any projective $KH$-module, quickly show that the projective cover $P_{26_a} = 26_a|26_a$ is uniserial, and similarly for $26_b$ and $26_c$; these involve no 13-dimensional factors, hence $P_{13}$ and $P_{13^{*}}$ have no 26-dimensional composition factors. Since $P_{1}$ is self-dual and involves only the modules $1$, $13$, $13^{*}$, dimension considerations imply that $P_1 = 1|(13 + 13^{*})|1$, and it follows that $P_{13} = 13|(1 + 13^{*})|13$. Thus neither $P_{13}$ nor $P_{13}^{*}$ has a uniserial quotient $13|1|13^{*}$ or $13^{*}|1|13$; a contradiction.

\subsection*{Case: $(G,H,p) = (E_8,\Alt_{10},3)$}
Here, every composition factor of the unique feasible character of $H$ on $L(G)$ (Table \ref{a10e8p3}) is self-dual, and all except one (that of dimension 84) have multiplicity 1. It follows that any irreducible $S$-submodule of dimension $\neq 84$ is in fact a direct summand. Thus if $L(G) \downarrow S$ has a reducible, indecomposable direct summand, say $W$, then $W$ is a quotient of the projective module $P_{84}$. The radical series of $P_{84}$ begins $84 | 1 + 41 + 84 | 34 + 41 + 84^{2} | ... $. Now if $N$ is a submodule of $P_{84}$ such that $N/\textup{Rad}(N) = 41$, then $N$ lies in the kernel of $P_{84} \to W$. Using Magma to facilitate calculations, we find that the quotient of $P_{84}$ by the sum of all such submodules is self-dual with shape $84 | (1 + 84) | 84$. It follows that either $L(G) \downarrow S$ is completely reducible, or
\[ L(G) \downarrow S = 1+9+34+36+(84|84). \]
In either case, $S$ fixes a 1-space on $L(G)$, and hence cannot be Lie primitive in $G$.

\subsection*{Case: $(G,H,p) = (E_8,\Alt_{10},2)$}
Here a Lie primitive subgroup $S \cong H$ of $G$ gives rise to Case 4) of Table \ref{a10e8p2}, so $L(G) \downarrow S = 1^{8}/8^{5}/26^{4}/48^{2}$. Hence the sum $\bigoplus H^{1}(S,W)$ over $S$-composition factors of $L(G)$ is 9-dimensional. Now, let $L(G) = W_1 \oplus \ldots \oplus W_t$ where each $W_i$ is indecomposable. If more than one $W_i$ had a trivial composition factor, then Proposition \ref{prop:substab}(iii) would apply to at least one such summand. Since $S$ fixes no nonzero vectors of $L(G)$, it follows that $S$ has an indecomposable summand $W$ on $L(G)$ containing all trivial composition factors, and also all $S$-composition factors with nonzero first cohomology group. Moreover $W$ is an image of $P_8^{2}+P_{26}^{2}+P_{48}$, by Lemma \ref{lem:projcover}. Now, let $Q_{8}$, $Q_{26}$ and $Q_{48}$ be the quotient of $P_{8}$, $P_{26}$ and $P_{48}$, respectively, by the sum of all submodules with a unique maximal submodule, and corresponding quotient not isomorphic to a member of $\{1,8,26,48\}$. Using Magma to facilitate calculations, we find that these quotients have the following structure:
\begin{align*}
Q_{8} &= 8 | 1 | 26 | 1 | (8 + 48) | (1 + 26) | (1 + 26) | (1 + 8) | (1 + 8),\\
Q_{26} &= 26 | (1 + 48) | (8 + 26) | 1^{2} | (8 + 26) | (1^{2} + 48) | 8 + 26 | 1 \\
Q_{48} &= 48 | 26 | 1 | 8 | (1 + 48) | 26 | 1
\end{align*}
Now, if $W/\textup{Rad}(W)$ has more than one composition factor `$8$' or `$26$', then $\textup{Rad}(W)$ satisfies Proposition \ref{prop:substab}(i) and has a trivial submodule. So we may assume that $W$ is an image of $Q_{48} \oplus Q_{26}$ or $Q_{48} \oplus Q_{8}$. Each of these modules has exactly nine trivial composition factors, so the kernel of the projection to $W$ can have at most one trivial factor. But also, both $Q_{48} \oplus Q_{26}$ and $Q_{48} \oplus Q_{8}$ contain a $2$-dimensional trivial submodule, and it follows that $S$ fixes a nonzero vector of $W \subseteq L(G)$. \qed

\section{Postponed Cases of Theorem \ref{THM:MAIN}} \label{sec:postponed}

Those pairs $(G,H)$ in Table \ref{tab:main} for which the above propositions do not apply are listed below. For these, a special argument is required, either because the groups themselves have have feasible characters with no trivial composition factors on $L(G)$ or $V_{\textup{min}}$, or because their representation theory allows for the existence of modules having appropriate composition factors but no fixed vectors. For these, we defer proving the conclusion of Theorem \ref{THM:MAIN} until Section \ref{sec:special}, where some ad-hoc arguments are applied.

\begin{table}[H]\small
\begin{tabular}{c|c}
$G$ & $(H,p)$ \\ \hline
$E_7$ & $\Alt_{10} \ (p = 5)$, $\Alt_{9} \ (p  \neq 2, 3)$ \\
$E_8$ & $\Alt_{16} \ (p = 2)$, $\Alt_{11}\ (p = 11)$, $\Alt_{10} \ (p > 3)$
\end{tabular}
\end{table}

%% file: AJL-stab.tex

\chapter{Normaliser Stability} \label{chap:stab}

In this chapter we complete the proof of Theorem \ref{THM:MAIN}. At this stage we have now shown that for $(G,H,p)$ as in Table \ref{tab:main}, with the exception of the `postponed cases' considered in Section \ref{sec:special}, every subgroup $S \cong H$ of $G$ lies in a proper connected subgroup of $G$. It remains to show the existence of a connected $N_{\textup{Aut}(G)}(S)$-stable subgroup. We apply a variety of techniques to achieve this; the following proposition summarises the results of this chapter.

\begin{proposition}
Let $G$ be an adjoint exceptional simple algebraic group in characteristic $p$, let $H$ be a finite simple group, not isomorphic to a member of $\textup{Lie}(p)$, and let $S \cong H$ be a subgroup of $G$. If $G$, $H$, $p$ appear in Table \ref{tab:main}, then one of the following holds:
\begin{itemize}
\item $S$ is not $G$-completely reducible, hence Lemma \ref{lem:gcrstab} applies;
\item $S$ is Lie primitive in a semisimple subgroup $X$ of $G$ as in Proposition \ref{prop:primstab};
\item $(G,H,p)$ appears in one of the tables \ref{tab:f4thm}--\ref{tab:e8thm}, hence Proposition \ref{prop:gcrstab} applies;
\item $(G,H,p)$ appears in Table \ref{tab:adhoc} in Section \ref{sec:special}.
\end{itemize}
Hence $S$ is contained in a proper, $N_{\textup{Aut}(G)}(S)$-stable connected subgroup of $G$.
\end{proposition}

\section{Complete Reducibility and Normaliser Stability}

Recall that a subgroup $S$ of a reductive group $G$ is called \emph{$G$-completely reducible} ($G$-cr) if whenever $S$ is contained in a parabolic subgroup $P$ of $G$, it is contained in a Levi subgroup of $P$. In \cite{MR2178661} it is shown that a subgroup of a reductive algebraic group $G$ is $G$-cr if and only if it is `strongly reductive' in the sense of Richardson \cite{MR952224}. A result of Liebeck, Martin and Shalev then states:
\begin{lemma}[\cite{MR2145743}*{Proposition 2.2 and Remark 2.4}]
\label{lem:gcrstab}
Let $S$ be a finite a finite subgroup of an adjoint simple algebraic group $G$. Then either $S$ is $G$-completely reducible, or $S$ is contained in a proper $N_{\textup{Aut}(G)}(S)$-stable parabolic subgroup of $G$.
\end{lemma}

It thus remains to prove Theorem \ref{tab:main} for $G$-cr subgroups $S \cong H$. In this case, let $L$ be minimal among Levi subgroups of $G$ containing $S$, so that $S$ lies in the semisimple subgroup $L'$ and is $L'$-irreducible. Then any connected subgroup of $L'$ containing $S$ is also $L'$-irreducible, and is therefore $G$-cr and semisimple. Thus there exists a proper semisimple subgroup $X = X_{1} \ldots X_{t}$ such that $S$ projects to an $X_{i}$-irreducible subgroup of each simple factor $X_{i}$. If we pick $X$ to be minimal among semisimple subgroups of $L'$ containing $S$, then the image of $S$ under projection to each simple factor is in fact Lie primitive in that factor.

\section{Subspace Stabilisers and Normaliser Stability}

We now assume that $S$ is $G$-completely reducible. Our strategy now is to construct a `small' $G$-cr semisimple subgroup $X$ of $G$ containing $S$. Informally, `small' means that the actions of $X$ and $S$ on $L(G)$ or $V_{\textup{min}}$ should be similar, which allows us to apply a number of results that we shall state in a moment. Note that this does not require that $S$ is Lie primitive in $X$, although this will be true in many cases. 

Our first result of interest is Proposition 1.12 of \cite{MR1458329}. Recall that if $M$ is a $G$-module with corresponding representation $\rho : G \to GL(M)$, the \emph{conjugate} $M^{\tau}$ of $M$ by an (abstract) automorphism $\tau$ of $G$ is the module corresponding to the representation $\tau \rho : G \to GL(M)$. If $G$ is an algebraic group, if $\tau$ is a morphism and $M$ is rational, then clearly $M^{\tau}$ is also rational.

Recall that a Suzuki-Ree group is the fixed-point subgroup $G_{\phi}$ when $\phi$ is an exceptional graph morphism of $G$. For a subspace $M$ of a $G$-module $V$, let $G_{M}$ denote the corresponding subspace stabiliser.

\begin{proposition}[\cite{MR1458329}*{Proposition 1.12}] \label{prop:onetwelve}
Let $G$ be a simple algebraic group over $K$, and let $\phi : G \to G$ be a morphism which is an automorphism of abstract groups.
\begin{itemize}
\item[\textup{(i)}] Suppose that $G_{\phi}$ is not a finite Suzuki or Ree group, and let $V$ be a $G$-composition factor of $L(G)$. If $M$ is a subspace of $V$, then $(G_M)^{\phi} = G_{M'}$ for some subspace $M'$ of $V$.
\item[\textup{(ii)}] Suppose $G_{\phi}$ is a finite Suzuki or Ree group, and let $V_1$, $V_2$ be the two $G$-composition factors of $L(G)$. If $M$ is a subspace of $V_i$ $(i = 1,2)$, then $G_M^{\phi} = G_{M'}$ for some subspace $M'$ of $V_{3-i}$.
\item[\textup{(iii)}] Let $S$ be a $\phi$-stable subgroup of $G$, and let $\mathcal{M}$ be the collection of all $S$-invariant subspaces of all $G$-composition factors of $L(G)$. Then the subgroup $\bigcap_{W \in \mathcal{M}} G_W$ of $G$ is $\phi$-stable.
\end{itemize}
\end{proposition}
Of particular interest to us here is part (iii). If $S$ is a finite subgroup of $G$, and if $X$ is a connected subgroup containing $S$ such that every $S$-submodule of every $G$-composition factor of $L(G)$ is an $X$-submodule, then the group $\bigcap_{W \in \mathcal{M}} G_W$ in (iii) contains $X$, and is therefore also of positive dimension. Applying this result for each morphism $\phi \in N_{\textup{Aut}(G)}(S)$, we obtain a positive-dimensional, $N_{\textup{Aut}(G)}(S)$-stable subgroup of $G$, whose identity component contains $X$ and therefore $S$.

It will be useful for us to extend the above result, since we will encounter cases when, for $X$ a minimal semisimple subgroup containing $S$, not every $S$-submodule of $L(G)$ is an $X$-submodule. We now do this by mimicking the proof given in \cite{MR1458329}.

\begin{proposition} \label{prop:newonetwelve}
Let $G$ be a simple algebraic group over $K$, let $\phi : G \to G$ be a morphism which is an automorphism of abstract groups.
\begin{itemize}
\item[\textup{(i)}] Let $V = \bigoplus V_G(\lambda_i)$ be a completely reducible $KG$-module such that the set $\{\lambda_i\}$ is stable under all graph morphisms of $G$. If $M$ is a subspace of $V$, then $(G_M)^{\phi} = G_{M\delta}$ where $\delta$ is an invertible semilinear transformation $V \to V$ depending on $\phi$ but not on $M$.
\item[\textup{(ii)}] If $S$ is a subgroup of $G$, then for each $KS$-submodule $M$ of $V$, the subspace $M\delta$ is a $KS^{\phi}$-submodule, of the same dimension as $M$, which is irreducible, indecomposable or completely reducible if and only if $M$ has the same property.
\end{itemize}
\end{proposition}

\proof (i) Let $V$ correspond to the representation $\rho : G \to GL(V)$. We may write $\phi = y\tau\sigma$ where $y$, $\tau$ and $\sigma$ are (possibly trivial) inner, graph and field morphisms of $G$, respectively. By assumption, the representations $\rho$ and $\tau\rho$ of $G$ are equivalent, since they are completely reducible with identical high weights. Hence if $\sigma$ is a $q$-power field automorphism, where $q = p^{e} \ge 1$, then the high weights of $\phi\rho$ are $\{q\lambda_i\}$. There is therefore a $q$-power field automorphism $\omega$ of $GL(V)$ such that $\phi\rho$ and $\rho\omega$ are equivalent. The automorphism $\omega$ is induced by a semilinear transformation $V \to V$ which we shall also denote by $\omega$. Then $y^{\omega} = \omega^{-1}y\omega$ for $y \in GL(V)$. Thus, identifying each $g \in G$ with its image $g\rho \in GL(V)$, there exists $x \in GL(V)$ such that $g^{\phi} = g^{\omega x} = x^{-1}\omega^{-1}g\omega x$, for all $g \in G$. Writing $\delta = \omega x$, this gives $\delta g^{\phi} = g\delta$ for all $g \in G$, and we have
\[ (v\delta)g^{\phi} = (vg)\delta \]
for all $v \in V$, $g \in G$. If $M$ is a subspace of $V$, and $m \in M$, $g \in G_M$, then $(m\delta)g^{\phi} = (mg)\delta \in M\delta$, and hence $g^{\phi} \in G_{M\delta}$. Therefore $(G_M)^{\phi} \le G_{M\delta}$. For the reverse inclusion, if $g \in G_{M\delta}$, then by the displayed equality above, for any $m \in M$ we have
\[ (m\delta)g = (m\delta).(g^{\phi^{-1}})^{\phi} = (mg^{\phi^{-1}})\delta = m'\delta \]
for some $m' \in M$. Therefore $mg^{\phi^{-1}} = m'$ and $g^{\phi^{-1}} \in G_M$, so $g \in (G_M)^{\phi}$ as required.

(ii) If $M$ is a $KS$-submodule of $V$, then the displayed equation above, applied to the elements of $S^{\phi}$, tells us that $S^{\phi}$ preserves the subspace $M\delta$ of $V$. It is clear that $M$ and $M\delta$ have the same dimension, since $\delta$ is invertible. If $W \subseteq M$ is a nonzero $KS$-submodule of $M$, then $W\delta$ is a nonzero $KS^{\phi}$-submodule of $M \delta$, and $M = M_1 + M_2$ as $KS$-modules if and only if $M\delta = M_1\delta + M_2\delta$, proving the final claim. \qed

Thus if $S = S^{\phi}$, then $\phi$ induces a permutation on the $S$-submodules of $L(G)$, and we immediately deduce:
\begin{corollary} \label{cor:phistab}
Let $S$ be a subgroup of $G$, and let $V$ be a $KG$-module as in Proposition \ref{prop:newonetwelve}\textup{(i)}, or the direct sum of the $G$-composition factors of $L(G)$ if $(G,p) = (F_4,2)$. Let $\phi$ be a morphism in $N_{\textup{Aut}(G)}(S)$, and let $\mathcal{M}$ be one of:
\begin{itemize}
\item The set of all $KS$-submodules of $V$, or all irreducible $KS$-submodules, or all indecomposable $KS$-submodules.
\item Those members of one of the above collections, with a prescribed set of composition factor dimensions.
\end{itemize}
Then the intersection $H \stackrel{\textup{def}}{=} \bigcap_{M \in \mathcal{M}} G_M$ is $\phi$-stable.

Further, if some member of $\mathcal{M}$ is not $G$-stable, then $H$ is proper. If $S$ lies in a positive-dimensional subgroup $X$ such that each member of $\mathcal{M}$ is $X$-invariant, then $H$ is a positive-dimensional. If $X$ is connected, then $S$ lies in $H^{\circ}$, which is connected and $\phi$-stable.
\end{corollary}

Thus with $S$ a finite simple subgroup of $G$, lying in the connected subgroup $X$, we are interested in techniques for spotting when $KS$-submodules of a given $KX$-module are $X$-invariant. The following result of Liebeck and Seitz provides such a method.
\begin{lemma}[\cite{MR1458329}*{Proposition 1.4}]
Let $X$ be an algebraic group over $K$ and let $S$ be a finite subgroup of $X$. Suppose $V$ is a finite-dimensional rational $KX$-module satisfying the following conditions:
\begin{itemize}
\item[\textup{(i)}] every $X$-composition factor of $V$ is $S$-irreducible,
\item[\textup{(ii)}] for any $X$-composition factors $M, N$ of $V$, the restriction map \\ ${\rm Ext}_{X}^{1}(M,N)$ $\to$ ${\rm Ext}_{S}^1(M,N)$ is injective,
\item[\textup{(iii)}] for any $X$-composition factors $M$, $N$ of $V$, if $M \downarrow S \cong N \downarrow S$, then $M \cong N$ as $X$-modules.
\end{itemize}
Then $X$ and $S$ fix precisely the same subspaces of $V$.
\label{lem:samesubs}
\end{lemma}
Conditions (i) and (iii) are straightforward to verify. Condition (ii) can often be checked by showing that the groups ${\rm Ext}_{X}^{1}(M,N)$ are trivial, for example using Lemmas \ref{lem:exthom} and \ref{lem:2step}.

The proof of Lemma \ref{lem:samesubs} given in \cite{MR1458329} uses condition (ii) only to deduce that an indecomposable $X$-module section of $V$ remains indecomposable as an $S$-module. This allows the following generalisation:
\begin{proposition} \label{prop:newsamesubs}
The conclusion of Lemma \ref{lem:samesubs} holds if we replace condition \textup{(ii)} with either:
\begin{itemize}
\item[\textup{(ii$'$)}] Each indecomposable $KX$-module section of $V$ is indecomposable as a $KS$-module.
\item[\textup{(ii$''$)}] As a $KX$-module, $V$ is completely reducible.
\end{itemize}
\end{proposition}

\proof Note that, assuming condition (i), we have implications (ii$''$) $\Rightarrow$ (ii$'$), and (ii) $\Rightarrow$ (ii$'$). Thus it suffices to assume that (i), (ii$'$) and (iii) hold. From here we proceed as in \cite{MR1458329}, by induction on $\textup{dim }V$, noting that the case $\textup{dim }V = 1$ is trivial.

For a contradiction, suppose that some $KS$-submodule of $V$ is not $X$-stable, and let $W$ be minimal among such submodules. If $W'$ is a proper $KS$-submodule of $W$, then $W'$ is $X$-invariant, and $W/W'$ is a $KS$-submodule of $V/W'$ which is not $X$-invariant. The inductive hypothesis thus forces $W' = 0$, so $W$ is irreducible.

Now let $U = \left<W^{X}\right>$, so $U \neq W$. If $U$ were irreducible for $X$, it would be irreducible for $S$ by (i), contradicting $U \neq W$. Thus there exists a proper, irreducible $X$-submodule $W_0$ of $W$.

Consider $V/W_0$. Then $S$ fixes the subspace $(W + W_0)/W_0$ of this, and by induction we deduce that $X$ fixes $W + W_0$. Hence $U = W + W_0$ (vector space direct sum). Since $W$ and $W_0$ are irreducible as $S$-modules, this is also a direct-sum decomposition of $U$ into $S$-submodules. Thus $U$ is not indecomposable as an $S$-module, and hence by (ii$'$) is also not indecomposable as an $X$-module. Hence there is an $X$-submodule $W_1$ of $U$ such that
\[ U = W + W_0 = W_1 + W_0 \]
where $W$ is $S$-isomorphic to $W_1$. If $W$ is not $S$-isomorphic to $W_0$, then $W$ and $W_0$ are the only irreducible $S$-submodules of $U$, so $W = W_1$ and $W$ is $X$-invariant, a contradiction. Thus $W$ is $S$-isomorphic to $W_0$, and thus $W_0$ and $W_1$ are $X$-isomorphic. Now $W \subseteq W_1 + W_0$, and we have
\[ W = \left\{ w + w\phi \ : \ w \in W_1 \right\} \]
for some $S$-isomorphism $\phi$ : $W_1 \to W_0$. But if $\alpha$ is any $X$-isomorphism : $W_1 \to W_0$, then $\alpha\phi^{-1}$ : $W_1 \to W_1$ is an $S$-isomorphism, hence by Schur's Lemma we have $\alpha \phi^{-1} = \lambda.\textup{id}_{W_1}$ for some $\lambda \in K^{*}$. Hence $\phi = \lambda.\alpha$ is an $X$-isomorphism, and $W$ is fixed by $X$, which is a contradiction. Therefore $W = U$, as required. \qed

\subsection{Restriction of $G$-modules to $G$-cr semisimple subgroups}

In order to compare the action of a finite simple subgroup of $G$ on various $G$-modules with the action of a $G$-cr semisimple subgroup $X$, we need to determine some details of how such a subgroup $X$ acts. To begin, let $L$ be minimal among Levi subgroups of $G$ containing $S$. We can thus assume that $S < X \le L'$. The action of $L'$ on $L(G)$ and $V_{\textup{min}}$ now follows from the known composition factors of $L'$, stated in \cite{MR1329942}*{Tables 8.1-8.7}, and Lemma \ref{lem:tilting}.

If we work with the simply-connected cover $\tilde{G}$ of $G$, and a minimal pre-image $\tilde{S}$ of $S$ in $\tilde{G}$, the derived subgroup of a Levi subgroup $L$ of $\tilde{G}$ is simply connected \cite{MR2850737}*{Proposition 12.14}, hence is a direct product of simply connected simple groups. The image of $\tilde{S}$ under projection to a simple factor $L_{0}$ of $L'$ is $L_{0}$-irreducible. If $L_{0}$ is classical, we can use Lemma \ref{lem:classicalparabs} to find a smaller connected subgroup of $L_{0}$ containing the image of $\tilde{S}$, such as the stabiliser of a direct-sum decomposition of the natural module. On the other hand, if $L_{0}$ is exceptional, then we can use Propositions \ref{prop:notprim} and \ref{prop:algorithm}, and the feasible characters in Chapter \ref{chap:thetables}, to find a smaller semisimple subgroup of $L_{0}$ containing the image of $S$ (if one exists).

This gives a `small' semisimple subgroup $X$ of $L'$ which contains $S$. The known action of $L'$ on the various $G$-modules is then usually enough information to determine the action of $X$. If $X$ is simple of rank $> \frac{1}{2}\textup{rank}(G)$, then $X$ is given up to conjugacy by \cite{MR1367085}*{Theorem 1}, and the restrictions of $L(G)$ and $V_{\textup{min}}$ to such a subgroup are given, at least up to composition factors, by \cite{MR1329942}*{Tables 8.1-8.7} or \cite{MR3075783}*{Chapter 5}. Once the composition factors are known, more precise information about the module structure can be inferred using Proposition \ref{prop:weyls} and Lemma \ref{lem:2step}.

The following summarises the module restrictions we need which are most difficult to verify from the above sources.

\begin{lemma} \label{lem:c4d4}
Let $G$ be a simply connected simple algebraic group of exceptional type in characteristic $p = 2$. If $X$ is a $G$-cr simple subgroup of type $C_4$ or $D_4$, then either $X = L'$ for a Levi subgroup $L$ of $G$, or $X$ is conjugate to a subgroup in Table \ref{tab:c4d4res}, acting on the $G$-module $V$ as stated.
\end{lemma}

\begin{table}[htbp]
\centering
\small
\caption{Non-Levi, $G$-cr simple subgroups of type $C_{4}$, $D_{4}$, $p = 2$}
\begin{tabular}{c|c|c|c}
$G$ & $X$ & $V$ & $V \downarrow X$ \\ \hline
$F_4$ & $C_4$ & $V_G(\lambda_1)$ & $0^{2}/V_X(\lambda_4)/V_X(2\lambda_1)$ \\
& & $V_G(\lambda_4)$ & $V_X(\lambda_2)$ \\
& $D_4$ (long) & $V_G(\lambda_1)$ & $V_{X}(\lambda_2)$ \\
& & $V_G(\lambda_4)$ & $0^{2} \oplus \lambda_1 \oplus \lambda_3 \oplus \lambda_4$ \\
& $D_4$ (short) & $V_G(\lambda_1)$ & $0^{2} \oplus V_X(2\lambda_1) \oplus V_X(2\lambda_3) \oplus V_X(2\lambda_4)$ \\
& & $V_G(\lambda_4)$ & $V_{X}(\lambda_2)$ \\
$E_6$ & $C_4 < F_4$ & $V_G(\lambda_1)$ & $0 \oplus V_X(\lambda_2)$ \\
& $D_4 < F_4$ (short) & $V_G(\lambda_1)$ & $0 \oplus V_X(\lambda_2)$ \\
$E_7$ & $C_4 < F_4$ & $V_G(\lambda_7)$ & $0^{4} \oplus V_X(\lambda_2)^{2}$ \\
& $C_{4} < A_{7}$ & $V_{G}(\lambda_7)$ & $(0|V_{X}(\lambda_2)|0)^{2}$ \\
& $D_4 < F_4$ (short) & $V_G(\lambda_7)$ & $0^{4} \oplus V_X(\lambda_2)^{2}$ \\
$E_{8}$ & $C_{4} < A_{7}$ Levi & $L(G)$ & $(\lambda_1 \otimes \lambda_1) \oplus \lambda_1^{4} \oplus (0|V_{X}(\lambda_2)|0)^{2} \oplus \lambda_3^{2}$ \\
& $C_{4} < F_{4} $ & $L(G)$ & $(0^{2}/V_{X}(\lambda_2)^{2}/V_{X}(\lambda_4)/V_{X}(2\lambda_1)) \oplus 0^{14} \oplus V_{X}(\lambda_2)^{6}$\\
& $D_{4} < A_{7}$ Levi & $L(G)$ & $(\lambda_1 \otimes \lambda_1) \oplus \lambda_1^{4} \oplus (0|V_{X}(\lambda_2)|0)^{2} \oplus V_X(\lambda_3 + \lambda_4)^{2}$ \\
& $D_4 < F_4$ (short) & $L(G)$ & $(0^{2}/V_{X}(\lambda_2)^{2}/V_{X}(2\lambda_1)/V_{X}(2\lambda_3)/V_{X}(2\lambda_4))$ \\ 
& & & ${}\oplus 0^{14} \oplus V_{X}(\lambda_2)^{6}$ \\
& $D_{4} < D_{4}D_{4}$ & & Infinitely many, $G$-irreducible, cf.\ \cite{MR3283715}*{Theorem 3} \\
\end{tabular}
\label{tab:c4d4res}
\end{table}

\proof Theorem 1 of \cite{MR1367085} lists all subgroups $C_{4}$ or $D_{4}$ when $G$ has type $F_{4}$, $E_{6}$ or $E_{7}$. The non-Levi subgroups $D_{4} < E_{7}$ given in part (IV) there are non-$E_7$-cr, except for the subgroup $D_{4} < F_{4}$ generated by short root subgroups. The $C_{4}$ and $D_{4}$ subgroups of $A_{7} < E_{7}$ are non-$E_{7}$-cr, by \cite{MR1274094}*{Lemma 4.9} (see also \cite{MR3283715}*{Lemma 6.1}).

For $G$ of type $E_{8}$, Theorem 3 of \cite{MR3283715} states that the only $G$-irreducible subgroups of type $C_{4}$ or $D_{4}$ are the subgroups $D_{4} < D_{4}D_{4}$ listed above. Each remaining $G$-cr subgroup lies in some Levi subgroup $L$ of $G$. Using the list of subgroups for $E_{6}$ and $E_{7}$, it follows that $L'$ is simple of type $A_{7}$ or $E_{6}$, which gives the remaining subgroups for $E_{8}$.

The composition factors of $V \downarrow X$ now follow from \cite{MR1329942}*{Table 8.1-8.7}. For $G$ of type $F_{4}$, the module structure of $V_G(\lambda_4) \downarrow X$ is stated in \cite{MR3075783}*{Chapter 5}. Since $V_G(\lambda_1)$ and $V_G(\lambda_4)$ are conjugate by an exceptional graph morphism of $G$, which swaps the long and short $D_4$ subgroups, a direct-sum decomposition of $V_G(\lambda_4)$ as a module for a long subgroup $D_4$, implies a decomposition of $V_G(\lambda_1)$ for a short $D_4$, and vice-versa. The given module structures for $G = F_4$ now follow.

For $G \neq F_{4}$ the stated module structures follow by first considering $V \downarrow F_4$ or $V \downarrow A_{7}$, which are straightforward to derive using the known composition factors and Lemmas \ref{lem:tilting} and \ref{lem:exthom}. Since $C_{4}$ and $D_{4}$ each support a unique nondegenerate bilinear form on their natural modules, there exists a unique nonzero $KX$-module homomorphism $\bigwedge^{2}(\lambda_1) \to K$ (up to scalars), and thus $V_{A_7}(\lambda_2) \downarrow X = V_{A_7}(\lambda_6) \downarrow X = 0 | V_{X}(\lambda_2)|0$ for each $X$. Finally, for $X = C_{4}$, $\bigwedge^{3}(\lambda_1) = \lambda_1 + \lambda_3$, by consideration of high weights and Lemma \ref{lem:exthom}. The remaining restrictions follow. \qed

\section{Proof of Theorem \ref{THM:MAIN}: Standard Cases}

In view of the above results, Propositions \ref{prop:primstab} and \ref{prop:gcrstab} below complete the proof of Theorem \ref{THM:MAIN} except for the `postponed cases' considered in Section \ref{sec:special}.

\begin{proposition} \label{prop:primstab}
Let $G$ be a simple algebraic group of exceptional type, let $S$ be a finite subgroup of $G$, and suppose that $S$ is Lie primitive in $G$-cr subgroup $X$ of $G$. If $X$ is the derived subgroup of a Levi subgroup of $G$, or if $X$ is as follows:
\begin{itemize}
\item $X = G_{2}$ in a Levi subgroup of type $D_{4}$ or $B_{3}$;
\item $X = B_{n-1}$ in a Levi subgroup of type $D_{n}$ $(n \ge 3)$;
\item $X < F_{4}$, and $X$ has type $A_{3}$, $B_{4}, C_{4}$, $D_{4}$ or $F_{4}$.
\end{itemize}
then $X$ is $N_{\textup{Aut}(G)}(S)$-stable.
\end{proposition}

\proof To begin, note that if $\textup{dim}(X) > \frac{1}{2}\textup{dim}(G)$, then $X \cap X^{\sigma}$ has positive dimension for any $\sigma \in \textup{Aut}(G)$. In particular, if $\sigma \in N_{\textup{Aut}(G)}(S)$ then $X \cap X^{\sigma} = X$ since $S$ is Lie primitive in $X$ and $X$ is connected. This gives the result when $G = F_{4}$ and $X$ is a proper subgroup of maximal rank listed above, so we now assume that this is not the case.

Let $V$ be the direct sum of non-trivial $G$-composition factors of either $L(G)$ or $V_{\textup{min}} \oplus V_{\textup{min}}^{*}$. For the remaining groups $X$, we prove that every fixed point of $S$ on $V$ is a fixed point of $X$ on $V$, and that $X$ has a nonzero fixed point on some such $V$. The desired conclusion then follows from Corollary \ref{cor:phistab}.

Firstly, the composition factors of $X$ on $V$ are known by \cite{MR1329942}*{Tables 8.1-8.7}. Then Lemmas \ref{lem:radfilts} and \ref{lem:exthom} show that $X$ fixes a nonzero vector on $V$. Moreover every composition factor of $V \downarrow X$ has dimension at most $\textup{dim}(X)$, and equality holds only if this composition factor is isomorphic to $L(X)$. Thus by Lemma \ref{lem:properembed}, since $S$ is Lie primitive in $X$ it cannot fix a nonzero vector in its action on any nontrivial $X$-composition factor of $V \downarrow X$.

Next, if $V \downarrow X$ has an indecomposable section of the form $K|W$, where $W$ is irreducible of dimension at most $\textup{dim}(X) - 2$, then this extension cannot split as an $S$-module (since the corresponding vector centraliser has positive dimension). It remains to show that this must also hold if $W$ instead has dimension $\ge \textup{dim}(X) - 1$. From the known action of $X$ on $V$, this can only occur in the following cases:

\begin{itemize}
\item $X = A_{n}$, $p \mid n-1$, $W = V_{X}(\lambda_1 + \lambda_n)$;
\item $X = D_{n}$ ($n$ odd) or $B_{n}$, $p = 2$, $W = V_{X}(\lambda_2)$;
\item $X = E_{7}$, $p = 2$, $W = V_{X}(\lambda_1)$;
\item $X = E_{6}$, $p = 3$, $W = V_{X}(\lambda_2)$.
\end{itemize}

In each case, $H^{1}(X,W)$ is $1$-dimensional, hence there is a unique indecomposable extension $K|W$ up to isomorphism. Moreover the representation $X \to GL(K|W)$ factors through the adjoint group $X_{\textup{ad}}$, and $K|W$ is isomorphic to $L(X_{\textup{ad}})$. In particular since $S$ is Lie primitive in $X$, its image in $X_{\textup{ad}}$ is also Lie primitive and so $S$ cannot fix a nonzero vector on $L(X_{\textup{ad}})$, by Lemma \ref{lem:properembed}. Thus every fixed point of $S$ on $W$ is fixed by $X$, as required. \qed

For reference, the following table lists the types $H$ of non-generic finite simple subgroup of $G$ such that each subgroup $S \cong H$ of $G$ is necessarily Lie primitive in some simple subgroup $X$ as in Proposition \ref{prop:primstab}. For each type $H$, the types of $X$ which may occur are straightforward to determine from Lemma \ref{lem:classicalparabs} and Propositions \ref{prop:notprim} and \ref{prop:algorithm}. For instance, when $p = 5$, $H = \Alt_{6}$ has irreducible modules of dimension $5$ and $8$, giving embeddings into $B_{2}$ and $D_{4}$. Since $\Alt_{6}$ has no nontrivial irreducible modules of dimension $4$ or less, it cannot occur as a subgroup of a (simply connected) subgroup of type $B_{2}$ or $A_{3}$ in $G = F_{4}$. Hence such a subgroup $\Alt_{6}$ of $G$ must be Lie primitive in a subgroup of type $D_{4}$.

\begin{table}[htbp]
\small \onehalfspacing
\caption{Types of subgroup necessarily in some $X$ satisfying Proposition \ref{prop:primstab}}
\begin{tabular}{c|c}
\hline $G$ & $H$ \\ \hline
$F_{4}$ & $\Alt_{7}$, $\Alt_{9-10}$, $M_{11}$, $J_{1}$, $J_{2}$, \\
& $\Alt_{6}$ $(p = 5)$, $L_{2}(17)$ $(p = 2)$, $U_{3}(3)$ $(p \neq 2,3,7)$ \\ \hline
$E_{6}$ & $\Alt_{9-12}$, $M_{22}$, $L_{2}(25)$, $L_{2}(27)$, $L_{4}(3)$, $U_{4}(3)$, ${}^{3}D_{4}(2)$, \\
& $M_{11}$ $(p \neq 3,5)$, $M_{12}$ $(p = 2)$, $L_{2}(11)$ $(p = 5)$, $L_{2}(17)$ $(p = 2)$, $L_{3}(3)$ $(p = 2)$ \\ \hline
$E_{7}$ & $\Alt_{11-13}$, $M_{11}$, $L_{2}(25)$, $L_{3}(3)$, $L_{4}(3)$, ${}^{3}D_{4}(2)$, ${}^{2}F_{4}(2)'$,\\
& $\Alt_{9-10}$ $(p = 2)$, $M_{12}$ $(p \neq 5)$ \\ \hline
$E_{8}$ & $\Alt_{12-15}$, $\Alt_{17}$, $M_{12}$, $L_{2}(27)$, $L_{2}(37)$, $L_{4}(3)$, $U_{3}(8)$, $\Omega_{8}^{+}(2)$, $G_{2}(3)$, \\
& $\Alt_{11}$ $(p = 2)$, $M_{11}$ $(p \neq 3,11)$, $^{2}F_{4}(2)'$ $(p \neq 3)$ \\ \hline
\end{tabular}
\label{tab:primstab}
\end{table}

\begin{lemma} \label{lem:nstab}
Let $S$ be a finite simple subgroup of a simple algebraic group $G$. If $S$ is contained in a semisimple subgroup $X$ of $G$, such that the following conditions all hold:
\begin{itemize}
\item[\textup{(i)}] $X$ is $G$-conjugate to $X^{\sigma}$ for all $\sigma \in N_{\textup{Aut}(G)}(S)$,
\item[\textup{(ii)}] If $g \in G$ and $S^{g} \le X$, then $S^{g} = S^{x}$ for some $x \in N_{G}(X)$,
\item[\textup{(iii)}] $N_{G}(S) \le N_{G}(X)$,
\end{itemize}
then $X$ is $N_{\textup{Aut}(G)}(S)$-stable.
\end{lemma}

\proof If $\sigma \in N_{\textup{Aut}(G)}(S)$, then $X^{\sigma} = X^{g}$ for some $g \in G$, by (i). Thus $S^{g^{-1}} \le X$, and $S^{g^{-1}} = S^{x}$ for some $x \in N_{G}(X)$, by (ii). Then $xg \in N_{G}(S) \le N_{G}(X)$, and $X^{\sigma} = X^{g} = X^{xg} = X$, as required. \qed

In view of the above, the following proposition now proves the conclusion of Theorem \ref{THM:MAIN} for those triples $(G,H,p)$ not postponed until Section \ref{sec:special}.

\begin{proposition} \label{prop:gcrstab}
Let $G$ be an adjoint exceptional simple algebraic group in characteristic $p$, and let $H \notin \textup{Lie}(p)$ be a non-abelian finite simple group, and let $S \cong H$ be $G$-cr subgroup of $G$. Let $\tilde{S}$ denote a minimal preimage of $S$ in the simply connected cover $\tilde{G}$ of $G$, and let $\tilde{H}$ denote the isomorphism type of $\tilde{S}$.

If $(G,\tilde{H},p)$ appears in one of the tables \ref{tab:f4thm} to \ref{tab:e8thm}, then $\tilde{S}$ is an $X$-irreducible subgroup of some subgroup $X$ listed there. Moreover one of the following holds:
\begin{itemize}
\item[\textup{(a)}] $S$ is Lie primitive in $X$ and Proposition \ref{prop:primstab} applies;
\item[\textup{(b)}] Lemma \ref{lem:nstab} applies to $X$;
\item[\textup{(c)}] $X$ has a submodule $W$ on $V$, such that every $S$-submodule of $W$ is $X$-stable, and the collection $\mathcal{M}$ of such submodules has the necessary form to apply Corollary \ref{cor:phistab}.
\end{itemize}
Thus $S$ lies in a proper connected $N_{\textup{Aut}(G)}(S)$-stable subgroup of $G$.
\end{proposition}

In Tables \ref{tab:f4thm}--\ref{tab:e8thm} we use the following notation. When $X$ is contained in a semisimple subgroup $Y$ with classical factors, we write `$X < Y$ via $\lambda$' to indicate the action of $X$ on the natural module for $Y$. If the simple factors of $Y$ are $Y_{i}$ $(i = 1,\ldots,r)$ we write $(V_1,V_2,\ldots,V_{r})$ to indicate a tensor product $V_{1} \otimes \ldots \otimes V_{r}$, where $V_{i}$ is a module for $Y_{i}$. Also $V^{[r]}$ denotes the conjugate of the module $V$ by a $p^{r}$-power Frobenius morphism. Finally, we write `$X$ (fpf)' to indicate that $S$ can be assumed to fix no nonzero vectors on any nontrivial $X$-composition factor of $V$ (otherwise $S$ lies in some other listed subgroup).

\begin{center}
\small \onehalfspacing
\begin{longtable}{c|c|c|c|c}
\caption{$G$-cr overgroups $X$: $G$ of type $F_4$, $p \neq 2$, $V = V_{G}(\lambda_4)$, $\delta = \delta_{p,3}$} \\
$\tilde{H}$ & $p$ & $X$ & $V \downarrow X$ & $W$ \\ \hline
$\Alt_{5}$ & $\ge 13$ & $A_{1}$ max $F_4$ & Lemma \ref{lem:nstab} applies \\
& $\neq 2,3,5$ & $A_{1} < B_{4}$ via $(1 \otimes 1^{[1]}) \oplus 4$ & $0 \oplus (1 \otimes 1^{[1]}) \oplus 4^{2} \oplus 2 \oplus (1^{[1]} \otimes 3)$ & $0$ \\
& & $A_{1} < B_{4}$ via $0 \oplus 2 \oplus 4$ & $0^{2} \oplus 2^{3} \oplus 4^{3} \oplus (1 \otimes 3)^{2}$ & $0^{2}$ \\
& $\neq 2,5$ & $A_{1}C_{3}$ (fpf) & Lemma \ref{lem:nstab} applies \\
& $\neq 2,3,5$ & $A_{1} < A_{2}^{2}$ via $(2,2)$ & $0^{2} \oplus 2^{3} \oplus 4^{3}$ & $2^{3}$ \\
& $p = 3$ & $A_{1} < A_{2}^{2}$ via $(2,2)$ & $(0|4|0)^{2} \oplus 2^{3} \oplus 4$ & $2^{3}$ \\
& $p \neq 2,5$ & $A_{1} < A_{2}^{2}$ via $(2^{[1]},2)$, & $(2^{[1]} \otimes 2)^{2} \oplus 2 \oplus V_{X}(4)$ & $2$ \\
& & $A_{1} < B_{3}$ via $(1 \otimes 1^{[1]}) \oplus 2$ & $0^{5-\delta} \oplus 2^{3} \oplus (1 \otimes 1^{[1]})^{3}$ & $V$ \\
& & $A_{1} < B_{3}$ via $0 \oplus 2 \oplus 2^{[1]}$ & $0^{4 - \delta} \oplus 2 \oplus 2^{[1]} \oplus (1 \otimes 1^{[1]})^{4}$ & $V$\\
& & $A_{1} <$ long $A_{2}$ Levi via $2$ & $0^{8 - \delta} \oplus 2^{6}$ & $V$ \\
& & $A_{1} <$ short $A_{2}$ Levi via $2$ & $2^{7} \oplus V_{A_1}(4)$ & $2^{7}$ \\
$L_{2}(7)$ & $3$ & $A_{2} < A_{2}\tilde{A}_{2}$ via $(10,10)$ & $(0|11|0)^{2} \oplus V_{X}(11)$ & $V$ \\
& & $A_{2} < A_{2}\tilde{A}_{2}$ via $(10,01)$ & $20 \oplus 02 \oplus 10 \oplus 01 \oplus V_{X}(11)$ & $V$ \\
& & $A_{2} < B_{3}$ via $V_{X}(11)$ & $0^{4} \oplus V_{X}(11)^{3}$ & $V$ \\
& & $A_{2}$ Levi & Prop.\ \ref{prop:primstab} applies
\label{tab:f4thm}
\end{longtable}

\begin{longtable}{c|c|c|c|c}
\caption{$G$-cr overgroups $X$: $G$ of type $E_{6}$, $V = V_G(\lambda_2)$} \\
$\tilde{H}$ & $p$ & $X$ & $V \downarrow X$ & $W$ \\ \hline
$3\cdot\Alt_{7}$ & $\neq 3,5$ & $A_5$ & Prop.\ \ref{prop:primstab} applies \\
$\Alt_{7}$ & $\neq 2,3,5,7$ & $A_3 < A_5$ & $0^{3} \oplus (\lambda_1+\lambda_3) \oplus 2\lambda_1^{2} \oplus 2\lambda_2 \oplus 2\lambda_3^{2}$ & $0^{3}$
\label{tab:e6thml2}
\end{longtable}

\begin{longtable}{c|c|c|c|c}
\caption{$G$-cr overgroups $X$: $G$ of type $E_6$, $V = V_G(\lambda_1) \oplus V_G(\lambda_6)$}\\
$\tilde{H}$ & $p$ & $X$ & $V \downarrow X$ & $W$ \\ \hline
$\Alt_{7}$ & $7$ & $D_5$ & Prop.\ \ref{prop:primstab} applies \\
& & $B_{2} < A_{4}$ & $0^{4} \oplus \lambda_1^{6} \oplus 2\lambda_{2}^{2}$ & $V$ \\
& $2$ & $A_3$ Levi & Prop.\ \ref{prop:primstab} applies \\
& & $A_3 < A_5$ & $\lambda_2^{4} \oplus (0^{2}/V_X(\lambda_1+\lambda_3)^{2})$ & $\lambda_2^{4}$ \\
$\Alt_{5}$ & $\neq 2,3,5$ & $A_{1} < A_{2}^{3}$ via $(2,2,2)$ & $0^{6} \oplus 2^{6} \oplus 4^{6}$ & $V$ \\
& & $A_{1} < A_{2}^{3}$ via $(2,2,2^{[1]})$ & $0^{2} \oplus 2^{2} \oplus 4^{2} \oplus (2 \otimes 2^{[1]})^{4}$ & $0^{2}$ \\
& & $A_{1} < A_{1}A_{5}$ via $(1,5)$ & $0^{2} \oplus (4/6/W_{X}(8)/W_{X}(4))^{2}$ & $0^{2}$ \\
& & $F_{4}$ or $A_{1}C_{3} < F_{4}$ (fpf) & $(V_{F_4}(\lambda_4)\downarrow X)^{2} \oplus 0^{2}$ & $0^{2}$ \\
& & $A_{1} < B_{4}$ via $4 \oplus (1 \otimes 1^{[1]})$ & $0^{4} \oplus 2^{2} \oplus 4^{4}$ & $0^{4}$ \\
& & & ${} \oplus (1 \otimes 1^{[1]})^{2} \oplus (1^{[1]} \otimes 3)^{2}$ & \\
& & $A_{1} < B_{3}$ via $2 \oplus (1 \otimes 1^{[1]})$ & $0^{12} \oplus 2^{6} \oplus (1 \otimes 1^{[1]})^{6}$ & $V$ \\
& & $A_{1} < D_{4}$ via $2 \oplus 4$ & $0^{6} \oplus 2^{6} \oplus 4^{6}$ & $V$ \\
& & $A_{1} < A_{4}$ via $4$ & $0^{4} \oplus 2^{2} \oplus 4^{6} \oplus 6^{2}$ & $0^{4}$ \\
& & $A_{1} < A_{3}$ Levi via $1 \otimes 1^{[1]}$ & $0^{10} \oplus 2^{2} \oplus (2^{[1]})^{2}$ & $V$ \\
& & & ${}\oplus (1 \otimes 1^{[1]})^{8}$ & \\
& & $A_{1} < A_{2}$ Levi via $2$ & $0^{18} \oplus 2^{12}$ & $V$ \\
& & $A_{1} < A_{2}^{2}$ Levi via $(2,2)$ & $0^{2} \oplus 2^{14} \oplus 4^{2}$ & $V$ \\
& & $A_{1} < A_{2}^{2}$ Levi via $(2,2^{[1]})$ & $(2 \oplus 2^{[1]})^{6} \oplus (2 \otimes 2^{[1]})^{2}$ & $(2 \oplus 2^{[1]})^{6}$ \\
$L_2(7)$ & $3$ & $A_{2} < A_{2}^{3}$ via $(10,10,10)$ & $(0|11|0)^{6}$ & $V$ \\
& & $A_{2} < A_{2}^{3}$ via $(10,10,01)$ & $(0|11|0)^{2} \oplus 10^{2} \oplus 01^{2}$ & $V$ \\
& & & ${} \oplus 20^{2} \oplus 02^{2}$ & \\
& & $A_{2} < A_{2}^{2}$ Levi via $(10,10)$ & $(0|11|0)^{2} \oplus 10^{6} \oplus 01^{6}$ & $V$ \\
& & $A_{2} < A_{2}^{2}$ Levi via $(10,01)$ & $20 \oplus 02 \oplus 10^{7} \oplus 01^{7}$ & $V$ \\
& & $A_{2}$ Levi & Prop.\ \ref{prop:primstab} applies \\
& & $A_{2} < G_{2} < D_{4}$ & Prop.\ \ref{prop:primstab} applies \\
& & $A_{3} < A_{5}$ & $\lambda_2^{2} \oplus (\lambda_1 + \lambda_3)$ & $\lambda_2^{2}$ \\
$L_{2}(8)$ & $7$ & $G_{2}$ max $F_{4}$ & $0^{2} \oplus V_{X}(20)^{2}$ & $0^{2}$ \\
& & $A_{2} < D_{4}$ via $V_{X}(11)$ & $0^{12} \oplus V_{X}(11)^{6}$ & $V$ \\
& & $G_{2}$ or $B_{3} < D_{4}$ & Prop.\ \ref{prop:primstab} applies \\
& & $D_{4}$, $F_{4}$ & Prop.\ \ref{prop:primstab} applies \\
$L_{2}(13)$ & $7$ & $G_{2}$ max $F_{4}$ & $0^{2} \oplus V_{X}(20)^{2}$ & $0^{2}$ \\
& & $G_{2} < D_{4}$ & Prop.\ \ref{prop:primstab} applies \\
& & $F_{4}$ & Prop.\ \ref{prop:primstab} applies \\
$U_{3}(3)$ & $7$ & $G_{2}$ max\ $F_{4}$ & $0^{2} \oplus V_{X}(20)	^{2}$ & $0^{2}$ \\
& & $C_{3} < A_{5}$ & $0 \oplus \lambda_1^{2} \oplus \lambda_2$ & $V$ \\
& & $G_{2} < D_{4}$ & Prop.\ \ref{prop:primstab} applies \\
& & $F_{4}$ & Prop.\ \ref{prop:primstab} applies \\
$U_4(2)$ & $\neq 2,3$ & $A_3 < A_5$ & $\lambda_2^{4} \oplus (\lambda_1 + \lambda_3)^{2}$ & $V$ \\
 & & $A_4$ & Prop.\ \ref{prop:primstab} applies
\label{tab:e6thml1l6}
\end{longtable}

\begin{longtable}{c|c|c|c|c}
\caption{$G$-cr overgroups $X$: $G$ of type $E_7$, $V = V_G(\lambda_7)$}\\
$\tilde{H}$ & $p$ & $X$ & $V \downarrow X$ & $W$ \\ \hline
$\Alt_{8}$ & $\neq 2,3,5$ & $B_3 < A_6$ & $\lambda_1^{2} \oplus \lambda_2^{2}$ & $V$ \\
$\Alt_{7}$ & $\neq 2,5,7$ & $A_{3} < A_{5}$ & $2\lambda_1 \oplus 2\lambda_3 \oplus \lambda_2^{6}$ & $V$ \\
& & $A_{3} < A_{5}'$ & $0^{2} \oplus \lambda_2^{4} \oplus (\lambda_1 + \lambda_4)^{2}$ & $V$ \\
& $7$ & $B_{2} < A_{4}$ & $0^{6} \oplus 10^{6} \oplus 02^{2}$ & $V$ \\
& & $D_{5}$ & Prop.\ \ref{prop:primstab} applies \\
& $2$ & $A_{3} < A_{5}$ & $\lambda_2^{6} \oplus (\lambda_2|(V_{X}(2\lambda_1) \oplus V_{X}(2\lambda_3))|\lambda_2)$ & $\lambda_2^{7}$ \\ 
& & $A_{3} < A_{5}'$ & $(0^{4}/V_{X}(\lambda_1 + \lambda_3)^{2}) \oplus \lambda_2^{4}$ & $\lambda_{2}^{4}$ \\
& & $A_{3}$ Levi & Prop.\ \ref{prop:primstab} applies \\
$2 \cdot \Alt_{7}$ & $3$ & $C_{3} < A_{5}$ & $\lambda_1^{7} \oplus \lambda_3$ & $\lambda_1^{7}$ \\
$J_1$ & $11$ & $G_2 < A_6$ & $10^{4} \oplus 01^{2}$ & $V$ \\
& & $G_{2}$ max $E_{6}$ & $0^{2} \oplus 20^{2}$ & $V$ \\
& & $G_{2} < D_{4}$ & Prop.\ \ref{prop:primstab} applies \\
& & $E_{6}$ & Prop.\ \ref{prop:primstab} applies \\
$2 \cdot J_2$ & $\neq 2$ & $C_3 < A_5$ & $\lambda_1^{7} \oplus \lambda_3$ & $V$ \\
$L_{2}(8)$ & $\neq 2,3,7$ & $G_{2}$ max $E_{6}$ & $0^{2} \oplus 20^{2}$ & $0^{2}$ \\
& & $F_{4}$, $E_{6}$ & Prop.\ \ref{prop:primstab} applies \\
$L_2(17)$ & $2$ & $C_{4} < F_{4}$ & Prop.\ \ref{prop:primstab} applies \\
& & $C_{4} < A_{7}$ & $(0|V_{X}(\lambda_2)|0)^{2}$ & $0^{2}$ \\
& & $B_{4} < D_{5}$ & Prop.\ \ref{prop:primstab} applies \\
& $\neq 2,17$ & $B_{4}$ & Prop.\ \ref{prop:primstab} applies \\
& & $C_{4}$ & $0^{2} \oplus \lambda_2^{2}$ & $0^{2}$ \\
& & $F_{4}$, $E_{6}$ & Prop.\ \ref{prop:primstab} applies \\
$U_4(2)$ & $\neq 2,3$ & $A_4$ & Prop.\ \ref{prop:primstab} applies \\
& & $A_3 < A_3 A_3$ & $0^{2} \oplus (\lambda_1 + \lambda_3)^{2} \oplus \lambda_2^{4}$ & $V$ \\
& & (2 classes) & $2\lambda_1 \oplus 2\lambda_3 \oplus \lambda_2^{6}$ & $V$ \\
$Sp_{6}(2)$ & $\neq 2$ & $B_3 < A_6$ & $\lambda_1^{2} \oplus \lambda_2$ & $V$
\label{tab:e7thml7}
\end{longtable}

\begin{longtable}{c|c|c|c|c}
\caption{$G$-cr overgroups $X$: $G$ of type $E_7$, $V = V_G(\lambda_1)$}\\
$\tilde{H}$ & $p$ & $X$ & $V \downarrow X$ & $W$ \\ \hline
$L_{2}(8)$ & $\neq 2,3,7$ & $B_{4}$, $F_{4}$ & Prop.\ \ref{prop:primstab} applies \\
& & $G_{2}$ or $B_{3} < D_{4}$ Levi & Prop.\ \ref{prop:primstab} applies \\
& & $B_{3} < A_{6}$ & $0 \oplus \lambda_1^{2} \oplus 2\lambda_1 \oplus \lambda_2 \oplus 2\lambda_3^{2}$ & $0$ \\
& & $G_{2} < A_{6}$ & $0^{3} \oplus 10^{5} \oplus 01 \oplus 20^{3}$ & $0^{3}$
\label{tab:e7thml1}
\end{longtable}

\begin{longtable}{c|c|c|c|c}
\caption{$G$-cr overgroups $X$: $G$ of type $E_8$, $V = L(G)$}\\
$\tilde{H}$ & $p$ & $X$ & $L(G) \downarrow X$ & $W$ \\ \hline
$\Alt_{10}$ & $3$ & $B_4 < D_{8}$ & $0 \oplus V_{X}(2\lambda_1) \oplus \lambda_2 \oplus \lambda_3^{2}$ & $\lambda_3^{2}$ \\
& 2 & $B_{4} < D_{8}$ via $\lambda_4$ & $(W_X(\lambda_2)/W_X(\lambda_3)) \oplus (\lambda_1 + \lambda_4)$ & $\lambda_1 + \lambda_4$ \\
& & $B_{4} < D_{5}$ & Prop.\ \ref{prop:primstab} applies \\
& & $C_{4} < A_{7}$ & $(\lambda_1 \otimes \lambda_1) \oplus \lambda_1^{4} \oplus (0|V_{X}(\lambda_2)|0)^{2} \oplus \lambda_3^{2}$ & $\lambda_3^{2}$ \\ 
& & $C_{4} < F_{4}$ & Prop.\ \ref{prop:primstab} applies \\
$\Alt_{9}$ & $\neq 2,3$ & $D_4 < A_7$ Levi & $0 \oplus \lambda_1^{2} \oplus \lambda_2^{3} \oplus 2\lambda_1 \oplus (\lambda_3 + \lambda_4)^{2}$ & $(\lambda_3 + \lambda_4)^{2}$ \\
$\Alt_{8}$ & $\neq 2$ & $B_{3} < A_{6}$ & $0^{4} \oplus \lambda_1^{6} \oplus \lambda_2^{5} \oplus 2\lambda_1 \oplus 2\lambda_3^{2}$ & $\lambda_2^{5}$ \\
& $3$, $5$ & $E_{7}$ & Prop.\ \ref{prop:primstab} applies \\
$J_1$ & $11$ & $G_2 < D_4$ Levi & Prop.\ \ref{prop:primstab} applies \\
& & $E_{6}$ & Prop.\ \ref{prop:primstab} applies \\
& & $G_2 < A_6$ & $0^{6} \oplus 10^{13} \oplus 01^{5} \oplus 20^{3}$ & $V$ \\
& & $G_2 < D_7$ via $01$ & $0 \oplus 01^{3} \oplus 30 \oplus 11^{2}$ & $V$ \\
& & $G_2 \textup{ max } E_6$ & $0^{8} \oplus 01 \oplus 20^{6} \oplus 11$ & $V$ \\
$J_{2}$ & $2$ & $E_{6}$ (fpf) & $0^{8} \oplus \lambda_2 \oplus \lambda_1^{3} \oplus \lambda_6^{3}$ & $0^{8}$ \\
& & $G_{2}$ max $E_{6}$ & $0^{8} \oplus 11 \oplus 01$ & $0^{8}$ \\
& & & ${} \oplus (01|20|00|10)^{3} \oplus (10|00|20|01)^{3}$ & \\
& & $E_{7}$ (fpf) & $0^{2} \oplus (0|V_{X}(\lambda_1)|0) \oplus \lambda_7^{2} \oplus \lambda_1$ & $0^{3}$ \\
& & $G_{2} < D_{7}$ via $01$ & $0 \oplus 11^{2} \oplus 01^{2} \oplus (30/01)$ & $11^{2}$ \\
& & $G_{2} < D_{4}$ Levi & Prop.\ \ref{prop:primstab} applies \\
& & $G_{2} < A_{5}$ & $0^{16} \oplus ((0 \oplus 01)|V_{X}(20)|(0 \oplus 01)) \oplus $ & $0^{17}$ \\
& & & ${} \oplus 01^{6} \oplus (V_{X}(10)|0|V_{X}(20)|0|V_{X}(10))^{2}$ & \\
$Sp_{6}(2)$ & $\neq 2$ & $B_3 < A_6$ & $(0^{4}/W_X(2\lambda_1)) \oplus \lambda_1^{6} \oplus \lambda_2^{5} \oplus 2\lambda_3^{2}$ & $\lambda_1^{6}$ \\
$U_{3}(3)$ & $\neq 2,3,7$ & $A_{6}$ & Prop.\ \ref{prop:primstab} applies \\
& & $E_{7}$ (fpf) & $0^{3} \oplus \lambda_1 \oplus \lambda_7^{2}$ & $0^{3}$ \\
& & $E_{6}$ & Prop.\ \ref{prop:primstab} applies \\
& & $G_{2}$ max $E_{6}$ & $0^{8} \oplus 20^{6} \oplus 01 \oplus 11$ & $01$ \\
& & $C_{3} < A_{5}$ & $0^{17} \oplus 2\lambda_1 \oplus \lambda_1^{14} \oplus \lambda_2^{7} \oplus \lambda_3^{2}$ & $0^{17}$ \\
& & $C_{3} < D_{7}$ & $0 \oplus 2\lambda_1 \oplus \lambda_2^{2} \oplus (\lambda_1 + \lambda_3) \oplus (\lambda_2 + \lambda_3)^{2}$ & $0$ \\
& & $G_{2} < D_{7}$ via $01$ & $0 \oplus 01^{3} \oplus 30 \oplus 11^{2}$ & $0$ \\
& & $G_{2} < D_{4}$ & Prop.\ \ref{prop:primstab} applies \\
& & $G_{2} < A_{6}$ & $0^{6} \oplus 10^{13} \oplus 01^{5} \oplus 20^{3}$ & $V$
\label{tab:e8thm}
\end{longtable}
\end{center}

\proof This is straightforward to verify on a case-by-case basis. We now give a general outline, and then illustrate by giving details in the most involved cases.

The groups $X$ are representatives of the $G$-cr semisimple subgroups which contain a copy of $\tilde{S}$ `minimally', in the sense that $\tilde{S}$ centralises no simple factor of $X$ and does not lie in a diagonal subgroup when $X$ has two or more isomorphic factors. Moreover if $X$ has a classical factor, then $\tilde{S}$ does not lie in the subgroup of $X$ corresponding to the stabiliser of a direct-sum decomposition, tensor-product decomposition or another bilinear or quadratic form on the natural module. If $X$ is exceptional, then the image of $\tilde{S}$ must correspond to a character marked `\possprim' in Chapter \ref{chap:thetables}.

Since $X$ is $G$-cr, it is $L'$-irreducible for some Levi subgroup $L$ of $G$. The structure of $V \downarrow L'$ is easy to determine from the composition factors given in \cite{MR1329942}*{Tables 8.1-8.7}, and Lemma \ref{lem:tilting}. This is enough information to determine $V \downarrow X$ in each case here. This restriction then limits the possible feasible characters of $\tilde{S}$ on $V$, and identifying the submodule $W$ is usually straightforward. In many cases, Proposition \ref{prop:newsamesubs} applies to the whole of $V$, and we take $W = V$.

The most complicated constructions of $X$ and $W$ are as follows.

\subsection*{Case: $(G,H,p) = (F_4,L_{2}(7),3)$ or $(E_{6},L_{2}(7),3)$} Here, $H$ has irreducible modules of dimension $3$, $6$ and $7$, giving irreducible embeddings into groups of type $A_{2}$, $A_{3}$ and $B_{3}$. Moreover an embedding of $H$ into $B_{3}$ lies in a subgroup $G_{2}$. Further, an embedding of $H$ into $G_{2}$ stabilises a vector on $V_{G_2}(01)$ and therefore lies in a subgroup $A_{2}$; this holds since $\bigwedge^{2}V_{G_2}(10) = V_{G_2}(10)^{2}/V_{G_2}(01)$ while the exterior square of the $7$-dimensional $H$-module has composition factors $7^{2}/3/3^{*}/1$, and $3$ and $3^{*}$ are projective since their dimension is the largest power of $3$ dividing $|H|$.

First suppose that $G = F_{4}$. Inspecting \cite{MR1329942}*{Table 8.4} we see that each simple subgroup $A_{3}$ of $G$ is simply connected. Since $S$ has no irreducible $4$-dimensional modules, $S$ cannot be contained irreducibly in such a subgroup $A_{3}$.

If $S$ is $G$-irreducible, then from Proposition \ref{prop:algorithm} we know that $S$ fixes a point on $L(G)$, and therefore lies in a proper subsystem subgroup of $G$. Now $S$ has no irreducible embeddings into a subgroup of type $B_{4}$ or $A_{1}C_{3}$ (Lemma \ref{lem:classicalparabs}), hence $S$ lies in a subgroup $A_{2}A_{2}$. This acts on $V = V_{G}(\lambda_4)$ as $(10,10) + (01,01) + (0,V_{A_2}(11))$. Since the $3$-dimensional irreducible $S$-modules are projective, so is any tensor product of them, hence the $A_{2}A_{2}$-modules $(01,10)$ and $(10,01)$ restrict to $S$ either as $3 + 6$ or $3^{*} + 6$, or as a uniserial projective module $1|7|1$. In any case, the embedding of $S$ factors through a diagonal subgroup $A_{2}$ and has a unique $7$-dimensional irreducible submodule on $V_{G}(\lambda_4)$, which is the restriction of the unique $A_{2}$-submodule $V_{A_2}(11)$.

If instead $S$ lies in a proper Levi subgroup $L$ of $G$, then $L'$ has type $A_{2}$ or $B_{3}$. In the former case, $S$ is Lie primitive in $X = L'$ and Corollary \ref{prop:primstab} applies. In the latter case, since $S$ is $L'$-irreducible and is not contained irreducibly in a subgroup $A_{3}$, it follows that $S$ is irreducible on $V_{B_3}(\lambda_1)$. Since the unique $7$-dimensional irreducible $S$-module is a section of $3 \otimes 3^{*}$, $S$ lies in an irreducible subgroup $X = A_{2}$ of $B_{3}$ acting as $V_{A_2}(11)$ on the natural module. The given action of $X$ follows and Proposition \ref{prop:newsamesubs} applies to the action of $S$ and $X$ on $V$.

Now suppose $G = E_{6}$ and $V = V_{G}(\lambda_1) \oplus V_{G}(\lambda_6)$. Then similar reasoning to the above applies, where we also note that a maximal connected subgroup $A_{2}^{3}$ acts on $V$ as a direct sum of tensor products of two $3$-dimensional modules (cf.\ \cite{MR1329942}*{Proposition 2.3}), so the same reasoning as for $S < A_{2}^{2}$ shows that $S$ lies in a diagonal subgroup $A_{2}$ stabilising every $S$-submodule of $V$.

\subsection*{Case: $(G,H) = (F_4,\Alt_{5})$} Here $p \neq 2,5$. In this case, every irreducible module for $S$ and its double cover $2 \cdot S \cong SL_{2}(5)$ is obtained as a symmetric power or tensor product of $2$-dimensional $SL_{2}(5)$-modules. Hence an embedding of $S$ into a classical simple algebraic group factors through an embedding of an adjoint simple subgroup of type $A_{1}$. Furthermore if $S$ lies in a subgroup of type $G_{2}$ then it stabilises a $3$- or $4$-dimensional subspace of the $7$-dimensional module, hence lies in a subgroup $A_{1}A_{1}$ or $A_{2}$ (cf.\ \cite{MR898346}*{Theorem 8}), and then in a further subgroup $A_{1}$.

We now break the proof up into several cases: (1) $S$ lies in no proper subsystem subgroup of $G$; (2) $S$ is $G$-irreducible and lies in a proper subsystem subgroup of $G$; (3) $S$ lies in a proper Levi subgroup of $G$. Furthermore we divide case (2) into the sub-cases where $S$ lies in a subgroup: (2a) $B_{4}$; (2b) $A_{1}C_{3}$; or (2c) $A_{2}A_{2}$.

(1) In this case $C_{G}(S) = 1$ and $S$ cannot fix a nonzero vector on $L(G)$ by Lemma \ref{lem:properembed}, and cannot fix a nonzero vector on $V_{\textup{min}}$ since the corresponding stabiliser has dimension at least $52 - 26$ and is either $G$-reducible or a subsystem subgroup $B_{4}$ or $D_{4}$. Thus $p \neq 3$ and $S$ corresponds to Case 1) of Table \ref{a5f4p0}. If $S$ is contained in a maximal subgroup $A_{1}G_{2}$ or $G_{2}$ $(p = 7)$, then the image of $S$ in $G_{2}$ lies in a subgroup $A_{1}A_{1}$ or $A_{2}$, hence centralises a non-trivial semisimple element of $G$, contradicting the fact that $S$ lies in no subsystem subgroup. Therefore $p \ge 13$ and $S$ lies in a maximal subgroup $X$ of type $A_{1}$. This is unique up to conjugacy in $G$, hence condition (i) of Lemma \ref{lem:nstab} holds. Moreover $SO_{3}(K)$ has a unique subgroup $\Alt_{5}$ up to conjugacy, hence condition (ii) of Lemma \ref{lem:nstab} also holds. Finally, the two $3$-dimensional irreducible $S$-modules are interchanged by an outer automorphism of $S$. Since these occur with different multiplicities as composition factors of $L(G)$, it follows that $G$ does not induce an outer automorphism on $S$, hence $N_{G}(S) = SC_{G}(S) = S$. Thus condition (iii) of Lemma \ref{lem:nstab} also holds, and $X$ is $N_{\textup{Aut}(G)}(S)$-stable.

(2a) If $S$ lies in $X = B_{4}$ and is $X$-irreducible then by Lemma \ref{lem:classicalparabs}, $S$ acts on $V_{B_4}(\lambda_1)$ with composition factor dimensions $1$, $3$ and $5$ or $4$ and $5$ (hence $p \neq 3$). The image of $S$ therefore lies in a subgroup $A_{1}$ acting as $0 \oplus 2 \oplus 4$ or $(1 \otimes 1^{[1]}) \oplus 4$. The restriction of $V_{G}(\lambda_4)$ to these respective subgroups follows from the restriction $V_{G}(\lambda_4) \downarrow B_{4} = 0 \oplus \lambda_1 \oplus \lambda_4$ and Lemma \ref{lem:spins}. Each nontrivial summand for such a subgroup $A_{1}$ restricts to $S$ with no trivial composition factors, and it follows that each trivial $S$-submodule of $V_{G}(\lambda_4)$ is a trivial submodule also for the relevant subgroup $A_{1}$.

(2b) Here $S$ lies in a subgroup $A_{1}C_{3}$; we show that Lemma \ref{lem:nstab} holds in this case. Firstly, since $G$ has a unique class of subgroups $A_{1}C_{3}$, condition (i) is immediate. Next, since the double cover $2\cdot S$ has a unique $6$-dimensional symplectic module and $SL_{2}(5)$ has a unique $2$-dimensional irreducible module up to conjugation by an outer automorphism, it follows that $A_{1}C_{3}$ has a unique class of subgroups $\Alt_{5}$ which centralise neither factor. Thus condition (ii) of Lemma \ref{lem:nstab} holds.

Finally, since $S$ lies in no Levi subgroup of $G$, if we assume that $S$ lies in no maximal subgroup $A_{2}A_{2}$ or $B_{4}$ of $G$ then the only nontrivial elements of $C_{G}(S)$ are involutions which are $G$-conjugate to the central involution $t$ of $A_{1}C_{3}$. Therefore $C_{G}(S)$ is commutative, and is therefore contained in $C_{G}(t) = A_{1}C_{3}$. Now $V_{G}(\lambda_4) \downarrow A_{1}C_{3} = (1,\lambda_1) \oplus (0,\lambda_2)$ restricts to $S$ as $(2 \otimes 6) \oplus \bigwedge^{2}(6)$, where `$2$' and `$6$' are irreducible modules for the double cover $2\cdot S$. The first summand has a unique $3$-dimensional composition factor, while the second has none. Thus the two $3$-dimensional $S$-modules occur with differing multiplicities, hence $N_{G}(S)$ does not induce an outer automorphism on $S$. This shows that $N_{G}(S) = S C_{G}(S) \le A_{1}C_{3} = N_{G}(A_{1}C_{3})$, so condition (iii) of Lemma \ref{lem:nstab} holds.

(2c) Here $S$ lies in a subgroup $A_{2}A_{2}$, acting irreducibly on the natural module for each factor. Since the two $3$-dimensional $S$-modules are the symmetric squares of the $2$-dimensional modules for $2 \cdot S$, and since these are conjugate under a Frobenius morphism of $SL_{2}(K)$, the image of $S$ lies in a diagonal subgroup $A_{1}$ acting as $2$ or $2^{[1]}$ on the natural module for each $A_{2}$ factor, giving two subgroups of type $A_{1}$, one of which contains $S$. Now, $V_{G}(\lambda_4) \downarrow A_{2}A_{2} = (10,01) + (01,10) + (0,V_{A_2}(11))$. The restrictions $V_{G}(\lambda_4) \downarrow A_{1}$ given in Table \ref{tab:f4thm} follow easily. Since $2^{[1]} \otimes 2$ restricts to $S$ as $3_{a} \otimes 3_{b} = 1/4^{2}$ $(p = 3)$ or $4/5$ $(p \neq 3)$, it follows that every $3$-dimensional $S$-submodule of $V_{G}(\lambda_4)$ is preserved by the subgroup $A_{1}$.

(3) Let $L$ be minimal among Levi subgroups of $G$ containing $S$. Since $G$ is simply connected, so is $L'$. Since $S$ itself has no nontrivial $2$-dimensional modules or irreducible symplectic modules, it follows that $L'$ has no factor $A_{1}$, $B_{2}$ or $C_{3}$. Thus $L'$ has type $B_{3}$ or $A_{2}$. If $L' = B_{3}$ then $S$ acts as $3_{a} + 4$ or $1 + 3_{a} + 3_{b}$ on the natural $L'$-module, hence lies in a subgroup $A_{1} < B_{3}$ via $(1 \otimes 1^{[1]}) + 2$ or $0 + 2 + 2^{[1]}$, as in Table \ref{tab:f4thm}. The restrictions for $X$ follow from $V_{G}(\lambda_4) \downarrow B_{3} = 0^{3} + \lambda_1 + \lambda_3^{2}$ and Lemma \ref{lem:spins}. If $L' = A_{2}$ then $S$ acts irreducibly on the natural $L'$-module. Since there are two $G$-classes of Levi subgroups $A_{2}$, this gives the final subgroups $A_{1} < A_{2}$ for $\Alt_{5}$ in Table \ref{tab:f4thm}.

\subsection*{Case: $(G,H) = (E_{6},\Alt_{5})$} Here $p \neq 2,3,5$. As argued for $G = F_{4}$ above, $S$ lies in an adjoint subgroup $A_{1}$ of $G$, and if $L'$ is minimal among Levi subgroups of $G$ containing $S$ then $L'$ is simply connected and therefore has no factors $A_{1}$ or $A_{5}$ as $S$ has no irreducible $2$- or $6$-dimensional modules. Thus $L'$ has type $E_{6}$, $D_{5}$, $D_{4}$, $A_{4}$, $A_{3}$, $A_{2}$ or $A_{2}^{2}$. In all but the first case, using Lemma \ref{lem:classicalparabs} and proceeding as in $F_{4}$ gives a list of possible subgroups of type $A_{1}$ containing $S$; the restrictions of the various $G$-modules are again straightforward to determine, and it is clear that the given summand $W$ satisfies the required property.

So now assume that $S$ is $G$-irreducible. If $S$ fixes a nonzero vector on $V_{G}(\lambda_1)$ then $S$ lies in a subgroup of dimension $\ge 78 - 27 = 51$, and such a subgroup is either simple of type $F_{4}$, or is $G$-reducible. If $S$ fixes no nonzero vector on $V_{G}(\lambda_1)$, then inspecting Table \ref{a5e6p0} we see that $S$ must fix a nonzero vector on $L(G)$, and therefore lies in a subgroup containing a maximal torus. This shows that $S$ lies in a subgroup $A_{2}^{3}$, $A_{1}A_{5}$ or $F_{4}$. In the first two cases, we again derive a list of possible subgroups $A_{1}$ containing $S$, which appear in Table \ref{tab:e6thml2} or \ref{tab:e6thml1l6}.

So now we assume that $S$ lies in a subgroup $F_{4}$ of $G$ and in no subsystem subgroup of $G$. Then $S$ must be fixed-point free on $V_{F_4}(\lambda_4)$, otherwise it lies in a subsystem subgroup of $F_{4}$, hence centralises a non-central semisimple element of $G$ and lies in the corresponding subsystem subgroup of $G$. Since $V_{G}(\lambda_1) \oplus V_{G}(\lambda_6) \downarrow F_{4} = 0^{2} \oplus V_{F_4}(\lambda_4)^{2}$, every trivial $S$-submodule of this is a trivial submodule for $F_{4}$.

\subsection*{Remaining cases} For the remaining triples $(G,H,p)$, the possible subgroups $X$ are straightforward to determine using the representation theory of $H$ and Propositions \ref{prop:notprim} and \ref{prop:algorithm}. Moreover comparing $V \downarrow X$ with the appropriate table of feasible characters in Chapter \ref{chap:thetables} shows that the conditions of Proposition \ref{prop:onetwelve} hold for the action of $X$ and $S$ on the submodule $W$. There are a few cases where further argument is required to show that the submodules of $W$ are of the form necessary to apply Corollary \ref{cor:phistab}, as follows:

\begin{itemize}
\item When $(G,p) = (E_7,7)$ and $H = L_{2}(8)$, $L_{2}(13)$ or $U_{3}(3)$, if $S \cong H$ lies in a subgroup $X = G_{2}$ such that $V \downarrow X$ is a sum of trivial modules and 26-dimensional modules $V_{G_2}(20)$, we claim that every trivial $S$-submodule of $V$ is an $X$-submodule. The symmetric square of each $7$-dimensional orthogonal $S$-module has a unique $1$-dimensional submodule and a unique $1$-dimensional quotient. Since $V_{G_2}(20)$ is a $26$-dimensional section of the 28-dimensional module $\bigwedge^{2}V_{G_2}(10)$, it follows that $S$ fixes no nonzero points on $V_{G_2}(20)$, and the claim follows.

\item Similarly if $S \cong L_{2}(17)$ lies in a subgroup $C_{4}$ when $p = 2$, we must show that $X$ fixes no nonzero vectors on $V_{C_4}(\lambda_2)$. This follows because this is a $26$-dimensional section of $\bigwedge^{2}V_{C_4}(\lambda_1)$, and the exterior square of the $8$-dimensional irreducible $H$-module has shape $1|8|1|8|1|8|1$.

\item If $S \cong \Alt_{7}$ is a subgroup of $G = E_{7}$ when $p = 2$, we must justify the stated structure of $V_{G}(\lambda_7) \downarrow X$ where $X = A_{3} < A_{5}$, given in Table \ref{tab:e7thml7}. This follows since $V_{G}(\lambda_7) \downarrow A_{5} = \lambda_1^{3} + \lambda_5^{3} + \lambda_3$; the factors $\lambda_1$ and $\lambda_5$ restrict to $X$ as $V_{X}(\lambda_2)$, while $\lambda_3$ restricts with high weights $\lambda_2^{2}$, $2\lambda_1$ and $2\lambda_3$. We verify computationally that the $\Alt_{7}$-module $\bigwedge^{3}(6)$ is indecomposable with shape $6|(4 + 4^{*}) | 6$.

\item If $(G,H,p) = (E_{8},J_{2},2)$ and $S \cong H$ lies in a subgroup $X = E_{6}$ or $E_{7}$ and in no proper subsystem subgroup of this, then as in the proof of Proposition \ref{prop:primstab} it follows that the only nonzero vectors of $L(G)$ which are fixed by $S$, must also be fixed by $X$, and so $X$ and $S$ fix exactly the same trivial submodules on $L(G)$, and $W$ is the sum of these.
\item Also with $(G,H,p) = (E_{8},J_{2},2)$, if $S \cong H$ lies either in a subgroup $G_{2}$ which is maximal in a subgroup $E_{6}$, or in a subgroup $G_{2} < A_{5}$, then the action of $G_{2}$ on $V_{E_6}(\lambda_1)$ is given by \cite{MR2044850}*{Table 10.2} or \cite{MR1329942}*{Table 8.1}. It follows that if $S$ fixes a vector which is not fixed by this $G_{2}$, then there is an indecomposable $G_{2}$-module $0|10$ or $0|20$ on which $S$ fixes a non-zero vector. This is impossible since this would place $S$ in the full centraliser of this vector, which is proper in $G_{2}$ and has dimension at least $7$.
\end{itemize}

\section{Proof of Theorem \ref{THM:MAIN}: Special Cases} \label{sec:special}

We now complete the proof of Theorem \ref{THM:MAIN} by applying ad-hoc arguments for the types $(G,H,p)$ not covered above. Recall from Section \ref{sec:postponed} that these remaining cases are as follows:

\begin{table}[H]\small
\begin{tabular}{c|c}
$G$ & $(H,p)$ \\ \hline
$E_7$ & $\Alt_{10} \ (p = 5)$, $\Alt_{9} \ (p  \neq 2, 3)$ \\
$E_8$ & $\Alt_{16} \ (p = 2)$, $\Alt_{11}\ (p = 11)$, $\Alt_{10} \ (p \neq 2, 3)$
\end{tabular}
\label{tab:adhoc}
\end{table}

In each case, we let $S \cong H$ be a subgroup of $G$, take a proper simple subgroup $A$ of $S$, construct a proper, connected subgroup $X$ of $G$ containing $A$, and show that $Y \stackrel{\textup{def}}{=} \left<S,X\right>^{\circ}$ is proper and stabilises an appropriate subset $\mathcal{M}$ of $Y$-submodules of $L(G)$ or $V_{\textup{min}}$, so that Corollary \ref{cor:phistab} applies and $Y$ is contained in the $N_{\textup{Aut}(G)}(S)$-stable proper connected subgroup $\left(\bigcap_{M \in \mathcal{M}} G_{M}\right)^{\circ}$. Since $Y \cap S$ is normal in $S$ and contains $A$, it follows that $S \le Y \le \left(\bigcap_{M \in \mathcal{M}} G_{M}\right)^{\circ}$, and the conclusion of Theorem \ref{THM:MAIN} holds for $(G,H,p)$.

\subsection{Remark} In order to construct an appropriate subgroup $X$ as above, we assume that the subgroup $A$ of $S$ is $G$-completely reducible. For this, we appeal to Theorem \ref{THM:NONGCR}, hence Theorem \ref{THM:MAIN} depends on Theorem \ref{THM:NONGCR} in these few cases. We take this opportunity to note that the proof of Theorem \ref{THM:NONGCR} given in Chapter \ref{chap:gcr} depends only on the feasible characters of simple subgroups (Theorem \ref{THM:FEASIBLES}), and not on Theorem \ref{THM:MAIN}.

\subsection*{Case: $G = E_7$, $H = \Alt_{10}$ or $\Alt_{9}$, $p \neq 2,3$} \label{sec:a10e7p5} \label{sec:a9e7p5} \label{sec:a9e7p7} \label{sec:a9e7p0}
Here the composition factors of $S$ on $L(G)$ and $V_{\textup{min}}$ are specified by one of the tables on pages \pageref{a10e7p5}-\pageref{a9e7p5}. Let $A \cong \Alt_{8}$ be a subgroup of $S$, lying in an intermediate subgroup isomorphic to $\Alt_{9}$. Then this subgroup $\Alt_{9}$ has an 8-dimensional composition factor on $L(G)$, whose restriction to $A$ contains a trivial submodule. Thus Proposition \ref{prop:substab} applies to $L(G) \downarrow A$, and $A$ is not Lie primitive in $G$. Since $p \neq 2,3$, Theorem \ref{THM:NONGCR} holds and $A$ is $G$-completely reducible. In addition, by Theorem \ref{thm:subtypes} there does not exist an embedding of $\Alt_{8}$ into a smaller exceptional algebraic group, and thus $A$ lies in a $G$-completely reducible semisimple subgroup, whose simple factors are all of classical type.

Now $\Alt_{8}$ has a $7$-dimensional irreducible module, giving an embedding into $B_3$, and no other non-trivial irreducible modules dimension $\le 12$. In addition, the only non-trivial faithful module for $2\cdot \Alt_{8}$ of dimension $\le 12$ is 8-dimensional, giving an embedding into $D_4$. If $2\cdot \Alt_{8}$ embeds into a simply connected group of type $D_4$ with its centre contained in $Z(D_4)$, then its centre acts trivially on one of the three $D_4$-modules $\lambda_1$, $\lambda_3$ or $\lambda_4$. Since these are self-dual, this restricts as a direct-sum $1 \oplus 7$, and therefore the quotient $\Alt_{8}$ in $D_4$ lies in a proper subgroup of type $B_3$.

Thus $A$ lies in a subgroup $B_3$ of $G$. The two conjugacy classes of such subgroups and their action on $V_G(\lambda_7)$ are given by \cite{MR1329942}*{Table 8.6}. Comparing these with the feasible characters of $A$ (Tables \ref{a8e7p5}, \ref{a8e7p7}, \ref{a8e7p0}), we see that $A$ lies in a subgroup $X = B_3 < A_6$ with $V_G(\lambda_7) \downarrow X = \lambda_1^{2}/\lambda_2^{2}$. By Lemma \ref{lem:exthom} and Proposition \ref{prop:weyls}, this is completely reducible. Thus Proposition \ref{prop:newsamesubs} applies and every $A$-submodule of $V_{G}(\lambda_7)$ is $X$-stable. In particular, every $S$-submodule of $V_G(\lambda_7)$ is preserved by $Y = \left<S,X\right>$, hence if $\mathcal{M}$ denotes the (non-empty) collection of proper $S$-submodules of $V_G(\lambda_7)$, we have $S < Y < \left(\bigcap_{M \in \mathcal{M}} G_M \right)^{\circ}$, and this latter group is proper, connected and $N_{\textup{Aut}(G)}(S)$-stable.

\subsection*{Case: $(G,H,p) = (E_8,\Alt_{16},2)$} 

Let $A \cong \Alt_{15}$ be a subgroup of $S$. Then Proposition \ref{prop:notprim} and Theorem \ref{THM:NONGCR} apply to $A$, and $A$ is not Lie primitive in $G$, and is $G$-cr. Since the smallest nontrivial $A$-module is $14$-dimensional and $H^{1}(A,14) = 0$, and since $A$ does not embed into a smaller exceptional algebraic group, we deduce that $A$ lies in a Levi subgroup of type $D_{7}$, call it $X$. From \cite{MR1329942}*{Table 8.1} it follows that
\[ L(G) \downarrow X = (0^{2}/V_X(\lambda_2)) \oplus \lambda_1^{2} \oplus \lambda_6 \oplus \lambda_7. \]
Comparing this with $L(G) \downarrow A$, whose factors are given by Table \ref{a15e8p2}, we find that $A$ must be irreducible on each $X$-composition factor. In particular every $14$- and $64$-dimensional $A$-submodule of $L(G)$ is the restriction of an $X$-submodule.

Now, if $S$ has no irreducible submodules of dimension $1$, $14$ or $64$ on $L(G)$, then $L(G) \downarrow S$ is an image of the projective cover of the 90-dimensional $S$-module. But this is absurd, since this composition factor is self-dual and occurs with multiplicity 1. Thus either $S$ has a nonzero trivial submodule on $L(G)$ or has irreducible submodules of dimension $14$ or $64$, in which case $\left<S,X\right>^{\circ}$ preserves every such submodule, and is therefore a proper, connected subgroup containing $S$.

So $S$ is not Lie primitive in $G$; hence by Lemma \ref{lem:gcrstab} we may assume that $S$ is $G$-cr. Since the smallest $S$-module is $14$-dimensional and supports a quadratic form, we deduce that $S$ lies in a subgroup $D_{7}$; this is the subgroup $X$ constructed above, and we deduce that every $14$- and $64$-dimensional $S$-composition factor on $L(G)$ is in fact an $X$-submodule. Thus $S < X < \left( \bigcap_{M} G_M \right)^{\circ}$, the intersection over $14$- and $64$-dimensional $S$-submodules, and by Corollary \ref{cor:phistab} this latter group is proper, connected and $N_{\textup{Aut}(G)}(S)$-stable.

\subsection*{Case: $(G,H,p) = (E_8,\Alt_{11},11)$ or $(E_8,\Alt_{10},p \neq 2,3,5)$} \label{sec:a11e8p11} \label{sec:a10e8p0} \label{sec:a10e8p7}
Assume first that $p \neq 7$. Note that each feasible character of $\Alt_{11}$ or $\Alt_{10}$ on $L(G)$ has two 84-dimensional composition factors (Table \ref{a11e8p11}). Let $S \cong H$ be a subgroup of $G$, and let $A \cong \Alt_{9}$ be a subgroup of $S$. Matching up Tables \ref{a11e8p11} and \ref{a9e8p0}, we see that
\[ L(G) \downarrow A = 1/8^{3}/27/28^{3}/56^{2}, \]
which is completely reducible since $p \nmid |A|$. Thus $A$ fixes a vector on $L(G)$ and lies in a proper subgroup of positive dimension, hence also in a proper, connected subgroup by Proposition \ref{prop:consub}. By Theorem \ref{THM:NONGCR}, $A$ also is $G$-completely reducible.

Inspecting \cite{MR1942256}*{Table 2} we see that the smallest irreducible modules for $A$ and its proper cover $2.A$ are of dimension $8$ and $21$, and the $8$-dimensional modules give an embedding into $D_{4}$. Since $\Alt_9$ admits no embeddings into a proper exceptional subgroup of $G$ by Theorem \ref{thm:subtypes}, it follows that a minimal connected reductive subgroup containing $A$ can only involve factors of type $D_4$. Comparing the above decomposition with \cite{MR1329942}*{Table 8.1}, we deduce that $A$ lies in a simple subgroup $X$ of type $D_4$, contained in a Levi subgroup $A_7$. Then
\begin{align*}
L(G) \downarrow X &= 0/\lambda_1^{2}/\lambda_2^{3}/2\lambda_1/(\lambda_3+\lambda_4)^{2},
\end{align*}
which is completely reducible by Corollary \ref{cor:2step} and Proposition \ref{prop:weyls}. Now, let $W$ be the span of all $28$-dimensional and $56$-dimensional $A$-submodules of $L(G)$. Then $W$ is the restriction of an $X$-submodule with factors $\lambda_2^{3}/(\lambda_3+\lambda_4)^{2}$, which contains every 28- and 56-dimensional $X$-composition factor of $L(G)$. Proposition \ref{prop:newsamesubs} then applies, so that every $A$-submodule of $W$ is $X$-invariant.

Now, as $L(G)$ is self-dual, $S$ has either a unique 84-dimensional submodule or has a summand of shape $84+84$. In either case, this is contained in $W$, since `84' must restrict to $A$ with composition factors $28/56$. Thus each 84-dimensional irreducible $S$-submodule of $L(G)$ is an $X$-submodule, hence preserved by $Y \stackrel{\textup{def}}{=} \left<S,X\right>$. Thus $S < Y \le (\bigcap G_M)^{\circ}$, the intersection over $84$-dimensional irreducible $S$-submodules of $L(G)$. This latter group is $N_{\textup{Aut}(G)}(S)$-stable by Corollary \ref{cor:phistab}.

If $p = 7$, then the same argument goes through with the caveat that $L(G) \downarrow A$ may no longer be completely reducible; however, it still has a trivial submodule since $H^{1}(\Alt_{9},M)$ vanishes for each composition factor $M$ in a feasible character (cf.\ Table \ref{tab:altreps}). The $28$- and $56$-dimensional $A$-modules are projective, since their dimension is divisible by the order of a Sylow subgroup of $A$, and therefore they are $A$-direct summands of $L(G)$.

\subsection*{Case: $(G,H,p) = (E_8,\Alt_{10},5)$} \label{sec:a10e8p5}
Here we have two feasible characters (Table \ref{a10e8p5}). The only $H$-module $W$ occurring in either such that $H^{1}(H,W) \neq \{0\}$ is the 8-dimensional module. Thus in Case 2), Proposition \ref{prop:substab} applies, and $G$ has no Lie primitive subgroups isomorphic to $H$ having this action on $L(G)$.

So let $S \cong H$ be a subgroup of $G$ giving rise to Case 1), and let $A \cong \Alt_{9}$ be a subgroup of $S$. Inspecting Table \ref{a9e8p5}, $A$ must act with composition factors $1/8^{3}/27/28^{3}/56^{2}$ on $L(G)$. Each factor has zero first cohomology group, and so $A$ fixes a nonzero vector, and is not Lie primitive in $G$. By Proposition \ref{prop:consub}, $A$ lies in a proper connected subgroup of $G$. It is also $G$-completely reducible by Proposition \ref{prop:levisub}. Now, $A$ does not lie in a subgroup $E_7$ (as this has three trivial composition factors on $L(G)$), and does not embed into a group of type $G_{2}$, $F_{4}$ or $E_{6}$, and so $A$ lies in a semisimple subgroup having only classical factors. The smallest non-trivial $A$-module has dimension $8$ and gives an embedding into a group of type $D_{4}$, hence $A$ lies in a subgroup $D_4$ of $G$. The conjugacy classes of these, and their action on $L(G)$, are given by \cite{MR1329942}*{Table 8.1}. Comparing this with the factors above, we deduce that $A$ lies in a subgroup $X$ of type $D_4$, such that
\[ L(G) \downarrow X = 0/\lambda_1^{2}/\lambda_2^{3}/2\lambda_1/(\lambda_3+\lambda_4)^2, \]
which is completely reducible by Corollary \ref{cor:2step} and Proposition \ref{prop:weyls}. Thus every irreducible $A$-submodule of dimension $1$, $28$ or $56$ is an $X$-submodule. In particular, if $S$ has any irreducible submodules of dimension $1$, $28$ or $56$ on $L(G)$, then $S < \left<S,X\right>^{\circ} \le \left(\bigcap_{M} G_M\right)^{\circ}$, the intersection over such $S$-submodules, and this is a proper, connected subgroup of $G$, which is $N_{\textup{Aut}(G)}(S)$-stable by Corollary \ref{cor:phistab}.

So suppose $L(G)$ has no irreducible $S$-submodules of dimension $1$, $28$ or $56$. Since $L(G)$ is self-dual it has no such irreducible quotients either, hence Lemma \ref{lem:projcover} implies that $L(G)$ is an image of the projective module $P_{8}^{2} \oplus P_{35}^{2}$. However, neither $P_8$ nor $P_{35}$ have a 56-dimensional composition factor; contradiction.

\subsection{Remark: Scope of Theorem \ref{THM:MAIN}} \label{rem:limitations}
The representation theory outlined here is not sufficient to prove results along the lines of Theorems \ref{THM:MAIN} and \ref{THM:NONGCR} for all the non-generic simple subgroups appearing in Table \ref{tab:subtypes} which are not known to occur as a Lie primitive subgroup of an exceptional group $G$.

For instance, when $p = 7$ the maximal subgroup $G_{2}$ of $G = F_{4}$ gives rise to a subgroup $U_{3}(3)$, which is fixed-point free on $L(G)$ and $V_{\textup{min}}$, but is not Lie primitive in $G$. In order to handle such subgroups in the manner of Theorems \ref{THM:MAIN}, \ref{THM:NONGCR} and their corollaries, it will be necessary to incorporate more information, such as the Lie algebra structure of $L(G)$.

\section{Proof of Theorem \ref{THM:FIN}} \label{sec:pfcorfin}

Having proved Theorem \ref{THM:MAIN}, the following result now implies Theorem \ref{THM:FIN}. The proof is similar to that of \cite{MR1458329}*{Theorem 6}.

\begin{proposition}
Let $G$ be an adjoint simple algebraic group and let $S$ be a non-abelian finite simple subgroup of $G$, not isomorphic to a member of $\textup{Lie}(p)$. Let $\sigma$ be a Frobenius morphism of $G$ such that $L = O^{p'}(G_\sigma)$ is simple, and let $L \le L_1 \le \textup{Aut}(L)$.

If there exists a proper closed, connected, $N_{\textup{Aut}(G)}(S)$-stable subgroup of $G$ containing $S$, then $L_{1}$ has no maximal subgroup with socle $S$.
\end{proposition}

\proof Suppose that $X < L_1$ is a maximal subgroup with socle $S$. It is well known that a simple group of Lie type has soluble outer automorphism group, hence the image of $S$ under the quotient map $L_{1} \to L_{1}/L$ is trivial. Thus $S \le L$ and $S$ is fixed point-wise by $\sigma$, and in particular $\sigma \in N_{\textup{Aut}(G)}(S)$. In addition, every automorphism of $L$ extends to a morphism $G \to G$ (cf.\ \cite{MR0407163}*{\S 12.2}). Since $S$ is normal in $X$, we can view $X\!\left<\sigma\right>$ as a subgroup of $N_{\textup{Aut}(G)}(S)$.

So let $\bar{S}$ be a proper connected, $N_{\textup{Aut}(G)}(S)$-stable subgroup of $G$ which contains $S$. Then $\bar{S}$ is $X\!\left<\sigma\right>$-stable. Let $Y$ be maximal among proper, connected, $X\!\left<\sigma\right>$-stable subgroups containing $\bar{S}$. Then $O^{p'}(Y_{\sigma})$ contains $S$ and is thus non-trivial. Now, we have containments
\[ X \le N_{L_1}(Y) \le L_1 \]
and if $N_{L_1}(Y) = L_1$ then $L$ normalises the non-trivial subgroup $O^{p'}(Y_{\sigma})$, which is a proper subgroup of $L$ since $Y$ is connected and proper in $G$, and this contradicts the simplicity of $L$. From the maximality of $X$ in $L_1$ it follows that $X = N_{L_1}(Y) \ge O^{p'}(Y_\sigma)$. In particular, if $Y$ is not reductive then $X$ normalises the non-trivial $p$-subgroup $R_u(Y)_\sigma$ of $L_1$, a contradiction. Hence $Y$ is reductive, and therefore $S = O^{p'}(Y_\sigma)$, contradicting $S \notin \textup{Lie}(p)$. \qed

%% file: AJL-gcr.tex

\chapter{Complete Reducibility} \label{chap:gcr}

In this chapter we prove Theorem \ref{THM:NONGCR}. We require an amount of background material concerning rational cohomology and complements in parabolic subgroups.

\section{Cohomology and Complements} \label{sec:complements}

Let $X$ be an algebraic group acting on a commutative algebraic group $V$, such that the action map $X \times V \to V$ is a morphism of varieties. We define $H^{1}(X,V)$ to be the quotient of the additive group $Z^{1}(X,V)$ of \emph{rational} 1-cocycles, by the subgroup $B^{1}(X,V)$ of rational 1-coboundaries. If $X$ is finite, then all cocycles are rational and we recover the usual first cohomology group.

This rational cohomology group parametrises conjugacy classes of closed complements to $V$ in $VX$, where a complement to $V$ (as an algebraic group) is now a subgroup $X'$ not only satisfying $VX' = VX$ and $V \cap X' = 1$, but also $L(V) \cap L(X') = 0$ (cf.\ \cite{MR2015057}*{I.7.9(2)}). Note that this latter condition is trivially satisfied when $X'$ is finite.

If $X$ is a closed subgroup of a parabolic subgroup $P$ of $G$ with $X \cap R_{u}(P) = 1$, then $R_{u}(P)$ admits a filtration by modules for the Levi factor, and we can use the cohomology groups of these to study complements to $R_{u}(P)$ in $R_{u}(P)X$, as follows.

\begin{proposition}
\label{prop:levisub}
Let $P = QL$ be a parabolic subgroup of $G$, and let $X$ be a finite subgroup of $P$ with $X \cap Q = 1$. Then the $KL$-modules $V_i = Q(i)/Q(i+1)$ also have the structure of $KX$-modules, and if $H^{1}(X,V_i) = \{0\}$ for all $i$, then all complements to $Q$ in $QX$ are $Q$-conjugate, and $X$ lies in a conjugate of $L$.
\end{proposition}

\proof The conjugation action of $X$ on $P$ induces an action on each $V_i$. If $\pi \ : \ P \twoheadrightarrow L$ is the natural quotient map, then we have $(qQ(i+1))^x = (qQ(i+1))^{\pi(x)}$ for all $q \in Q(i)$ and $x \in X$, since $Q(i)/Q(i+1)$ is central in $Q/Q(i+1)$. Thus the (linear) action of $L$ on each $V_i$ gives rise to linear action of $X$.

To prove that all closed complements to $Q$ in $QX$ are $Q$-conjugate, we work by induction on $i$, proving that all copies of $X$ in $(Q/Q(i))X$ are $Q/Q(i)$-conjugate. When $i = 1$ we have $Q/Q(i) = V_1$ and the vanishing cohomology group gives the result. Now assume this holds for some $i \ge 1$. If $Q(i) = \{0\}$ then we are done, so suppose not and let $Y$ be a complement to $Q/Q(i+1)$ in $(Q/Q(i+1))X$. Now, consider the projection $(Q/Q(i+1))X \to (Q/Q(i))X$. By the inductive hypothesis, we may replace $Y$ by a conjugate whose image under this projection is $X$. Then we have
\[ Y = \{ \phi(x).x \ : \ x \in X \} \]
for some rational map $\phi \ : \ X \to Q/Q(i+1)$, whose image lies in the kernel of the projection $Q/Q(i+1) \to Q/Q(i)$, which is $Q(i)/Q(i+1)$. Hence $\phi \in Z^1(X,V_i) = B^{1}(X,V_i)$ and $Y$ is $(Q/Q(i+1))$-conjugate to $X$, as required.

Finally, since $Q \cap X = 1$, the projection $P \twoheadrightarrow L$ restricts to an isomorphism of $X$ onto its image $\bar{X} \le L$. This is a complement to $Q$ in $QX$, and is thus conjugate to $X$ by the above. \qed

Next, we note that if $Y$ is a complement to $X$ in $QX$, then the composition factors of $Y$ and of $X$ in the filtration of $Q$ correspond. This is proved for reductive $X$ in \cite{MR3075783}*{Lemma 3.4.3}; the proof here is identical.

\begin{lemma} \label{lem:complement_factors}
Let $Q$ be a unipotent algebraic group over $K$, and let $X$ be a finite group, with no non-trivial normal $p$-subgroups if $K$ has characteristic $p > 0$. If $Y$ is a complement to $Q$ in the semidirect product $QX$, and if $V$ is a rational $QX$-module, then the composition factors of $V \downarrow X$ correspond to the composition factors of $V \downarrow Y$ under an isomorphism $X \to Y$.
\end{lemma}

\proof Without loss of generality we can assume that $V$ is an irreducible $QX$-module. Since $Q$ is unipotent, the fixed-point space of $Q$ in $V$ is non-trivial \cite{MR0396773}*{Theorem 17.5}, and since $Q$ is normal in $QX$, the space of $Q$-fixed points is $QX$-invariant, hence equal to $V$. Thus $Q$ acts trivially on $V$, hence the representation of $QX$ on $V$ factors through the projection $QX \to X$, and the composed map $Y \hookrightarrow QX \twoheadrightarrow X$ is the required isomorphism. \qed

\begin{corollary} \label{cor:levisub}
If $X$ is a non-abelian finite simple subgroup of $G$ which is not $G$-completely reducible, then $L(G) \downarrow X$ has a trivial composition factor, as well as a factor $W$ such that
\begin{itemize}
\item $H^{1}(X,W) \neq \{0\}$,
\item Either $W$ has multiplicity $\ge 2$, or $W \ncong W^{*}$,
\item $W$ has dimension at most $14$, $20$, $35$, $64$ when $G$ is respectively of type $F_4$, $E_6$, $E_7$, $E_8$.
\end{itemize}
\end{corollary}

\proof Let $P = QL$ be a parabolic subgroup of $G$ containing $X$, such that $X$ is not contained in a conjugate of the Levi factor $L$. The torus $Z(L)$ gives rise to a 1-dimensional $X$-composition factor on $L(G)$, which must be trivial. 

By Proposition \ref{prop:levisub}, some $X$-composition factor $W$ exists in the filtration of $Q$ by $X$-modules, such that $H^{1}(X,W) \neq \{0\}$, which is then an $X$-composition factor of $L(P) \subset L(G)$. Further, $W^{*}$ occurs as a composition factor of $L(G)/L(P) \cong L(Q^{\textup{op}})$, where $Q^{\textup{op}}$ is the unipotent radical of the opposite parabolic subgroup (see \cite{MR1047327}*{Remark 6, p.\ 561}), hence either $W$ has multiplicity $\ge 2$ as a factor of $L(G)$, or $W \ncong W^{*}$.

Finally, Lemma \ref{lem:radfilts} gives us the high weights of each simple factor of $L$ on the modules occurring in the filtration of $Q$, which allows us to determine the largest dimensions of a module occurring for some $L$. For instance, if $G$ has type $F_{4}$ then the types of Levi subgroup and the highest dimension of an irreducible module in the filtration of $Q$ is as follows:

\begin{center}
\begin{tabular}{c|ccccccc}
Levi subgroup & $B_{3}$ & $C_{3}$ & $A_{1} A_{2}$ & $A_{2}$ & $B_{2}$ & $A_{1} A_{1}$ & $A_{1}$ \\ \hline
High weight & $\lambda_{3}$ & $\lambda_{3}$ & $2 \otimes \lambda_{1}$ & $2\lambda_{1}$ & $\lambda_{1}$ & $1 \otimes 2$ & $2$ \\
Dimension & $8$ & $14$ & $12$ & $9$ & $5$ & $6$ & $3$
\end{tabular}
\end{center}
where we note that the high weights $2$ for a factor of type $A_{1}$, and $2\lambda_{1}$, $2\lambda_{2}$ for type $A_{2}$, can only occur when this factor contains short root subgroups; hence no module $2 \otimes 2$ can occur for $L$ of type $A_{1} A_{1}$, and no module $2 \otimes 2\lambda_{i}$ for $i \in \{ 1,2 \}$ can occur for $L$ of type $A_{1} A_{2}$. Hence the largest module dimension occurring is $14$.

Similarly, for $G$ of type $E_{6}$, the highest possible dimension of a module for a Levi subgroup with weights as in Lemma \ref{lem:radfilts} are the 20-dimensional module $\lambda_{3}$ for $L$ of type $A_{5}$, and the 20-dimensional module $\lambda_{1} \otimes \lambda_{2}$ for $L$ of type $A_{1}A_{4}$.

For $G$ of type $E_{7}$, and $L$ of type $A_{6}$ we get a 35-dimensional module $\lambda_{3}$. The only higher-dimensional module for a Levi subgroup with weights as in Lemma \ref{lem:radfilts} is the $40$-dimensional module $1 \otimes \lambda_{3}$ for $L$ of type $A_{1}A_{5}$. However, $G$ has a unique standard parabolic subgroup of this type, and we check directly that the high weights of modules occurring in the unipotent radical are $1 \otimes \lambda_{2}$, $0 \otimes \lambda_{4}$ and $1 \otimes 0$, and all such modules have dimension less than $35$.

For $G$ of type $E_{8}$ and $L$ of type $D_{7}$ we get a 64-dimensional module $\lambda_{7}$. The only modules of larger dimension for a Levi factor with weights as in Lemma \ref{lem:radfilts} are the 70-dimensional modules $1 \otimes \lambda_{3}$ and $1 \otimes \lambda_{4}$ when $L$ has type $A_{1} A_{6}$, and again we check directly that these do not occur in standard parabolic with this Levi subgroup type. \qed

\section{Completely Reducible Subgroups}

Corollary \ref{cor:levisub} above places restrictions on the feasible characters potentially arising from a non-$G$-cr subgroup of an adjoint exceptional simple algebraic group $G$.

Let $S \cong H$ be a non-$G$-cr subgroup of $G$, and let $\tilde{S}$ be a minimal preimage of $S$ in the simply connected group $\tilde{G}$. Then $\tilde{S}$ lies in a proper parabolic subgroup $P$ of $G$, with Levi decomposition $P = QL$. Since $\tilde{G}$ is simply connected, so is $L'$, hence $L'$ is a direct product of simply connected simple groups. In addition the image of $\tilde{S}$ under the projection $P \twoheadrightarrow L$ lies in $L'$, since $\tilde{S}$ is perfect. Thus if $L_{0}$ is a simple factor of $L'$, then $\tilde{S}$ acts on the natural module $V_{L_{0}}(\lambda_1)$ if $L_{0}$ is classical, and if $L_{0}$ is exceptional then the image of $\tilde{S}$ under $L' \twoheadrightarrow L_{0}$ corresponds to a feasible character in Chapter \ref{chap:thetables}.

Thus it is straightforward to determine whether $\tilde{S}$ admits an embedding into a proper Levi subgroup of $\tilde{G}$. In Chapter \ref{chap:thetables} we have labelled `\nongcr' those feasible characters which satisfy the conclusion of Corollary \ref{cor:levisub}, for subgroup types of $\tilde{S}$ which admit an embedding into a proper Levi subgroup of $G$. In Table \ref{tab:candidates}, we collect together all triples $(G,H,p)$ with such a feasible character. Thus if $G$ has a non-$G$-cr subgroup isomorphic to the group $H \notin \textup{Lie}(p)$, then $(G,H,p)$ appears in Table \ref{tab:candidates}.

Recall that the groups $\Alt_5 \cong L_2(4) \cong L_2(5)$, $\Alt_6 \cong L_2(9) \cong Sp_4(2)'$, $\Alt_8 \cong L_4(2)$, $L_2(7) \cong L_3(2)$, $U_4(2) \cong PSp_4(3)$ and $U_3(3) \cong G_2(2)'$ are considered to be of Lie type in each corresponding characteristic, hence do not appear in those characteristics in Table \ref{tab:candidates}. 

Note also that if $G$ is of type $G_2$, then according to Table \ref{tab:onlyprim}, the only non-generic finite simple subgroups of $G$ which are not Lie primitive in $G$, are isomorphic to $\Alt_{5}$ or $L_2(7)$, in which case $p \neq 2$. But neither $\Alt_{5}$ nor $L_2(7)$ has a non-trivial 2-dimensional module for $p \neq 2$, and so these groups have no embeddings into a Levi subgroup of $G$, hence all non-generic finite simple subgroups of $G$ are $G$-irreducible in this case.

\begin{table}[H]\small
\caption{Candidate Non-$G$-cr subgroup types}
\begin{tabularx}{.9\linewidth}{c|X}
$G$ & \multicolumn{1}{c}{$(H,p)$} \\ \hline
$F_4$ & $(\Alt_{5},3)$, $(J_2,2)$, $(L_2(7),3)$, $(L_2(8),7)$, $(L_2(8),3)$, $(L_2(13),2)$ \\ \hline

$E_6$ & $(\Alt_{10},2)$, $(\Alt_{9},2)$, $(\Alt_{7},7)$, $(\Alt_{7},3)$, $(\Alt_{7},2)$, $(\Alt_{6},5)$, $(\Alt_{5},3)$, $(M_{11},5)$, $(M_{11},3)$, $(M_{11},2)$, $(M_{22},2)^{\dag}$, $(J_2,2)$, $(L_2(7),3)$, $(L_2(8),7)$, $(L_2(8),3)$, $(L_2(11),5)$, $(L_2(11),3)$, $(L_2(11),2)$, $(L_2(13),7)$, $(L_2(13),2)$, $(L_2(17),3)$, $(L_2(17),2)$, $(U_4(3),2)^{\dag}$ \\ \hline

$E_7$ & $(\Alt_{12},2)$, $(\Alt_{10},2)$, $(\Alt_{9},3)$, $(\Alt_9,2)$, $(\Alt_{8},3)$, $(\Alt_{7},7)$, $(\Alt_{7},5)$, $(\Alt_{7},3)$, $(\Alt_{7},2)$, $(\Alt_{6},5)$, $(\Alt_{5},3)$, $(M_{11},11)$, $(M_{11},5)$, $(M_{11},3)$, $(M_{11},2)$, $(M_{12},3)^{\dag}$, $(M_{12},2)$, $(J_2,3)^{\dag}$, $(J_2,2)$, $(L_2(7),3)$, $(L_2(8),7)$, $(L_2(8),3)$, $(L_2(11),5)$, $(L_2(11),3)$, $(L_2(11),2)$, $(L_2(13),7)$, $(L_2(13),3)$, $(L_2(13),2)$, $(L_2(17),3)$, $(L_2(17),2)$, $(L_2(19),5)^{\dag}$, $(L_2(19),3)$, $(L_2(19),2)$, $(L_2(25),3)$, $(L_2(25),2)$, $(L_2(27),13)^{\dag}$, $(L_2(27),7)$, $(L_2(27),2)$, $(L_3(3),13)$, $(L_3(3),2)$, $(L_{3}(4),3)^{\dag}$, $(U_3(3),7)$, $(^{3}D_4(2),3)$, $(^{2}F_4(2)',5)$ \\ \hline

$E_8$ & $(\Alt_{16},2)$, $(\Alt_{14},2)$, $(\Alt_{12},2)$, $(\Alt_{10},5)$, $(\Alt_{10},2)$, $(\Alt_{9},3)$, $(\Alt_{9},2)$, $(\Alt_{8},7)$, $(\Alt_{8},3)$, $(\Alt_{7},7)$, $(\Alt_{7},5)$, $(\Alt_{7},3)$, $(\Alt_{7},2)$, $(\Alt_{6},5)$, $(\Alt_{5},3)$, $(M_{11},11)$, $(M_{11},5)$, $(M_{11},3)$, $(M_{11},2)$, $(M_{12},2)$, $(J_2,2)$, $(L_2(7),3)$, $(L_2(8),7)$, $(L_2(8),3)$, $(L_2(11),5)$, $(L_2(11),3)$, $(L_2(11),2)$, $(L_2(13),7)$, $(L_2(13),3)$, $(L_2(13),2)$, $(L_2(17),3)$, $(L_2(17),2)$, $(L_2(19),5)$, $(L_2(19),3)$, $(L_2(19),2)$, $(L_2(25),3)$, $(L_2(25),2)$, $(L_2(27),7)$, $(L_2(27),2)$, $(L_2(29),2)$, $(L_2(37),2)$, $(L_3(3),13)$, $(L_3(3),2)$, $(U_3(3),7)$, $(U_3(8),3)$, $(U_4(2),5)$, $(PSp_4(5),2)$, $(^{3}D_4(2),3)$, $(^{2}F_4(2)',5)$
\end{tabularx}\\
If a subgroup type is marked with $^{\dag}$, then any non-$G$-cr subgroup of this type lifts to a proper cover in the simply connected cover of $G$.
\label{tab:candidates}
\end{table}

\section{Proof of Theorem \ref{THM:NONGCR}}

We now prove that certain triples $(G,H,p)$ appearing in Table \ref{tab:candidates} do not give rise to non-$G$-cr subgroups. This proves Theorem \ref{THM:NONGCR} and Corollary \ref{cor:nongcr}.

\begin{proposition} \label{prop:moregcrtypes}
If $G$ is an adjoint exceptional simple algebraic group in characteristic $p$, if $H \notin \textup{Lie}(p)$, and if $(G,H,p)$ appears in Table \ref{tab:onlycgr}, then all subgroups of $G$ isomorphic to $H$ are $G$-completely reducible.
\end{proposition}

\begin{table}[H]\small
\caption{Candidates giving rise only to $G$-cr subgroups}
\begin{tabularx}{.9\linewidth}{c|X}
$G$ & $(H,p)$ \\ \hline
$E_6$ & $(\Alt_{10},2)$, $(\Alt_{9},2)$, $(\Alt_{7},2)$, $(M_{11},2)$, $(J_2,2)$, $(L_2(13),7)$ \\
$E_7$ & $(J_{2},3)$, $(L_2(19),5)$, $(L_2(19),3)$, $(L_2(25),2)$, $(L_2(27),13)$, $(L_2(27),7)$, $(L_{3}(4),3)$ \\
$E_8$ & $(\Alt_{14},2)$, $(\Alt_{8},7)$, $(L_2(27),7)$, $(L_2(29),2)$
\end{tabularx}
\label{tab:onlycgr}
\end{table}

\proof In each case, let $S \cong H$ be a subgroup of $G$, let $\tilde{S}$ be a preimage of $S$ in the simply connected cover $\tilde{G}$ of $G$, and suppose that $P$ is a proper parabolic subgroup, with Levi decomposition $P = QL$, which is minimal among parabolic subgroups of $\tilde{G}$ containing $\tilde{S}$, so that the image of $\tilde{S}$ in $L$ is an $L'$-irreducible subgroup of $L'$.

If $G$ is of type $E_7$, note that no cover of any of the corresponding groups $L_{2}(q)$ in Table \ref{tab:onlycgr} has a faithful, irreducible orthogonal module of dimension $\le 12$, a symplectic module of dimension $\le 10$, or any faithful module of dimension $\le 7$ (cf.\ \cite{MR1835851}*{Table 2}). In particular, this means that such a group cannot have a non-trivial homomorphism into a simple factor of $L$ of classical type. Thus $L'$ is simple of type $E_6$. The corresponding unipotent radical is abelian and irreducible as an $L'$-module, with high weight $\lambda_1$. Since $V_{E_7}(\lambda_7) \downarrow L' = 0^{2}/\lambda_1/\lambda_6$ by \cite{MR1329942}*{Table 8.2}, we deduce that a non-$G$-cr subgroup isomorphic to $H$ must have a composition factor on both $L(G)$ and $V_{E_7}(\lambda_7)$ with a nonzero first cohomology group. Inspecting the relevant tables (\ref{l219e7p5}, \ref{l219e7p3}, \ref{l225e7p2}, \ref{l227e7p13}, \ref{l227e7p7}), we see that this never occurs.

If $(G,H,p) = (E_{7},J_{2},3)$, the double cover $2 \cdot J_{2}$ has a faithful 6-dimensional symplectic irreducible module, giving an irreducible embedding into a Levi subgroup of type $A_{5}$. Now $H^{1}(2 \cdot J_2,6_{a}) = H^{1}(2 \cdot J_2,6_{b}) = 0$, and also $J_{2}$ does not embed into a smaller exceptional group when $p = 3$ by Theorem \ref{thm:subtypes}. It follows that $L'$ is simple of type $A_{5}$. Comparing Table \ref{2j2e7p3} with the composition factors of each on $L(G)$ (given by \cite{MR1329942}*{Table 8.2}), we deduce that $L(G) \downarrow L' = 0^{8}/\lambda_2^{3}/\lambda_4^{3}/(\lambda_1+\lambda_5)$, and in particular, the non-trivial irreducible $L'$-modules occurring in the filtration of $Q$ each have high weight $\lambda_2$ or $\lambda_4$. The corresponding Weyl modules are irreducible and are isomorphic to the alternating squares of the modules $\lambda_1$ and $\lambda_5$, respectively, which restrict to irreducible 6-dimensional $\tilde{S}$-modules. Using Magma to help with calculations, we find that $\bigwedge^{2}(6_{a})$ and $\bigwedge^{2}(6_{b})$ are each self-dual and uniserial with composition factor dimensions $1$, $13$, $1$. Since $H^{1}(J_2,13_{a}) \cong H^{1}(J_2,13_b)$ is $1$-dimensional, from Proposition \ref{prop:substab}(ii) it follows that $H^{1}(\tilde{S},\lambda_2) = H^{1}(\tilde{S},\lambda_4) = 0$, and so Proposition \ref{prop:levisub} applies and $\tilde{S}$ is $G$-completely reducible.

If $(G,H,p) = (E_{7},L_{3}(4),3)$, the double cover $2 \cdot L_{3}(4)$ of $H$ has a faithful 6-dimensional irreducible module, giving an embedding into a subgroup of type $A_{5}$. Since the image of $S$ in $L'$ is $L'$-irreducible, and since $L_{3}(4)$ does not embed into a smaller exceptional group by Theorem \ref{thm:subtypes}, we deduce that $L'$ is simple of type $A_{5}$. But by \cite{MR1329942}*{Table 8.6}, each subgroup of $G$ of type $A_{5}$ has at least four six-dimensional composition factors on $V_{56}$. Thus the embedding of $H$ must correspond to Case 3) of Table \ref{2l34e7p3}, which cannot come from a non-$G$-cr subgroup by Corollary \ref{cor:levisub}.

If $(G,H,p) = (E_6,\Alt_{10},2)$, we must have $L'$ simple of type $D_5$, as no other Levi subgroup admits an embedding of $H$. Since $\Alt_{10}$ does not preserve a nondegenerate quadratic form on its $8$-dimensional irreducible module, and has no other non-trivial irreducible modules of dimension $\le 10$, a subgroup $S \cong H$ of $L'$ must act uniserially on the natural 10-dimensional $L'$-module, with shape $1|8|1$. Inspecting Table \ref{a10e6p2} we see that the $L'$-module of high weight $\lambda_5$ occurring in the filtration $Q$ restricts to $S$ as an irreducible 16-dimensional module. Applying Proposition \ref{prop:substab}(ii) to $1|8|1$, in each case deduce that $H^{1}(S,V) = \{0\}$, and so all subgroups $S \cong H$ of $G$ are $G$-completely reducible.

If $(G,H,p) = (E_6,\Alt_{9},2)$, then $L'$ is simple of type $D_4$. Every non-trivial $L'$-module occurring in the filtration of the unipotent radical has dimension $8$, hence restricts to a subgroup $S \cong H$ with only trivial or 8-dimensional composition factors. The first cohomology group vanishes for each such $S$-module, hence by \ref{cor:levisub} all subgroups $S \cong H$ of $G$ are $G$-completely reducible.

If $(G,H,p) = (E_6,\Alt_{7},2)$, a non-$G$-cr subgroup isomorphic to $H$ must correspond to Case 1) or 3) of Table \ref{a7e6p2}, or to Case 1) of Table \ref{3a7e6p2}. The Levi subgroups of $G$ containing a copy of $H$ are of type $A_3$ and $A_5$. The only factors in the feasible character with nonzero first cohomology group are $14-$ or $20-$dimensional. By Lemma \ref{lem:radfilts}, the non-trivial modules occurring in the filtration of the unipotent radical of an $A_{3}$-parabolic subgroup containing $S$ are have high weight $\lambda_1$, $\lambda_2$ or $\lambda_3$, hence dimension 4 or 6, and hence have trivial first cohomology group, so such a parabolic cannot contain a non-$G$-cr subgroup isomorphic to $H$.

On the other hand, an $A_5$ Levi subgroup of a parabolic $P$ acts on $R_u(P)$ with high weights $\lambda_3$ and $0$. It is routine to check that the module $\lambda_3 \downarrow S = \bigwedge^{3}(\lambda_1) \downarrow S = \bigwedge^{3}(6_a)$ or $\bigwedge^{3}(6_b)$ has no 14- or 20-dimensional composition factors, when $S \cong 3.\Alt_{7}$ is a subgroup of $A_5$. Thus no non-$G$-cr subgroups isomorphic to $H$ occur here, either.

If $(G,H,p) = (E_6,M_{11},2)$, the only parabolic subgroups of $G$ admitting an embedding of $H$ have Levi factor $L$ with $L'$ simple of type $D_5$. Then $Q$ is a 16-dimensional irreducible $L'$-module, which by comparison with Table \ref{m11e6p2} must restrict irreducibly to a subgroup $S \cong H$. Thus $H^{1}(H,Q) = \{0\}$ and so all subgroups $S \cong H$ of $G$ are $G$-completely reducible.

If $(G,H,p) = (E_6,J_{2},2)$, a parabolic subgroup $P = QL$ containing an $L'$-irreducible copy of $H$ must have $L'$ of type $D_4$ or $A_5$. In the first case, the unipotent radical has a filtration by 8-dimensional $D_4$-modules which each admit a nondegenerate quadratic form. Since $H$ does not preserve a nonzero quadratic form on its 6-dimensional modules, and has no other faithful modules of dimension $\le 8$, each of these $D_4$-modules must restrict to $H$ with shape $1|6_x|1$, where $x \in \{a,b\}$. In particular, $H^{1}(H,V)$ vanishes for each such module (Proposition \ref{prop:substab}(ii)), hence $H^{1}(H,Q) = \{0\}$ and no non-$G$-cr subgroups arise here. If instead $L'$ is of type $A_5$, then the only non-trivial $L'$-module occurring in $Q$ has high weight $\lambda_3$. We find that $\bigwedge^{3}(6_a)$ is uniserial of shape $6_a|1|6_b|1|6_a$, and similarly for $6_b$. Now, if $H^{1}(H,\bigwedge^{3}(6_a)) \neq \{0\}$, let $E$ be a non-split extension of $\bigwedge^{3}(6_a)$ by 1-dimensional trivial $KS$-module. Then $E$ has socle $6_a$, and so $E^{*}$ is an image of the projective cover $P_{6_a}$. Using Magma to help with calculations, we find that the largest quotient of $P_{6_a}$ having only 1- and 6-dimensional composition factors, has only two trivial composition factors, and is thus not isomorphic to $E^{*}$. It follows that $H^{1}(H,\bigwedge^{3}(6_a))$ vanishes, and similarly $H^{1}(H,\bigwedge^{3}(6_a))$ vanishes, so all subgroups $S \cong H$ of $G$ are $G$-completely reducible.

If $(G,H,p) = (E_6,L_2(13),7)$ we must have $L'$ of type $D_{4}$. Then $V_{27} \downarrow L'$ respectively has composition factor dimensions $8,8,8,1,1,1$ or $10,16,1$. Now, if there exists a non-$G$-cr subgroup $S \cong H$, then as well as a trivial composition factor, $S$ must have at least two composition factors of dimension 12, since these are self-dual and are the only factors with nonzero first cohomology group. Thus we are in Case 2) of Table \ref{l213e6p7}, and so $V_{27} \downarrow S = 1/12/14_b$. This is not compatible with the $L'$-composition factors of $L(G)$, hence all subgroups $S \cong H$ of $G$ are $G$-completely reducible.

For $(G,H,p) = (E_8,\Alt_{14},2)$, noting that $H$ does not preserve a nonzero quadratic form on its 12-dimensional irreducible module, we necessarily have $L'$ simple of type $D_7$, corresponding to a self-dual uniserial module of shape $1|12|1$. By Table \ref{a14e8p2}, we have $L(G) \downarrow S = 1^{8}/12^{4}/64_a/64_b^{2}$. Now, $R_u(P)$ has two levels, which are irreducible $L'$-modules of high weight $\lambda_6$ and $\lambda_1$. These restrict to $S$ as $64_b$ and a uniserial module of shape $1|12|1$, respectively. In each case we have $H^{1}(H,V) = \{0\}$, and hence $S$ is $G$-completely reducible.

For $(G,H,p) = (E_8,\Alt_{8},7)$, suppose first that $L'$ has type $E_{7}$. Then the non-trivial $L'$-modules occurring in the corresponding unipotent radical have high weight $\lambda_7$, and by Table \ref{a8e7p7} these we have $H^{1}(S,M) = 0$ for each $S$-composition factor $M$ of such module, hence $S$ is $G$-completely reducible in this case.

So now assume that $L'$ has no exceptional factors. Then $L'$ is simple of type $D_{7}$, $A_{6}$ or $D_{4}$, since the smallest non-trivial irreducible $H$-modules have dimensions $7$ and $14$. If $L'$ has type $D_{7}$ then $V_{D_7}(\lambda_1)$ occurs with multiplicity 2 as a composition factor of $L(G)$ (cf. \cite{MR1329942}*{Table 8.1}). Since $S$ is irreducible on this module, this contradicts the fact that $L(G) \downarrow S$ has at most one 14-dimensional composition factor (Table \ref{a8e8p7}). Similarly if $L'$ has type $D_{4}$ then $L'$ has at least $28$ trivial composition factors on $L(G)$, contradicting the feasible characters in Table \ref{a8e8p7}. Hence $L'$ has type $A_{6}$. By Lemma \ref{lem:radfilts}, the $L$-modules occurring in the corresponding unipotent radical $Q$ have high weights $\lambda_1$, $\lambda_2$, $\ldots$, $\lambda_6$, and these are respectively isomorphic to $\bigwedge^{1}(V_{A_6}(\lambda_1))$, $\bigwedge^{2}(V_{A_6}(\lambda_1))$, $\ldots$, $\bigwedge^{6}(V_{A_6}(\lambda_1))$. Using Magma to help with calculations, we find that the $S$-modules $\bigwedge^{1}(7)$, $\bigwedge^{2}(7)$, $\ldots$, $\bigwedge^{6}(7)$ are each irreducible, and in particular none of these involve a 19-dimensional composition factor. Hence $H^{1}(S,M) = 0$ for each $S$-composition factor of each $L'$-module occurring in the filtration of $Q$, and so $S$ is $G$-completely reducible.

For $(G,H,p) = (E_8,L_2(27),7)$, the only Levi subgroups of $G$ admitting an embedding of $H$ have derived subgroup $E_6$ or $E_7$. In the corresponding parabolic subgroups, the unipotent radicals have filtrations by modules of high weight $0$, $\lambda_1$ or $\lambda_6$ for $E_6$, or $\lambda_7$ for $E_7$. On the other hand, the composition factors of a subgroup $S \cong H$ on these modules are given by Tables \ref{l227e6p7} and \ref{l227e7p7}; they are always 1- or 13-dimensional, and hence their first cohomology groups vanish. Thus Proposition \ref{prop:levisub} applies to any such subgroup $S$, so all subgroups $S \cong H$ of $G$ are $G$-completely reducible.

For $(G,H,p) = (E_8,L_2(29),2)$, if $P = QL$ is a parabolic subgroup containing a copy of $H$, then $L'$ has type $E_7$. From Table \ref{l229e7p2}, the 56-dimensional module $V_{L'}(\lambda_7)$ restricts to a subgroup isomorphic to $H$ with 28-dimensional composition factor dimensions, and $H^{1}(H,V)$ vanishes for each. Now, $Q$ has two levels, which are irreducible $L'$-modules of respective high weights $\lambda_7$ and $0$. Thus $H^{1}(H,Q)$ vanishes and all subgroups $S \cong H$ of $G$ are $G$-completely reducible. \qed

\subsection{Remark: Existence of Non-Completely Reducible subgroups}
\label{rem:existence}
At this stage, we do not attempt the converse problem of classifying non-$G$-cr subgroups of each candidate type. In many cases, finding a non-$G$-cr subgroup isomorphic to $H$ is straightforward using the existence of indecomposable $H$-modules; for instance, since $\Alt_{5}$ has an indecomposable module of shape $1|4$ when $p = 3$, this gives a non-completely reducible embedding into $SL_5(K)$. A result of Serre \cite{MR2167207}*{Proposition 3.2} then implies that the image of $\Alt_{5}$ is non-$G$-completely reducible whenever $SL_5(K)$ is embedded as a Levi subgroup into $G = E_6(K)$, $E_7(K)$ or $E_8(K)$.

In general, however, classifying non-$G$-cr subgroups requires determining properties of the \emph{non-abelian cohomology set} $H^{1}(S,Q)$ when $S$ is a finite subgroup of $G$, and $Q$ is the unipotent radical of a parabolic subgroup containing $S$. When $Q$ is not abelian, this is not a group, only a pointed set, and its exact structure is not straightforward to determine. For instance, consider the case $(G,H,p) = (E_6,M_{22},2)$, where $3.M_{22}$ admits an embedding into $L' = SL_6(K)$ and is $L'$-irreducible, where $L$ a Levi subgroup of $G$. The corresponding unipotent radical has two levels, and $H^{1}(M_{22},Q)$ fits into an exact sequence of pointed sets:
\[ \{0\} \to H^{1}(M_{22},Q) \to H^{1}(M_{22},10|10^{*}) \to H^{2}(M_{22},K) \]
where $H^{1}(M_{22},10|10^{*}) \cong H^{2}(M_{22},K) \cong K$.

%% file: AJL-thetables.tex

\chapter{Tables of Feasible Characters} \label{chap:thetables}

\addtocontents{toc}{\setcounter{tocdepth}{-1}}
\section*{Notes on the Tables}
\addtocontents{toc}{\setcounter{tocdepth}{1}}
If $G$ is a simply connected exceptional simple algebraic group, and $H$ is quasisimple finite group which embeds into $G$, and $H/Z(H) \notin \textup{Lie}(p)$, then the following tables give all feasible characters of $H$ on the adjoint and minimal modules for $G$, such that $Z(H)$ acts as a group of scalars on $V_{\textup{min}}$. This therefore contains the composition factors of the possible restrictions $L(G) \downarrow \tilde{S}$ and $V_{\textup{min}} \downarrow \tilde{S}$, whenever $S$ is a simple subgroup of $G/Z(G)$ and $\tilde{S}$ is a minimal preimage of $S$ in $G$.

The relevant $KG$-modules are denoted as follows:
\begin{center}
\begin{tabular}{c|c|c}
$G$ & $L(G)$ & $V_{\textup{min}}$ \\ \hline
$E_8$ & $W(\lambda_1) = V_{248}$ & $V_{248}$ \\
$E_7$ & $W(\lambda_1) = V_{133}$ & $W(\lambda_7) = V_{56}$ \\
$E_6$ & $W(\lambda_2) = V_{78}$ & $W(\lambda_1) = V_{27}$ \\
$F_4$ & $W(\lambda_1) = V_{52}$ & $W(\lambda_4) = V_{26}$ \\
\end{tabular}
\end{center}

In calculating these tables, we have used information on elements of $\tilde{S}$ and $\tilde{G}$ of orders $2$ to $37$. A very small number of characters given here may be ruled out as occurring via an embedding $\tilde{S} \to \tilde{G}$ by consideration of elements of higher order. For example, $PSL_2(61)$ has two 30-dimensional irreducible modules in characteristic 2, whose Brauer characters differ only on elements of order 61.

If two sets of feasible characters differ only by permuting the module isomorphism types, we list only one member of each orbit. For example, line 1) of the first table for $\Alt_{6}$ overleaf corresponds to \emph{two} sets of compatible feasible characters of $\Alt_6$ on $L(G)$ and $V_G(\lambda_4)$, when $G = F_4$ and $p = 0$ or $p > 5$; these respectively have factors $8_b^{4}/10^{2}$ and $8_a^{4}/10^{2}$ on $L(G)$ (and similarly for $V_{\textup{min}}$). When we have shortened a table in this way, the permutations used will be noted underneath the table.

For $G$ not of type $E_6$, each irreducible $G$-module is self-dual, and hence each irreducible factor of a feasible character occurs with the same multiplicity as its dual. For $G$ of type $E_6$, note that a feasible character on $V_G(\lambda_1)$ gives rise to a feasible character on $V_G(\lambda_6) \cong V_G(\lambda_1)^{*}$ by the taking of duals. We therefore omit characters which arise by taking duals.

Finally, since several of the results of this paper follow by inspecting each of these tables, for convenience we attach labels to feasible characters satisfying certain conditions, as follows. Let $L$ and $V$ respectively be the sum of non-trivial $G$-composition factors of $L(G)$ and $V_{\textup{min}}$.

\begin{center}
\begin{tabularx}{\textwidth}{c|X}
\possprim & Possibly Lie \textbf{P}rimitive. A row has this label if and only if a subgroup $\tilde{S}$ of $G$ with these composition factors on $L(G)$ and $V_{\textup{min}}$ fails to satisfy the conditions of Proposition \ref{prop:substab}(i) and (iii) in its action on both $L$ and $V$. Every non-generic Lie primitive finite simple subgroup of $G$ must give rise to feasible characters with this label, though the converse is not true. \\ \hline
\nongcr & Possibly \textbf{N}on-$G$-cr. A feasible character has this label if and only if both:\\
& \ $\bullet$ \ $H$ admits an embedding into a proper Levi subgroup of $G$, and\\
& \ $\bullet$ \ The feasible character of $H$ on $L$ satisfies the conclusion of Corollary \ref{cor:levisub}, that is, it contains a trivial factor, as well as a factor $W$ such that $H^{1}(H,W)$ is nonzero, and either $W$ occurs with multiplicity at least two, or $W \ncong W^{*}$. Furthermore, $W$ has dimension at least $14$, $20$, $35$ or $64$, if the type of $G$ is respectively $F_{4}$, $E_{6}$, $E_{7}$ or $E_{8}$.\\
& Every non-$G$-cr non-generic finite simple subgroup of $G$ must give rise to feasible characters with this label, though the converse is not true.
\end{tabularx}
\end{center}
Thus a triple $(G,H,p)$ appears in Table \ref{tab:substabprop} (page \pageref{tab:substabprop}) if and only if none of the corresponding feasible characters are marked with \possprim, and the triple appears in Table \ref{tab:candidates} (page \pageref{tab:candidates}) if and only if some corresponding feasible character is marked with \nongcr.

\section{$F_{4}$} \label{sec:F4tabs}

\subsection{Alternating Groups} \ 

\begin{table}[H]\small
\caption{Alt$_{10} < F_4$, $p = 2$}
\begin{tabular}{r|cccc|cccc}
& \multicolumn{4}{c|}{$V_{52}$} & \multicolumn{4}{c}{$V_{26}$} \\
& 1 & 8 & 16 & 26 & 1 & 8 & 16 & 26 \\ \hline
1) & 2 & 1 & 1 & 1 & 0 & 0 & 0 & 1 \\
2) & 2 & 1 & 1 & 1 & 2 & 1 & 1 & 0
\end{tabular}
\label{a10f4p2}
\end{table}

\begin{table}[H]\small
\caption{Alt$_{9} < F_4$, $p = 2$}
\begin{tabular}{r|ccccc|ccccc}
& \multicolumn{5}{c|}{$V_{52}$} & \multicolumn{5}{c}{$V_{26}$} \\
& 1 & 8$_a$ & 8$_b$ & 8$_c$ & 26 & 1 & 8$_a$ & 8$_b$ & 8$_c$ & 26 \\ \hline
1) & 2 & 1 & 1 & 1 & 1 & 0 & 0 & 0 & 0 & 1 \\
2) & 2 & 1 & 1 & 1 & 1 & 2 & 1 & 1 & 1 & 0
\end{tabular}\\
Permutations: $(8_b,8_c)$.
\label{a9f4p2}
\end{table}

\begin{table}[H]\small
\caption{Alt$_{7} < F_4$, $p = 5$}
\begin{tabular}{r|cccc|cccc}
& \multicolumn{4}{c|}{$V_{52}$} & \multicolumn{4}{c}{$V_{26}$} \\
& 1 & 8 & 10 & 10$^{*}$ & 1 & 8 & 10 & 10$^{*}$  \\ \hline
1) & 0 & 4 & 1 & 1 & 2 & 3 & 0 & 0
\end{tabular}
\label{a7f4p5}
\end{table}

\begin{table}[H]\small
\caption{Alt$_{7} < F_4$, $p = 2$}
\begin{tabular}{r|cccccc|cccccc}
& \multicolumn{6}{c|}{$V_{52}$} & \multicolumn{6}{c}{$V_{26}$} \\
& 1 & 4 & 4$^{*}$ & 6 & 14 & 20 & 1 & 4 & 4$^{*}$ & 6 & 14 & 20 \\ \hline
1) & 4 & 2 & 2 & 2 & 0 & 1 & 0 & 0 & 0 & 1 & 0 & 1 \\
2) & 4 & 2 & 2 & 3 & 1 & 0 & 0 & 0 & 0 & 2 & 1 & 0 \\
3) & 4 & 2 & 2 & 2 & 0 & 1 & 4 & 2 & 2 & 1 & 0 & 0 \\
4) & 4 & 2 & 2 & 3 & 1 & 0 & 4 & 2 & 2 & 1 & 0 & 0 \\
5) & 4 & 2 & 2 & 3 & 1 & 0 & 0 & 0 & 0 & 1 & 0 & 1
\end{tabular}
\label{a7f4p2}
\end{table}

\begin{table}[H]\small
\caption{Alt$_{6} < F_4$, $p = 0$ or $p \nmid |H|$}
\begin{tabular}{rc|ccccccc|ccccccc}
& & \multicolumn{7}{c|}{$V_{52}$} & \multicolumn{7}{c}{$V_{26}$} \\
& & 1 & 5$_a$ & 5$_b$ & 8$_a$ & 8$_b$ & 9 & 10 & 1 & 5$_a$ & 5$_b$ & 8$_a$ & 8$_b$ & 9 & 10 \\ \hline
1) & & 0 & 0 & 0 & 0 & 4 & 0 & 2 & 2 & 0 & 0 & 0 & 3 & 0 & 0 \\
2) & \possprim & 0 & 0 & 0 & 1 & 3 & 0 & 2 & 0 & 0 & 0 & 1 & 0 & 2 & 0 \\
3) & & 0 & 0 & 0 & 2 & 2 & 0 & 2 & 1 & 0 & 0 & 1 & 1 & 1 & 0 \\
4) & & 2 & 1 & 1 & 0 & 0 & 0 & 4 & 1 & 0 & 0 & 1 & 1 & 1 & 0 \\
5) & & 3 & 0 & 0 & 0 & 0 & 1 & 4 & 0 & 1 & 1 & 1 & 1 & 0 & 0 \\
6) & & 3 & 0 & 0 & 0 & 0 & 1 & 4 & 1 & 0 & 3 & 0 & 0 & 0 & 1
\end{tabular}\\
Permutations: $(5_a,5_b)$, $(8_a,8_b)$.
\label{a6f4p0}
\end{table}

\begin{table}[H]\small
\caption{Alt$_{6} < F_4$, $p = 5$}
\begin{tabular}{rc|ccccc|ccccc}
& & \multicolumn{5}{c|}{$V_{52}$} & \multicolumn{5}{c}{$V_{26}$} \\
& & 1 & 5$_a$ & 5$_b$ & 8 & 10 & 1 & 5$_a$ & 5$_b$ & 8 & 10 \\ \hline
1) & \possprim & 0 & 0 & 0 & 4 & 2 & 2 & 0 & 0 & 3 & 0 \\
2) & & 2 & 1 & 1 & 0 & 4 & 2 & 0 & 0 & 3 & 0 \\
3) & & 4 & 0 & 0 & 1 & 4 & 0 & 1 & 1 & 2 & 0 \\
4) & & 4 & 0 & 0 & 1 & 4 & 1 & 0 & 3 & 0 & 1
\end{tabular}\\
Permutations: $(5_a,5_b)$.
\label{a6f4p5}
\end{table}

\begin{table}[H]\small
\caption{Alt$_{5} < F_4$, $p = 0$ or $p \nmid |H|$}
\begin{tabular}{rc|ccccc|ccccc}
& & \multicolumn{5}{c|}{$V_{52}$} & \multicolumn{5}{c}{$V_{26}$} \\
& & $1$ & $3_a$ & $3_b$ & $4$ & $5$ & $1$ & $3_a$ & $3_b$ & $4$ & $5$\\ \hline
1) & \possprim & 0 & 3 & 5 & 2 & 4 & 0 & 1 & 0 & 2 & 3\\
2) & & 3 & 3 & 5 & 5 & 1 & 0 & 1 & 0 & 2 & 3\\
3) & & 14& 1 & 0 & 0 & 7 & 0 & 7 & 0 & 0 & 1\\
4) & & 0 & 4 & 4 & 2 & 4 & 1 & 1 & 1 & 1 & 3\\
5) & & 3 & 4 & 4 & 5 & 1 & 1 & 1 & 1 & 1 & 3\\
6) & & 0 & 6 & 2 & 2 & 4 & 2 & 3 & 0 & 0 & 3\\
7) & & 3 & 6 & 2 & 5 & 1 & 2 & 3 & 0 & 0 & 3\\
8) & & 3 & 3 & 5 & 5 & 1 & 3 & 1 & 0 & 5 & 0\\
9) & & 3 & 4 & 4 & 5 & 1 & 4 & 1 & 1 & 4 & 0\\
10)& & 3 & 6 & 2 & 5 & 1 & 5 & 3 & 0 & 3 & 0\\
11)& & 8 & 0 & 13& 0 & 1 & 8 & 0 & 6 & 0 & 0
\end{tabular}\\
Permutations: $(3_a,3_b)$.
\label{a5f4p0}
\end{table}

\begin{table}[H]\small
\caption{Alt$_{5} < F_4$, $p = 3$}
\begin{tabular}{rc|cccc|cccc}
& & \multicolumn{4}{c|}{$V_{52}$} & \multicolumn{4}{c}{$V_{26}$}\\
& & 1 & 3$_a$ & 3$_b$ & 4 & 1 & 3$_a$ & 3$_b$ & 4 \\ \hline
1) & \nongcr & 21 & 1 & 0 & 7 & 1 & 7 & 0 & 1\\
2) & \possprim, \nongcr & 4 & 3 & 5 & 6 & 3 & 1 & 0 & 5\\
3) & \possprim, \nongcr & 4 & 4 & 4 & 6 & 4 & 1 & 1 & 4\\
4) & \nongcr & 4 & 6 & 2 & 6 & 5 & 3 & 0 & 3\\
5) & & 9 & 0 & 13 & 1 & 8 & 0 & 6 & 0 
\end{tabular}\\
Permutations: $(3_a,3_b)$.
\label{a5f4p3}
\end{table}

\subsection{Sporadic Groups} \ 

\begin{table}[H]\small
\caption{$M_{11} < F_4$, $p = 11$}
\begin{tabular}{r|ccc|ccc}
& \multicolumn{3}{c|}{$V_{52}$} & \multicolumn{3}{c}{$V_{26}$} \\
& 10 & $10^{*}$ & 16 & 1 & 9 & 16 \\ \hline
1) & 1 & 1 & 2 & 1 & 1 & 1
\end{tabular}
\label{m11f4p11}
\end{table}

\begin{table}[H]\small
\caption{$J_{1} < F_4$, $p = 11$}
\begin{tabular}{r|ccc|cc}
& \multicolumn{3}{c|}{$V_{52}$} & \multicolumn{2}{c}{$V_{26}$} \\
& 1 & 7 & 14 & 1 & 7 \\ \hline
1) & 3 & 5 & 1 & 5 & 3
\end{tabular}
\label{j1f4p11}
\end{table}

\begin{table}[H]\small
\caption{$J_{2} < F_4$, $p = 2$}
\begin{tabular}{rc|ccccc|ccccc}
& & \multicolumn{5}{c|}{$V_{52}$} & \multicolumn{5}{c}{$V_{26}$} \\
& &  1 & $6_a$ & $6_b$ & $14_a$ & $14_b$ & 1 & $6_a$ & $6_b$ & $14_a$ & $14_b$ \\ \hline
1) & \nongcr & 8 & 2 & 3 & 1 & 0 & 0 & 2 & 0 & 1 & 0 \\
2) & \nongcr & 8 & 5 & 0 & 1 & 0 & 8 & 3 & 0 & 0 & 0 
\end{tabular}
\\
Permutations: $(6_a,6_b)(14_a,14_b)$.
\label{j2f4p2}
\end{table}

\subsection{Cross-characteristic Groups $L_2(q)$ $(q \neq 4$, $5$, $9)$} \ 

\begin{table}[H]\small
\caption{$L_2(7) < F_4$, $p = 0$ or $p \nmid |H|$}
\begin{tabular}{rc|cccccc|ccccccc}
& & \multicolumn{6}{c|}{$V_{52}$} & \multicolumn{6}{c}{$V_{26}$} \\
& & 1 & 3 & $3^{*}$ & 6 & 7 & 8 & 1 & 3 & $3^{*}$ & 6 & 7 & 8 \\ \hline
1) & \possprim & 0 & 1 & 1 & 0 & 2 & 4 & 0 & 1 & 1 & 2 & 0 & 1 \\
2) & & 0 & 1 & 1 & 0 & 2 & 4 & 2 & 0 & 0 & 0 & 0 & 3 \\
3) & & 3 & 1 & 1 & 0 & 5 & 1 & 0 & 1 & 1 & 2 & 0 & 1 \\
4) & & 3 & 1 & 1 & 0 & 5 & 1 & 2 & 0 & 0 & 0 & 0 & 3 \\
5) & & 3 & 1 & 1 & 0 & 5 & 1 & 5 & 0 & 0 & 0 & 3 & 0 \\
6) & & 8 & 0 & 0 & 6 & 0 & 1 & 0 & 3 & 3 & 0 & 0 & 1 \\
7) & & 8 & 6 & 6 & 0 & 0 & 1 & 8 & 3 & 3 & 0 & 0 & 0
\end{tabular}
\label{l27f4p0}
\end{table}

\begin{table}[H]\small
\caption{$L_2(7) < F_4$, $p = 3$}
\begin{tabular}{rc|ccccc|cccccc}
& & \multicolumn{5}{c|}{$V_{52}$} & \multicolumn{5}{c}{$V_{26}$} \\
& & 1 & 3 & $3^{*}$ & $6_a$ & 7 & 1 & 3 & $3^{*}$ & $6_a$ & 7 \\ \hline
1) & \possprim, \nongcr & 4 & 1 & 1 & 0 & 6 & 1 & 1 & 1 & 2 & 1 \\
2) & \nongcr & 4 & 1 & 1 & 0 & 6 & 5 & 0 & 0 & 0 & 3 \\
3) & \nongcr & 5 & 2 & 2 & 0 & 5 & 0 & 0 & 0 & 2 & 2 \\
4) & & 9 & 0 & 0 & 6 & 1 & 1 & 3 & 3 & 0 & 1 \\
5) & & 9 & 6 & 6 & 0 & 1 & 8 & 3 & 3 & 0 & 0
\end{tabular}
\label{l27f4p3}
\end{table}

\begin{table}[H]\small
\caption{$L_2(8) < F_4$, $p = 0$ or $p \nmid |H|$}
\begin{tabular}{rc|*{9}{c}|*{9}{c}}
& & \multicolumn{9}{c|}{$V_{52}$} & \multicolumn{9}{c}{$V_{26}$} \\
& & 1 & $7_a$ & $7_b$ & $7_c$ & $7_d$ & 8 & $9_a$ & $9_b$ & $9_c$ & 1 & $7_a$ & $7_b$ & $7_c$ & $7_d$ & 8 & $9_a$ & $9_b$ & $9_c$ \\ \hline
1) & \possprim & 0 & 2 & 1 & 1 & 1 & 1 & 0 & 0 & 1 & 0 & 0 & 0 & 0 & 0 & 1 & 1 & 0 & 1 \\
2) & & 0 & 2 & 1 & 1 & 1 & 1 & 0 & 0 & 1 & 1 & 1 & 0 & 0 & 0 & 0 & 1 & 0 & 1 \\
3) & & 1 & 2 & 1 & 1 & 1 & 2 & 0 & 0 & 0 & 2 & 0 & 0 & 0 & 0 & 3 & 0 & 0 & 0 \\
4) & & 1 & 2 & 1 & 1 & 1 & 2 & 0 & 0 & 0 & 3 & 1 & 0 & 0 & 0 & 2 & 0 & 0 & 0 \\
5) & & 3 & 1 & 0 & 5 & 1 & 0 & 0 & 0 & 0 & 5 & 0 & 0 & 3 & 0 & 0 & 0 & 0 & 0 
\end{tabular}\\
Permutations: $(7_c,7_c,7_d)(9_a,9_b,9_c)$.
\label{l28f4p0}
\end{table}

\begin{table}[H]\small
\caption{$L_2(8) < F_4$, $p = 7$}
\begin{tabular}{rc|cccccc|cccccc}
& & \multicolumn{6}{c|}{$V_{52}$} & \multicolumn{6}{c}{$V_{26}$} \\
& & 1 & $7_a$ & $7_b$ & $7_c$ & $7_d$ & 8 & 1 & $7_a$ & $7_b$ & $7_c$ & $7_d$ & 8 \\ \hline
1) & \possprim, \nongcr & 1 & 2 & 1 & 1 & 1 & 2 & 2 & 0 & 0 & 0 & 0 & 3 \\
2) & \nongcr & 1 & 2 & 1 & 1 & 1 & 2 & 3 & 1 & 0 & 0 & 0 & 2 \\
3) & & 3 & 1 & 0 & 1 & 5 & 0 & 5 & 0 & 0 & 0 & 3 & 0
\end{tabular}\\
Permutations: $(7_b,7_c,7_d)$.
\label{l28f4p7}
\end{table}

\begin{table}[H]\small
\caption{$L_2(8) < F_4$, $p = 3$}
\begin{tabular}{rc|ccccc|ccccc}
& & \multicolumn{5}{c|}{$V_{52}$} & \multicolumn{5}{c}{$V_{26}$} \\
& & 1 & 7 & $9_a$ & $9_b$ & $9_c$ & 1 & 7 & $9_a$ & $9_b$ & $9_c$ \\ \hline
1) & \possprim, \nongcr & 1 & 6 & 0 & 0 & 1 & 1 & 1 & 0 & 1 & 1 \\
2) & \nongcr & 3 & 7 & 0 & 0 & 0 & 5 & 3 & 0 & 0 & 0
\end{tabular}\\
Permutations: $(9_a,9_b,9_c)$.
\label{l28f4p3}
\end{table}

\begin{table}[H]\small
\caption{$L_2(13) < F_4$, $p = 0$ or $p \nmid |H|$}
\begin{tabular}{rc|*{7}{c}|*{7}{c}}
& & \multicolumn{7}{c|}{$V_{52}$} & \multicolumn{7}{c}{$V_{26}$} \\
& & 1 & $7_a$ & $7_b$ & $12_a$ & $12_b$ & $12_c$ & $14_a$ & 1 & $7_a$ & $7_b$ & $12_a$ & $12_b$ & $12_c$ & $14_b$ \\ \hline
1) & \possprim & 0 & 0 & 0 & 0 & 1 & 1 & 2 & 0 & 0 & 0 & 0 & 1 & 0 & 1 \\
2) & & 3 & 0 & 5 & 0 & 0 & 0 & 1 & 5 & 0 & 3 & 0 & 0 & 0 & 0
\end{tabular}\\
Permutations: $(7_a,7_b)$, $(12_{a},12_{b},12_{c})$.
\label{l213f4p0}
\end{table}

\begin{table}[H]\small
\caption{$L_2(13) < F_4$, $p = 7$}
\begin{tabular}{rc|*{6}{c}|*{6}{c}}
& & \multicolumn{6}{c|}{$V_{52}$} & \multicolumn{6}{c}{$V_{26}$} \\
& & 1 & $7_a$ & $7_b$ & $12$ & $14_a$ & $14_b$ & 1 & $7_a$ & $7_b$ & $12$ & $14_a$ & $14_b$ \\ \hline
1) & \possprim & 0 & 0 & 0 & 2 & 2 & 0 & 0 & 0 & 0 & 1 & 0 & 1 \\
2) & & 3 & 0 & 5 & 0 & 1 & 0 & 5 & 0 & 3 & 0 & 0 & 0
\end{tabular}\\
Permutations: $(7_a,7_b)$.
\label{l213f4p7}
\end{table}

\begin{table}[H]\small
\caption{$L_2(13) < F_4$, $p = 3$}
\begin{tabular}{rc|*{7}{c}|*{7}{c}}
& & \multicolumn{7}{c|}{$V_{52}$} & \multicolumn{7}{c}{$V_{26}$} \\
& & 1 & $7_a$ & $7_b$ & $12_a$ & $12_b$ & $12_c$ & $13$ & 1 & $7_a$ & $7_b$ & $12_a$ & $12_b$ & $12_c$ & $13$ \\ \hline
1) & \possprim & 0 & 2 & 2 & 0 & 1 & 1 & 0 & 1 & 0 & 0 & 0 & 1 & 0 & 1 \\
2) & & 3 & 1 & 6 & 0 & 0 & 0 & 0 & 5 & 0 & 3 & 0 & 0 & 0 & 0
\end{tabular}\\
Permutations: $(7_a,7_b)$, $(12_{a},12_{b},12_{c})$.
\label{l213f4p3}
\end{table}

\begin{table}[H]\small
\caption{$L_2(13) < F_4$, $p = 2$}
\begin{tabular}{rc|*{7}{c}|*{7}{c}}
& & \multicolumn{7}{c|}{$V_{52}$} & \multicolumn{7}{c}{$V_{26}$} \\
& & 1 & $6_a$ & $6_b$ & $12_a$ & $12_b$ & $12_c$ & 14 & 1 & $6_a$ & $6_b$ & $12_a$ & $12_b$ & $12_c$ & 14 \\ \hline
1) & \possprim & 0 & 0 & 0 & 0 & 1 & 1 & 2 & 0 & 0 & 0 & 0 & 0 & 1 & 1 \\
2) & \nongcr & 8 & 0 & 5 & 0 & 0 & 0 & 1 & 8 & 0 & 3 & 0 & 0 & 0 & 0 \\
3) & \nongcr & 8 & 2 & 3 & 0 & 0 & 0 & 1 & 0 & 2 & 0 & 0 & 0 & 0 & 1
\end{tabular}\\
Permutations: $(6_a,6_b)$, $(12_{a},12_{b},12_{c})$.
\label{l213f4p2}
\end{table}

\begin{table}[H]\small
\caption{$L_2(17) < F_4$, $p = 0$ or $p \nmid |H|$}
\begin{tabular}{rc|*{3}{c}|*{5}{c}}
& & \multicolumn{3}{c|}{$V_{52}$} & \multicolumn{5}{c}{$V_{26}$} \\
& & $16_a$ & $18_b$ & $18_c$ & 1 & $9_a$ & $9_b$ & $16_a$ & $17$ \\ \hline
1) & \possprim & 1 & 1 & 1 & 0 & 0 & 1 & 0 & 1 \\
2) & & 1 & 1 & 1 & 1 & 0 & 1 & 1 & 0
\end{tabular}\\
Permutations: $(9_a,9_b)$.
\label{l217f4p0}
\end{table}

\begin{table}[H]\small
\caption{$L_2(17) < F_4$, $p = 3$}
\begin{tabular}{rc|*{4}{c}|*{4}{c}}
& & \multicolumn{4}{c|}{$V_{52}$} & \multicolumn{4}{c}{$V_{26}$} \\
& & 16 & $18_a$ & $18_b$ & $18_c$ & 1 & $9_a$ & $9_b$ & 16 \\ \hline
1) & \possprim & 1 & 0 & 1 & 1 & 1 & 0 & 1 & 1
\end{tabular}\\
Permutations: $(9_a,9_b)$.
\label{l217f4p3}
\end{table}

\begin{table}[H]\small
\caption{$L_2(17) < F_4$, $p = 2$}
\begin{tabular}{r|*{4}{c}|*{4}{c}}
& \multicolumn{4}{c|}{$V_{52}$} & \multicolumn{4}{c}{$V_{26}$} \\
& 1 & $8_a$ & $8_b$ & $16_a$ & 1 & $8_a$ & $8_b$ & $16_a$ \\ \hline
1) & 4 & 2 & 2 & 1 & 2 & 0 & 1 & 1 \\
2) & 4 & 2 & 2 & 1 & 2 & 1 & 2 & 0
\end{tabular}\\
Permutations: $(8_a,8_b)$.
\label{l217f4p2}
\end{table}

\begin{table}[H]\small
\caption{$L_2(25) < F_4$, $p \neq 2,3,5$}
\begin{tabular}{rc|*{2}{c}|*{1}{c}}
& & \multicolumn{2}{c|}{$V_{52}$} & \multicolumn{1}{c}{$V_{26}$} \\
& & $26_b$ & $26_c$ & $26_a$ \\ \hline
1) & \possprim & 1 & 1 & 1
\end{tabular}\\
$26_b$, $26_c$ are $\textup{Aut}(L_2(25))$-conjugate.
\label{l225f4p0}
\label{l225f4p13}
\end{table}

\begin{table}[H]\small
\caption{$L_2(25) < F_4$, $p = 3$}
\begin{tabular}{rc|*{1}{c}|*{2}{c}}
& & \multicolumn{1}{c|}{$V_{52}$} & \multicolumn{2}{c}{$V_{26}$} \\
& & 26 & 1 & 25 \\ \hline
1) & \possprim & 2 & 1 & 1
\end{tabular}
\label{l225f4p3}
\end{table}

\begin{table}[H]\small
\caption{$L_2(25) < F_4$, $p = 2$}
\begin{tabular}{rc|*{1}{c}|*{1}{c}}
& & \multicolumn{1}{c|}{$V_{52}$} & \multicolumn{1}{c}{$V_{26}$} \\
& & 26 & 26 \\ \hline
1) & \possprim & 2 & 1
\end{tabular}
\label{l225f4p2}
\end{table}

\begin{table}[H]\small
\caption{$L_2(27) < F_4$, $p \neq 2,3,7$}
\begin{tabular}{rc|*{3}{c}|*{3}{c}}
& & \multicolumn{3}{c|}{$V_{52}$} & \multicolumn{3}{c}{$V_{26}$} \\
& & $26_a$ & $26_b$ & $26_c$ & $26_d$ & $26_e$ & $26_f$ \\ \hline
1) & \possprim & 0 & 1 & 1 & 1 & 0 & 0
\end{tabular}\\
Permutations: $(26_a,26_b,26_c)(26_d,26_e,26_f)$.
\label{l227f4p0}
\label{l227f4p13}
\end{table}

\begin{table}[H]\small
\caption{$L_2(27) < F_4$, $p = 7$}
\begin{tabular}{rc|c|cc}
& & \multicolumn{1}{c|}{$V_{52}$} & \multicolumn{2}{c}{$V_{26}$} \\
& & 26 & $13_a$ & $13_b$ \\ \hline
1) & \possprim & 2 & 1 & 1
\end{tabular}
\label{l227f4p7}
\end{table}

\begin{table}[H]\small
\caption{$L_2(27) < F_4$, $p = 2$}
\begin{tabular}{rc|*{2}{c}|*{1}{c}}
& & \multicolumn{2}{c|}{$V_{52}$} & \multicolumn{1}{c}{$V_{26}$} \\
& & $26_a$ & $26_b$ & $26_a$ \\ \hline
1) & \possprim & 1 & 1 & 1 \\
\end{tabular}\\
Permutations: $(26_a,26_b,26_c)$.
\label{l227f4p2}
\end{table}

\subsection{Cross-characteristic Groups $\ncong L_2(q)$} \ 

\begin{table}[H]\small
\caption{$L_3(3) < F_4$, $p \neq 2, 3$}
\begin{tabular}{rc|*{2}{c}|*{1}{c}}
& & \multicolumn{2}{c|}{$V_{52}$} & \multicolumn{1}{c}{$V_{26}$} \\
& & $26_b$ & $(26_b)^{*}$ & $26_a$ \\ \hline
1) & \possprim & 1 & 1 & 1
\end{tabular}
\label{l33f4p13}
\label{l33f4p0}
\end{table}

\begin{table}[H]\small
\caption{$L_3(3) < F_4$, $p = 2$}
\begin{tabular}{rc|c|c}
& & \multicolumn{1}{c|}{$V_{52}$} & \multicolumn{1}{c}{$V_{26}$} \\
& & 26 & 26 \\ \hline
1) & \possprim & 2 & 1
\end{tabular}
\label{l33f4p2}
\end{table}

\begin{table}[H]\small
\caption{$L_4(3) < F_4$, $p = 2$}
\begin{tabular}{rc|cc|c}
& & \multicolumn{2}{c|}{$V_{52}$} & \multicolumn{1}{c}{$V_{26}$} \\
& & $26_a$ & $26_b$ & $26_a$ \\ \hline
1) & \possprim & 1 & 1 & 1
\end{tabular}\\
Permutations: $(26_a,26_b)$.
\label{l43f4p2}
\end{table}

\begin{table}[H]\small
\caption{$U_3(3) < F_4$, $p = 0$ or $p \nmid |H|$}
\begin{tabular}{r|*{7}{c}|*{4}{c}}
& \multicolumn{7}{c|}{$V_{52}$} & \multicolumn{4}{c}{$V_{26}$} \\
& 1 & 6 & $7_a$ & $(7_a)^{*}$ & $7_b$ & 14 & $21_b$ & 1 & 6 & $7_b$ & 14 \\ \hline
1) & 3 & 0 & 0 & 0 & 5 & 1 & 0 & 5 & 0 & 3 & 0 \\
2) & 3 & 0 & 2 & 2 & 0 & 0 & 1 & 0 & 2 & 0 & 1
\end{tabular}
\label{u33f4p0}
\end{table}

\begin{table}[H]\small
\caption{$U_3(3) < F_4$, $p = 7$}
\begin{tabular}{rc|*{8}{c}|*{5}{c}}
& & \multicolumn{8}{c|}{$V_{52}$} & \multicolumn{5}{c}{$V_{26}$} \\
& & 1 & 6 & $7_a$ & $7_b$ & $(7_b)^{*}$ & 14 & $21_a$ & 26 & 1 & 6 & $7_a$ & 14 & 26 \\ \hline
1) & \possprim & 0 & 2 & 0 & 0 & 0 & 1 & 0 & 1 & 0 & 0 & 0 & 0 & 1 \\
2) & & 3 & 0 & 0 & 2 & 2 & 0 & 1 & 0 & 0 & 2 & 0 & 1 & 0 \\
3) & & 3 & 0 & 5 & 0 & 0 & 1 & 0 & 0 & 5 & 0 & 3 & 0 & 0
\end{tabular}
\label{u33f4p7}
\end{table}

\addtocounter{table}{1}
\begin{table}[H]\small
\begin{center}
\arabic{chapter}.\arabic{section}.\arabic{table}: $^{3}D_4(2) < F_4$, $p \neq 2, 3$. Irreducible on $V_{52}$ and $V_{26}$. \possprim
\end{center}
\label{3d4f4p7}
\label{3d4f4p13}
\label{3d4f4p0}
\end{table}

\begin{table}[H]\small
\caption{$^{3}D_4(2) < F_4$, $p = 3$}
\begin{tabular}{rc|*{1}{c}|*{2}{c}}
& & \multicolumn{1}{c|}{$V_{52}$} & \multicolumn{2}{c}{$V_{26}$} \\
& & 52 & 1 & 25 \\ \hline
1) & \possprim & 1 & 1 & 1
\end{tabular}
\label{3d4f4p3}
\end{table}

\section{$E_6$} \label{sec:E6tabs}

\subsection{Alternating Groups} \ 

\begin{table}[H]\small
\caption{Alt$_{12} < E_6$, $p = 2$}
\begin{tabular}{r|*{4}{c}|*{4}{c}}
& \multicolumn{4}{c|}{$V_{78}$} & \multicolumn{4}{c}{$V_{27}$} \\
& 1 & 16 & 16$^{*}$ & 44 & 1 & 10 & 16 & 16$^{*}$\\ \hline
1) & 2 & 1 & 1 & 1 & 1 & 1 & 0 & 1
\end{tabular}
\label{a12e6p2}
\end{table}

\begin{table}[H]\small
\caption{Alt$_{11} < E_6$, $p = 2$}
\begin{tabular}{r|cccc|cccc}
& \multicolumn{4}{c|}{$V_{78}$} & \multicolumn{4}{c}{$V_{27}$} \\
& 1 & 16 & 16$^{*}$ & 44 & 1 & 10 & 16 & 16$^{*}$ \\ \hline
1) & 2 & 1 & 1 & 1 & 1 & 1 & 1 & 0
\end{tabular}
\label{a11e6p2}
\end{table}

\begin{table}[H]\small
\caption{Alt$_{10} < E_6$, $p = 2$}
\begin{tabular}{r|cccc|cccc}
& \multicolumn{4}{c|}{$V_{78}$} & \multicolumn{4}{c}{$V_{27}$} \\
& 1 & 8 & 16 & 26 & 1 & 8 & 16 & 26 \\ \hline
1) & 2 & 1 & 1 & 2 & 1 & 0 & 0 & 1\\
2) & 4 & 2 & 2 & 1 & 3 & 1 & 1 & 0
\end{tabular}
\label{a10e6p2}
\end{table}

\begin{table}[H]\small
\caption{Alt$_{9} < E_6$, $p = 2$}
\begin{tabular}{rc|ccccc|ccccc}
& & \multicolumn{5}{c|}{$V_{78}$} & \multicolumn{5}{c}{$V_{27}$} \\
& & 1 & 8$_a$ & 8$_b$ & 8$_c$ & 26 & 1 & 8$_a$ & 8$_b$ & 8$_c$ & 26 \\ \hline
1) & \possprim & 2 & 1 & 1 & 1 & 2 & 1 & 0 & 0 & 0 & 1\\
2) & & 4 & 2 & 2 & 2 & 1 & 3 & 1 & 1 & 1 & 0
\end{tabular}\\
Permutations: $(8_b,8_c)$.
\label{a9e6p2}
\end{table}

\begin{table}[H]\small
\caption{Alt$_{7} < E_6$, $p = 0$ or $p \nmid |H|$}
\begin{tabular}{r|*{7}{c}|*{3}{c}}
& \multicolumn{7}{c|}{$V_{78}$} & \multicolumn{3}{c}{$V_{27}$} \\
& 1 & $6_a$ & 10 & 10$^{*}$ & 14$_a$ & 14$_b$ & $15_a$ & 1 & $6_a$ & $15_a$ \\ \hline
1) & 3 & 1 & 2 & 2 & 1 & 0 & 1 & 0 & 2 & 1
\end{tabular}\\
Permutations: $(10,10^{*})$. \\
$14_a$ is a section of $6_a \otimes 6_a$.
\label{a7e6p0}
\end{table}

\begin{table}[H]\small
\caption{$3\cdot$Alt$_{7} < E_6$, $p = 0$ or $p \nmid |H|$}
\begin{tabular}{r|*{8}{c}|*{2}{c}}
& \multicolumn{8}{c|}{$V_{78}$} & \multicolumn{2}{c}{$V_{27}$} \\
& 1 & 6$_a$ & 10 & 10$^{*}$ & 14$_a$ & 14$_b$ & 15$_a$ & 21$_a$ & 6$_b$ & $15_c$ \\ \hline
1) & 3 & 0 & 2 & 2 & 1 & 0 & 0 & 1 & 2 & 1
\end{tabular}\\
$14_a$ is a section of $6_a \otimes 6_a$.
\label{3a7e6p0}
\end{table}

\begin{table}[H]\small
\caption{Alt$_{7} < E_6$, $p = 7$}
\begin{tabular}{rc|ccccc|ccc}
& & \multicolumn{5}{c|}{$V_{78}$} & \multicolumn{3}{c}{$V_{27}$} \\
& & 1 & 5 & 10 & 14$_a$ & 14$_b$ & 1 & 5 & 10 \\ \hline
1) & \nongcr & 4 & 2 & 5 & 1 & 0 & 2 & 3 & 1
\end{tabular}\\
$14_a$ is a section of $5 \otimes 5$.
\label{a7e6p7}
\end{table}

\begin{table}[H]\small
\caption{$3\cdot$Alt$_{7} < E_6$, $p = 7$}
\begin{tabular}{r|ccccccc|*{6}{c}}
& \multicolumn{7}{c|}{$V_{78}$} & \multicolumn{6}{c}{$V_{27}$} \\
& 1 & 5 & 10 & 14$_a$ & 14$_b$ & 21 & 35 & 6 & 6$^{*}$ & 9 & 9$^{*}$ & 15 & 15$^{*}$\\ \hline
1) & 3 & 0 & 4 & 1 & 0 & 1 & 0 & 0 & 2 & 0 & 0 & 0 & 1
\end{tabular}\\
$14_a$ is a section of $5 \otimes 5$.
\label{3a7e6p7}
\end{table}

\begin{table}[H]\small
\caption{Alt$_{7} < E_6$, $p = 5$}
\begin{tabular}{rc|ccccccc|ccccc}
& & \multicolumn{7}{c|}{$V_{78}$} & \multicolumn{5}{c}{$V_{27}$} \\
& & 1 & 6 & 8 & 10 & 10$^{*}$ & 15 & 35 & 1 & 6 & 8 & 13 & 15 \\ \hline
1) & \possprim & 0 & 0 & 1 & 0 & 0 & 0 & 2 & 0 & 1 & 1 & 1 & 0 \\
2) & & 3 & 2 & 1 & 2 & 2 & 1 & 0 & 0 & 2 & 0 & 0 & 1 \\
3) & & 2 & 0 & 7 & 1 & 1 & 0 & 0 & 3 & 0 & 3 & 0 & 0 
\end{tabular}
\label{a7e6p5}
\end{table}

\begin{table}[H]\small
\caption{$3\cdot$Alt$_{7} < E_6$, $p = 5$}
\begin{tabular}{rc|cccccccc|*{5}{c}}
& & \multicolumn{8}{c|}{$V_{78}$} & \multicolumn{5}{c}{$V_{27}$} \\
& & 1 & 6$_a$ & 8 & 10 & 10$^{*}$ & 13 & 15$_a$ & 35 & $3$ & 6$_b$ & $(6_b)^{*}$ & 15$_b$ & 21 \\ \hline
1) & \possprim & 0 & 0 & 1 & 0 & 0 & 0 & 0 & 2 & 0 & 0 & 1 & 0 & 1 \\
2) & & 3 & 1 & 2 & 2 & 2 & 1 & 0 & 0 & 0 & 0 & 2 & 1 & 0 \\
3) & & 14 & 0 & 8 & 0 & 0 & 0 & 0 & 0 & 7 & 1 & 0 & 0 & 0
\end{tabular}\\
$15_b$ is a section of $6_a \otimes 6_b$.
\label{3a7e6p5}
\end{table}

\begin{table}[H]\small
\caption{Alt$_{7} < E_6$, $p = 3$}
\begin{tabular}{rc|cccccc|cccccc}
& & \multicolumn{6}{c|}{$V_{78}$} & \multicolumn{6}{c}{$V_{27}$} \\
& & 1 & 6 & 10 & 10$^{*}$ & 13 & 15 & 1 & 6 & 10 & 10$^{*}$ & 13 & 15 \\ \hline
1) & \possprim, \nongcr & 4 & 1 & 2 & 2 & 1 & 1 & 0 & 2 & 0 & 0 & 0 & 1
\end{tabular}
\label{a7e6p3}
\end{table}

\begin{table}[H]\small
\caption{Alt$_{7} < E_6$, $p = 2$}
\begin{tabular}{rc|cccccc|cccccc}
& & \multicolumn{6}{c|}{$V_{78}$} & \multicolumn{6}{c}{$V_{27}$} \\
& & 1 & 4 & 4$^{*}$ & $6_a$ & 14 & 20 & 1 & 4 & 4$^{*}$ & $6_a$ & 14 & 20 \\ \hline
1) & \nongcr & 4 & 2 & 2 & 3 & 0 & 2 & 1 & 0 & 0 & 1 & 0 & 1\\
2) & & 4 & 2 & 2 & 4 & 1 & 1 & 1 & 0 & 0 & 1 & 0 & 1\\
3) & \nongcr & 4 & 2 & 2 & 5 & 2 & 0 & 1 & 0 & 0 & 2 & 1 & 0\\
4) & & 8 & 4 & 4 & 3 & 0 & 1 & 5 & 2 & 2 & 1 & 0 & 0\\
5) & & 8 & 4 & 4 & 4 & 1 & 0 & 5 & 2 & 2 & 1 & 0 & 0
\end{tabular}
\label{a7e6p2}
\end{table}

\begin{table}[H]\small
\caption{$3\cdot$Alt$_{7} < E_6$, $p = 2$}
\begin{tabular}{rc|cccccc|cccc}
& & \multicolumn{6}{c|}{$V_{78}$} & \multicolumn{4}{c}{$V_{27}$} \\
& & 1 & 4 & 4$^{*}$ & 6$_a$ & 14 & 20 & 6$_b$ & $(6_b)^{*}$ & 15 & $15^{*}$ \\ \hline
1) & \nongcr & 4 & 2 & 2 & 3 & 0 & 2 & 0 & 2 & 0 & 1\\
2) & & 4 & 2 & 2 & 4 & 1 & 1 & 0 & 2 & 0 & 1
\end{tabular}
\label{3a7e6p2}
\end{table}

\begin{table}[H]\small
\caption{Alt$_{6} < E_6$, $p = 0$ or $p \nmid |H|$}
\begin{tabular}{rc|ccccccc|ccccccc}
& & \multicolumn{7}{c|}{$V_{78}$} & \multicolumn{7}{c}{$V_{27}$} \\
& & 1 & 5$_a$ & 5$_b$ & 8$_a$ & 8$_b$ & 9 & 10 & 1 & 5$_a$ & 5$_b$ & 8$_a$ & 8$_b$ & 9 & 10 \\ \hline
1) & \possprim & 0 & 0 & 0 & 2 & 3 & 2 & 2 & 0 & 1 & 1 & 1 & 0 & 1 & 0\\
2) & & 2 & 1 & 1 & 0 & 1 & 2 & 4 & 0 & 1 & 1 & 1 & 0 & 1 & 0\\
3) & & 0 & 0 & 0 & 2 & 3 & 2 & 2 & 1 & 0 & 0 & 1 & 0 & 2 & 0\\
4) & & 2 & 1 & 1 & 0 & 1 & 2 & 4 & 1 & 0 & 0 & 1 & 0 & 2 & 0\\
5) & & 1 & 0 & 0 & 3 & 3 & 1 & 2 & 1 & 1 & 1 & 1 & 1 & 0 & 0\\
6) & & 3 & 1 & 1 & 1 & 1 & 1 & 4 & 1 & 1 & 1 & 1 & 1 & 0 & 0\\
7) & & 1 & 0 & 0 & 3 & 3 & 1 & 2 & 2 & 0 & 0 & 1 & 1 & 1 & 0\\
8) & & 3 & 1 & 1 & 1 & 1 & 1 & 4 & 2 & 0 & 0 & 1 & 1 & 1 & 0\\
9) & & 4 & 0 & 3 & 0 & 0 & 1 & 5 & 2 & 0 & 3 & 0 & 0 & 0 & 1\\
10) & & 2 & 0 & 0 & 7 & 0 & 0 & 2 & 3 & 0 & 0 & 3 & 0 & 0 & 0
\end{tabular}\\
Permutations: $(5_a,5_b)(8_a,8_b)$.
\label{a6e6p0}
\end{table}

\begin{table}[H]\small
\caption{$3\cdot$Alt$_{6} < E_6$, $p = 0$ or $p \nmid |H|$}
\begin{tabular}{rc|ccccccc|*{6}{c}}
& & \multicolumn{7}{c|}{$V_{78}$} & \multicolumn{6}{c}{$V_ {27}$} \\
& & 1 & 5$_a$ & 5$_b$ & 8$_a$ & 8$_b$ & 9$_a$ & 10 & $3_a$ & $3_b$ & 6 & 6$^{*}$ & 9$_b$ & 15 \\ \hline
1) &  & 1  & 0 & 0 & 3 & 3 & 1 & 2 & 0 & 0 & 0 & 2 & 0 & 1\\
2) &  & 3  & 1 & 1 & 1 & 1 & 1 & 4 & 0 & 0 & 0 & 2 & 0 & 1\\
3) & \possprim & 0  & 0 & 0 & 2 & 3 & 2 & 2 & 0 & 1 & 1 & 0 & 2 & 0\\
4) &  & 2  & 1 & 1 & 0 & 1 & 2 & 4 & 0 & 1 & 1 & 0 & 2 & 0\\
5) &  & 2  & 0 & 0 & 7 & 0 & 0 & 2 & 0 & 3 & 3 & 0 & 0 & 0\\
6) &  & 14 & 0 & 0 & 8 & 0 & 0 & 0 & 0 & 7 & 1 & 0 & 0 & 0\\
7) &  & 1  & 0 & 0 & 3 & 3 & 1 & 2 & 1 & 1 & 2 & 0 & 1 & 0\\
8) &  & 3  & 1 & 1 & 1 & 1 & 1 & 4 & 1 & 1 & 2 & 0 & 1 & 0\\
9) &  & 8  & 0 & 0 & 1 & 1 & 6 & 0 & 3 & 3 & 0 & 0 & 1 & 0
\end{tabular}\\
Permutations: $(3_a,3_b)(8_a,8_b)$.
\label{3a6e6p0}
\end{table}

\begin{table}[H]\small
\caption{Alt$_{6} < E_6$, $p = 5$}
\begin{tabular}{rc|ccccc|ccccc}
& & \multicolumn{5}{c|}{$V_{78}$} & \multicolumn{5}{c}{$V_{27}$} \\
& & 1 & 5$_a$ & 5$_b$ & 8 & 10 & 1 & 5$_a$ & 5$_b$ & 8 & 10 \\ \hline
1) & \possprim, \nongcr & 2 & 0 & 0 & 7 & 2 & 1 & 1 & 1 & 2 & 0\\
2) & \nongcr & 4 & 1 & 1 & 3 & 4 & 1 & 1 & 1 & 2 & 0\\
3) & & 5 & 0 & 3 & 1 & 5 & 2 & 0 & 3 & 0 & 1\\
4) & \nongcr & 2 & 0 & 0 & 7 & 2 & 3 & 0 & 0 & 3 & 0\\
5) & \nongcr & 4 & 1 & 1 & 3 & 4 & 3 & 0 & 0 & 3 & 0
\end{tabular}\\
Permutations: $(5_a,5_b)$.
\label{a6e6p5}
\end{table}

\begin{table}[H]\small
\caption{$3\cdot$Alt$_{6} < E_6$, $p = 5$}
\begin{tabular}{rc|ccccc|*{6}{c}}
& & \multicolumn{5}{c|}{$V_{78}$} & \multicolumn{6}{c}{$V_{27}$} \\
& & 1 & 5$_a$ & 5$_b$ & 8 & 10 & 3 & 3$^{*}$ & 6 & 6$^{*}$ & 15 & 15$^{*}$ \\ \hline
1) & \possprim, \nongcr & 2  & 0 & 0 & 7 & 2 & 0 & 0 & 0 & 2 & 1 & 0 \\
2) & \nongcr & 4  & 1 & 1 & 3 & 4 & 0 & 0 & 0 & 2 & 1 & 0 \\
3) & \possprim, \nongcr & 2  & 0 & 0 & 7 & 2 & 0 & 3 & 3 & 0 & 0 & 0 \\
4) & \nongcr & 4  & 1 & 1 & 3 & 4 & 0 & 3 & 3 & 0 & 0 & 0 \\
5) & \nongcr & 14 & 0 & 0 & 8 & 0 & 0 & 7 & 1 & 0 & 0 & 0
\end{tabular}
\label{3a6e6p5}
\end{table}

\begin{table}[H]\small
\caption{Alt$_{5} < E_6$, $p = 0$ or $p \nmid |H|$}
\begin{tabular}{r|ccccc|ccccc}
& \multicolumn{5}{c|}{$V_{78}$} & \multicolumn{5}{c}{$V_{27}$} \\
& 1 & 3$_a$ & 3$_b$ & 4 & 5 & 1 & 3$_a$ & 3$_b$ & 4 & 5 \\ \hline
1) & 5 & 1 & 1 & 3 & 11 & 0 & 3 & 3 & 1 & 1\\
2) & 8 & 1 & 1 & 6 & 8 & 0 & 3 & 3 & 1 & 1\\
3) & 0 & 4 & 5 & 4 & 7 & 1 & 1 & 0 & 2 & 3\\
4) & 3 & 4 & 5 & 7 & 4 & 1 & 1 & 0 & 2 & 3\\
5) & 14 & 8 & 0 & 0 & 8 & 1 & 7 & 0 & 0 & 1\\
6) & 1 & 5 & 5 & 3 & 7 & 2 & 1 & 1 & 1 & 3\\
7) & 4 & 5 & 5 & 6 & 4 & 2 & 1 & 1 & 1 & 3\\
8) & 2 & 9 & 2 & 2 & 7 & 3 & 3 & 0 & 0 & 3\\
9) & 5 & 9 & 2 & 5 & 4 & 3 & 3 & 0 & 0 & 3\\
10) & 6 & 4 & 5 & 10 & 1 & 4 & 1 & 0 & 5 & 0\\
11) & 7 & 5 & 5 & 9 & 1 & 5 & 1 & 1 & 4 & 0\\
12) & 8 & 9 & 2 & 8 & 1 & 6 & 3 & 0 & 3 & 0\\
13) & 16 & 0 & 19 & 0 & 1 & 9 & 0 & 6 & 0 & 0
\end{tabular}\\
Permutations: $(3_a,3_b)$.
\label{a5e6p0}
\end{table}

\begin{table}[H]\small
\caption{Alt$_{5} < E_6$, $p = 3$}
\begin{tabular}{rc|cccc|cccc}
& & \multicolumn{4}{c|}{$V_{78}$} & \multicolumn{4}{c}{$V_{27}$}\\
& & 1 & 3$_a$ & 3$_b$ & 4 & 1 & 3$_a$ & 3$_b$ & 4 \\ \hline
1) & \possprim, \nongcr & 15 & 0 & 1 & 15 & 0 & 3 & 2 & 3 \\
2) & \nongcr & 16 & 1 & 1 & 14 & 1 & 3 & 3 & 2\\
3) & \nongcr & 22 & 8 & 0 & 8 & 2 & 7 & 0 & 1\\
4) & \possprim, \nongcr & 7 & 4 & 5 & 11 & 4 & 1 & 0 & 5\\
5) & \nongcr & 8 & 5 & 5 & 10 & 5 & 1 & 1 & 4\\
6) & \nongcr & 9 & 9 & 2 & 9 & 6 & 3 & 0 & 3\\
7) & & 17 & 0 & 19 & 1 & 9 & 0 & 6 & 0
\end{tabular}\\
Permutations: $(3_a,3_b)$.
\label{a5e6p3}
\end{table}

\subsection{Sporadic Groups} \ 

\begin{table}[H]\small
\caption{$M_{11} < E_6$, $p = 0$ or $p \nmid |H|$}
\begin{tabular}{r|cccc|ccccc}
& \multicolumn{4}{c|}{$V_{78}$} & \multicolumn{5}{c}{$V_{27}$} \\
& 1 & $16$ & $16^{*}$ & 45 & 1 & 10$_a$ & 11 & $16$ & $16^{*}$ \\ \hline
1) & 1 & 1 & 1 & 1 & 0 & 0 & 1 & 0 & 1 \\
2) & 1 & 1 & 1 & 1 & 1 & 1 & 0 & 0 & 1
\end{tabular}\\
$10_a$ is self-dual.
\label{m11e6p0}
\end{table}

\begin{table}[H]\small
\caption{$M_{11} < E_6$, $p = 11$}
\begin{tabular}{r|ccccc|cccc}
& \multicolumn{5}{c|}{$V_{78}$} & \multicolumn{4}{c}{$V_{27}$} \\
& 1 & 9 & 10 & $10^{*}$ & 16 & 1 & 9 & 11 & 16 \\ \hline
1) & 1 & 1 & 1 & 1 & 3 & 0 & 0 & 1 & 1 \\
2) & 1 & 1 & 1 & 1 & 3 & 2 & 1 & 0 & 1 
\end{tabular}
\label{m11e6p11}
\end{table}

\begin{table}[H]\small
\caption{$M_{11} < E_6$, $p = 5$}
\begin{tabular}{rc|cccc|ccccc}
& & \multicolumn{4}{c|}{$V_{78}$} & \multicolumn{5}{c}{$V_{27}$} \\
& & 1 & 16 & $16^{*}$ & 45 & 1 & $10_a$ & 11 & 16 & $16^{*}$ \\ \hline
1) & \possprim, \nongcr & 1 & 1 & 1 & 1 & 0 & 0 & 1 & 0 & 1\\
2) & \possprim, \nongcr & 1 & 1 & 1 & 1 & 0 & 0 & 1 & 1 & 0\\
3) & \nongcr & 1 & 1 & 1 & 1 & 1 & 1 & 0 & 0 & 1\\
4) & \nongcr & 1 & 1 & 1 & 1 & 1 & 1 & 0 & 1 & 0
\end{tabular}\\
$10_a$ is self-dual.
\label{m11e6p5}
\end{table}

\begin{table}[H]\small
\caption{$M_{11} < E_6$, $p = 3$}
\begin{tabular}{rc|cccccccc|cccccc}
& & \multicolumn{8}{c|}{$V_{78}$} & \multicolumn{6}{c}{$V_{27}$} \\
& & 1 & 5 & $5^{*}$ & 10$_a$ & 10$_b$ & $(10_b)^{*}$ & 24 & 45 & 1 & 5 & $5^{*}$ & 10$_a$ & 10$_b$ & $(10_b)^{*}$ \\ \hline
1) & \possprim, \nongcr & 3 & 1 & 1 & 0 & 1 & 1 & 0 & 1  & 2 & 0 & 1 & 1 & 0 & 1\\
2) & \possprim, \nongcr & 3 & 1 & 1 & 0 & 1 & 1 & 0 & 1  & 2 & 1 & 2 & 0 & 0 & 1\\
3) & \possprim, \nongcr & 4 & 1 & 1 & 0 & 2 & 2 & 1 & 0  & 2 & 1 & 2 & 0 & 1 & 0 
\end{tabular}
\label{m11e6p3}
\end{table}

\begin{table}[H]\small
\caption{$M_{11} < E_6$, $p = 2$}
\begin{tabular}{rc|ccccc|cccc}
& & \multicolumn{5}{c|}{$V_{78}$} & \multicolumn{4}{c}{$V_{27}$} \\
& & 1 & 10 & 16 & $16^{*}$ & 44 & 1 & 10 & 16 & $16^{*}$ \\ \hline
1) & & 2 & 0 & 1 & 1 & 1 & 1 & 1 & 0 & 1 \\
2) & \nongcr & 4 & 3 & 0 & 0 & 1 & 1 & 1 & 0 & 1
\end{tabular}
\label{m11e6p2}
\end{table}

\begin{table}[H]\small
\caption{$M_{12} < E_6$, $p = 5$}
\begin{tabular}{rc|c|*{3}{c}}
& & \multicolumn{1}{c|}{$V_{78}$} & \multicolumn{3}{c}{$V_{27}$} \\
& & 78 & $11_a$ & $16$ & $16^{*}$ \\ \hline
1) & \possprim & 1 & 1 & 1 & 0
\end{tabular}\\
Permutations: $(11_a,11_b)$, $(16,16^{*})$.
\label{m12e6p5}
\end{table}

\begin{table}[H]\small
\caption{$M_{12} < E_6$, $p = 2$}
\begin{tabular}{r|ccccc|*{4}{c}}
& \multicolumn{5}{c|}{$V_{78}$} & \multicolumn{4}{c}{$V_{27}$} \\
& 1 & 10 & 16 & $16^{*}$ & 44 & 1 & 10 & 16 & $16^{*}$ \\ \hline
1) & 2 & 0 & 1 & 1 & 1 & 1 & 1 & 0 & 1
\end{tabular}
\label{m12e6p2}
\end{table}

\begin{table}[H]\small
\caption{$3\cdot M_{22} < E_6$, $p = 2$}
\begin{tabular}{rc|cccc|cc}
& & \multicolumn{4}{c|}{$V_{78}$} & \multicolumn{2}{c}{$V_{27}$} \\
& & 1 & 10 & $10^{*}$ & 34 & $6$ & $15$ \\ \hline
1) & \nongcr & 4 & 2 & 2 & 1 & 2 & 1
\end{tabular}\\
Permutations: $(6,6^{*})(15,15^{*})$. \\
$15 = \bigwedge^{2}6^{*}$.
\label{3m22e6p2}
\end{table}

\begin{table}[H]\small
\caption{$J_{1} < E_6$, $p = 11$}
\begin{tabular}{rc|cccc|ccc}
& & \multicolumn{4}{c|}{$V_{78}$} & \multicolumn{3}{c}{$V_{27}$} \\
& & 1 & 7 & 14 & 64 & 1 & 7 & 27 \\ \hline
1) & \possprim & 0 & 0 & 1 & 1 & 0 & 0 & 1 \\
2) & & 8 & 8 & 1 & 0 & 6 & 3 & 0
\end{tabular}
\label{j1e6p11}
\end{table}

\begin{table}[H]\small
\caption{$J_{2} < E_6$, $p = 2$}
\begin{tabular}{rc|cccccccc|ccccc}
& & \multicolumn{8}{c|}{$V_{78}$} & \multicolumn{5}{c}{$V_{27}$} \\
& & 1 & $6_a$ & $6_b$ & $14_a$ & $14_b$ & 36 & $64_a$ & $64_b$ & 1 & $6_a$ & $6_b$ & $14_a$ & $14_b$ \\ \hline
1) & \nongcr & 8 & 3 & 4 & 0 & 2 & 0 & 0 & 0 & 1 & 0 & 2 & 0 & 1 \\
2) & \possprim & 0 & 0 & 0 & 1 & 0 & 0 & 1 & 0 & 1 & 1 & 1 & 1 & 0 \\
3) & \nongcr & 16 & 8 & 0 & 1 & 0 & 0 & 0 & 0 & 9 & 3 & 0 & 0 & 0
\end{tabular}\\
Permutations: $(6_a,6_b)(14_a,14_b)(64_a,64_b)$.
\label{j2e6p2}
\end{table}

\begin{table}[H]\small
\caption{$3 \cdot J_{3} < E_6$, $p = 2$}
\begin{tabular}{rc|cc|cccccc}
& & \multicolumn{2}{c|}{$V_{78}$} & \multicolumn{6}{c}{$V_{27}$} \\
& & $78_a$ & $78_b$ & $9$ & $9^{*}$ & $18_a$ & $(18_a)^{*}$ & $18_b$ & $(18_b)^{*}$ \\ \hline
1) & \possprim & 1 & 0 & 0 & 1 & 0 & 0 & 0 & 1 \\
2) & \possprim & 0 & 1 & 0 & 1 & 1 & 0 & 0 & 0 \\
3) & \possprim & 0 & 1 & 1 & 0 & 0 & 0 & 1 & 0 \\
4) & \possprim & 1 & 0 & 1 & 0 & 0 & 1 & 0 & 0
\end{tabular}
\label{3j3e6p2}
\end{table}

\addtocounter{table}{1}
\begin{table}[H]\small
\begin{center}
\arabic{chapter}.\arabic{section}.\arabic{table}: $3 \cdot Fi_{22} < E_6$, $p = 2$. Irreducible on $V_{78}$ and $V_{27}$. \possprim
\end{center}
\label{3fi22e6p2}
\end{table}

\subsection{Cross-characteristic Groups $L_2(q)$ $(q \neq 4$, $5$, $9)$} \ 

\begin{table}[H]\small
\caption{$L_2(7) < E_6$, $p = 0$ or $p \nmid |H|$}
\begin{tabular}{rc|cccccc|ccccccc}
& & \multicolumn{6}{c|}{$V_{78}$} & \multicolumn{6}{c}{$V_{27}$} \\
& & 1 & 3 & $3^{*}$ & 6 & 7 & 8 & 1 & 3 & $3^{*}$ & 6 & 7 & 8 \\ \hline
1) & \possprim & 0 & 2 & 2 & 2 & 2 & 5 & 0 & 0 & 0 & 2 & 1 & 1 \\
2) & & 0 & 2 & 2 & 2 & 2 & 5 & 1 & 1 & 1 & 2 & 0 & 1 \\
3) & & 2 & 1 & 1 & 0 & 2 & 7 & 3 & 0 & 0 & 0 & 0 & 3 \\
4) & & 3 & 2 & 2 & 2 & 5 & 2 & 0 & 0 & 0 & 2 & 1 & 1 \\
5) & & 3 & 2 & 2 & 2 & 5 & 2 & 1 & 1 & 1 & 2 & 0 & 1 \\
6) & & 5 & 1 & 1 & 0 & 5 & 4 & 3 & 0 & 0 & 0 & 0 & 3 \\
7) & & 8 & 1 & 1 & 0 & 8 & 1 & 6 & 0 & 0 & 0 & 3 & 0 \\
8) & & 8 & 3 & 3 & 6 & 0 & 2 & 1 & 3 & 3 & 0 & 0 & 1 \\
9) & & 14 & 0 & 0 & 0 & 0 & 8 & 0 & 0 & 7 & 1 & 0 & 0 \\
10) & & 14 & 0 & 0 & 0 & 0 & 8 & 0 & 7 & 0 & 1 & 0 & 0 \\
11) & & 16 & 9 & 9 & 0 & 0 & 1 & 9 & 3 & 3 & 0 & 0 & 0
\end{tabular}
\label{l27e6p0}
\end{table}

\begin{table}[H]\small
\caption{$L_2(7) < E_6$, $p = 3$}
\begin{tabular}{rc|ccccc|cccccc}
& & \multicolumn{5}{c|}{$V_{78}$} & \multicolumn{5}{c}{$V_{27}$} \\
& & 1 & 3 & $3^{*}$ & 6 & 7 & 1 & 3 & $3^{*}$ & 6 & 7 \\ \hline
1) & \possprim, \nongcr & 5 & 2 & 2 & 2 & 7 & 1 & 0 & 0 & 2 & 2 \\
2) & \nongcr & 5 & 2 & 2 & 2 & 7 & 2 & 1 & 1 & 2 & 1 \\
3) & \nongcr & 9 & 1 & 1 & 0 & 9 & 6 & 0 & 0 & 0 & 3 \\
4) & \nongcr & 10 & 3 & 3 & 6 & 2 & 2 & 3 & 3 & 0 & 1 \\
5) & & 17 & 9 & 9 & 0 & 1 & 9 & 3 & 3 & 0 & 0 \\
6) & \nongcr & 22 & 0 & 0 & 0 & 8 & 0 & 0 & 7 & 1 & 0 \\
7) & \nongcr & 22 & 0 & 0 & 0 & 8 & 0 & 7 & 0 & 1 & 0
\end{tabular}
\label{l27e6p3}
\end{table}

\begin{table}[H]\small
\caption{$L_2(8) < E_6$, $p = 0$ or $p \nmid |H|$}
\begin{tabular}{rc|*{9}{c}|*{9}{c}}
& & \multicolumn{9}{c|}{$V_{78}$} & \multicolumn{9}{c}{$V_{27}$} \\
& & 1 & $7_a$ & $7_b$ & $7_c$ & $7_d$ & 8 & $9_a$ & $9_b$ & $9_c$ & 1 & $7_a$ & $7_b$ & $7_c$ & $7_d$ & 8 & $9_a$ & $9_b$ & $9_c$ \\ \hline
1) & & 0 & 2 & 1 & 1 & 1 & 2 & 0 & 1 & 2 & 1 & 0 & 0 & 0 & 0 & 1 & 0 & 1 & 1 \\
2) & \possprim & 0 & 2 & 1 & 1 & 1 & 2 & 1 & 1 & 1 & 0 & 0 & 0 & 0 & 0 & 0 & 1 & 1 & 1 \\
3) & & 1 & 3 & 1 & 1 & 1 & 1 & 0 & 1 & 2 & 2 & 1 & 0 & 0 & 0 & 0 & 0 & 1 & 1 \\
4) & & 3 & 2 & 1 & 1 & 1 & 5 & 0 & 0 & 0 & 3 & 0 & 0 & 0 & 0 & 3 & 0 & 0 & 0 \\
5) & & 4 & 3 & 1 & 1 & 1 & 4 & 0 & 0 & 0 & 4 & 1 & 0 & 0 & 0 & 2 & 0 & 0 & 0 \\
6) & & 8 & 1 & 0 & 1 & 8 & 0 & 0 & 0 & 0 & 6 & 0 & 0 & 0 & 3 & 0 & 0 & 0 & 0
\end{tabular}\\
Permutations: $(7_b,7_c,7_d)(9_a,9_b,9_c)$.
\label{l28e6p0}
\end{table}

\begin{table}[H]\small
\caption{$L_2(8) < E_6$, $p = 7$}
\begin{tabular}{rc|cccccc|cccccc}
& & \multicolumn{6}{c|}{$V_{78}$} & \multicolumn{6}{c}{$V_{27}$} \\
& & 1 & $7_a$ & $7_b$ & $7_c$ & $7_d$ & 8 & 1 & $7_a$ & $7_b$ & $7_c$ & $7_d$ & 8 \\ \hline
1) & \nongcr & 3 & 2 & 1 & 1 & 1 & 5 & 3 & 0 & 0 & 0 & 0 & 3 \\
2) & \nongcr & 4 & 3 & 1 & 1 & 1 & 4 & 4 & 1 & 0 & 0 & 0 & 2 \\
3) & & 8 & 1 & 0 & 1 & 8 & 0 & 6 & 0 & 0 & 0 & 3 & 0
\end{tabular}\\
Permutations: $(7_b,7_c,7_d)$.
\label{l28e6p7}
\end{table}

\begin{table}[H]\small
\caption{$L_2(8) < E_6$, $p = 3$}
\begin{tabular}{rc|ccccc|ccccc}
& & \multicolumn{5}{c|}{$V_{78}$} & \multicolumn{5}{c}{$V_{27}$} \\
& & 1 & 7 & $9_a$ & $9_b$ & $9_c$ & 1 & 7 & $9_a$ & $9_b$ & $9_c$ \\ \hline
1) & \nongcr & 2 & 7 & 0 & 1 & 2 & 2 & 1 & 0 & 1 & 1 \\
2) & \possprim, \nongcr & 2 & 7 & 1 & 1 & 1 & 0 & 0 & 1 & 1 & 1 \\
3) & \nongcr & 8 & 10 & 0 & 0 & 0 & 6 & 3 & 0 & 0 & 0
\end{tabular}\\
Permutations: $(9_a,9_b,9_c)$.
\label{l28e6p3}
\end{table}

\begin{table}[H]\small
\caption{$L_2(11) < E_6$, $p = 0$ or $p \nmid |H|$}
\begin{tabular}{rc|*{7}{c}|*{8}{c}}
& & \multicolumn{7}{c|}{$V_{78}$} & \multicolumn{8}{c}{$V_{27}$} \\
& & 1 & 5 & $5^{*}$ & $10_b$ & 11 & $12_a$ & $12_b$ & 1 & 5 & $5^{*}$ & $10_a$ & $10_b$ & 11 & $12_a$ & $12_b$ \\ \hline
1) & \possprim & 0 & 1 & 1 & 1 & 2 & 1 & 2 & 0 & 0 & 1 & 1 & 0 & 0 & 1 & 0 \\
2) & & 1 & 1 & 1 & 1 & 3 & 1 & 1 & 1 & 0 & 1 & 1 & 0 & 1 & 0 & 0 \\
3) & & 4 & 1 & 1 & 4 & 0 & 1 & 1 & 2 & 1 & 2 & 0 & 1 & 0 & 0 & 0 \\
\end{tabular}\\
Permutations: $(12_a,12_b)$.
\label{l211e6p0}
\end{table}

\begin{table}[H]\small
\caption{$L_2(11) < E_6$, $p = 5$}
\begin{tabular}{rc|*{6}{c}|*{6}{c}}
& & \multicolumn{6}{c|}{$V_{78}$} & \multicolumn{6}{c}{$V_{27}$} \\
& & 1 & 5 & $5^{*}$ & $10_a$ & $10_b$ & $11$ & 1 & 5 & $5^{*}$ & $10_a$ & $10_b$ & $11$ \\ \hline
1) & \nongcr & 3 & 1 & 1 & 0 & 1 & 5 & 1 & 0 & 1 & 1 & 0 & 1 \\
2) & \nongcr & 6 & 1 & 1 & 0 & 4 & 2 & 2 & 2 & 1 & 0 & 1 & 0
\end{tabular}
\label{l211e6p5}
\end{table}

\begin{table}[H]\small
\caption{$L_2(11) < E_6$, $p = 3$}
\begin{tabular}{rc|*{6}{c}|*{6}{c}}
& & \multicolumn{6}{c|}{$V_{78}$} & \multicolumn{6}{c}{$V_{27}$} \\
& & 1 & 5 & $5^{*}$ & 10 & $12_a$ & $12_b$ & 1 & 5 & $5^{*}$ & 10 & $12_a$ & $12_b$ \\ \hline
1) & \possprim, \nongcr & 2 & 1 & 1 & 3 & 1 & 2 & 0 & 1 & 2 & 0 & 1 & 0 \\
2) & \nongcr & 4 & 1 & 1 & 4 & 1 & 1 & 2 & 2 & 1 & 1 & 0 & 0
\end{tabular}\\
Permutations: $(12_a,12_b)$.
\label{l211e6p3}
\end{table}

\begin{table}[H]\small
\caption{$L_2(11) < E_6$, $p = 2$}
\begin{tabular}{rc|*{6}{c}|*{6}{c}}
& & \multicolumn{6}{c|}{$V_{78}$} & \multicolumn{6}{c}{$V_{27}$} \\
& & 1 & 5 & $5^{*}$ & 10 & $12_a$ & $12_b$ & 1 & 5 & $5^{*}$ & 10 & $12_a$ & $12_b$ \\ \hline
1) & & 2 & 0 & 0 & 4 & 1 & 2 & 0 & 0 & 1 & 1 & 1 & 0 \\
2) & \possprim, \nongcr & 2 & 3 & 3 & 1 & 1 & 2 & 0 & 0 & 1 & 1 & 1 & 0 \\\
3) & \nongcr & 4 & 1 & 1 & 4 & 1 & 1 & 2 & 1 & 2 & 1 & 0 & 0 \\
4) & \possprim, \nongcr & 4 & 4 & 4 & 1 & 1 & 1 & 2 & 1 & 2 & 1 & 0 & 0
\end{tabular}\\
Permutations: $(12_a,12_b)$.
\label{l211e6p2}
\end{table}

\begin{table}[H]\small
\caption{$L_2(13) < E_6$, $p = 0$ or $p \nmid |H|$}
\begin{tabular}{rc|*{9}{c}|*{6}{c}}
& & \multicolumn{9}{c|}{$V_{78}$} & \multicolumn{6}{c}{$V_{27}$} \\
& & 1 & $7_a$ & $7_b$ & $12_a$ & $12_b$ & $12_c$ & 13 & $14_a$ & $14_b$ & 1 & $7_a$ & $7_b$ & $12_a$ & 13 & $14_b$  \\ \hline
1) & \possprim & 0 & 0 & 0 & 1 & 1 & 1 & 0 & 2 & 1 & 0 & 0 & 0 & 0 & 1 & 1 \\
2) & & 0 & 0 & 0 & 2 & 1 & 0 & 0 & 2 & 1 & 1 & 0 & 0 & 1 & 0 & 1 \\
3) & & 1 & 2 & 2 & 1 & 1 & 1 & 1 & 0 & 0 & 0 & 0 & 0 & 0 & 1 & 1 \\
4) & & 1 & 2 & 2 & 2 & 1 & 0 & 1 & 0 & 0 & 1 & 0 & 0 & 1 & 0 & 1 \\
5) & & 8 & 0 & 8 & 0 & 0 & 0 & 0 & 1 & 0 & 6 & 0 & 3 & 0 & 0 & 0
\end{tabular}\\
Permutations: $(7_a,7_b)$, $(12_a,12_b,12_c)$.
\label{l213e6p0}
\end{table}

\begin{table}[H]\small
\caption{$L_2(13) < E_6$, $p = 7$}
\begin{tabular}{rc|*{6}{c}|*{6}{c}}
& & \multicolumn{6}{c|}{$V_{78}$} & \multicolumn{6}{c}{$V_{27}$} \\
& & 1 & $7_a$ & $7_b$ & $12$ & $14_a$ & $14_b$ & 1 & $7_a$ & $7_b$ & $12$ & $14_a$ & $14_b$ \\ \hline
1) & & 0 & 0 & 0 & 3 & 2 & 1 & 1 & 0 & 0 & 1 & 0 & 1 \\
2) & \nongcr & 2 & 2 & 2 & 4 & 0 & 0 & 1 & 0 & 0 & 1 & 0 & 1 \\
3) & \nongcr & 8 & 0 & 8 & 0 & 1 & 0 & 6 & 0 & 3 & 0 & 0 & 0
\end{tabular}\\
Permutations: $(7_a,7_b)$.
\label{l213e6p7}
\end{table}

\begin{table}[H]\small
\caption{$L_2(13) < E_6$, $p = 3$}
\begin{tabular}{rc|*{7}{c}|*{7}{c}}
& & \multicolumn{7}{c|}{$V_{78}$} & \multicolumn{7}{c}{$V_{27}$} \\
& & 1 & $7_a$ & $7_b$ & $12_a$ & $12_b$ & $12_c$ & $13$ & 1 & $7_a$ & $7_b$ & $12_a$ & $12_b$ & $12_c$ & $13$ \\ \hline
1) & \possprim & 1 & 2 & 2 & 1 & 1 & 1 & 1 & 1 & 0 & 0 & 0 & 0 & 0 & 2 \\
2) & & 1 & 2 & 2 & 2 & 1 & 0 & 1 & 2 & 0 & 0 & 1 & 0 & 0 & 1 \\
3) & & 8 & 1 & 9 & 0 & 0 & 0 & 0 & 6 & 0 & 3 & 0 & 0 & 0 & 0
\end{tabular}\\
Permutations: $(7_a,7_b)$, $(12_a,12_b,12_c)$.
\label{l213e6p3}
\end{table}

\begin{table}[H]\small
\caption{$L_2(13) < E_6$, $p = 2$}
\begin{tabular}{rc|*{7}{c}|*{7}{c}}
& & \multicolumn{7}{c|}{$V_{78}$} & \multicolumn{7}{c}{$V_{27}$} \\
& & 1 & $6_a$ & $6_b$ & $12_a$ & $12_b$ & $12_c$ & 14 & 1 & $6_a$ & $6_b$ & $12_a$ & $12_b$ & $12_c$ & 14 \\ \hline
1) & \possprim & 0 & 0 & 0 & 1 & 1 & 1 & 3 & 1 & 1 & 1 & 0 & 0 & 0 & 1 \\
2) & & 0 & 0 & 0 & 2 & 0 & 1 & 3 & 1 & 0 & 0 & 1 & 0 & 0 & 1 \\
3) & & 2 & 0 & 1 & 0 & 0 & 0 & 5 & 1 & 0 & 2 & 0 & 0 & 0 & 1 \\
4) & \nongcr & 6 & 3 & 3 & 1 & 1 & 1 & 0 & 1 & 1 & 1 & 0 & 0 & 0 & 1 \\
5) & \nongcr & 6 & 3 & 3 & 2 & 0 & 1 & 0 & 1 & 0 & 0 & 1 & 0 & 0 & 1 \\
6) & \nongcr & 8 & 3 & 4 & 0 & 0 & 0 & 2 & 1 & 0 & 2 & 0 & 0 & 0 & 1 \\
7) & \nongcr & 16 & 0 & 8 & 0 & 0 & 0 & 1 & 9 & 0 & 3 & 0 & 0 & 0 & 0
\end{tabular}\\
Permutations: $(6_a,6_b)$, $(12_a,12_b,12_c)$.
\label{l213e6p2}
\end{table}

\begin{table}[H]\small
\caption{$L_2(17) < E_6$, $p = 0$ or $p \nmid |H|$}
\begin{tabular}{rc|*{7}{c}|*{6}{c}}
& & \multicolumn{7}{c|}{$V_{78}$} & \multicolumn{6}{c}{$V_{27}$} \\
& & 1 & $9_a$ & $9_b$ & $16_d$ & $17$ & $18_b$ & $18_c$ & 1 & $9_a$ & $9_b$ & $16_d$ & $17$ & $18_a$ \\ \hline
1) & \possprim & 0 & 0 & 1 & 1 & 1 & 1 & 1 & 0 & 0 & 1 & 0 & 0 & 1 \\
2) & & 0 & 0 & 1 & 1 & 1 & 1 & 1 & 1 & 0 & 1 & 0 & 1 & 0 \\
3) & & 1 & 0 & 1 & 2 & 0 & 1 & 1 & 2 & 0 & 1 & 1 & 0 & 0
\end{tabular}\\
Permutations: $(9_a,9_b)$.
\label{l217e6p0}
\end{table}

\begin{table}[H]\small
\caption{$L_2(17) < E_6$, $p = 3$}
\begin{tabular}{rc|*{6}{c}|*{5}{c}}
& & \multicolumn{6}{c|}{$V_{78}$} & \multicolumn{5}{c}{$V_{27}$} \\
& & 1 & $9_a$ & $9_b$ & 16 & $18_b$ & $18_c$ & 1 & $9_a$ & $9_b$ & 16 & $18_a$ \\ \hline
1) & \possprim, \nongcr & 1 & 0 & 1 & 2 & 1 & 1 & 0 & 0 & 1 & 0 & 1 \\
2) & \nongcr & 1 & 0 & 1 & 2 & 1 & 1 & 2 & 0 & 1 & 1 & 0 
\end{tabular}\\
Permutations: $(9_a,9_b)$.
\label{l217e6p3}
\end{table}

\begin{table}[H]\small
\caption{$L_2(17) < E_6$, $p = 2$}
\begin{tabular}{rc|*{4}{c}|*{4}{c}}
& & \multicolumn{4}{c|}{$V_{78}$} & \multicolumn{4}{c}{$V_{27}$} \\
& & 1 & $8_a$ & $8_b$ & $16_a$ & 1 & $8_a$ & $8_b$ & $16_a$ \\ \hline
1) & \nongcr & 6 & 2 & 3 & 2 & 3 & 0 & 1 & 1 \\
2) & \nongcr & 6 & 3 & 4 & 1 & 3 & 1 & 2 & 0
\end{tabular}\\
Permutations: $(8_a,8_b)$.
\label{l217e6p2}
\end{table}

\begin{table}[H]\small
\caption{$L_2(19) < E_6$, $p = 0$ or $p \nmid |H|$}
\begin{tabular}{rc|*{5}{c}|*{4}{c}}
& & \multicolumn{5}{c|}{$V_{78}$} & \multicolumn{4}{c}{$V_{27}$} \\
& & $18_c$ & $18_d$ & $20_b$ & $20_c$ & $20_d$ & 9 & $9^{*}$ & $18_a$ & $18_b$ \\ \hline
1) & \possprim & 0 & 1 & 1 & 1 & 1 & 0 & 1 & 1 & 0
\end{tabular}\\
Permutations: $(18_a,18_b)(18_c,18_d)$.
\label{l219e6p0}
\end{table}

\begin{table}[H]\small
\caption{$L_2(19) < E_6$, $p = 5$}
\begin{tabular}{rc|*{4}{c}|*{2}{c}}
& & \multicolumn{4}{c|}{$V_{78}$} & \multicolumn{2}{c}{$V_{27}$} \\
& & 18 & $20_b$ & $20_c$ & $20_d$ & 9 & $9^{*}$ \\ \hline
1) & \possprim & 1 & 1 & 1 & 1 & 1 & 2 
\end{tabular}
\label{l219e6p5}
\end{table}

\begin{table}[H]\small
\caption{$L_2(19) < E_6$, $p = 3$}
\begin{tabular}{rc|*{4}{c}|*{4}{c}}
& & \multicolumn{4}{c|}{$V_{78}$} & \multicolumn{4}{c}{$V_{27}$} \\
& & 1 & $18_c$ & $18_d$ & 19 & 9 & $9^{*}$ & $18_a$ & $18_b$ \\ \hline
1) & \possprim & 3 & 0 & 1 & 3 & 0 & 1 & 1 & 0 \\
\end{tabular}\\
Permutations: $(18_a,18_b)(18_c,18_d)$.
\label{l219e6p3}
\end{table}

\begin{table}[H]\small
\caption{$L_2(19) < E_6$, $p = 2$}
\begin{tabular}{rc|*{5}{c}|*{4}{c}}
& & \multicolumn{5}{c|}{$V_{78}$} & \multicolumn{4}{c}{$V_{27}$} \\
& & $18_a$ & $18_b$ & $20_b$ & $20_c$ & $20_d$ & 9 & $9^{*}$ & $18_a$ & $18_b$ \\ \hline
1) & \possprim & 0 & 1 & 1 & 1 & 1 & 0 & 1 & 1 & 0
\end{tabular}\\
Permutations: $(18_a,18_b)$.
\label{l219e6p2}
\end{table}

\begin{table}[H]\small
\caption{$L_2(25) < E_6$, $p = 0$ or $p \nmid |H|$}
\begin{tabular}{r|*{6}{c}|*{3}{c}}
& \multicolumn{6}{c|}{$V_{78}$} & \multicolumn{3}{c}{$V_{27}$} \\
& 1 & 25 & $26_a$ & $26_c$ & $26_d$ & $26_e$ & 1 & $26_b$ & $26_c$ \\ \hline
1) & 0 & 0 & 1 & 0 & 1 & 1 & 1 & 1 & 0 \\
2) & 1 & 1 & 0 & 2 & 0 & 0 & 1 & 0 & 1
\end{tabular}\\
$26_b$, $26_c$ are $\textup{Aut}(L_2(25))$-conjugate.
\label{l225e6p0}
\end{table}

\begin{table}[H]\small
\caption{$L_2(25) < E_6$, $p = 13$}
\begin{tabular}{r|*{6}{c}|*{3}{c}}
& \multicolumn{6}{c|}{$V_{78}$} & \multicolumn{3}{c}{$V_{27}$} \\
& 1 & 24 & $26_a$ & $26_b$ & $26_c$ & $26_d$ & 1 & $26_a$ & $26_e$ \\ \hline
1) & 0 & 0 & 1 & 0 & 1 & 1 & 1 & 1 & 0 \\
2) & 2 & 1 & 0 & 2 & 0 & 0 & 1 & 0 & 1 \\
3) & 2 & 1 & 0 & 2 & 0 & 0 & 1 & 1 & 0 
\end{tabular}\\
$26_b$, $26_c$ are $\textup{Aut}(L_2(25))$-conjugate.
\label{l225e6p13}
\end{table}

\begin{table}[H]\small
\caption{$L_2(25) < E_6$, $p = 3$}
\begin{tabular}{r|*{3}{c}|*{4}{c}}
& \multicolumn{3}{c|}{$V_{78}$} & \multicolumn{4}{c}{$V_{27}$} \\
& 1 & 25 & 26 & 1 & $13_a$ & $13_b$ & 25 \\ \hline
1) & 1 & 1 & 2 & 1 & 1 & 1 & 0 \\
2) & 1 & 1 & 2 & 2 & 0 & 0 & 1
\end{tabular}
\label{l225e6p3}
\end{table}

\begin{table}[H]\small
\caption{$L_2(25) < E_6$, $p = 2$}
\begin{tabular}{r|*{4}{c}|*{2}{c}}
& \multicolumn{4}{c|}{$V_{78}$} & \multicolumn{2}{c}{$V_{27}$} \\
& 1 & $12_a$ & $12_b$ & 26 & 1 & 26 \\ \hline
1) & 0 & 0 & 0 & 3 & 1 & 1 \\
2) & 6 & 3 & 3 & 0 & 1 & 1
\end{tabular}
\label{l225e6p2}
\end{table}

\begin{table}[H]\small
\caption{$L_2(27) < E_6$, $p = 0$ or $p \nmid |H|$}
\begin{tabular}{r|*{6}{c}|*{4}{c}}
& \multicolumn{6}{c|}{$V_{78}$} & \multicolumn{4}{c}{$V_{27}$} \\
& $26_a$ & $26_b$ & $26_c$ & $26_d$ & $26_e$ & $26_f$ & 1 & $26_d$ & $26_e$ & $26_f$ \\ \hline
1) & 0 & 1 & 1 & 0 & 1 & 0 & 1 & 0 & 1 & 0
\end{tabular}\\
Permutations: $(26_a,26_b,26_c)(26_d,26_e,26_f)$.
\label{l227e6p0}
\end{table}

\begin{table}[H]\small
\caption{$L_2(27) < E_6$, $p = 13$}
\begin{tabular}{r|*{6}{c}|*{4}{c}}
& \multicolumn{6}{c|}{$V_{78}$} & \multicolumn{4}{c}{$V_{27}$} \\
& $26_a$ & $26_b$ & $26_c$ & $26_d$ & $26_e$ & $26_f$ & 1 & $26_d$ & $26_e$ & $26_f$ \\ \hline
1) & 0 & 1 & 1 & 0 & 0 & 1 & 1 & 0 & 0 & 1
\end{tabular}\\
Permutations: $(26_a,26_b,26_c)(26_d,26_e,26_f)$.
\label{l227e6p13}
\end{table}

\begin{table}[H]\small
\caption{$L_2(27) < E_6$, $p = 7$}
\begin{tabular}{r|*{3}{c}|*{3}{c}}
& \multicolumn{3}{c|}{$V_{78}$} & \multicolumn{3}{c}{$V_{27}$} \\
& $13_a$ & $13_b$ & 26 & 1 & $13_a$ & $13_b$ \\ \hline
1) & 1 & 1 & 2 & 1 & 1 & 1
\end{tabular}
\label{l227e6p7}
\end{table}

\begin{table}[H]\small
\caption{$L_2(27) < E_6$, $p = 2$}
\begin{tabular}{r|*{3}{c}|*{4}{c}}
& \multicolumn{3}{c|}{$V_{78}$} & \multicolumn{4}{c}{$V_{27}$} \\
& $26_a$ & $26_b$ & $26_c$ & 1 & 13 & $13^{*}$ & $26_a$ \\ \hline
1) & 1 & 1 & 1 & 1 & 1 & 1 & 0 \\
2) & 2 & 0 & 1 & 1 & 0 & 0 & 1
\end{tabular}\\
Permutations: $(26_a,26_b,26_c)$.
\label{l227e6p2}
\end{table}

\subsection{Cross-characteristic Groups $\ncong L_2(q)$} \ 

\begin{table}[H]\small
\caption{$L_3(3) < E_6$, $p = 0$ or $p \nmid |H|$}
\begin{tabular}{rc|*{9}{c}|*{3}{c}}
& & \multicolumn{9}{c|}{$V_{78}$} & \multicolumn{3}{c}{$V_{27}$} \\
& & 1 & 13 & $16_a$ & $(16_a^{*})$ & $16_b$ & $(16_b)^{*}$ & $26_a$ & $26_b$ & $(26_b)^{*}$ & 1 & $26_a$ & 27 \\ \hline
1) & \possprim & 0 & 0 & 0 & 0 & 0 & 0 & 1 & 1 & 1 & 0 & 0 & 1 \\
2) & & 0 & 0 & 0 & 0 & 0 & 0 & 1 & 1 & 1 & 1 & 1 & 0 \\
3) & & 1 & 1 & 1 & 1 & 1 & 1 & 0 & 0 & 0 & 0 & 0 & 1 \\
4) & & 1 & 1 & 1 & 1 & 1 & 1 & 0 & 0 & 0 & 1 & 1 & 0
\end{tabular}
\label{l33e6p0}
\end{table}

\begin{table}[H]\small
\caption{$L_3(3) < E_6$, $p = 13$}
\begin{tabular}{rc|*{6}{c}|*{4}{c}}
& & \multicolumn{6}{c|}{$V_{78}$} & \multicolumn{4}{c}{$V_{27}$} \\
& & 1 & 13 & 16 & $26_a$ & $26_b$ & $(26_b)^{*}$ & 1 & 11 & 16 & $26_a$ \\ \hline
1) & \possprim & 0 & 0 & 0 & 1 & 1 & 1 & 0 & 1 & 1 & 0 \\
2) & & 0 & 0 & 0 & 1 & 1 & 1 & 1 & 0 & 0 & 1 \\
3) & & 1 & 1 & 4 & 0 & 0 & 0 & 0 & 1 & 1 & 0 \\
4) & & 1 & 1 & 4 & 0 & 0 & 0 & 1 & 0 & 0 & 1 
\end{tabular}
\label{l33e6p13}
\end{table}

\begin{table}[H]\small
\caption{$L_3(3) < E_6$, $p = 2$}
\begin{tabular}{r|*{7}{c}|*{2}{c}}
& \multicolumn{7}{c|}{$V_{78}$} & \multicolumn{2}{c}{$V_{27}$} \\
& 1 & 12 & $16_a$ & $(16_a)^{*}$ & $16_b$ & $(16_b)^{*}$ & 26 & 1 & 26 \\ \hline
1) & 0 & 0 & 0 & 0 & 0 & 0 & 3 & 1 & 1 \\
2) & 2 & 1 & 1 & 1 & 1 & 1 & 0 & 1 & 1
\end{tabular}
\label{l33e6p2}
\end{table}

\begin{table}[H]\small
\caption{$L_4(3) < E_6$, $p = 2$}
\begin{tabular}{r|*{2}{c}|*{2}{c}}
& \multicolumn{2}{c|}{$V_{78}$} & \multicolumn{2}{c}{$V_{27}$} \\
& $26_a$ & $26_b$ & 1 & $26_a$ \\ \hline
1) & 2 & 1 & 1 & 1
\end{tabular}\\
Permutations: $(26_a,26_b)$.
\label{l43e6p2}
\end{table}

\begin{table}[H]\small
\caption{$U_3(3) < E_6$, $p = 0$ or $p \nmid |H|$}
\begin{tabular}{rc|*{11}{c}|*{5}{c}}
& & \multicolumn{11}{c|}{$V_{78}$} & \multicolumn{5}{c}{$V_{27}$} \\
& & 1 & 6 & $7_a$ & $7_b$ & $(7_b)^{*}$ & 14 & $21_a$ & 28 & $ 28^{*}$ & $32$ & $32^{*}$ & 1 & 6 & $7_a$ & 14 & 27 \\ \hline
1) & \possprim & 0 & 0 & 0 & 0 & 0 & 1 & 0 & 0 & 0 & 1 & 1 & 0 & 0 & 0 & 0 & 1 \\
2) & & 1 & 0 & 1 & 1 & 1 & 0 & 0 & 1 & 1 & 0 & 0 & 0 & 0 & 0 & 0 & 1 \\
3) & & 3 & 2 & 0 & 2 & 2 & 1 & 1 & 0 & 0 & 0 & 0 & 1 & 2 & 0 & 1 & 0 \\
4) & & 8 & 0 & 8 & 0 & 0 & 1 & 0 & 0 & 0 & 0 & 0 & 6 & 0 & 3 & 0 & 0
\end{tabular}
\label{u33e6p0}
\end{table}

\begin{table}[H]\small
\caption{$U_3(3) < E_6$, $p = 7$}
\begin{tabular}{r|*{10}{c}|*{5}{c}}
& \multicolumn{10}{c|}{$V_{78}$} & \multicolumn{5}{c}{$V_{27}$} \\
& 1 & 6 & $7_a$ & $7_b$ & $(7_b)^{*}$ & 14 & $21_a$ & 26 & 28 & $ 28^{*}$ & 1 & 6 & $7_a$ & 14 & 26 \\ \hline
1) & 0 & 2 & 0 & 0 & 0 & 1 & 0 & 2 & 0 & 0 & 1 & 0 & 0 & 0 & 1 \\
2) & 1 & 0 & 1 & 1 & 1 & 0 & 0 & 0 & 1 & 1 & 1 & 0 & 0 & 0 & 1 \\
3) & 3 & 2 & 0 & 2 & 2 & 1 & 1 & 0 & 0 & 0 & 1 & 2 & 0 & 1 & 0 \\
4) & 8 & 0 & 8 & 0 & 0 & 1 & 0 & 0 & 0 & 0 & 6 & 0 & 3 & 0 & 0
\end{tabular}
\label{u33e6p7}
\end{table}

\begin{table}[H]\small
\caption{$U_4(2) < E_6$, $p = 0$ or $p \nmid |H|$}
\begin{tabular}{r|*{8}{c}|*{6}{c}}
& \multicolumn{8}{c|}{$V_{78}$} & \multicolumn{6}{c}{$V_{27}$} \\
& 1 & 5 & $5^{*}$ & 10 & $10^{*}$ & $15_b$ & 20 & 24 & 1 & 5 & $5^{*}$ & 6 & 10 & $15_b$  \\ \hline
1) & 3 & 0 & 0 & 2 & 2 & 1 & 1 & 0 & 0 & 0 & 0 & 2 & 0 & 1 \\
2) & 4 & 1 & 1 & 2 & 2 & 0 & 0 & 1 & 2 & 1 & 2 & 0 & 1 & 0
\end{tabular}
\label{u42e6p0}
\end{table}

\begin{table}[H]\small
\caption{$U_4(2) < E_6$, $p = 5$}
\begin{tabular}{r|*{8}{c}|*{6}{c}}
& \multicolumn{8}{c|}{$V_{78}$} & \multicolumn{6}{c}{$V_{27}$} \\
& 1 & 5 & $5^{*}$ & 10 & $10^{*}$ & $15_b$ & 20 & 23 & 1 & 5 & $5^{*}$ & 6 & 10 & $15_b$ \\ \hline
1) & 3 & 0 & 0 & 2 & 2 & 1 & 1 & 0 & 0 & 0 & 0 & 2 & 0 & 1 \\
2) & 5 & 1 & 1 & 2 & 2 & 0 & 0 & 1 & 2 & 1 & 2 & 0 & 1 & 0
\end{tabular}
\label{u42e6p5}
\end{table}

\begin{table}[H]\small
\caption{$3 \cdot U_4(3) < E_6$, $p = 2$}
\begin{tabular}{rc|*{3}{c}|*{4}{c}}
& & \multicolumn{3}{c|}{$V_{78}$} & \multicolumn{4}{c}{$V_{27}$} \\
& & 1 & 20 & $34_a$ & 6 & $6^{*}$ & 15 & $15^{*}$ \\ \hline
1) & \nongcr & 4 & 2 & 1 & 0 & 2 & 0 & 1 \\
2) & \nongcr & 4 & 2 & 1 & 2 & 0 & 1 & 0
\end{tabular}\\
There are two triple covers of $U_4(3)$ up to isomorphism, however one of these has no 27-dimensional faithful modules, hence has no feasible characters here.
\label{u43e6p2}
\end{table}

\addtocounter{table}{1}
\begin{table}[H]\small
\thetable: $3 \cdot \Omega_7(3) < E_6$, $p = 2$. Irreducible on $V_{78}$ and $V_{27}$. \possprim
\label{omega73e6p2}
\end{table}

\addtocounter{table}{1}
\begin{table}[H]\small
\thetable: $3 \cdot G_2(3) < E_6$, $p = 2$. Irreducible on $V_{78}$ and $V_{27}$. \possprim
\label{3g23e6p2}
\end{table}

\begin{table}[H]\small
\caption{$^{3}D_4(2) < E_6$, $p \neq 2,3$}
\begin{tabular}{r|*{2}{c}|*{2}{c}}
& \multicolumn{2}{c|}{$V_{78}$} & \multicolumn{2}{c}{$V_{27}$} \\
& 26 & 52 & 1 & 26 \\ \hline
1) & 1 & 1 & 1 & 1
\end{tabular}
\label{3d4e6p7}
\label{3d4e6p13}
\label{3d4e6p0}
\end{table}

\begin{table}[H]\small
\caption{$^{3}D_4(2) < E_6$, $p = 3$}
\begin{tabular}{r|*{3}{c}|*{2}{c}}
& \multicolumn{3}{c|}{$V_{78}$} & \multicolumn{2}{c}{$V_{27}$} \\
& 1 & 25 & 52 & 1 & 25 \\ \hline
1) & 1 & 1 & 1 & 2 & 1
\end{tabular}
\label{3d42e6p3}
\end{table}

\addtocounter{table}{1}
\begin{table}[H]\small
\thetable: $^{2}F_4(2)' < E_6$, $p \neq 2$, $3$. Irreducible on $V_{78}$ and $V_{27}$. \possprim
\label{2f42e6p5}
\label{2f42e6p13}
\label{2f42e6p0}
\end{table}

\begin{table}[H]\small
\caption{$^{2}F_4(2)' < E_6$, $p = 3$}
\begin{tabular}{rc|*{2}{c}|*{1}{c}}
& & \multicolumn{2}{c|}{$V_{78}$} & \multicolumn{1}{c}{$V_{27}$} \\
& & 1 & 77 & 27 \\ \hline
1) & \possprim & 1 & 1 & 1
\end{tabular}
\label{2f42e6p3}
\end{table}

\section{$E_7$}
\label{sec:E7tabs}

\subsection{Alternating Groups} \ 

\begin{table}[H]\small
\caption{Alt$_{13} < E_7$, $p = 2$}
\begin{tabular}{r|cccc|ccc}
& \multicolumn{4}{c|}{$V_{133}$} & \multicolumn{3}{c}{$V_{56}$} \\
& 1 & 32$_a$ & 32$_b$ & 64 & 12 & 32$_a$ & 32$_b$ \\ \hline
1) & 5 & 0 & 2 & 1 & 2 & 0 & 1
\end{tabular}\\
Permutations: $(32_a,32_b)$.
\label{a13e7p2}
\end{table}

\begin{table}[H]\small
\caption{Alt$_{12} < E_7$, $p = 2$}
\begin{tabular}{rc|ccccc|cccc}
& & \multicolumn{5}{c|}{$V_{133}$} & \multicolumn{4}{c}{$V_{56}$} \\
& & 1 & 10 & 16 & 16$^{*}$ & 44 & 1 & 10 & 16 & 16$^{*}$ \\ \hline
1) & \nongcr & 5 & 2 & 2 & 2 & 1 & 4 & 2 & 1 & 1
\end{tabular}
\label{a12e7p2}
\end{table}

\begin{table}[H]\small
\caption{Alt$_{11} < E_7$, $p = 2$}
\begin{tabular}{r|ccccc|cccc}
& \multicolumn{5}{c|}{$V_{133}$} & \multicolumn{4}{c}{$V_{56}$} \\
& 1 & 10 & 16 & 16$^{*}$ & 44 & 1 & 10 & 16 & 16$^{*}$ \\ \hline
1) & 5 & 2 & 2 & 2 & 1 & 4 & 2 & 1 & 1 
\end{tabular}
\label{a11e7p2}
\end{table}

\begin{table}[H]\small
\caption{Alt$_{10} < E_7$, $p = 5$}
\begin{tabular}{rc|cccc|c}
& & \multicolumn{4}{c|}{$V_{133}$} & \multicolumn{1}{c}{$V_{56}$} \\
& & 28 & 35$_a$ & 35$_b$ & 35$_c$ & 28 \\ \hline
1) & \possprim & 1 & 1 & 1 & 1 & 2
\end{tabular}
\label{a10e7p5}
\end{table}

\begin{table}[H]\small
\caption{Alt$_{10} < E_7$, $p = 2$}
\begin{tabular}{rc|cccc|cccc}
& & \multicolumn{4}{c|}{$V_{133}$} & \multicolumn{4}{c}{$V_{56}$} \\
& & 1 & 8 & 16 & 26 & 1 & 8 & 16 & 26 \\ \hline
1) & \nongcr & 5 & 1 & 1 & 4 & 4 & 0 & 0 & 2 \\
2) & \nongcr & 11 & 4 & 4 & 1 & 8 & 2 & 2 & 0
\end{tabular}
\label{a10e7p2}
\end{table}

\begin{table}[H]\small
\caption{Alt$_{9} < E_7$, $p = 0$ or $p \nmid |H|$}
\begin{tabular}{rc|ccccc|c}
& & \multicolumn{5}{c|}{$V_{133}$} & \multicolumn{1}{c}{$V_{56}$} \\
& & 8 & 27 & 28 & 35$_a$ & 35$_b$ & 28 \\ \hline
1) & \possprim & 1 & 1 & 1 & 1 & 1 & 2
\end{tabular}
\label{a9e7p0}
\end{table}

\begin{table}[H]\small
\caption{Alt$_{9} < E_7$, $p = 7$}
\begin{tabular}{rc|ccccc|c}
& & \multicolumn{5}{c|}{$V_{133}$} & \multicolumn{1}{c}{$V_{56}$} \\
& & 8 & 19 & 28 & 35$_a$ & 35$_b$ & 28 \\ \hline
1) & \possprim & 2 & 1 & 1 & 1 & 1 & 2
\end{tabular}
\label{a9e7p7}
\end{table}

\begin{table}[H]\small
\caption{Alt$_{9} < E_7$, $p = 5$}
\begin{tabular}{rc|ccccc|c}
& & \multicolumn{5}{c|}{$V_{133}$} & \multicolumn{1}{c}{$V_{56}$} \\
& & 8 & 27 & 28 & 35$_a$ & 35$_b$ & 28 \\ \hline
1) & \possprim & 1 & 1 & 1 & 1 & 1 & 2
\end{tabular}
\label{a9e7p5}
\end{table}

\begin{table}[H]\small
\caption{Alt$_{9} < E_7$, $p = 3$}
\begin{tabular}{rc|ccccc|cc}
& & \multicolumn{5}{c|}{$V_{133}$} & \multicolumn{2}{c}{$V_{56}$} \\
& & 1 & 7 & 21 & 27 & 35 & 7 & 21 \\ \hline
1) & \possprim, \nongcr & 1 & 2 & 1 & 1 & 2 & 2 & 2
\end{tabular}
\label{a9e7p3}
\end{table}

\begin{table}[H]\small
\caption{Alt$_{9} < E_7$, $p = 2$}
\begin{tabular}{rc|ccccc|ccccc}
& & \multicolumn{5}{c|}{$V_{133}$} & \multicolumn{5}{c}{$V_{56}$} \\
& & 1 & 8$_a$ & 8$_b$ & 8$_c$ & 26 & 1 & 8$_a$ & 8$_b$ & 8$_c$ & 26 \\ \hline
1) & \nongcr & 5 & 1 & 1 & 1 & 4 & 4 & 0 & 0 & 0 & 2 \\
2) & & 11 & 4 & 4 & 4 & 1 & 8 & 2 & 2 & 2 & 0
\end{tabular}
\label{a9e7p2}
\end{table}

\begin{table}[H]\small
\caption{Alt$_{8} < E_7$, $p = 0$ or $p \nmid |H|$}
\begin{tabular}{r|ccccc|cc}
& \multicolumn{5}{c|}{$V_{133}$} & \multicolumn{2}{c}{$V_{56}$} \\
& 1 & 7 & 20 & 21$_a$ & 35 & 7 & 21$_a$ \\ \hline
1) & 1 & 3 & 1 & 1 & 2 & 2 & 2
\end{tabular}
\label{a8e7p0}
\end{table}

\begin{table}[H]\small
\caption{Alt$_{8} < E_7$, $p = 7$} 
\begin{tabular}{r|ccccc|cc}
& \multicolumn{5}{c|}{$V_{133}$} & \multicolumn{2}{c}{$V_{56}$} \\
& 1 & 7 & 19 & 21$_a$ & 35 & 7 & 21$_a$ \\ \hline
1) & 2 & 3 & 1 & 1 & 2 & 2 & 2
\end{tabular}
\label{a8e7p7}
\end{table}

\begin{table}[H]\small
\caption{Alt$_{8} < E_7$, $p = 5$}
\begin{tabular}{rc|cccccccc|ccc}
& & \multicolumn{8}{c|}{$V_{133}$} & \multicolumn{3}{c}{$V_{56}$} \\
& & 1 & 7 & 13 & 20 & 21$_a$ & 35 & 43 & 70 & 7 & 21$_a$ & 21$_b$\\ \hline
1) & \possprim & 0 & 0 & 0 & 1 & 0 & 0 & 1 & 1 & 2 & 0 & 2\\
2) & & 1 & 2 & 1 & 0 & 0 & 1 & 0 & 1 & 2 & 0 & 2\\
3) & & 1 & 3 & 0 & 1 & 1 & 2 & 0 & 0 & 2 & 2 & 0
\end{tabular}
\label{a8e7p5}
\end{table}

\begin{table}[H]\small
\caption{Alt$_{8} < E_7$, $p = 3$}
\begin{tabular}{rc|*{5}{c}|cc}
& & \multicolumn{5}{c|}{$V_{133}$} & \multicolumn{2}{c}{$V_{56}$} \\
& & 1 & 7 & 13 & 21 & 35 & 7 & 21 \\ \hline
1) & \possprim, \nongcr & 1 & 4 & 1 & 1 & 2 & 2 & 2
\end{tabular}
\label{a8e7p3}
\end{table}

\begin{table}[H]\small
\caption{Alt$_{7} < E_7$, $p = 0$ or $p \nmid |H|$}
\begin{tabular}{r|*{6}{c}|*{5}{c}}
& \multicolumn{6}{c|}{$V_{133}$} & \multicolumn{5}{c}{$V_{56}$} \\
& 1 & 6 & 10 & 10$^{*}$ & 14$_a$ & 15 & 1 & 6 & 10 & 10$^{*}$ & 15 \\ \hline
1) & 8 & 1 & 0 & 0 & 1 & 7 & 0 & 6 & 1 & 1 & 0\\
2) & 4 & 5 & 2 & 2 & 1 & 3 & 2 & 4 & 0 & 0 & 2
\end{tabular}\\
$14_a$ is a section of $5 \otimes 5$.
\label{a7e7p0}
\end{table}

\begin{table}[H]\small
\caption{Alt$_{7} < E_7$, $p = 7$}
\begin{tabular}{rc|*{4}{c}|*{3}{c}}
& & \multicolumn{4}{c|}{$V_{133}$} & \multicolumn{3}{c}{$V_{56}$} \\
& & 1 & 5 & 10 & 14$_a$ & 1 & 5 & 10 \\ \hline
1) & \nongcr & 9 & 8 & 7 & 1 & 6 & 6 & 2
\end{tabular}\\
$14_a$ is a section of $5 \otimes 5$.
\label{a7e7p7}
\end{table}

\begin{table}[H]\small
\caption{Alt$_{7} < E_7$, $p = 5$}
\begin{tabular}{rc|cccccccc|ccccccc}
& & \multicolumn{8}{c|}{$V_{133}$} & \multicolumn{7}{c}{$V_{56}$} \\
& & 1 & 6 & 8 & 10 & 10$^{*}$ & 13 & 15 & 35 & 1 & 6 & 8 & 10 & 10$^{*}$ & 13 & 15 \\ \hline
1) & \possprim & 0 & 3 & 7  & 1 & 1 & 3 & 0 & 0 & 0 & 0 & 2 & 2 & 2 & 0 & 0\\
2) & & 8 & 2 & 1  & 0 & 0 & 0 & 7 & 0 & 0 & 6 & 0 & 1 & 1 & 0 & 0\\
3) & \nongcr & 1 & 2 & 3  & 0 & 0 & 2 & 0 & 2 & 2 & 2 & 2 & 0 & 0 & 2 & 0\\
4) & & 1 & 3 & 2  & 0 & 0 & 1 & 1 & 2 & 2 & 2 & 2 & 0 & 0 & 2 & 0\\
5) & \nongcr & 4 & 1 & 2  & 1 & 1 & 4 & 0 & 1 & 2 & 2 & 2 & 0 & 0 & 2 & 0\\
6) & \nongcr & 4 & 2 & 1  & 1 & 1 & 3 & 1 & 1 & 2 & 2 & 2 & 0 & 0 & 2 & 0\\
7) & & 4 & 6 & 1  & 2 & 2 & 0 & 3 & 0 & 2 & 4 & 0 & 0 & 0 & 0 & 2\\
8) & & 9 & 0 & 13 & 1 & 1 & 0 & 0 & 0 & 8 & 0 & 6 & 0 & 0 & 0 & 0
\end{tabular}
\label{a7e7p5}
\end{table}

\begin{table}[H]\small
\caption{$2\cdot$Alt$_{7} < E_7$, $p = 5$}
\begin{tabular}{rc|*{4}{c}|*{5}{c}}
& & \multicolumn{4}{c|}{$V_{133}$} & \multicolumn{5}{c}{$V_{56}$} \\
& & 8 & 10 & 10$^{*}$ & 35 & 4 & 4$^{*}$ & 14$_a$ & 14$_b$ & 20$_a$ \\ \hline
1) & \possprim & 1 & 1 & 1 & 3 & 1 & 1 & 0 & 2 & 1
\end{tabular}\\
Permutations: $(14_a,14_b)$.
\label{2a7e7p5}
\end{table}

\begin{table}[H]\small
\caption{Alt$_{7} < E_7$, $p = 3$}
\begin{tabular}{rc|cccccc|cccccc}
& & \multicolumn{6}{c|}{$V_{133}$} & \multicolumn{6}{c}{$V_{56}$} \\
& & 1 & 6 & 10 & 10$^{*}$ & 13 & 15 & 1 & 6 & 10 & 10$^{*}$ & 13 & 15 \\ \hline
1) & & 9 & 1 & 0 & 0 & 1 & 7 & 0 & 6 & 1 & 1 & 0 & 0\\
2) & \nongcr & 5 & 5 & 2 & 2 & 1 & 3 & 2 & 4 & 0 & 0 & 0 & 2
\end{tabular}
\label{a7e7p3}
\end{table}

\begin{table}[H]\small
\caption{$2\cdot$Alt$_{7} < E_7$, $p = 3$}
\begin{tabular}{rc|*{5}{c}|*{4}{c}}
& & \multicolumn{5}{c|}{$V_{133}$} & \multicolumn{4}{c}{$V_{56}$} \\
& & 1 & 10 & 10$^{*}$ & 13 & 15 & 4 & 4$^{*}$ & 6$_b$ & 6$_c$ \\ \hline
1) & \nongcr & 22 & 1 & 1 & 7 & 0 & 1 & 1 & 1 & 7
\end{tabular}\\
Permutations: $(6_b,6_c)$.
\label{2a7e7p3}
\end{table}

\begin{table}[H]\small
\caption{Alt$_{7} < E_7$, $p = 2$}
\begin{tabular}{rc|cccccc|cccccc}
& & \multicolumn{6}{c|}{$V_{133}$} & \multicolumn{6}{c}{$V_{56}$} \\
& & 1 & 4 & 4$^{*}$ & 6 & 14 & 20 & 1 & 4 & 4$^{*}$ & 6 & 14 & 20 \\ \hline
1) & \nongcr & 15 & 0 & 0 & 1 & 8 & 0 & 0 & 1 & 1 & 8 & 0 & 0\\
2) & \nongcr & 7 & 2 & 2 & 5 & 0 & 4 & 4 & 0 & 0 & 2 & 0 & 2\\
3) & \nongcr & 7 & 2 & 2 & 6 & 1 & 3 & 4 & 0 & 0 & 2 & 0 & 2\\
4) & \nongcr & 7 & 2 & 2 & 9 & 4 & 0 & 4 & 0 & 0 & 4 & 2 & 0\\
5) & & 19 & 8 & 8 & 5 & 0 & 1 & 12 & 4 & 4 & 2 & 0 & 0\\
6) & & 19 & 8 & 8 & 6 & 1 & 0 & 12 & 4 & 4 & 2 & 0 & 0
\end{tabular}
\label{a7e7p2}
\end{table}

{\small
\begin{longtable}{rc|ccccccc|ccccccc}
\caption{Alt$_{6} < E_7$, $p = 0$ or $p \nmid |H|$} \\
& & \multicolumn{7}{c|}{$V_{133}$} & \multicolumn{7}{c}{$V_{56}$} \\
& & 1 & 5$_a$ & 5$_b$ & 8$_a$ & $8_b$ & 9 & 10 & 1 & 5$_a$ & 5$_b$ & 8$_a$ & $8_b$ & 9 & 10 \\ \hline
1) & \possprim & 0 & 3 & 3 & 3 & 4  & 3 & 2 & 0 & 0 & 0 & 0 & 2 & 0 & 4 \\
2) & & 1 & 2 & 5 & 2 & 3  & 3 & 3 & 0 & 0 & 0 & 0 & 2 & 0 & 4 \\
3) & & 2 & 1 & 1 & 3 & 4  & 5 & 2 & 0 & 0 & 0 & 0 & 2 & 0 & 4 \\
4) & & 2 & 4 & 4 & 1 & 2  & 3 & 4 & 0 & 0 & 0 & 0 & 2 & 0 & 4 \\
5) & & 3 & 0 & 3 & 2 & 3  & 5 & 3 & 0 & 0 & 0 & 0 & 2 & 0 & 4 \\
6) & & 4 & 2 & 2 & 1 & 2  & 5 & 4 & 0 & 0 & 0 & 0 & 2 & 0 & 4 \\
7) & & 3 & 0 & 0 & 3 & 4  & 6 & 2 & 0 & 4 & 4 & 0 & 2 & 0 & 0 \\
8) & & 5 & 1 & 1 & 1 & 2  & 6 & 4 & 0 & 4 & 4 & 0 & 2 & 0 & 0 \\
9) & & 2 & 2 & 2 & 3 & 5  & 3 & 2 & 1 & 2 & 2 & 1 & 0 & 3 & 0 \\
10)& & 3 & 1 & 4 & 2 & 4  & 3 & 3 & 1 & 2 & 2 & 1 & 0 & 3 & 0 \\
11)& & 4 & 3 & 3 & 1 & 3  & 3 & 4 & 1 & 2 & 2 & 1 & 0 & 3 & 0 \\
12)& & 1 & 2 & 2 & 3 & 4  & 4 & 2 & 2 & 2 & 2 & 0 & 2 & 2 & 0 \\
13)& & 2 & 1 & 4 & 2 & 3  & 4 & 3 & 2 & 2 & 2 & 0 & 2 & 2 & 0 \\
14)& & 3 & 3 & 3 & 1 & 2  & 4 & 4 & 2 & 2 & 2 & 0 & 2 & 2 & 0 \\
15)& & 4 & 0 & 0 & 3 & 5  & 5 & 2 & 3 & 0 & 0 & 1 & 0 & 5 & 0 \\
16)& & 6 & 1 & 1 & 1 & 3  & 5 & 4 & 3 & 0 & 0 & 1 & 0 & 5 & 0 \\
17)& & 3 & 0 & 0 & 3 & 4  & 6 & 2 & 4 & 0 & 0 & 0 & 2 & 4 & 0 \\
18)& & 5 & 1 & 1 & 1 & 2  & 6 & 4 & 4 & 0 & 0 & 0 & 2 & 4 & 0 \\
19)& & 4 & 2 & 2 & 5 & 5  & 1 & 2 & 4 & 2 & 2 & 2 & 2 & 0 & 0 \\
20)& & 5 & 1 & 4 & 4 & 4  & 1 & 3 & 4 & 2 & 2 & 2 & 2 & 0 & 0 \\
21)& & 6 & 3 & 3 & 3 & 3  & 1 & 4 & 4 & 2 & 2 & 2 & 2 & 0 & 0 \\
22)& & 6 & 0 & 0 & 5 & 5  & 3 & 2 & 6 & 0 & 0 & 2 & 2 & 2 & 0 \\
23)& & 8 & 1 & 1 & 3 & 3  & 3 & 4 & 6 & 0 & 0 & 2 & 2 & 2 & 0 \\
24)& & 9 & 0 & 9 & 0 & 0  & 1 & 7 & 6 & 0 & 6 & 0 & 0 & 0 & 2 \\
25)& & 9 & 0 & 0 & 0 & 13 & 0 & 2 & 8 & 0 & 0 & 0 & 6 & 0 & 0 \\
\multicolumn{16}{c}{Permutations: $(5_a,5_b)$, $(8_a,8_b)$.}
\label{a6e7p0}
\end{longtable}
}

{\small
\begin{longtable}{rc|ccccccc|ccccccc}
\caption{$2\cdot$Alt$_{6} < E_7$, $p = 0$ or $p \nmid |H|$} \\
& & \multicolumn{7}{c|}{$V_{133}$} & \multicolumn{6}{c}{$V_{56}$} \\
& & 1 & 5$_a$ & 5$_b$ & 8$_a$ & 8$_b$ & 9 & 10$_a$ & 4$_a$ & 4$_b$ & 8$_c$ & 8$_d$ & 10$_b$ & 10$_c$ \\ \cline{1-14}
1) & \possprim & 0 & 0 & 0 & 4 & 3 & 3 & 5 & 0 & 0 & 0 & 2 & 1 & 3 \\
2) & & 2 & 1 & 1 & 2 & 1 & 3 & 7 & 0 & 0 & 0 & 2 & 1 & 3 \\
3) & \possprim & 0 & 0 & 0 & 4 & 3 & 3 & 5 & 0 & 0 & 0 & 2 & 3 & 1 \\
4) & & 2 & 1 & 1 & 2 & 1 & 3 & 7 & 0 & 0 & 0 & 2 & 3 & 1 \\
5) & & 3 & 1 & 1 & 3 & 1 & 2 & 7 & 0 & 3 & 2 & 1 & 0 & 2 \\
6) & & 3 & 1 & 1 & 3 & 1 & 2 & 7 & 0 & 3 & 2 & 1 & 2 & 0 \\
7) & & 1 & 0 & 0 & 5 & 3 & 2 & 5 & 1 & 1 & 1 & 0 & 1 & 3 \\
8) & & 3 & 1 & 1 & 3 & 1 & 2 & 7 & 1 & 1 & 1 & 0 & 1 & 3 \\
9) & & 1 & 0 & 0 & 5 & 3 & 2 & 5 & 1 & 1 & 1 & 0 & 3 & 1 \\
10)& & 3 & 1 & 1 & 3 & 1 & 2 & 7 & 1 & 1 & 1 & 0 & 3 & 1 \\
11)& & 3 & 1 & 1 & 3 & 1 & 2 & 7 & 3 & 0 & 2 & 1 & 0 & 2 \\
12)& & 3 & 1 & 1 & 3 & 1 & 2 & 7 & 3 & 0 & 2 & 1 & 2 & 0 \\
\multicolumn{16}{c}{Permutations: $(4_a,4_b)(5_a,5_b)$, $(8_a,8_b)(8_c,8_d)(10_b,10_c)$.}
\label{2a6e7p0}
\end{longtable}
}

\begin{table}[H]\small
\caption{Alt$_{6} < E_7$, $p = 5$}
\begin{tabular}{rc|ccccc|ccccc}
& & \multicolumn{5}{c|}{$V_{133}$} & \multicolumn{5}{c}{$V_{56}$} \\
& & 1 & 5$_a$ & 5$_b$ & 8 & 10 & 1 & 5$_a$ & 5$_b$ & 8 & 10 \\ \hline
1) & \possprim, \nongcr & 3 & 3 & 3 & 10 & 2 & 0 & 0 & 0 & 2 & 4\\
2) & \possprim, \nongcr & 4 & 2 & 5 & 8 & 3 & 0 & 0 & 0 & 2 & 4\\
3) & \possprim, \nongcr & 5 & 4 & 4 & 6 & 4 & 0 & 0 & 0 & 2 & 4\\
4) & \possprim, \nongcr & 7 & 1 & 1 & 12 & 2 & 0 & 0 & 0 & 2 & 4\\
5) & \possprim, \nongcr & 8 & 0 & 3 & 10 & 3 & 0 & 0 & 0 & 2 & 4\\
6) & \nongcr & 9 & 2 & 2 & 8 & 4 & 0 & 0 & 0 & 2 & 4\\
7) & \possprim, \nongcr & 9 & 0 & 0 & 13 & 2 & 0 & 4 & 4 & 2 & 0\\
8) & \nongcr & 11 & 1 & 1 & 9 & 4 & 0 & 4 & 4 & 2 & 0\\
9) & \nongcr & 5 & 2 & 2 & 11 & 2 & 4 & 2 & 2 & 4 & 0\\
10) & \nongcr & 6 & 4 & 1 & 9 & 3 & 4 & 2 & 2 & 4 & 0\\
11) & \nongcr & 7 & 3 & 3 & 7 & 4 & 4 & 2 & 2 & 4 & 0\\
12) & & 10 & 0 & 9 & 1 & 7 & 6 & 0 & 6 & 0 & 2\\
13) & \nongcr & 9 & 0 & 0 & 13 & 2 & 8 & 0 & 0 & 6 & 0\\
14) & \nongcr & 11 & 1 & 1 & 9 & 4 & 8 & 0 & 0 & 6 & 0
\label{a6e7p5}
\end{tabular}\\
Permutations: $(5_a,5_b)$.
\end{table}

\begin{table}[H]\small
\caption{$2\cdot$Alt$_{6} < E_7$, $p = 5$}
%
}

\begin{table}[H]\small
\caption{$2\cdot$Alt$_{5} < E_7$, $p = 0$ or $p \nmid |H|$}
%
\\
Permutations: $(2_a,2_b)(3_b,3_b)$. \\
$4_b$ is faithful for $2\cdot\Alt_5$, $4_a$ is not.
\label{2a5e7p0}
\end{table}

\begin{table}[H]\small
\caption{Alt$_{5} < E_7$, $p = 3$}
\begin{tabular}{rc|cccc|cccc}
& & \multicolumn{4}{c|}{$V_{133}$} & \multicolumn{4}{c}{$V_{56}$}\\
& & 1 & 3$_a$ & 3$_b$ & 4 & 1 & 3$_a$ & 3$_b$ & 4\\ \hline
1) & \possprim, \nongcr & 17 & 5 & 7 & 20 & 1 & 5 & 4 & 7\\
2) & \nongcr & 21 & 6 & 10 & 16 & 2 & 1 & 9 & 6\\
3) & \possprim, \nongcr & 16 & 5 & 6 & 21 & 2 & 4 & 6 & 6\\
4) & \nongcr & 19 & 7 & 7 & 18 & 4 & 6 & 6 & 4\\
5) & \nongcr & 27 & 0 & 22 & 10 & 6 & 0 & 14 & 2\\
6) & \nongcr & 22 & 2 & 15 & 15 & 6 & 4 & 10 & 2\\
7) & \possprim, \nongcr & 17 & 5 & 7 & 20 & 9 & 1 & 0 & 11\\
8) & \nongcr & 16 & 5 & 6 & 21 & 10 & 0 & 2 & 10\\
9) & \nongcr & 19 & 7 & 7 & 18 & 12 & 2 & 2 & 8\\
10) & \nongcr & 22 & 2 & 15 & 15 & 14 & 0 & 6 & 6\\
11) & \nongcr & 36 & 0 & 31 & 1 & 20 & 0 & 12 & 0
\end{tabular}\\
Permutations: $(3_a,3_b)$.
\label{a5e7p3}
\end{table}

\begin{table}[H]\small
\caption{$2\cdot$Alt$_{5} < E_7$, $p = 3$}
%
\\
Permutations: $(2_a,2_b)(3_a,3_b)$.
\label{2a5e7p3}
\end{table}

\subsection{Sporadic Groups} \ 

\begin{table}[H]\small
\caption{$M_{11} < E_7$, $p = 0$ or $p \nmid |H|$}
%
\label{m11e7p0}
\end{table}

\begin{table}[H]\small
\caption{$M_{11} < E_7$, $p = 11$}
%
\label{m11e7p11}
\end{table}

\begin{table}[H]\small
\caption{$M_{11} < E_7$, $p = 5$}
%
\label{m11e7p5}
\end{table}

\begin{table}[H]\small
\caption{$M_{11} < E_7$, $p = 3$}
%
\label{m11e7p2}
\end{table}

\begin{table}[H]\small
\caption{$2\cdot M_{12} < E_7$, $p = 0$ or $p \nmid |H|$}
%
\label{2m12e7p0}
\end{table}

\begin{table}[H]\small
\caption{$2\cdot M_{12} < E_7$, $p = 11$}
%
\label{2m12e7p11}
\end{table}

\begin{table}[H]\small
\caption{$M_{12} < E_7$, $p = 5$}
%
\\
Permutations: $(11_a,11_b)$.
\label{m12e7p5}
\end{table}

\begin{table}[H]\small
\caption{$2\cdot M_{12} < E_7$, $p = 5$}
%
\\
$55_c$ is fixed by $\textup{Out}(G)$.
\label{2m12e7p5}
\end{table}

\begin{table}[H]\small
\caption{$2\cdot M_{12} < E_7$, $p = 3$}
%
\label{2m12e7p3}
\end{table}

\begin{table}[H]\small
\caption{$M_{12} < E_7$, $p = 2$}
%
\label{m12e7p2}
\end{table}

\begin{table}[H]\small
\caption{$2\cdot M_{22} < E_7$, $p = 5$}
%
\label{2m22e7p5}
\end{table}

\begin{table}[H]\small
\caption{$J_{1} < E_7$, $p = 11$}
%
\label{j1e7p11}
\end{table}

\begin{table}[H]\small
\caption{$2 \cdot J_{2} < E_7$, $p = 0$ or $p \nmid |H|$}
%
\\
Permutations: $(6_a,6_b)(14_a,14_b)(21_a,21_b)$.
\label{2j2e7p0}
\end{table}

\begin{table}[H]\small
\caption{$2 \cdot J_{2} < E_7$, $p = 7$}
%
\\
Permutations: $(6_a,6_b)(14_a,14_b)(21_a,21_b)$.
\label{2j2e7p7}
\end{table}

\begin{table}[H]\small
\caption{$2 \cdot J_{2} < E_7$, $p = 5$}
%
\\
Outer automorphism fixes these isomorphism types.
\label{2j2e7p5}
\end{table}

\begin{table}[H]\small
\caption{$2 \cdot J_{2} < E_7$, $p = 3$}
%
\\
Permutations: $(6_a,6_b)(13_a,13_b)(21_a,21_b)$.
\label{2j2e7p3}
\end{table}

\begin{table}[H]\small
\caption{$J_{2} < E_7$, $p = 2$}
%
\\
Permutations: $(6_a,6_b)(14_a,14_b)(64_a,64_b)$.
\label{j2e7p2}
\end{table}

\begin{table}[H]\small
\caption{$2 \cdot Ru < E_7$, $p = 5$}
%
\label{2Rue7p5}
\end{table}

\begin{table}[H]\small
\caption{$2 \cdot HS < E_7$, $p = 5$}
%
\label{2HSe7p5}
\end{table}

\subsection{Cross-characteristic Groups $L_2(q)$ $(q \neq 4$, $5$, $9)$} \ 

{\small
%
}

\begin{table}[H]\small
\caption{$SL_2(7) < E_7$, $p = 0$ or $p \nmid |H|$}
%
\\
Permutations: $(6_{b},6_{c})$
\label{2l27e7p0}
\end{table}

\begin{table}[H]\small
\caption{$L_2(7) < E_7$, $p = 3$}
%
\label{l27e7p3}
\end{table}

\begin{table}[H]\small
\caption{$SL_2(7) < E_7$, $p = 3$}
%
\\
Permutations: $(6_{b},6_{c})$
\label{2l27e7p3}
\end{table}

{\small
\begin{table}[H]\small
\caption{$L_2(8) < E_7$, $p = 0$ or $p \nmid |H|$}
%
\\
Permutations: $(7_b,7_c,7_d)(9_a,9_b,9_c)$.
\label{l28e7p0}
\end{table}
}

\begin{table}[H]\small
\caption{$L_2(8) < E_7$, $p = 7$}
%
\\
Permutations: $(7_b,7_c,7_d)$.
\label{l28e7p7}
\end{table}

\begin{table}[H]\small
\caption{$L_2(8) < E_7$, $p = 3$}
%
\\
Permutations: $(9_a,9_b,9_c)$.
\label{l28e7p3}
\end{table}

\begin{table}[H]\small
\caption{$L_2(11) < E_7$, $p = 0$ or $p \nmid |H|$}
%
\\
Permutations: $(12_{a},12_{b})$.
\label{l211e7p0}
\end{table}

\begin{table}[H]\small
\caption{$SL_2(11) < E_7$, $p = 0$ or $p \nmid |H|$}
%
\\
Permutations: $(12_{a},12_{b})(12_{c},12_{d})$.
\label{2l211e7p0}
\end{table}

\begin{table}[H]\small
\caption{$L_2(11) < E_7$, $p = 5$}
%
\label{l211e7p5}
\end{table}

\begin{table}[H]\small
\caption{$SL_2(11) < E_7$, $p = 5$}
%
\label{2l211e7p5}
\end{table}

\begin{table}[H]\small
\caption{$L_2(11) < E_7$, $p = 3$}
%
\label{l211e7p3}
\end{table}

\begin{table}[H]\small
\caption{$SL_2(11) < E_7$, $p = 3$}
%
\\
Permutations: $(12_a,12_b)(12_c,12_d)$.
\label{2l211e7p3}j
\end{table}

\begin{table}[H]\small
\caption{$L_2(11) < E_7$, $p = 2$}
%
\\
Permutations: $(12_a,12_b)$.
\label{l211e7p2}
\end{table}

{\small
\begin{table}[H]\small
\caption{$L_2(13) < E_7$, $p = 0$ or $p \nmid |H|$}
%
\\
Permutations: $(7_a,7_b)$, $(12_a,12_b,12_c)$.
\label{l213e7p0}
\end{table}
}

{\small
\begin{table}[H]\small
\caption{$SL_2(13) < E_7$, $p = 0$ or $p \nmid |H|$}
%
\\
Permutations: $(6_a,6_b)(7_a,7_b)$, $(12_a,12_b,12_c)(12_d,12_e,12_f)$.\\
$12_a$ and $12_b$ occur with equal multiplicities, as do $14_d$ and $14_e$.
\label{2l213e7p0}
\end{table}
}

\begin{table}[H]\small
\caption{$L_2(13) < E_7$, $p = 7$}
%
\\
Permutations: $(7_a,7_b)$.
\label{l213e7p7}
\end{table}

\begin{table}[H]\small
\caption{$SL_2(13) < E_7$, $p = 7$}
%
\\
Permutations: $(7_{a},7_{b})(6_{a},6_{b})$.
\label{2l213e7p7}
\end{table}

\begin{table}[H]\small
\caption{$L_2(13) < E_7$, $p = 3$}
%
\\
Permutations: $(7_a,7_b)$, $(12_a,12_b,12_c)$.
\label{l213e7p3}
\end{table}

\begin{table}[H]\small
\caption{$SL_2(13) < E_7$, $p = 3$}
%
\\
Permutations: $(6_a,6_b)(7_a,7_b)$, $(12_a,12_b,12_c)(12_d,12_e,12_f)$.
\label{2l213e7p3}
\end{table}

\begin{table}[H]\small
\caption{$L_2(13) < E_7$, $p = 2$}
%
\\
Permutations: $(6_a,6_b)$, $(12_a,12_b,12_c)$.
\label{l213e7p2}
\end{table}

\begin{table}[H]\small
\caption{$L_2(17) < E_7$, $p = 0$ or $p \nmid |H|$}
%
\\
Permutations: $(9_a,9_b)$.
\label{l217e7p0}
\end{table}

\begin{table}[H]\small
\caption{$L_2(17) < E_7$, $p = 3$}
%
\\
Permutations: $(9_a,9_b)$.
\label{l217e7p3}
\end{table}

\begin{table}[H]\small
\caption{$L_2(17) < E_7$, $p = 2$}
%
\\
Permutations: $(8_a,8_b)$.
\label{l217e7p2}
\end{table}

\begin{table}[H]\small
\caption{$L_2(19) < E_7$, $p = 0$ or $p \nmid |H|$}
%
\label{l219e7p0}
\end{table}

\begin{table}[H]\small
\caption{$SL_2(19) < E_7$, $p = 0$ or $p \nmid |H|$}
%
\\
Permutation: $(18_b,18_c)(18_f,18_g)(18_h,18_i)$.
\label{2l219e7p0}
\end{table}

\begin{table}[H]\small
\caption{$L_2(19) < E_7$, $p = 5$}
%
\label{2l219e7p5}
\end{table}

\begin{table}[H]\small
\caption{$L_2(19) < E_7$, $p = 3$}
%
\label{2l219e7p3}
\end{table}

\begin{table}[H]\small
\caption{$L_2(19) < E_7$, $p = 2$}
%
\label{l219e7p2}
\end{table}

\begin{table}[H]\small
\caption{$L_2(25) < E_7$, $p = 0$ or $p \nmid |H|$}
%
\label{l225e7p0}
\end{table}

\begin{table}[H]\small
\caption{$L_2(25) < E_7$, $p = 13$}
%
\label{l225e7p13}
\end{table}

\begin{table}[H]\small
\caption{$L_2(25) < E_7$, $p = 3$}
%
\label{l225e7p3}
\end{table}

\begin{table}[H]\small
\caption{$L_2(25) < E_7$, $p = 2$}
%
\label{l225e7p2}
\end{table}

\begin{table}[H]\small
\caption{$L_2(27) < E_7$, $p = 0$ or $p \nmid |H|$}
%
\label{l227e7p0}
\end{table}

\begin{table}[H]\small
\caption{$SL_2(27) < E_7$, $p = 0$ or $p \nmid |H|$}
%
\\
Permutations: $(28_a,28_b,\ldots,28_f)(28_g,28_h,\ldots,28_l)$.
\label{2l227e7p0}
\end{table}

\begin{table}[H]\small
\caption{$L_2(27) < E_7$, $p = 13$}
%
\\
Permutations: $(26_a,26_b,26_c)(26_d,26_e,26_f)$.
\label{l227e7p13}
\end{table}

\begin{table}[H]\small
\caption{$SL_2(27) < E_7$, $p = 13$}
%
\\
Permutations: $(28_a,28_b,\ldots,28_f)(28_g,28_h,\ldots,28_l)$. \\
\label{2l227e7p7}
\end{table}

\begin{table}[H]\small
\caption{$L_2(27) < E_7$, $p = 2$}
%
\\
Permutations: $(26_a,26_b,26_c)$, $(28_a,28_b,28_c,28_d,28_e,28_f)$.
\label{l227e7p2}
\end{table}

\begin{table}[H]\small
\caption{$SL_2(29) < E_7$, $p = 0$ or $p \nmid |H|$}
%
\\
Permutations: $(15_a,15_b)$, $(28_a,28_b)(28_i,28_j)(28_k,28_l)$.
\label{2l229e7p0}
\end{table}

\begin{table}[H]\small
\caption{$SL_2(29) < E_7$, $p = 7$}
%
\\
Permutations: $(15_a,15_b)$, $(28_k,28_l)(28_m,28_n)$.
\label{2l229e7p7}
\end{table}

\begin{table}[H]\small
\caption{$SL_2(29) < E_7$, $p = 5$}
%
\\
Permutations: $(15_a,15_b)$.
\label{2l229e7p5}
\end{table}

\begin{table}[H]\small
\caption{$SL_2(29) < E_7$, $p = 3$}
%
\\
Permutations: $(15_a,15_b)$, $(28_b,28_c)(28_d,28_e)$.
\label{2l229e7p3}
\end{table}	

\begin{table}[H]\small
\caption{$L_2(29) < E_7$, $p = 2$}
%
\\
Permutations: $(14_a,14_b)$, $(28_b,28_c)(28_d,28_e)(28_f,28_g)$.
\label{l229e7p2}
\end{table}

\begin{table}[H]\small
\caption{$SL_2(37) < E_7$, $p = 0$ or $p \nmid |H|$}
%
\\
Permutations: $(18_a,18_b)$, $(19_a,19_b)$, $(38_i,38_j,38_k)$.
\label{2l237e7p0}
\end{table}

\begin{table}[H]\small
\caption{$SL_2(37) < E_7$, $p = 19$}
%
\\
Permutations: $(18_a,18_b)$, $(19_a,19_b)$, $(38_i,38_l,38_m)$.
\label{2l237e7p19}
\end{table}

\begin{table}[H]\small
\caption{$SL_2(37) < E_7$, $p = 3$}
%
\\
where $r = 1,\ldots,7$.\\
Permutations: $(18_a,18_b)$.
\label{2l237e7p3}
\end{table}

\begin{table}[H]\small
\caption{$L_2(37) < E_7$, $p = 2$}
%
\\
Permutations: $(18_a,18_b)$.
\label{l237e7p2}
\end{table}

\subsection{Cross-characteristic Groups $\ncong L_2(q)$} \ 

\begin{table}[H]\small
%
\\
$16_a^{*}$, $16_b^{*}$ and $26_b^{*}$ occur with the same multiplicities as their duals.\\
Permutations: $(16_a,16_b)(16_a^{*},16_b^{*})$.
\label{l33e7p0}
\end{table}

\begin{table}[H]\small
\caption{$L_3(3) < E_7$, $p = 13$}
%
\\
Permutations: $(16_a,16_b)(16_a^{*},16_b^{*})$.
\label{l33e7p2}
\end{table}

\begin{table}[H]\small
\caption{$2 \cdot L_3(4) < E_7$, $p = 0$ or $p \nmid |H|$}
Note that although $L_3(4)$ has Schur multiplier $C_3 \times C_4 \times C_4$, all double-covers of $L_3(4)$ are isomorphic.\\
%
\label{2l34e7p0}
\end{table}

\begin{table}[H]\small
\caption{$2 \cdot L_3(4) < E_7$, $p = 7$}
%
\label{2l34e7p7}
\end{table}

\begin{table}[H]\small
\caption{$2 \cdot L_3(4) < E_7$, $p = 5$}
%
\label{2l34e7p5}
\end{table}

\begin{table}[H]\small
\caption{$2 \cdot L_3(4) < E_7$, $p = 3$}
%
\\
$15_{b} = \bigwedge^{2}(6)$.
\label{2l34e7p3}
\end{table}

\begin{table}[H]\small
\caption{$L_4(3) < E_7$, $p = 2$}
%
\label{l43e7p2}
\end{table}

\begin{table}[H]\small
\caption{$U_3(3) < E_7$, $p = 0$ or $p \nmid |H|$}
Here, the duals of $7_b$, $21_b$, 28 and 32 also occur, with the same multiplicity. \\
%
\label{u33e7p0}
\end{table}

\begin{table}[H]\small
\caption{$U_3(3) < E_7$, $p = 7$}
Here, the duals of $7_b$, $21_b$ and $28$ also occur, with the same multiplicity.
%
\label{u33e7p7}
\end{table}

\begin{table}[H]\small
\caption{$U_3(8) < E_7$, $p \neq 2$}
%
\\
Permutations: $(133_a,133_b,133_c)$.
\label{u38e7p3}
\label{u38e7p7}
\label{u38e7p19}
\label{u38e7p0}
\end{table}

\begin{table}[H]\small
\caption{$U_4(2) < E_7$, $p = 0$ or $p \nmid |H|$}
%
\label{u42e7p0}
\end{table}

\begin{table}[H]\small
\caption{$U_4(2) < E_7$, $p = 5$}
%
\label{u42e7p5}
\end{table}

\begin{table}[H]\small
\caption{$Sp_6(2) < E_7$, $p = 0$ or $p \nmid |H|$}
%
\label{sp62e7p0}
\end{table}

\begin{table}[H]\small
\caption{$Sp_6(2) < E_7$, $p = 7$}
%
\label{sp62e7p7}
\end{table}

\begin{table}[H]\small
\caption{$Sp_6(2) < E_7$, $p = 5$}
%
\label{sp62e7p5}
\end{table}

\begin{table}[H]\small
\caption{$Sp_6(2) < E_7$, $p = 3$}
%
\label{sp62e7p3}
\end{table}

\begin{table}[H]\small
\caption{$\Omega_8^{+}(2) < E_7$, $p \neq 2$}
%
\label{omega8plus2e7p3}
\label{omega8plus2e7p5}
\label{omega8plus2e7p7}
\label{omega8plus2e7p30}
\end{table}

\begin{table}[H]\small
\caption{$^{3}D_4(2) < E_7$, $p \neq 2, 3$}
%
\label{3d42e7p7}
\label{3d42e7p13}
\label{3d42e7p0}
\end{table}

\begin{table}[H]\small
\caption{$^{3}D_4(2) < E_7$, $p = 3$}
%
\label{3d42e7p3}
\end{table}

\begin{table}[H]\small
\caption{$^{2}F_4(2)' < E_7$, $p \neq 2, 3, 5$}
%
\label{2f42e7p13}
\label{2f42e7p0}
\end{table}

\begin{table}[H]\small
\caption{$^{2}F_4(2)' < E_7$, $p = 5$}
%
\label{2f42e7p5}
\end{table}

\begin{table}[H]\small
\caption{$^{2}F_4(2)' < E_7$, $p = 3$}
%
\label{2f42e7p3}
\end{table}

\begin{table}[H]\small
\caption{$^{2}B_2(8) < E_7$, $p = 5$}
%
\end{table}

\section{$E_8$} \label{sec:E8tabs}

\subsection{Alternating Groups} \ 

Each group in the following table has a unique feasible character on $L(E_8(K))$:
\begin{table}[H]\small
\caption{Alt$_{n} < E_8$, $n \ge 10$, unique feasible character}
%
\label{a17e8p2}
\label{a16e8p2}
\label{a15e8p2}
\label{a14e8p2}
\label{a13e8p2}
\label{a12e8p2}
\label{a11e8p11}
\label{a11e8p2}
\label{a10e8p0}
\label{a10e8p7}
\label{a10e8p3}
\end{table}

\begin{table}[H]\small
\caption{Alt$_{10} < E_8$, $p = 5$}
%
\\
Permutations: $(64_a,64_b)$.
\label{a10e8p2}
\end{table}

\begin{table}[H]\small
\caption{Alt$_9 < E_8$, $p = 0$ or $p \nmid |H|$}
%
\label{a9e8p0}
\end{table}

\begin{table}[H]\small
\caption{Alt$_9 < E_8$, $p = 7$}
%
\label{a9e8p7}
\end{table}

\begin{table}[H]\small
\caption{Alt$_9 < E_8$, $p = 5$}
%
\label{a9e8p5}
\end{table}

\begin{table}[H]\small
\caption{Alt$_9 < E_8$, $p = 3$}
%
\\
Permutations: $(8_b,8_c)$.
\label{a9e8p2}
\end{table}

\begin{table}[H]\small
\caption{Alt$_8 < E_8$, $p = 0$ or $p \nmid |H|$}
%
\\
$21_a = \bigwedge^{2}7$.
\label{a8e8p0}
\end{table}

\begin{table}[H]\small
\caption{Alt$_8 < E_8$, $p = 7$}
%
\\
$21_a = \bigwedge^{2}7$.
\label{a8e8p7}
\end{table}

\begin{table}[H]\small
\caption{Alt$_8 < E_8$, $p = 5$}
%
\\
$21_a = \bigwedge^{2}7$.
\label{a8e8p5}
\end{table}

\begin{table}[H]\small
\caption{Alt$_8 < E_8$, $p = 3$}
%
\label{a8e8p3}
\end{table}

{\small
\begin{table}[H]\small
\caption{Alt$_7 < E_8$, $p = 0$ or $p \nmid |H|$}
%
\\
$14_a$ is a section of $6 \otimes 6$.
\label{a7e8p0}
\end{table}
}

\begin{table}[H]\small
\caption{Alt$_7 < E_8$, $p = 7$}
%
\\
$14_a$ is a section of $5 \otimes 5$.
\label{a7e8p7}
\end{table}

\begin{table}[H]\small
\caption{Alt$_7 < E_8$, $p = 5$}
%
\label{a7e8p2}
\end{table}

\begin{table}[H]\small
\caption{Alt$_6 < E_8$, $p = 0$ or $p \nmid |H|$}
%
Permutations: $(5_a,5_b)$, $(8_a,8_b)$.
\label{a6e8p0}
\end{table}

\begin{table}[H]\small
\caption{Alt$_6 < E_8$, $p = 5$}
%
Permutations: $(5_a,5_b)$.
\label{a6e8p5}
\end{table}

\begin{table}[H]\small
\caption{Alt$_5 < E_8$, $p = 0$ or $p \nmid |H|$}
%
\\
Permutations: $(3_a,3_b)$.
\label{a5e8p0}
\end{table}

\begin{table}[H]\small
\caption{Alt$_5 < E_8$, $p = 3$}
%
\\
Permutations: $(3_a,3_b)$.
\label{a5e8p3}
\end{table}

\subsection{Sporadic Groups} \ 

\begin{table}[H]\small
\caption{$M_{11} < E_8$, $p = 0$ or $p \nmid |H|$}
%
\\
Permutations: $(11_a,11_b)$.
\label{m12e8p5}
\end{table}

\begin{table}[H]\small
\caption{$M_{12} < E_8$, $p = 2$}
%
\label{m12e8p2}
\end{table}

\begin{table}[H]\small
\caption{$J_{1} < E_8$, $p = 11$}
%
\\
$77_a$ is a section of $\bigwedge^{2}14$.
\label{j1e8p11}
\end{table}

\begin{table}[H]\small
\caption{$J_{2} < E_8$, $p = 2$}
%
\\
Permutations: $(6_a,6_b)(14_a,14_b)(64_a,64_b)$.
\label{j2e8p2}
\end{table}

\begin{table}[H]\small
\caption{$J_{3} < E_8$, $p = 2$}
%
\label{j3e8p2}
\end{table}

\addtocounter{table}{1}
\begin{table}[H]\small
\thetable: $Th < E_8$, $p = 3$. Irreducible on $V_{248}$. \possprim
\label{the8p3}
\end{table}

\subsection{Cross-characteristic Groups $L_2(q)$ $(q \neq 4$, $5$, $9)$} \ 

\begin{table}[H]\small
\caption{$L_2(7) < E_8$, $p = 0$ or $p \nmid |H|$}
%

\begin{table}[H]\small
\caption{$L_2(8) < E_8$, $p = 7$}
%
\\
Permutations: $(7_b,7_c,7_d)$.
\label{l28e8p7}
\end{table}

\begin{table}[H]\small
\caption{$L_2(8) < E_8$, $p = 3$}
%
\\
Permutations: $(9_a,9_b,9_c)$.
\label{l28e8p3}
\end{table}

\begin{table}[H]\small
\caption{$L_2(11) < E_8$, $p = 0$ or $p \nmid |H|$}
%
\\
Permutations: $(12_a,12_b)$.
\label{l211e8p0}
\end{table}

\begin{table}[H]\small
\caption{$L_2(11) < E_8$, $p = 5$}
%
}

\begin{table}[H]\small
\caption{$L_2(13) < E_8$, $p = 7$}
%
\\
Permutations: $(7_{a},7_{b})$, $(12_a,12_b,12_c)$.
\label{l213e8p3}
\end{table}

\begin{table}[H]\small
\caption{$L_2(13) < E_8$, $p = 2$}
%
\\
Permutations: $(6_{a},6_{b})$, $(12_a,12_b,12_c)$.
\label{l213e8p2}
\end{table}

\begin{table}[H]\small
\caption{$L_2(16) < E_8$, $p = 0$ or $p \nmid |H|$}
%
\\
Permutations: $(15_a,15_b,\ldots,15_h)$. \\
$17_a$, $17_b$ and $17_c$ occur with equal multiplicities, as do $17_d$, $17_e$, $17_f$ and $17_g$.
\label{l216e8p0}
\end{table}

\begin{table}[H]\small
\caption{$L_2(16) < E_8$, $p = 17$}
%
\\
Permutations: $(17_b,17_c)(17_d,17_e,17_f,17_g)$.
\label{l216e8p17}
\end{table}

\begin{table}[H]\small
\caption{$L_2(16) < E_8$, $p = 5$}
%
\\
Permutations: $(15_a,15_b,\ldots,15_h)$.
\label{l216e8p5}
\end{table}

\begin{table}[H]\small
\caption{$L_2(16) < E_8$, $p = 3$}
%
\\
Permutations: $(15_a,15_b,\ldots,15_h)$.
\label{l216e8p3}
\end{table}

\begin{table}[H]\small
\caption{$L_2(17) < E_8$, $p = 0$ or $p \nmid |H|$}
%
\\
Permutations: $(9_a,9_b)$, $(16_b,16_c,16_d)$.
\label{l217e8p0}
\end{table}

\begin{table}[H]\small
\caption{$L_2(17) < E_8$, $p = 3$}
%
\\
Permutations: $(8_a,8_b)$, $(16_b,16_c,16_d)$.
\label{l217e8p2}
\end{table}

\begin{table}[H]\small
\caption{$L_2(19) < E_8$, $p = 0$ or $p \nmid |H|$}
%
\\
Permutations: $(18_a,18_b)(18_c,18_d)$, $(20_b,20_c,20_d)$.
\label{l219e8p0}
\end{table}

\begin{table}[H]\small
\caption{$L_2(19) < E_8$, $p = 5$}
%
\\
Permutations: $(20_b,20_c,20_d)$.
\label{l219e8p5}
\end{table}

\begin{table}[H]\small
\caption{$L_2(19) < E_8$, $p = 3$}
%
\\
Permutations: $(18_a,18_b)(18_c,18_d)$.
\label{l219e8p3}
\end{table}

\begin{table}[H]\small
\caption{$L_2(19) < E_8$, $p = 2$}
%
\\
Permutations: $(18_a,18_b)$, $(20_b,20_c,20_d)$.
\label{l219e8p2}
\end{table}

\begin{table}[H]\small
\caption{$L_2(25) < E_8$, $p = 0$ or $p \nmid |H|$}
%
\\
Permutations: $(24_a,24_b,\ldots,24_f)$.
\label{l225e8p3}
\end{table}

\begin{table}[H]\small
\caption{$L_2(25) < E_8$, $p = 2$}
%
\\
Permutations: $(12_a, 12_b)$, $(24_a,24_b, \ldots, 24_f)$.
\label{l225e8p2}
\end{table}

\begin{table}[H]\small
\caption{$L_2(27) < E_8$, $p = 0$, $13$ or $p \nmid |H|$}
%
\\
Permutations: $(26_a,26_b,26_c)(26_d,26_e,26_f)$.
\label{l227e8p0}
\label{l227e8p13}
\end{table}

\begin{table}[H]\small
\caption{$L_2(27) < E_8$, $p = 7$}
%
\label{l227e8p7}
\end{table}

\begin{table}\small
\caption{$L_2(27) < E_8$, $p = 2$}
%
Permutations: $(26_a,26_b,26_c)$, $(28_a,28_b, \ldots, 28_f)$.
\label{l227e8p2}
\end{table}

\begin{table}[H]\small
\caption{$L_2(29) < E_8$, $p = 0$ or $p \nmid |H|$}
%
\\
Permutations: $(15_a,15_b)$, $(30_a,30_b,30_c)$, $(15_a,15_e)(15_b,15_f)(15_c,15_g)(15_d,15_h)$.
\label{l229e8p0}
\end{table}

\begin{table}[H]\small
\caption{$L_2(29) < E_8$, $p = 7$}
%
\\
Permutations: $(15_a,15_b)$, $(30_a,30_b,30_c)$.
\label{l229e8p5}
\end{table}

\begin{table}[H]\small
\caption{$L_2(29) < E_8$, $p = 3$}
%
\\
Permutations: $(15_a,15_b)$, $(30_a,30_b,30_c)$.
\label{l229e8p3}
\end{table}

\begin{table}[H]\small
\caption{$L_2(29) < E_8$, $p = 2$}
%
\\
Permutations: $(14_a,14_b)$, $(30_a,30_b,30_c)$, $(28_b,28_c)(28_d,28_e,28_f,28_g)$.
\label{l229e8p2}
\end{table}

\begin{table}[H]\small
\caption{$L_2(31) < E_8$, $p = 0$ or $p \nmid |H|$}
%
\label{l231e8p2}
\end{table}

\begin{table}[H]\small
\caption{$L_2(32) < E_8$, $p = 0$, $31$ or $p \nmid |H|$}
%
\\
Permutations: $(31_a,31_b,\ldots,31_e)(31_f,31_g,\ldots,31_j)(31_k,31_l,\ldots,31_o)$.\\
Outer automorphism group induces this permutation. Elements of order 3 have Brauer character value $-2$ on $31_a,\ldots,31_e$, and $1$ on the rest.
\label{l232e8p0}
\label{l232e8p31}
\end{table}

\begin{table}[H]\small
\caption{$L_2(32) < E_8$, $p = 11$}
%
\\
Elements of order 3 have Brauer character value $-2$ on $31_a$, and $1$ on $31_b$.
\label{l232e8p11}
\end{table}

\begin{table}[H]\small
\caption{$L_2(32) < E_8$, $p = 3$}
%
\\
Permutations: $(31_a,31_b,31_c,31_d,31_e)$.\\
These modules are all conjugate by an outer automorphism.
\label{l232e8p3}
\end{table}

\begin{table}[H]\small
\caption{$L_2(37) < E_8$, $p = 2$}
%
\\
Permutations: $(18_a,18_b)$, $(38_b,38_c,38_d)$.
\label{l237e8p2}
\end{table}

\begin{table}[H]\small
\caption{$L_2(41) < E_8$, $p = 0$, $3$ or $p \nmid |H|$}
%
\\
Permutations: $(40_b,40_c,40_d)$.
\label{l241e8p0}
\label{l241e8p3}
\end{table}

\begin{table}[H]\small
\caption{$L_2(41) < E_8$, $p = 7$}
%
\\
Permutations: $(20_a,20_b)$, $(40_b,40_c,40_d)$, $(40_e,40_f)(40_g,40_h)(40_i,40_j)$.
\label{l241e8p2}
\end{table}

\begin{table}[H]\small
\caption{$L_2(49) < E_8$, $p = 0$ or $p \nmid |H|$}
%
\\
Permutations: $(48_a,48_b)$.
\label{l249e8p2}
\end{table}

\begin{table}[H]\small
\caption{$L_2(61) < E_8$, $p = 0$, $31$ or $p \nmid |H|$}
%
\label{l261e8p0}
\label{l261e8p31}
\end{table}

\begin{table}[H]\small
\caption{$L_2(61) < E_8$, $p = 5$}
%
\end{table}

\begin{table}[H]\small
\caption{$L_2(61) < E_8$, $p = 2$}
%
\label{l261e8p2}
\end{table}

\subsection{Cross-characteristic Groups $\ncong L_2(q)$} \ 

\begin{table}[H]\small
\caption{$L_3(3) < E_8$, $p = 0$ or $p \nmid |H|$}
%
\label{l33e8p2}
\end{table}

\begin{table}[H]\small
\caption{$L_3(5) < E_8$, $p \neq 2$, $5$}
%
\label{l35e8p3}
\label{l35e8p31}
\label{l35e8p0}
\end{table}

\begin{table}[H]\small
\caption{$L_3(5) < E_8$, $p = 2$}
%
\label{2b232e8p5}
\end{table}

%% file: AJL-auxiliary.tex

\appendix

\chapter{Auxiliary Data} \label{chap:auxiliary}

\label{sec:altreps}
\label{sec:sporadicreps}
\label{sec:l2qreps}
\label{sec:ccreps}

We collect here various data calculated for use in Chapters \ref{chap:disproving}--\ref{chap:gcr}, as well as references to data already in the literature. All of the following information is found either in the Atlas of Finite Groups \cite{MR827219}, the modular Atlas \cite{MR1367961}, the list of Hiss and Malle \cites{MR1835851,MR1942256}, or is straightforward to calculate using well-known routines which have been implemented in Magma \cite{MR1484478}.

In calculating the tables of feasible characters, we have made use of Brauer characters for various finite groups in Table \ref{tab:subtypes} and their proper covers. Of those we require, those which do not appear elsewhere in the literature are those of the alternating groups $\Alt_{n}$ with $13 \le n \le 17$ in characteristic $2$. Here, we give enough character values to verify that each feasible character appears in the relevant table of Chapter \ref{chap:thetables}.

\begin{table}[H]\small
\caption{Alt$_{n}$ Brauer Character Values, Degree $\le 248$, $p = 2$}
\begin{tabular}{c|*{12}{r}}
& \multicolumn{12}{c}{Cycle type/Character Value} \\ \hline
$H = \Alt_{17}$ & $3$ & $3^{2}$ & $3^{3}$ & $3^{4}$ & $3^{5}$ & $5$ & $5^{2}$ & $5^{3}$ & $7$ & $7^{2}$ & $11$ & $13$ \\ \hline
$16$ & $13$ & $10$ & $7$ & $4$ & $1$ & $11$ & $6$ & $1$ & $9$ & $2$ & $5$ & $3$ \\
$118$ & $76$ & $43$ & $19$ & $4$ & $-2$ & $53$ & $13$ & $-2$ & $34$ & $-1$ & $8$ & $1$ \\
$128_a$, $128_b$ & $-64$ & $32$ & $-16$ & $8$ & $-4$ & $-32$ & $8$ & $-2$ & $16$ & $2$ & $-4$ & $-2$ \\ \hline \hline
$H = \Alt_{16}$ & $3$ & $3^{2}$ & $3^{3}$ & $3^{4}$ & $3^{5}$ & $5$ & $5^{2}$ & $5^{3}$ & $7$ & $7^{2}$ & $11$ & $13$ \\ \hline
$14$ & $11$ & $8$ & $5$ & $2$ & $-1$ & $9$ & $4$ & $-1$ & $7$ & $0$ & $3$ & $1$ \\
$64_a$, $64_b$ & $-32$ & $16$ & $-8$ & $4$ & $-2$ & $-16$ & $4$ & $-1$ & $8$ & $1$ & $-2$ & $-1$ \\
$90$ & $54$ & $27$ & $9$ & $0$ & $0$ & $35$ & $5$ & $0$ & $20$ & $-1$ & $2$ & $-1$ \\ \hline \hline
$H = \Alt_{15}$ & \multicolumn{12}{c}{Same as $\Alt_{16}$} \\ \hline \hline
$H = \Alt_{14}$ & $3$ & $3^{2}$ & $3^{3}$ & $3^{4}$ & $5$ & $5^{2}$ & $7$ & $7^{2}$ & $11$ & $13$ \\ \cline{1-11}
$12$ & $9$ & $6$ &$3$ &$0$& $7$& $2$ &$5$& $-2$ &$1$ &$-1$ \\
$64_a$ & $34$ & $13$ &$1$ &$-2$ &$19$ &$-1$ &$8$ &$1$ &$-2$ &$-1$ \\
$64_b$ & $-32$ & $16$ &$-8$& $4$& $-16$& $4$ &$8$& $1$& $-2$ &$-1$ \\
$208$ & $76$ &$16$& $1$& $4$& $28$& $-2$& $5$& $-2$& $-1$& $0$ \\ \hline \hline
$H = \Alt_{13}$ & $3$ & $3^{2}$ & $3^{3}$ & $3^{4}$ & $5$ & $5^{2}$ & $7$ & $11$ & \multicolumn{2}{c}{$13_a$} & \multicolumn{2}{c}{$13_b$} \\ \cline{1-13}
$12$ & $9$ & $6$ & $3$ & $0$ & $7$ & $2$ & $5$ & $1$ & \multicolumn{2}{c}{$-1$} & \multicolumn{2}{c}{$-1$} \\
$32_a$ & $-16$ & $8$ & $-4$ & $2$ & $-8$ & $2$ & $4$ & $-1$ & \multicolumn{2}{c}{$\frac{-1 + \sqrt{13}}{2}$} & \multicolumn{2}{c}{$\frac{-1 - \sqrt{13}}{2}$} \\
$32_b$ & $-16$ & $8$ & $-4$ & $2$ & $-8$ & $2$ & $4$ & $-1$ & \multicolumn{2}{c}{$\frac{-1 - \sqrt{13}}{2}$} & \multicolumn{2}{c}{$\frac{-1 + \sqrt{13}}{2}$} \\
$64$ & $34$ & $13$ & $1$ & $-2$ & $19$ & $-1$ & $8$ & $-2$ & \multicolumn{2}{c}{$-1$} & \multicolumn{2}{c}{$-1$} \\
$144_a$, $144_b$ & $-48$ & $12$ & $0$ & $-3$ & $-16$ & $-1$ & $4$ & $1$ & \multicolumn{2}{c}{$1$} & \multicolumn{2}{c}{$1$} \\
$208$ & $76$ & $16$ & $1$ & $4$ & $28$ & $-2$ & $5$ & $-1$ & \multicolumn{2}{c}{$0$} & \multicolumn{2}{c}{$0$}
\end{tabular}
\end{table}

In addition to the Brauer characters, for many of the irreducible $KH$-modules encountered, we have made use of the Frobenius-Schur indicator, as well as the dimension of the first the cohomology group.

In the following tables we list all $KH$-modules of dimension at most 248, where $H$ is a finite simple group which embeds into an adjoint exceptional simple algebraic group over the algebraically closed field $K$. Frobenius-Schur indicators are taken from \cite{MR1942256}, and the dimension of the group $H^{1}(H,V)$, when given, has been calculated using Magma. Note that the list of \cite{MR1942256} does not distinguish between non-isomorphic $KH$-modules of the same dimension which have the same Frobenius-Schur indicator, however we list all modules in such a case, to clarify when cohomology groups differ.

Whenever $V \ncong V^{*}$ but $H^{1}(H,V) \cong H^{1}(H,V^{*})$, we omit $V^{*}$. Subscripts denote a collection of modules having identical properties, for example `($15_{a-h}$, $+$, 0)' denotes eight pairwise non-isomorphic modules of dimension 15, such that $H^{1}(H,V)$ vanishes for each such module $V$. 

\begin{longtable}{c|c|l}
\caption{$H$ alternating}\\
\hline $n$ & char $K = p$ & \multicolumn{1}{c}{($V$, ind($V$), dim $H^{1}(H,V)$)} \\ \hline
5 & 2 & ($2_{a,b}$, $-$, 1), (4, $+$, 0) \\
& 3 & ($3_{a,b}$, $+$, 0), (4, $+$, 1) \\
& 5 & (3, $+$, 1), (5, $+$, 0) \\
& $p \neq 2, 3, 5$ & ($3_{a,b}$, $+$), ($4_a$, $+$), (5, $+$) \\ \hline
6 & 2 & ($4_{a,b}$, $-$, 1), ($8_{a,b}$, $+$, 0), \\
& 3 & ($3_{a,b}$, $+$, 0), (4, $+$, 2), (9, $+$, 0) \\
& 5 & ($5_{a,b}$, $+$, 0), (8, $+$, 1), ($10_a$, $+$, 0) \\
& $p \neq 2, 3, 5$ & ($5_{a,b}$, $+$), ($8_{a,b}$, $+$), (9, $+$), (10, $+$) \\ \hline
7 & 2 & ($4$, $\circ$, 0), ($6_a$, $+$, 0), ($14$, $+$, 1), (20, $-$, 1) \\
& 3 & (6, $+$, 0), (10, $\circ$, 1), (13, $+$, 2), (15, $+$, 0) \\
& 5 & ($6_a$, $+$, 0), (8, $+$, 0), (10, $\circ$, 0), (13, $+$, 1), ($15_a$, $+$, 0),\\
& & (35, $+$, 0) \\
& 7 & (5, $+$, 1), (10, $+$, 0), ($14_{a,b}$, $+$, 0), ($21_a$, $+$, 0), (35, $+$, 0) \\
& $p \neq 2,3,5,7$ & ($6_a$, $+$), (10, $\circ$), ($14_{a,b}$, $+$), ($15_a$, $+$), ($21_a$, $+$), (35, $+$) \\ \hline
8 & 2 & ($4$, $\circ$, 0), ($6$, $+$, 1), ($14$, $+$, 1), (20, $\circ$, 1), (64, $+$, 0) \\
& 3 & (7, $+$, 0), (13, $+$, 1), (21, $+$, 0), (28, $+$, 0), (35, $+$, 1), \\
& & (45, $\circ$, 0) \\
& 5 & (7, $+$, 0), (13, $+$, 1), (20, $+$, 0), ($21_{a,b}$, $+$, 0), (35, $+$, 0), \\
& & (43, $+$, 0), (45, $\circ$, 0), (70, $+$, 0) \\
& 7 & (7, $+$, 0), (14, $+$, 0), (19, $+$, 1), ($21_a$, $+$, 0), ($21_b$, $\circ$, 0), \\
& &  (28, $+$, 0), (35, $+$, 0), (45, $+$, 0), (56, $+$, 0), (70, $+$, 0) \\
& $p \neq 2, 3, 5,7$ & (7, $+$), (14, $+$), (20, $+$), ($21_a$, $+$), ($21_b$, $\circ$), (28, $+$), \\
& & (35, $+$), (45, $\circ$), ($56_a$, $+$), ($64_a$, $+$), (70, $+$) \\ \hline
9 & 2 & ($8_{a,b,c}$, $+$, 0), (20, $\circ$, 1), (26, $+$, 2), (48, $+$, 0), (78, $+$, 1), \\
& & (160, $+$, 0) \\
& 3 & (7, $+$, 1), (21, $+$, 0), (27, $+$, 0), (35, $+$, 1), (41, $+$, 1), \\
& & (162, $+$, 0), (189, $+$, 0) \\
& 5 & ($8_a$, $+$, 0), (21, $+$, 0), (27, $+$, 0), (28, $+$, 0), (34, $+$, 0), ($35_{a,b}$, $+$, 0), \\
& & ($56_a$, $+$, 0), (83, $+$, 1), (105, $+$, 0), (120, $+$, 0), (133, $+$, 0), \\
& & (134, $+$, 0) \\
& 7 & ($8_a$, $+$, 0), (19, $+$, 0), (21, $\circ$, 0), (28, $+$, 0), ($35_{a,b}$, $+$, 0),\\
& & (42, $+$, 0), (47, $+$, 1), ($56_a$, $+$, 0), (84, $+$, 0), (101, $+$, 0),\\
& & (105, $+$, 0), (115, $+$, 0), (168, $+$, 0), (189, $+$, 0) \\
& $p \neq 2,3,5,7$ & ($8_a$, $+$), (21, $\circ$), (27, $+$), (28, $+$), ($35_{a,b}$, $+$),\\ 
& & (42, $+$), ($48_a$, $+$), ($56_a$, $+$), (84, $+$), (105, $+$), \\
& & ($120_a$, $+$), (162, $+$), (168, $+$), (189, $+$), (216, $+$) \\ \hline
10 & 2 & (8, $-$, 1), (16, $+$, 0), (26, $+$, 1), (48, $+$, 0), ($64_{a,b}$, $+$, 0), \\
& & (160, $+$, 0), (198, $+$, 2), (200, $+$, 1) \\
& 3 & (9, $+$, 0), (34, $+$, 1), (36, $+$, 0), (41, $+$, 1), (84, $+$, 1), \\
& & (90, $+$, 0), (126, $+$, 0), (224, $+$, 1) \\
& 5 & ($8_a$, $+$, 1), (28, $+$, 0), (34, $+$, 0), ($35_{a,b,c}$, $+$, 0), (55, $+$, 0), \\
& & ($56_a$, $+$, 0), (75, $+$, 0), ($133_{a,b}$, $+$, 0), \\
& & (155, $+$, 0), ($160_a$, $+$, 0), (217, $+$, 1), (225, $+$, 0) \\
& 7 & (9, $+$, 0), (35, $+$, 0), (36, $+$, 0), (42, $+$, 0), (66, $+$, 0), \\
& & (84, $+$, 0), (89, $+$, 1), (101, $+$, 0), (124, $+$, 0), (126, $+$, 0), \\
& & (199, $+$, 0), (210, $+$, 0), ($224_{a,b}$, $+$, 0) \\
& $p \neq 2, 3, 5, 7$ & (9, $+$), (35, $+$), (36, $+$), (42, $+$), (75, $+$), (84, $+$), (90, $+$), \\
& & (126, $+$), (160, $+$), (210, $+$), ($224_{a,b}$, $+$), (225, $+$) \\ \hline
11 & 2 & (10, $+$, 0), (16, $\circ$, 0), (44, $+$, 1), (100, $+$, 0), (144, $+$, 0),\\
& & (164, $-$, 1), (186, $+$, 1), (198, $+$, 2) \\
& 11 & (9, $+$, 1), (36, $+$, 0), (44, $+$, 0), (84, $+$, 0), (110, $+$, 0), \\
& & (126, $+$, 0), (132, $+$, 0), (165, $+$, 0), (231, $+$, 0) \\ \hline
12 & 2 & (10, $+$, 1), (16, $\circ$, 0), (44, $+$, 1), (100, $+$, 0), (144, $\circ$, 0),\\
& & (164, $-$, 1) \\
13 & 2 &(12, $+$, 0), (32$_{a,b}$, $+$, 0), (64, $-$, 2), (144, $\circ$, 0), (208, $+$, 0) \\
14 & 2 & (12, $-$, 1), (64$_a$, $-$, 1), (64$_b$, $+$, 0), (208, $+$) \\
15 & 2 & (14, $+$, 0), (64, $\circ$, 0), (90, $+$, 1) \\
16 & 2 & (14, $+$, 1), (64, $\circ$, 0), (90, $+$, 1) \\
17 & 2 & (16, $+$, 0), (118, $+$, 2), ($128_{a,b}$, $+$, 0)
\label{tab:altreps}
\end{longtable}

\begin{longtable}{c|c|l}
\caption{$H$ sporadic}\\
\hline $H$ & char $K = p$ & \multicolumn{1}{c}{($V$, ind($V$), dim $H^{1}(H,V))$} \\ \hline
\endfirsthead
$M_{11}$ & 2 & (10, $+$, 1), (16, $\circ$, 0), (44, $+$, 1) \\
& 3 & (5, $\circ$, 0), ($5^{*}$, $\circ$, 1), ($10_a$, $+$, 0), ($10_b$, $\circ$, 1), ($(10_b)^{*}$, $\circ$, 0),\\
& & (24, $+$, 0), (45, $+$, 0) \\
& 5 & ($10_a$, $+$, 0), ($10_b$, $\circ$, 0), (11, $+$, 0), (16, $\circ$, 0), ($16^{*}$, $\circ$, 1), \\
& & (45, $+$, 0), (55, $+$, 0) \\
& 11 & (9, $+$, 1), (10, $\circ$, 0), (11, $+$, 0), (16, $+$, 0), (44, $+$, 0), \\
& & (55, $+$, 0) \\
& $p \neq 2,3,5,11$ & ($10_a$, $+$), ($10_b$, $\circ$), (11, $+$), (16, $\circ$), (44, $+$), (45, $+$),\\
& & (55, $+$) \\ \hline
$M_{12}$ & 2 & (10, $+$, 2), (16, $\circ$, 0), (44, $+$, 1), (144, $+$, 0) \\
& 3 & ($10_{a,b}$, $+$, 1), (15, $\circ$, 1), (34, $+$, 1), ($45_{a-c}$, $+$, 0),\\
& & ($54$, $+$, 0), ($99$, $+$, 0) \\
& 5 & ($11_{a,b}$, $+$, 0), (16, $\circ$, 0), (45, $+$, 0), ($55_{a-c}$, $+$, 0),\\
& & (66, $+$, 0), (78, $+$, 0), (98, $+$, 1), ($120_a$, $+$, 0) \\
& 11 & ($11_{a,b}$, $+$, 0), (16, $+$, 0), (29, $+$, 0), (53, $+$, 1), \\
& & ($55_{a-c}$, $+$, 0), (66, $+$, 0), (91, $+$, 0), (99, $+$, 0), (176, $+$, 0) \\
& $p\neq 2,3,5,11$ & ($11_{a,b}$, $+$), (16, $\circ$), (44, $+$), (45, $+$), (54, $+$), ($55_{a-c}$, $+$),\\
& & (66, $+$), (99, $+$), ($120_a$, $+$), (144, $+$), (176, $+$)\\ \hline
$M_{22}$ & 2 & (10, $\circ$, 1), ($10^{*}$, $\circ$, 0), (34, $-$, 1), (70, $\circ$, 1), (98, $-$, 1) \\
& 5 & ($21_a$, $+$, 0), ($45_a$, $\circ$, 0), (55, $+$, 0), (98, $+$, 1), (133, $+$, 0),\\
& & ($210_a$, $+$, 0) \\ \hline
$J_1$ & 11 & (7, $+$, 0), (14, $+$, 0), (27, $+$, 0), (49, $+$, 0), (56, $+$, 0),\\
& & (64, $+$, 0), (69, $+$, 0), ($77_{a-c}$, $+$, 0), (106, $+$, 0), \\
& & (119, $+$, 1), (209, $+$, 0) \\ \hline
$J_2$ & 2 & ($6_{a,b}$, $-$, 1), ($14_{a,b}$, $+$, 0), (36, $+$, 0), ($64_{a,b}$, $+$, 0),\\
& & (84, $+$, 1), (160, $+$, 0) \\
& 3 & ($13_{a,b}$, $+$, 1), ($21_{a,b}$, $+$, 0), ($36_a$, $+$, 0), ($57_{a,b}$, $+$, 0),\\
& & (63, $+$, 0), (90, $+$, 0), (133, $+$, 0), ($189_{a,b}$, $+$, 0),\\
& & (225, $+$, 0) \\
& 5 & ($14_a$, $+$, 0), (21, $+$, 1), (41, $+$, 0), (70, $+$, 0), (85, $+$, 0), \\
& & (90, $+$, 0), (175, $+$, 0), (189, $+$, 0), (225, $+$, 0) \\
& 7 & ($14_{a,b}$, $+$, 0), ($21_{a,b}$, $+$, 0), (36, $+$, 0), (63, $+$, 0), \\
& & ($70_{a,b}$, $+$, 0), (89, $+$, 1), (101, $+$, 0), (124, $+$, 0),\\
& & (126, $+$, 0), (175, $+$, 0), ($189_{a,b}$, $+$, 0), (199, $+$, 0), \\
& & ($224_{a,b}$, $+$, 0) \\
& $p \neq 2,3,5,7$ & ($14_{a,b}$, $+$), ($21_{a,b}$, $+$), (36, $+$), (63, $+$), $(70_{a,b},+)$, (90, $+$),\\
& & ($126$, $+$), (160, $+$), (175, $+$), $(189_{a,b},+)$, \\
& & ($224_{a,b}$, $+$), ($225$, $+$)\\ \hline
$J_3$ & 2 & ($78_{a,b}$, $+$, 0), (80, $+$, 0), (84, $\circ$, 1), (244, $+$, 1)
\label{tab:sporadicreps}
\end{longtable}

\begin{longtable}{c|c|l}
\caption{$H \cong L_{2}(q)$}\\
\hline $q$ & char $K = p$ & \multicolumn{1}{c}{($V$, ind($V$), dim $H^{1}(H,V)$)} \\ \hline
\endfirsthead
7 & 3 & ($3$, $\circ$, 0), ($6_a$, $+$, 0), (7, $+$, 1) \\
	& $p \neq 2,3,7$ & (3, $\circ$), ($6_a$, $+$), (7, $+$), (8, $+$) \\ \hline
8 & 3 & (7, $+$, 1), ($9_{a-c}$, $+$, 0), \\
	& 7 & ($7_{a-d}$, $+$, 0), (8, $+$, 1) \\
	& $p \neq 2,3,7$ & ($7_{a-d}$, $+$), (8, $+$), $(9_{a-c}$, $+$) \\ \hline
11 & 2 & (5, $\circ$, 1), (10, $+$, 0), ($12_{a,b}$, $+$, 0) \\
	& 3 & (5, $\circ$, 0), ($10_a$, $+$, 1), ($12_{a,b}$, $+$, 0) \\
	& 5 & (5, $\circ$, 0), ($10_{a,b}$, $+$, 0), (11, $+$, 1) \\
	& $p \neq 2,3,5,11$ & (5, $\circ$), ($10_{a,b}$, $+$), (11, $+$), ($12_{a,b}$, $+$) \\ \hline
13 & 2 & ($6_{a,b}$, $-$, 1), ($12_{a-c}$, $+$, 0), (14, $+$, 0) \\
	& 3 & ($7_{a,b}$, $+$, 0), ($12_{a-c}$, $+$, 0), (13, $+$, 1) \\
	& 7 & ($7_{a,b}$, $+$, 0), ($12$, $+$, 1), ($14_{a,b}$, $+$, 0) \\
	& $p \neq 2,3,7,13$ & ($7_{a,b}$, $+$), ($12_{a-c}$, $+$), (13, $+$), ($14_{a,b}$, $+$) \\ \hline
16 & 3 & ($15_{a-h}$, $+$, 0), (16, $+$, 1), ($17_{a,b}$, $+$, 0) \\
	& 5 & ($15_{a-h}$, $+$, 0), (16, $+$, 1), (17, $+$, 0) \\
	& 17 & (15, $+$, 1), ($17_{a-g}$, $+$, 0) \\
	& $p \neq 2,3,5,17$ & ($15_{a-h}$, $+$), ($16_a$, $+$), ($17_{a-g}$, $+$) \\ \hline
17 & 2 & ($8_{a,b}$, $-$, 1), ($16_{a-d}$, $+$, 0) \\
	& 3 & ($9_{a,b}$, $+$, 0), (16, $+$, 1), ($18_{a-c}$, $+$, 0) \\
	& $p \neq 2,3,17$ & ($9_{a,b}$, $+$), ($16_{a-d}$, $+$), (17, $+$), ($18_{a-c}$, $+$) \\ \hline
19 & 2 & (9, $\circ$, 1), ($18_{a,b}$, $+$, 0), ($20_{a-d}$, $+$, 0) \\
	& 3 & (9, $\circ$, 0), ($18_{a-d}$, $+$, 0), (19, $+$, 1) \\
	& 5 & (9, $\circ$, 0), ($18_a$, $+$, 1), ($20_{a-d}$, $+$, 0) \\
	& $p \neq 2,3,5,19$ & (9, $\circ$), ($18_{a-d}$, $+$), (19, $+$), ($20_{a-d}$, $+$) \\ \hline
25 & 2 & ($12_{a,b}$, $-$, 1), ($24_{a-f}$, $+$, 0), (26, $+$, 0) \\
	& 3 & ($13_{a,b}$, $+$, 0), ($24_{a-f}$, $+$, 0), (25, $+$, 1), ($26_a$, $+$, 0) \\
	& 13 & ($13_{a,b}$, $+$, 0), (24, $+$, 1), ($26_{a-e}$, $+$, 0) \\
	& $p \neq 2,3,5,13$ & ($13_{a,b}$, $+$), ($24_{a-f}$, $+$), (25, $+$), ($26_{a-e}$, $+$) \\ \hline
27 & 2 & (13, $\circ$, 1), ($26_{a-c}$, $+$, 0), ($28_{a-f}$, $+$, 0) \\
	& 7 & (13, $\circ$, 0), ($26_a$, $+$, 1), ($28_{a-f}$, $+$, 0) \\
	& 13 & (13, $\circ$, 0), ($26_{a-f}$, $+$, 0), (27, $+$, 1) \\
	& $p \neq 2,3,7,13$ & (13, $\circ$), ($26_{a-f}$, $+$), (27, $+$), ($28_{a-f}$, $+$) \\ \hline
29 & 2 & ($14_{a,b}$, $-$, 1), ($28_{a-g}$, $+$, 0), ($30_{a-c}$, $+$, 0) \\
	& 3 & ($15_{a,b}$, $+$, 0), ($28_a$, $+$, 1), ($28_{b,c}$, $+$, 0), ($30_{a-f}$, $+$, 0) \\
	& 5 & ($15_{a,b}$, $+$, 0), ($28_a$, $+$, 1), ($28_b$, $+$, 0), ($30_{a-f}$, $+$, 0) \\
	& 7 & ($15_{a,b}$, $+$, 0), ($28_{a-g}$, $+$, 0), (29, $+$, 1) \\
	& $p \neq 2,3,5,7,29$ & ($15_{a,b}$, $+$), ($28_{a-g}$, $+$), (29, $+$), ($30_{a-f}$, $+$) \\ \hline
31 & 2 & (15, $\circ$, 1), ($32_{a-g}$, $+$, 0) \\
	& 3 & (15, $\circ$, 0), ($30_{a-g}$, $+$, 0), (31, $+$, 1), ($32_{a,b}$, $+$, 0) \\
	& 5 & (15, $\circ$, 0), ($30_{a-g}$, $+$, 0), (31, $+$, 1), (32, $+$, 0) \\
	& $p \neq 2,3,5,31$ & (15, $\circ$), ($30_{a-g}$, $+$), (31, $+$), ($32_{a-g}$, $+$) \\ \hline
32 & 3 & ($31_a$, $+$, 1), ($31_{b-f}$, $+$, 0), ($33_{a-o}$, $+$, 0) \\
	& 11 & ($31_a$, $+$, 1), ($31_{b}$, $+$, 0), ($33_{a-o}$, $+$, 0) \\
	& 31 & ($31_{a-p}$, $+$, 0), (32, $+$, 1) \\
	& $p \neq 2,3,11,31$ & ($31_{a-p}$, $+$), (32, $+$), ($33_{a-o}$, $+$) \\ \hline
37 & 2 & ($18_{a,b}$, $-$, 1), ($36_{a-i}$, $+$, 0), ($38_{a-d}$, $+$, 0) \\
	& 3 & ($19_{a,b}$, $+$, 0), ($36_{a-i}$, $+$, 0), (37, $+$, 1) \\
	& 19 & ($19_{a,b}$, $+$, 0), (36, $+$, 1), ($38_{a-h}$, $+$, 0) \\
	& $p \neq 2,3,19,37$ & ($19_{a,b}$, $+$), ($36_{a-i}$, $+$), (37, $+$), ($38_{a-h}$, $+$) \\ \hline
41 & 2 & ($20_{a,b}$, $-$, 1), ($40_{a-j}$, $+$, 0), ($42_{a,b}$, $+$, 0) \\
	& 3 & ($21_{a,b}$, $+$, 0), ($40_a$, $+$, 1), ($40_{b-d}$, $+$, 0), ($42_{a-i}$, $+$, 0) \\
	& 5 & ($21_{a,b}$, $+$, 0), ($40_{a-j}$, $+$, 0), (41, $+$, 1), ($42_a$, $+$, 0) \\
	& 7 & ($21_{a,b}$, $+$, 0), ($40_a$, $+$, 1), ($40_{b}$, $+$, 0), ($42_{a-i}$, $+$, 0) \\
	& $p \neq 2,3,5,7,41$ & ($21_{a,b}$, $+$), ($40_{a-j}$, $+$), (41, $+$), ($42_{a-i}$, $+$) \\ \hline
49 & 2 & ($24_{a,b}$, $-$, 1), ($48_{a-l}$, $+$, 0), (50, $+$, 0) \\
	& 3 & ($25_{a,b}$, $+$, 0), ($48_{a-l}$, $+$, 0), (49, $+$, 1), ($50_{a-c}$, $+$, 0) \\
	& 5 & ($25_{a,b}$, $+$, 0), (48, $+$, 1), ($50_{a-k}$, $+$, 0) \\
	& $p \neq 2,3,5,7$ & ($25_{a,b}$, $+$), ($48_{a-l}$, $+$), (49, $+$), ($50_{a-k}$, $+$) \\ \hline
61 & 2 & ($30_{a,b}$, $-$, 1), ($60_{a-o}$, $+$, 0), ($62_{a-g}$, $+$, 0) \\
	& 3 & ($31_{a,b}$, $+$, 0), ($60_{a-o}$, $+$, 0), (61, $+$, 1), ($62_{a-d}$, $+$, 0) \\
	& 5 & ($31_{a,b}$, $+$, 0), ($60_{a-o}$, $+$, 0), (61, $+$, 1), ($62_{a,b}$, $+$, 0) \\
	& 31 & ($31_{a,b}$, $+$, 0), (60, $+$, 1), ($62_{a-n}$, $+$, 0) \\
	& $p \neq 2,3,5,31,61$ & ($31_{a,b}$, $+$), ($60_{a-o}$, $+$), (61, $+$), ($62_{a-n}$, $+$)
\label{tab:l2qreps}
\end{longtable}

\begin{longtable}{c|c|>{\raggedright\arraybackslash}p{0.65\textwidth}}
\caption{$H$ of Lie type, $H \ncong L_{2}(q)$}\\
\hline $H$ & char $K = p$ & \multicolumn{1}{c}{($V$, ind($V$), dim $H^{1}(H,V)$)} \\ \hline
\endfirsthead
$L_3(3)$ & 2 & $(12,+,1)$, $(16_a,\circ,0)$, $(16_b,\circ,0)$, $(26,+,1)$ \\ 
 & 13 & $(11,+,1)$, $(13,+,0)$, $(16,+,0)$, $(26_a,+,0)$, $(26_b,\circ,0)$, $(39,+,0)$ \\
 & $p \neq 2,3,13$ & $(12,+)$, $(13,+)$, $(16_{a,b},\circ)$, $(26_a,+)$, $(26_b,\circ)$, $(27,+)$, $(39,+)$ \\ \hline
$L_3(4)$ & 3 & $(15_{a-c},+,0)$, $(19,+,2)$, $(45,\circ,0)$, $(63_{a,b},+,0)$ \\
 & 5 & $(20,+,0)$, $(35_{a-c},+,0)$, $(45_a,\circ,0)$, $(63,+,1)$ \\
 & 7 & $(19,+,1)$, $(35_{a-c},+,0)$, $(45,+,0)$, $(63_{a,b},+,0)$ \\
 & $p \neq 2,3,5,7$ & $(20,+)$, $(35_{a-c},+)$, $(45_a,\circ)$, $(63_{a,b},+)$, $(64,+)$ \\ \hline
$L_3(5)$ & 2 & $(30,+,1)$, $(96_{a-e},\circ,0)$, $(124_a,+,0)$, $(124_b,-,1)$ \\
 & 3 & $(30,+,0)$, $(31_a,+,0)$, $(31_b,\circ,0)$, $(96_{a-e},\circ,0)$, $(124_a,+,1)$, $(124_b,+,0)$, $(124_{c,d},\circ,0)$, $(186,+,0)$ \\
 & 31 & $(29,+,1)$, $(31_a,+,0)$, $(31_b,\circ,0)$, $(96,+,0)$, $(124_{a,b},+,0)$, $(124_{c-f},\circ,0)$, $(155_a,+,0)$, $(155_b,\circ,0)$, $(186,+,0)$ \\
 & $p \neq 2,3,5,31$ & $(30,+)$, $(31_a,+)$, $(31_b,\circ)$, $(96_{a-e},\circ)$, $(124_{a,b},+)$, $(124_{c-f},\circ)$, $(125,+)$, $(155_a,+)$, $(155_b,\circ)$, $(186,+)$ \\ \hline
$L_4(3)$ & 2 & $(26_{a,b},+,0)$, $(38,+,2)$, $(208_{a,b},+,0)$ \\ \hline
$L_4(5)$ & 2 & (154,+,2), $(248_{a,b},+,0)$ \\ \hline
$U_3(3)$ & 7 & $(6,-,0)$, $(7_a,+,0)$, $(7_b,\circ,0)$, $(14,+,0)$, $(21_a,+,0)$, $(21_b,\circ,0)$, $(26,+,1)$, $(28,\circ,0)$ \\
 & $p \neq 2,3,7$ & $(6,-)$, $(7_a,+)$, $(7_b,\circ)$, $(14,+)$, $(21_a,+)$, $(21_b,\circ)$, $(27,+)$, $(28,\circ)$, $(32,\circ)$ \\ \hline
$U_3(8)$ & 3 & $(56,-,1)$, $(133_{a-c},+,0)$ \\
 & 7 & $(56,-,0)$, $(57,\circ,0)$, $(133_{a-c},+,0)$ \\
 & 19 & $(56,-,0)$, $(57,\circ,0)$, $(133_{a-c},+,0)$ \\
 & $p \neq 2,3,7,19$ & $(57,\circ)$, $(133_{a-c},+)$ \\ \hline
$U_4(2)$ & 5 & $(5,\circ,0)$, $(6,+,0)$, $(10,\circ,0)$, $(15_{a,b},+,0)$, $(20_a,+,0)$, $(23,+,1)$, $(30_a,+,0)$, $(30_b,\circ,0)$, $(40,\circ,0)$, $(45,\circ,0)$, $(58,+,0)$, $(60_a,+,0)$ \\
 & $p \neq 2,3,5$ & $(5,\circ)$, $(6,+)$, $(10,\circ)$, $(15_{a,b},+)$, $(20_a,+)$, $(24,+)$, $(30_a,+)$, $(30_b,\circ)$, $(40,\circ)$, $(45,\circ)$, $(60_a,+)$, $(64_a,+)$, $(81,+)$ \\ \hline
$U_4(3)$ & 2 & $(20,+,1)$, $(34_{a,b},-,1)$, $(70_{a,b},\circ,0)$, $(120,+,0)$ \\ \hline
$PSp_4(5)$ & 2 & $(12_{a,b},-,1)$, $(40,+,0)$, $(64,-,1)$, $(104_{a,b},+,0)$, $(208_{a,b},+,0)$, $(248_{a,b},+,0)$ \\ \hline
$Sp_6(2)$ & 3 & $(7,+,0)$, $(14,+,1)$, $(21,+,0)$, $(27,+,0)$, $(34,+,1)$, $(35,+,0)$, $(49,+,0)$, $(91,+,0)$, $(98,+,1)$, $(189_{a-c},+,0)$, $(196,+,0)$ \\
 & 5 & $(7,+,0)$, $(15,+,0)$, $(21_{a,b},+,0)$, $(27,+,0)$, $(35_{a,b},+,0)$, $(56,+,0)$, $(70,+,0)$, $(83,+,1)$, $(105_{a-c},+,0)$, $(120,+,0)$, $(133,+,0)$, $(141,+,0)$, $(168_{a,b},+,0)$, $(210_{a,b},+,0)$ \\
 & 7 & $(7,+,0)$, $(15,+,0)$, $(21_{a,b},+,0)$, $(26,+,1)$, $(35_{a,b},+,0)$, $(56,+,0)$, $(70,+,0)$, $(84,+,0)$, $(94,+,0)$, $(105_{a-c},+,0)$, $(168,+,0)$, $(189_{a-c},+,0)$, $(201,+,0)$, $(210_{a,b},+,0)$ \\
& $p \neq 2,3,5,7$ & $(7,+)$, $(15,+)$, $(21_{a,b},+)$, $(27,+)$, $(35_{a,b},+)$, $(56,+)$, $(70,+$), $(84,+)$, $(105_{a-c},+)$, $(120,+)$, $(168,+)$, $(189_{a-c},+)$, $(210_{a,c},+)$, $(216,+)$ \\ \hline
$\Omega_8^{+}(2)$ & 3 & $(28,+,0)$, $(35_{a-c},+,0)$, $(48,+,2)$, $(147,+,0)$ \\
& 5 & $(28,+,0)$, $(35_{a-c},+,0)$, $(50,+,0)$, $(83_{a-c},+,1)$, $(175,+,0)$, $(210_{a-c},+,0)$ \\
& $p \neq 2,3,5$ & $(28,+,0)$, $(35_{a-c},+,0)$, $(50,+,0)$, $(84_{a-c},+,0)$, $(175,+,0)$, $(210_{a-c},+,0)$ \\ \hline
$G_2(3)$ & 2 & $(14,+,0)$, $(64,\circ,0)$, $(78,+,0)$, $(90_{a-c},+,1)$ \\
& 7 & $(14,+,0)$, $(64,\circ,0)$, $(78,+,0)$, $(91_{a-c},+,0)$, $(103,+,1)$, $(168,+,0)$, $(182_{a,b},+,0)$ \\
& 13 & $(14,+,0)$, $(64,\circ,0)$, $(78,+,0)$, $(91_{a-c},+,0)$, $(104,+,0)$, $(167,+,1)$, $(182_{a,b},+,0)$ \\
& $p \neq 2,3,7,13$ & $(14,+)$, $(64,\circ)$, $(78,+)$, $(91_{a-c},+)$, $(104,+)$, $(168,+)$, $(182_{a,b},+)$ \\ \hline
$^{3}D_4(2)$ & 3 & $(25,+,1)$, $(52,+,0)$, $(196,+,0)$ \\
& $p \neq 2,3$ & $(26,+,0)$, $(52,+,0)$, $(196,+,0)$ \\ \hline
$^{2}F_4(2)'$ & 3 & $(26,\circ,0)$, $(27,\circ,0)$, $(77,\circ,1)$, $(124_{a,b},+,0)$ \\
& 5 & $(26,\circ,0)$, $(27,\circ,1)$, $(27^{*},\circ,0)$, $(78,+,0)$, $(109_{a,b},+,0)$ \\
& $p \neq 2,3,5$ & $(26,\circ,0)$, $(27,\circ,0)$, $(78,+,0)$ \\ \hline
$^{2}B_2(8)$ & 5 & $(14,\circ,0)$, $(35_{a-c},+,0)$, $(63,+,1)$, $(65_{a-c},+,0)$ \\
& 7 & $(14,\circ,0)$, $(35_{a-c},+,0)$, $(64,+,1)$, $(91,+,0)$\\
& 13 & $(14,\circ,0)$, $(14^{*},\circ,1)$, $(35,+,0)$, $(65_{a-c},+,0)$, $(91,+,0)$ \\
& $p \neq 5,7,13$ & $(14,\circ)$, $(35_{a-c},\circ)$, $(64,+)$, $(65_{a-c},+)$, $(91,+)$ \\ \hline
$^{2}B_2(32)$ & 5 & $(124,\circ)$
\label{tab:ccreps}
\end{longtable}